\documentclass[12pt]{amsart}
\usepackage[all]{xy}
\usepackage[backend=bibtex,giveninits=true,maxbibnames=99,maxcitenames=99,maxalphanames=99,sorting=nyt,style=alphabetic,useprefix=true]{biblatex}
\usepackage{amsthm,amsmath,amsfonts,amssymb,mathrsfs,graphicx,subcaption,float,indentfirst,listings,verbatim,ifthen,tikz,tikz-cd,calc,enumitem,gensymb,multicol,lscape,ltablex,hyperref}
\usepackage[makeroom]{cancel}
\usepackage[capitalize]{cleveref}
\usepackage[margin=1in]{geometry}
\renewbibmacro{begentry}{\midsentence}
\usetikzlibrary{cd,decorations.pathreplacing,patterns}
\DeclareMathOperator{\sign}{sign}
\DeclareMathOperator{\spn}{span}
\DeclareMathOperator{\support}{supp}

\DeclareMathOperator{\kernel}{ker}
\DeclareMathOperator{\cokernel}{coker}
\DeclareMathOperator{\image}{im}
\DeclareMathOperator{\coimage}{coim}

\DeclareMathOperator{\diameter}{diam}
\DeclareMathOperator{\diagonal}{diag}
\DeclareMathOperator{\evaluation}{ev}
\DeclareMathOperator{\thickening}{th}
\DeclareMathOperator{\CF}{CF}
\DeclareNameAlias{author}{family-given}
\DeclareSymbolFont{bbold}{U}{bbold}{m}{n}
\DeclareSymbolFontAlphabet{\mathbbold}{bbold}
\DeclareRobustCommand{\gobblefive}[5]{}

\tikzcdset{scale cd/.style={every label/.append style={scale=#1},
    cells={nodes={scale=#1}}}}

\title{\(K\)-theory of two-dimensional substitution tiling spaces from \(AF\)-algebras}
\author{Jianlong Liu}

\theoremstyle{theorem}
\newtheorem{theorem}{Theorem}[section]
\newtheorem{proposition}[theorem]{Proposition}
\newtheorem{lemma}[theorem]{Lemma}
\newtheorem{corollary}[theorem]{Corollary}
\theoremstyle{definition}
\newtheorem{definition}[theorem]{Definition}
\newtheorem{example}[theorem]{Example}
\theoremstyle{remark}
\newtheorem{remark}[theorem]{Remark}
\newtheorem{warning}[theorem]{Warning}

\begin{document}
\sloppy

\begin{abstract}
Given a two-dimensional substitution tiling space, we show that, under some reasonable assumptions, the \(K\)-theory of the groupoid \(C^\ast\)-algebra of its unstable groupoid can be explicitly reconstructed from the \(K\)-theory of the \(AF\)-algebras of the substitution rule and its analogue on the \(1\)-skeleton. We prove this by generalizing the calculations done for the chair tiling in \cite{juliensavinien16} using relative \(K\)-theory and excision, and packaging the result into an exact sequence purely in topology. From this exact sequence, it appears that one cannot use only ordinary \(K\)-theory to compute using the dimension-filtration on the unstable groupoid. Several examples are computed using Sage and the results are compiled in a table.
\end{abstract}

\maketitle
\tableofcontents

\section{Introduction}
\indent Under suitable conditions, a \(d\)-dimensional substitution tiling space \(\Omega_T\) has its substitution rule as a self-homeomorphism. With this self-homeomorphism, the tiling space carries a Smale space structure, and an associated Smale bracket returns an unstable and a stable coordinate. The unstable coordinate turns out to correspond to the \(\mathbb{R}^d\) direction that is translations, and the stable coordinate the Cantor set direction that is substitutions, giving us the classical picture that locally, the tiling space is a product of the two. From the orbits returned of each of the two types of actions, one obtains a groupoid and a groupoid \(C^\ast\)-algebra (\cite{putnam96}). Due to its characterization, it was shown in \cite{goncalves11} and \cite{goncalvesramirezsolano17} that, for \(d\leq 2\), one can directly calculate the \(K\)-theory of the \(C^\ast\)-algebra of the stable groupoid using the substitution rule. Can one do the same for the unstable groupoid?\\
\indent To be more explicit, by identifying each prototile as a copy of \(\mathbb{C}\), the substitution rule induces a natural inclusion, up to unitary equivalence, of a direct sum of matrix algebras into another, each of the summands being a matrix algebra whose diagonal represents the tiles inside a given supertile. Its limit completes to an \(AF\)-algebra that includes into the groupoid \(C^\ast\)-algebra of the unstable groupoid, and this inclusion induces a map between their \(K\)-theories. We would like to see if this induced map (and perhaps its lower-dimensional analogues) is sufficient in reconstructing the \(K\)-theory of the unstable groupoid.\\
\indent While one can compute the \(K\)-theory of the unstable groupoid directly (\cite{savinien08}, \cite{savinienbellissard09}), one motivation for studying it from the perspective of \(AF\)-algebras is due to Elliott's classification of \(AF\)-algebras (\cite{elliott76}), where they are completely classified by their ordered \(K_0\)-groups (\(K_1\) being trivial). For \(d=1\), by using the order inherent on \(\mathbb{R}\), one can form an associated Bratteli--Vershik system whose underlying Bratteli diagram is generated from the substitution rule, giving us that the \(C^\ast\)-algebra of the unstable groupoid associated to the Bratteli--Vershik system is a nontrivial extension of the \(AF\)-algebra whose presence in the unstable groupoid is as an orbit-breaking subalgebra. It was computed in \cite{putnam89} that the \(K_0\)-group is generally a quotient of the \(K_0\)-group of the \(AF\)-algebra, and \(K_1\) is \(\mathbb{Z}\), realized as a generator of the translation action. In fact, by \cite{durandhostskau99}, one can reencode the substitution so that it becomes proper, i.e. every supertile begins with the same prototile and ends with the same prototile. This reencoding process leaves the tiling space unmodified, but the \(K_0\)-group of the unstable groupoid becomes order isomorphic to the \(K_0\)-group of this new \(AF\)-algebra. In higher dimensions, it was shown that there exists an order structure (\cite{kellendonk97} and \cite{putnam00} via traces; \cite{ormesradinsadun02} via the Perron--Frobenius eigenvector).\\
\indent There is also a physical reason for analyzing this induced inclusion map, at least for the \(K_0\)-group of the unstable groupoid. With the discovery of quasicrystals in \cite{shechtmanblechgratiascahn84}, it was realized that punctured tiling spaces can be used as models in their study, with one of the two main constructions of tiling spaces being those arising from substitution rules. One then attempts to characterize the spectra of the associated Schr\"odinger operators through gap labelling, which was discovered to be intricately related to the \(K_0\)-group, specifically the image of the \(K_0\)-group of the \(AF\)-algebra under the induced inclusion map (\cite{bellissard86}, \cite{bellissard92}, \cite{kellendonk95}, \cite{kellendonk97}). In fact, it was shown in \cite{vanelst94} that for \(d=2\), a quotient of the \(K_0\)-group of the \(AF\)-algebra recovers the \(K_0\)-group of the unstable groupoid, up to a direct summand of \(\mathbb{Z}\).\\
\indent In \cite{juliensavinien16}, it was demonstrated for the chair tiling that the \(K_0\)-group of the \(AF\)-algebra arising from the induced substitution rule on the boundaries of prototiles is required in reconstructing the \(K_1\)-group of the unstable groupoid. However, the authors there relied on splitting each of the supertiles into three squares and reencoding the substitution based on these squares, thus obtaining two induced substitution rules, one on the horizontal boundaries, the other on the vertical, both of which are stationary and one-dimensional. They then use a more general form of \cite{putnam89} in \cite{putnam97} and \cite{putnam98} to compute the \(K_1\)-group from these one-dimensional substitutions. It is unclear that if one attempts the same procedure on a general substitution rule, the induced substitution rules on the boundaries of the newly-split supertiles would remain stationary.\\
\indent For \(d\leq 2\), using the six-term sequence in relative \(K\)-theory (\cite{haslehurst21}) and an observation that, frequently, the \(K_0\)-groups of the \(AF\)-algebra obtained from inducing the substitution rule on lower-dimensional boundaries coincide exactly with the cochain group of the same dimension, we show that there is an isomorphism of exact sequences from the six-term sequence arising from the inclusion map of the \(AF\)-algebra into the \(C^\ast\)-algebra of the unstable groupoid to an exact sequence in topology. Using the observation, we recover the fact that a quotient (with a direct sum of \(\mathbb{Z}\) for \(d=2\)) of the \(K_0\)-group of the \(AF\)-algebra gives the \(K_0\)-group of the unstable groupoid, and the \(K_1\)-group is a subgroup of a quotient of the \(K_0\)-group of the \(AF\)-algebra of the one-dimensional substitution rule.
\begin{theorem}
\label{theorem:isomorphism-d-1}
For \(d=1\), we have an isomorphism of exact sequences
\[
\begin{tikzcd}
0\arrow[r]&\check{H}^0(\Omega_T)\arrow[r,hook]\arrow[d]&C^0\arrow[r,"\delta^0"]\arrow[d,"\thickening"]&C^1\arrow[r,two heads]\arrow[d,"\CF"]&\check{H}^1(\Omega_T)\arrow[r]\arrow[d]&0\\
0\arrow[r]&K_{1,u}\arrow[r,hook]&K_{0,AF;u}\arrow[r,"\evaluation"]&K_{0,AF}\arrow[r,"\iota_\ast"]&K_{0,u}\arrow[r]&0
\end{tikzcd}
\]
\noindent where \(K_{i,u}=K_i(C_r^\ast(\dot{G}_u))\), \(K_{i,AF}=K_i(C_r^\ast(\dot{G}_{AF}))\), and \(K_{i,AF;u}=K_i(C_r^\ast(\dot{G}_{AF});C_r^\ast(\dot{G}_u))\). Thus
\begin{align*}
K_0(C_r^\ast(\dot{G}_u))&\cong K_0(C_r^\ast(\dot{G}_{AF}))/\image\evaluation\\
K_1(C_r^\ast(\dot{G}_u))&\cong\mathbb{Z}.
\end{align*}
\end{theorem}
\begin{theorem}
\label{theorem:isomorphism-d-2}
For \(d=2\), we have an isomorphism of exact sequences
\[
\begin{tikzcd}
0\arrow[r]&\check{H}^1(\Omega_T)\arrow[r,hook]\arrow[ddd]&C^1/\image\delta^0\arrow[r,"\delta^1"]\arrow[ddd,"\thickening"]&C^2\arrow[r,two heads]&\check{H}^2(\Omega_T)\arrow[r]&0&\\[-28pt]
&&&\oplus&\oplus&\oplus&\\[-28pt]
&&\vphantom{C^1/\image\delta^0}&0\arrow[r]\arrow[d,"\CF"]&\check{H}^0(\Omega_T)\arrow[r,hook,two heads]\arrow[d]&\check{H}^0(\Omega_T)\arrow[r]\arrow[d]&0\\
0\arrow[r]&K_{1,u}\arrow[r,hook]&K_{0,AF;u}\arrow[r,"\evaluation"]&K_{0,AF}\arrow[r,"\iota_\ast"]&K_{0,u}\arrow[r,two heads]&K_{1,AF;u}\arrow[r]&0,
\end{tikzcd}
\]
\noindent where \(K_{i,u}=K_i(C_r^\ast(\dot{G}_u))\), \(K_{i,AF}=K_i(C_r^\ast(\dot{G}_{AF}))\), and \(K_{i,AF;u}=K_i(C_r^\ast(\dot{G}_{AF});C_r^\ast(\dot{G}_u))\). Thus, if \(T\) satisfies the boundary hyperplane condition,
\begin{align*}
K_0(C_r^\ast(\dot{G}_u))&\cong K_0(C_r^\ast(\dot{G}_{AF}))/\image\evaluation\oplus\mathbb{Z}\\
K_1(C_r^\ast(\dot{G}_u))&\cong\kernel[K_0(C_r^\ast(\dot{G}_{AF}^{(1)}))/\thickening(K_0(C_r^\ast(\dot{G}_{AF}^{(0)})))\xrightarrow{\evaluation} K_0(C_r^\ast(\dot{G}_{AF}))]
\end{align*}
\noindent where \(\thickening:K_0(C_r^\ast(\dot{G}_{AF}^{(0)}))\rightarrow K_0(C_r^\ast(\dot{G}_{AF}^{(1)}))\) returns an alternating sum of elements of \(K_0(C_r^\ast(\dot{G}_{AF}^{(1)}))\) that coincides with \(\delta^0\).
\end{theorem}
\noindent We then use this to compute the \(K\)-theory (in fact the six-term sequence in relative \(K\)-theory) for several examples, the two-dimensional ones being done via Sage and summarized in \cref{table:computation-d-2}. Many of these computations are not in the literature.\\
\indent While the statement and the proof of this result can be given very succinctly by starting from punctures, we have elected to take a more scenic route by starting from the unstable groupoid as defined on the tiling space, carrying it over to the Anderson--Putnam inverse limit by the Robinson map, constructing various subgroupoids, then restricting to the punctures. This allows us to very easily observe the existence \(AF\)-algebras for the boundaries of each of the intermediate dimensions, and how one can inductively construct the unstable groupoid from these \(AF\)-algebras, thereby recovering weaker forms of \cite[Theorem 3.22]{bellissardjuliensavinien10} and \cite[Theorem 4.9]{juliensavinien10} that are only as groupoid equivalences rather than homeomorphisms.\\
\indent Furthermore, this isomorphism of exact sequences relies on the classical fact that in low dimensions, \(K\)-theory is isomorphic to direct sums of \v Cech cohomology groups of the same parity, with a dimension shift. While there are a variety of proofs of this fact (e.g. the Chern and Connes--Thom isomorphisms, and the Pimsner--Voiculescu exact sequence), using any of them requires unwrapping said isomorphisms in order to provide the explicit correspondences between the generators and relations of each group. Since our isomorphism of exact sequences uses relative \(K\)-theory, we believe it to be in the spirit of the paper to reprove this fact purely using the six-term sequence in relative \(K\)-theory (and excision, \cite{putnam97}, \cite{putnam21}, \cite{putnam98}) by adapting the computations performed in \cite{juliensavinien16} for the chair tiling.\\
\indent In \cref{section:tiling-spaces}, we recall the basic definitions of substitution tiling spaces and sketch the topological conjugacy to the Anderson--Putnam inverse limit. In \cref{section:groupoids}, we construct various groupoids on the inverse limit and define their (reduced) \(C^\ast\)-algebras. In \cref{section:k-theory}, we briefly recall \(K\)-theory and state the facts we need from relative \(K\)-theory and excision. In \cref{section:k-theory-cohomology}, we first state and prove the existence of a map of exact sequences between topology and the six-term sequence in relative \(K\)-theory, then show that it is an isomorphism for \(d=1\). For \(d=2\), we reprove the isomorphism between \(K\)-theory and \v Cech cohomology, then derive our desired isomorphism as a corollary by realizing that this isomorphism coincides with the maps introduced at the beginning of this section. In \cref{section:computations}, for \(d=1\), we manually churn through the computations for each of the \(K\)-groups (including relative \(K_0\)) for the Fibonacci and the Silver Mean substitutions. We then illustrate the difference in computations for the two-dimensional dyadic solenoid, and leave the computations for the less trivial substitution tiling spaces to be done using a script in Sage, whose results are presented in a table. Other than the latter half of \cref{section:k-theory-cohomology} and \cref{section:computations}, our statements hold for tilings of arbitrary dimensions.

\section*{Acknowledgments}
\indent I am very thankful for my Ph.D. advisor, Rodrigo Trevi\~no, for proposing this problem, and the innumerous helpful discussions and guidance he has provided. This was supported in part by NSF grant DMS-\(2154762\) and Simons Collaboration Grant \(712227\).

\section{Substitution tiling spaces and inverse limits}
\label{section:tiling-spaces}
\indent In this section, we recall the necessary definitions and theorems of substitution tiling spaces and their inverse limit structure.

\subsection{Substitution tiling spaces}
\indent A \emph{tiling} is a partition of \(\mathbb{R}^d\) into compact sets with disjoint interiors (possibly with labels), each of which is homeomorphic to the closed ball in \(\mathbb{R}^d\). Each element in this partition is a \emph{tile}. The set of \emph{prototiles} is the quotient of the set of tiles by translations in \(\mathbb{R}^d\). A \emph{patch} is a finite collection of tiles whose union is homeomorphic to the closed ball. If \(P\) is a patch, denote \(\support P\), the \emph{support} of the patch, the union of the tiles in \(P\).\\
\indent Given a tiling \(T\) and an \(x\in\mathbb{R}^d\), denote \(T-x\) the tiling where each tile is translated by \(-x\). Alternatively, this is the tiling with the origin shifted to \(x\). Two tilings \(T_1\) and \(T_2\) are less than \(\epsilon\)-apart if there exist \(x_1,x_2\in B_\epsilon(0)\) such that the tilings \(T_1-x_1\) and \(T_2-x_2\) agree on \(B_{1/\epsilon}(0)\). This is the \emph{tiling metric}. The \emph{tiling space} associated to \(T\), \(\Omega_T\), is the closure of all of its translations under the tiling metric, that is
\[
\Omega_T=\overline{\{T-x:x\in\mathbb{R}^d\}}^\textnormal{tiling metric}.
\]
\noindent This metric yields a basis of \emph{cylinder sets} for the topology, defined as
\[
C(P,U)=\{T'\in\Omega_T:\exists x\in U\textnormal{ such that }P-x\subseteq T'\},
\]
\noindent where \(P\) is a patch and \(U\subseteq\support P\) is an open ball in \(\mathbb{R}^d\). If \(T\) arises from a substitution rule, we can refine this basis further to only use (unions of) supertiles (of the same level), since there is a mutual refining of topologies generated by both types of cylinder sets.\\
\indent We shall assume that
\begin{itemize}
\item
\(T\) is (strongly) \emph{aperiodic}, or \(\Omega_T\) contains no periodic components, and
\item
\(T\) has \emph{finite local complexity}, or given any \(r>0\), there exist finitely many patches, up to translation, of radius \(r\).
\end{itemize}%

\indent There is an important recurring property called \emph{extension} that we first state in terms of cylinder sets that allows us to extend (or reduce) the patch of a cylinder set.
\begin{proposition}[Extension]
Let \(P\) be a patch. If \(Q\) is a patch that contains \(P\) so that each occurrence of \(P\) in \(T\) induces the same inclusion of \(P\) in \(Q\) up to translation, then for any open \(U\subseteq\support P\), \(C(P,U)=C(Q,U)\). More generally, if \(\{Q_i\}_i\) is a collection of patches, each containing \(P\) so that each occurrence \(P\) in \(T\) is contained in some \(Q_i\), and each occurrence of \(Q_i\) in \(T\) contains \(P\) in the same relative location, then \(C(P,U)=\bigcup_iC(Q_i,U)\).
\begin{proof}
The inclusion \(C(Q,U)\hookrightarrow C(P,U)\) and the restriction \(C(P,U)|_{C(Q,U)}\) maps are well-defined, and are inverses to each other. The second statement is a straightforward generalization by partitioning the occurrences of \(P\) correctly.
\end{proof}
\end{proposition}
\indent Since translations are used to form the tiling space, they induce an action \(\mathbb{R}^d\ \rotatebox[origin=c]{-90}{\(\circlearrowright\)}\ \Omega_T\) by translations as well, by
\[
x\cdot T'=T'-x.
\]
\indent There is a notion of a \emph{punctured prototile}, consisting of a prototile with a choice of a point in its interior. One can then form a \emph{punctured tiling} by replacing each tile with a punctured tile and requiring the origin to be situated on a puncture, and a \emph{punctured tiling space}, where the only allowed translations are those moving between punctured tilings.\\
\indent A \emph{substitution rule}, \(\varsigma\), is a map from the set of tiles to the set of patches, formed from the composition of two operations
\begin{enumerate}
\item
\emph{Inflation}, where a tile (with its position relative to the origin) is expanded by a scalar greater than \(1\), and
\item
\emph{Subdivision}, where the inflated tile is subdivided into a set of tiles.
\end{enumerate}%

\noindent Substitution rules extend to patches. The patch that is the result of applying the substitution rule \(n\)-times to a single tile is a \emph{level-\(n\) supertile}, denoted \(\varsigma^n(t)\) with \(t\) a tile. A \emph{supertile} is a level-\(1\) supertile, and a level-\(0\) supertile is a tile. We may similarly consider patches and supertiles up to translation, which (unfortunately) are again called patches and supertiles. We will let context dictate which type we are referring to. A tiling \emph{arises from a substitution rule} if every patch is contained in a level-\(n\) supertile, for some sufficiently large \(n\).

\subsection{The Anderson--Putnam inverse limit}
\begin{figure}[t]
\centering
\begin{tikzpicture}
\fill[lightgray](1.25,0.5)circle[radius=0.1];
\draw(0,0)--(2,0)--(2,1)--(1,1)--(1,2)--(0,2)--(0,0);
\draw(1,0)--(1,0.5)--(0.5,0.5)--(0.5,1)--(0,1);
\draw(0.5,1)--(0.5,1.5)--(1,1.5);
\draw(1,0.5)--(1.5,0.5)--(1.5,1);
\draw[->](2.25,1)--(3.5,2.25)node[midway,above,rotate=45]{\(\pi_0\)};
\fill[lightgray](4.6,1.75)arc(0:180:0.1);
\draw(3.75,1.75)--(4.75,1.75)--(4.75,2.25)--(4.25,2.25)--(4.25,2.75)--(3.75,2.75)--(3.75,1.75);
\fill[lightgray](5.15,2.25)arc(180:360:0.1);
\draw(5,1.75)--(6,1.75)--(6,2.75)--(5.5,2.75)--(5.5,2.25)--(5,2.25)--(5,1.75);
\draw[->](2.25,1)--(3.5,-0.25)node[midway,below,rotate=-45]{\(\pi_1\)};
\fill[lightgray](5.25,-0.75)circle[radius=0.1];
\draw(3.75,-1.25)--(5.75,-1.25)--(5.75,-0.25)--(4.75,-0.25)--(4.75,0.75)--(3.75,0.75)--(3.75,-1.25);
\end{tikzpicture}
\caption{A cylinder set in a substitution tiling space where we have fixed the (level-\(1\)) supertile and a small ball of origins (light gray), as projected onto the first two coordinates of the inverse system.}
\label{figure:inverse-limit-projection}
\end{figure}
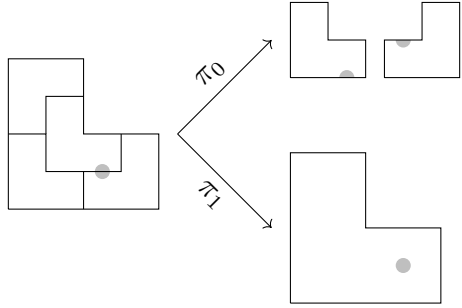%

\indent If \(T\) arises from a substitution rule, then one can define an inverse limit space using supertiles of all levels, each viewed as the support of the corresponding patch. More precisely, given any \(n\in\mathbb{N}\), let the \emph{level-\(n\) Anderson--Putnam complex}, denoted \(AP_n\), be formed from the quotient of the set of (supports of) level-\(n\) supertiles (up to translation) by the identifications given in the tiling \(T\). In other words, we identify the boundary of the supertiles to yield a continuous projection map \(T\twoheadrightarrow AP_n\) that extends to
\[
\pi_n:\Omega_T\twoheadrightarrow AP_n
\]
\noindent where a tiling is sent to the point in \(AP_n\) corresponding to the placement of the origin relative to the level-\(n\) supertiles it is situated in (\cref{figure:inverse-limit-projection}). The substitution rule induces a map
\[
\sigma:AP_n\rightarrow AP_{n-1}
\]
by sending each \(d\)-cell that is a level-\(n\) supertile to a collection of \(d\)-cells obtained by subdividing the level-\(n\) supertile into level-\(n-1\) supertiles. This gives us an inverse system of \(CW\)-complexes, more particularly of branched manifolds, whose inverse limit is denoted \(\varprojlim_n(AP_n,\sigma)\), and, since each \(\pi_n\) is compatible with \(\sigma\), a canonical projection map
\[
\pi:\Omega_T\twoheadrightarrow\varprojlim_n(AP_n,\sigma),
\]
 called the \emph{Robinson map} that was first defined in \cite{kellendonk95}. Letting \(AP_n^{(k)}\) denote the \(k\)-skeleton of the corresponding \(CW\)-complex, we similarly get an inverse limit \(\varprojlim_n(AP_n^{(k)},\sigma)\) for each \(0\leq k\leq d\), that we call the \emph{\(k\)-skeleton} of \(\varprojlim_n(AP_n,\sigma)\).\\
\indent Hereafter, we shall assume that our tiling space \(\Omega_T\) arises from a substitution rule \(\sigma\). Furthermore, we will assume that
\begin{itemize}
\item
\(\sigma\) is \emph{primitive}, or there exists an \(n\in\mathbb{N}\) such that each level-\(n\) supertile contains all of the prototiles, and
\item
\(T\), and therefore each of the \(AP\)-complexes, is \emph{collared}, or the supertiles in \(T\) are further labelled by their incident same-level supertiles.
\end{itemize}%

\noindent The substitution map induces a map on collared \(AP\)-complexes, and one can form a corresponding inverse limit and obtain an analogous Robinson map.\\
\indent We have \cite[Theorem 4.3]{andersonputnam98}.
\begin{theorem}[Anderson--Putnam]
\label{theorem:anderson-putnam}
The Robinson map is a homeomorphism.
\end{theorem}
\noindent We will only provide explanation and intuition to convince the reader of why the Robinson map is injective after collaring. From collaring, we will then see that there exists an action by translations by \(\mathbb{R}^d\) on \(\varprojlim_n(AP_n,\sigma)\), even though an action does not appear to exist on a finite level if \(AP_n\) is branched. Under this action, the Robinson map turns into a topological conjugacy.\\
\indent Suppose that we are given a tile whose boundary contains the origin. Substituting infinitely many times may only yield a partial tiling of \(\mathbb{R}^d\), particularly if the tile has its support on one side of a hyperplane intersecting the origin. There may be multiple tilings that contain this partial tiling. Applying the Robinson map gives that all such tilings map to the same point, and thus the map cannot be injective. Phrased differently, the inverse limit is one of branched manifolds, and it may be that the limit still results in at least one branching point. This branching point gives multiple tilings that project to it.\\
\indent One way to resolve this is with collaring, since placing the origin on the boundary of a collared tile still results in the origin being completely surrounded by tiles coming from the boundary tile labels, and substituting infinitely always results in a complete tiling, even though, formally, we lack the knowledge of the pattern on all of \(\mathbb{R}^d\). Therefore, the moral is that \emph{under the collaring process, branches eventually disappear.} Let us phrase this more precisely with a definition and a proposition.
\begin{definition}
For \(k\geq 1\), a set is
\begin{itemize}
\item
\emph{Contractible in \(AP_n^{(k)}\)} if its restriction and quotient to each \(AP_n^{(\ell)}/AP_n^{(\ell-1)}\) is contractible for all \(1\leq\ell\leq k\), and
\item
\emph{Unbranched in \(AP_n^{(k)}\)} if it is the homeomorphic image of an open ball in \(\mathbb{R}^k\).
\end{itemize}%

\end{definition}
\noindent Part of the goal of this definition of contractibility is to prevent an unbranched set from intersecting too many supertiles of a certain level, because such sets should be thought of as belonging to a higher-level \(AP\)-complex.
\begin{proposition}[Border-forcing]
\label{proposition:ap-open-flat-branched}
There exists an \(N\in\mathbb{N}\) such that for any open ball \(U\subseteq AP_N\) with \(\sigma^N(U)\subseteq AP_0\) contractible, \(\sigma^N(U)\) is unbranched. For such an \(N\) and for all \(n\geq N\), this holds between \(AP_n\) and \(AP_{n-N}\).
\begin{proof}
Let \(c(t)\) denote the patch formed from \(t\) and its collars. Let \(N\) be the level at which for all collared prototiles \(t\), for all \(t'\in\varsigma^N(c(t)\backslash t)\) so that \(t'\) intersects the boundary of \(\support\varsigma^N(t)\), \(c(t')\subseteq\varsigma^N(c(t))\). In other words, \(\varsigma^N(c(t))\) contains the neighbors of neighbors of \(\varsigma^N(t)\). This \(N\) exists since the substitution rule is expansive.\\
\indent With this choice of \(N\), the boundary collared tiles of \(\varsigma^N(t)\) are explicitly determined. More precisely, given a collection \(\{t_i\}_i\) of all collared tiles so that each \(t_i\) intersects \(t\) at the same (branched) boundary and each forms a patch with \(t\) along that boundary in \(T\), the collared tiles in \(\varsigma^N(t_i)\) that intersect the boundary of \(\support\varsigma^N(t)\) are identical among all \(i\). Picking a \(t_i\) then repeating the same argument inductively until all supertiles intersecting this boundary has been exhausted gives that the projection \(\pi_N\) to \(AP_0\) of the small neighborhood of this boundary in \(AP_N\) is unbranched.\\
\indent For the first statement, \(U\subseteq AP_N\) being contractible in \(AP_0\) means that for an appropriate choice of boundary, it is contained inside this small neighborhood.\\
\indent The second statement immediately follows from self-similarity of the substitution rule.
\end{proof}
\end{proposition}
\noindent Using terminology from \cite{kellendonk95}, the minimum such \(N\) is called the level at which the substitution rule \emph{forces the border}, and the boundary tiles are the \emph{forced border}. For \(n\geq N\), we take the forced border to be the patch obtained from applying \(\varsigma^{n-N}\) to the forced border at level-\(N\)\footnote{Strictly speaking, the forced border can be taken to be larger, since the border itself may yield additional forced border.}. Hereafter, when we refer to a (super)tile, we will always assume that it is collared.\\
\indent We have the following corollary due to \cite{kellendonk95} as an immediate consequence to the proof.
\begin{corollary}[Extension; Kellendonk]
Let \(N\) be the level at which the substitution rule forces the border, and let \(n\geq N\). If \(\varsigma^n(t)\) is a supertile and \(P\) is the patch formed from the supertile together with its forced border, then \(C(\varsigma^n(t),U)=C(P,U)\).
\begin{proof}
The forced border is unique and appears with every occurrence of \(\varsigma^n(t)\). An application of extension gives equality between the two cylinder sets.
\end{proof}
\end{corollary}
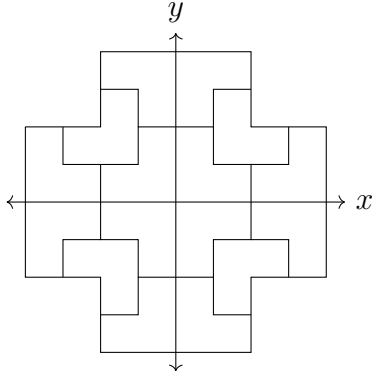
\begin{figure}[t]
\centering
\begin{tikzpicture}
\draw[<->](-2.25,0)--(2.25,0);
\draw[<->](0,-2.25)--(0,2.25);
\draw(2.25,0)node[right]{\(x\)};
\draw(0,2.25)node[above]{\(y\)};
\draw(-2,1)--(-2,-1)--(-1,-1)--(-1,-2)--(1,-2)--(1,-1)--(2,-1)--(2,1)--(1,1)--(1,2)--(-1,2)--(-1,1)--(-2,1);
\draw(-1.0,.5)--(-1.0,-.5)--(-.5,-.5)--(-.5,-1.0)--(.5,-1.0)--(.5,-.5)--(1.0,-.5)--(1.0,.5)--(.5,.5)--(.5,1.0)--(-.5,1.0)--(-.5,.5)--(-1.0,.5);
\draw(-1.5,1)--(-1.5,0.5)--(-1,0.5);
\draw(-1.5,-1)--(-1.5,-0.5)--(-1,-0.5);
\draw(1.5,1)--(1.5,0.5)--(1,0.5);
\draw(1.5,-1)--(1.5,-0.5)--(1,-0.5);
\draw(-1,1.5)--(-0.5,1.5)--(-0.5,1);
\draw(-1,-1.5)--(-0.5,-1.5)--(-0.5,-1);
\draw(1,1.5)--(0.5,1.5)--(0.5,1);
\draw(1,-1.5)--(0.5,-1.5)--(0.5,-1);
\end{tikzpicture}
\caption{The substitution rule sends this particular chair tiling to itself. Observing that the \(1\)-skeleton contains the \(x\)- and \(y\)-axes and is sent to itself under the induced substitution rule, the origin remains branched in \(\protect\varprojlim_n(AP_n^{(1)},\sigma)\).}
\label{figure:branched-1-skeleton}
\end{figure}%

\begin{warning}
The same cannot be said about the lower-dimensional skeleta! For example, one can observe this failure in the chair tiling on the \(1\)-skeleton (\cref{figure:branched-1-skeleton}). More generally, this exists in the \(1\)-skeleton for all two-dimensional substitution tiling spaces.
\end{warning}
\begin{figure}[t]
\centering
\begin{tikzpicture}
\draw[->](0.5176,1.9319)--(0.2588,0.9659);
\draw(0.2588,0.9659)--(0,0);
\path(0,0)--(0.5176,1.9319)node[midway,above,rotate=75]{\((a)b(b)\)};
\draw[->](0,0)--(0.9659,0.2588);
\draw(0.9659,0.2588)--(1.9319,0.5176);
\path(0,0)--(1.9319,0.5176)node[midway,below,rotate=15]{\((b)b(a)\)};
\draw[->](1.9319,0.5176)--(1.2247,1.2247);
\draw(1.2247,1.2247)--(0.5176,1.9319);
\path(0.5176,1.9319)--(1.9319,0.5176)node[midway,above,rotate=-45]{\((b)a(b)\)};
\draw[->](0,0)arc(45:225:0.3183);
\draw(-0.4502,-0.4502)arc(-135:45:0.3183);
\path(-0.9003,0)--(0,-0.9003)node[midway,below,rotate=-45]{\((b)b(b)\)};
\end{tikzpicture}
\caption{The \(AP\)-complex of the collared Silver Mean substitution, \(a\mapsto b\), \(b\mapsto bab\). The prototile \((b)b(b)\) has its endpoints attached, but \((b)b.b(b)\) is not a patch in the tiling.}
\label{figure:silver-mean-allowed}
\end{figure}
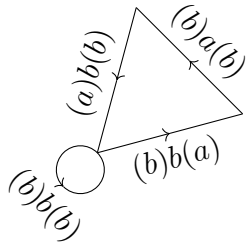%

\indent Even after collaring, not all contractible, unbranched \(U\subseteq AP_n\) can be realized as the Robinson map applied to some patch. In the proof of the above proposition, this is the reason for referring back to patches in the tiling and inducing on the pairs of branches that are images of patches under the Robinson map. See the Silver Mean substitution, for example (\cref{figure:silver-mean-allowed}). This forces us make the following definition intrinsic to the inverse system that we will show is the analogue of a patch in the tiling.
\begin{definition}
By an \emph{allowed} \(U\subseteq AP_n\), we mean one so that there exists an \(m\geq n\) and an open ball \(V\subseteq AP_m\) in the subspace topology so that \(\sigma^{m-n}(V)=U\).
\end{definition}
\noindent Let us abbreviate ``allowed, contractible, and unbranched'' by \hypertarget{acu}{\emph{acu}}.
\begin{proposition}
Acu sets correspond to patches in the tiling.
\begin{proof}
Given an acu \(U\subseteq AP_n\) with a preimage \(V\subseteq AP_m\) that is an open ball, since the Robinson map is a homeomorphism, there exist collections of level-\(m\) supertiles each of whose union forms a patch under the inverse of the Robinson map and whose intersection with \(V\) projects to \(U\) via \(\sigma^{m-n}\). This makes \(\pi^{-1}(U)\) subsets of each of the patches, and \(\pi^{-1}(U)\) itself is therefore a collection of patches. In fact, existence of a single such preimage implies that for sufficiently large \(m\), each of the preimages \(\sigma^{-(m-n)}(U)\) will be an open ball, since a patch forms a cylinder set, and its image under the Robinson map is open.\\
\indent Conversely, given any patch \(P\), we consider the level \(n\) at which \(\pi_n(P)\) is contractible and a collection of open balls. This exists since the substitution rule is expansive. Each open ball can be further written as a union of acu sets.
\end{proof}
\end{proposition}
\noindent These two propositions tell us that \hyperlink{acu}{acu} sets and the standard basis for the subspace topology in the \(AP\)-complexes mutually refine each other in \(\varprojlim_n(AP_n,\sigma)\). In light of this, let us restrict the cylinder sets in \(\Omega_T\) to those that have the second parameter satisfying the contractibility condition under \(\pi_0\).\\
\indent We can now give a nice description of the cylinder sets in \(\varprojlim_n(AP_n,\sigma)\), which have the form
\[
C(U_m)=\{(p_0,\ldots)\in\varprojlim_n(AP_n,\sigma):p_m\in U_m\}
\]
\noindent where \(U_m\subseteq AP_m\) is \hyperlink{acu}{acu}. While it appears to not be open in the subspace topology of the \(AP\)-complexes, border-forcing says that one can look at sufficiently-high preimages to obtain (unions of) open balls in the correct topology.\\
\indent Notice that, on taking the inverse system at face value, it appears impossible for there to be a translation action, since it may be that the translation vector is large enough to cross from one supertile into another, thus cross a branch, where it becomes ill-defined. In fact, this is possible had we not collared. Border-forcing resolves this issue of crossing a boundary.
\begin{figure}[t]
\centering
\begin{tikzpicture}
\draw[lightgray,line width=0.1cm](0.3586,0.6414)--(1.6414,-0.6414);
\draw[lightgray,line width=0.1cm](0.3586,-0.6414)--(1.6414,0.6414);
\draw[->](0,1.0)--(.25,.75);
\draw(.25,.75)--(1.0,0);
\draw[->](0,-1.0)--(.25,-.75);
\draw(.25,-.75)--(1.0,0);
\draw[->](1.0,0)--(1.75,.75);
\draw(1.75,.75)--(2.0,1.0);
\draw[->](1.0,0)--(1.75,-.75);
\draw(1.75,-.75)--(2.0,-1.0);
\fill(0.5,0.5)circle[radius=0.05];
\draw(1,-1)node[below]{\(AP_{m+N}\)};
\draw[->](2.25,0)--(2.75,0)node[midway,above]{\(\sigma^N\)};
\draw[lightgray,line width=0.1cm](3.3586,0.6414)--(4,0);
\draw[lightgray,line width=0.1cm](4,0)--(4.6414,0.6414);
\draw[->](3,1.0)--(3.25,.75);
\draw(3.25,.75)--(4.0,0);
\draw[->](3,-1.0)--(3.25,-.75);
\draw(3.25,-.75)--(4.0,0);
\draw[->](4.0,0)--(4.75,.75);
\draw(4.75,.75)--(5.0,1.0);
\draw[->](4.0,0)--(4.75,-.75);
\draw(4.75,-.75)--(5.0,-1.0);
\fill(3.5,0.5)circle[radius=0.05];
\draw(4,-1)node[below]{\(AP_m\)};
\end{tikzpicture}
\caption{We want to define a translation on \(AP_m\) on a point that crosses over a branching point (right). By collaring (and border-forcing), any open ball in the branched manifold topology that contains both the point and the branching point (left, light gray) has its image under \(\sigma^N\) unbranched (right, light gray). The choice of the branch associated to the translation is given by this unbranched image.}
\label{figure:inverse-limit-action}
\end{figure}
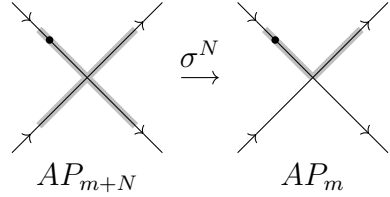%

\indent Suppose that we are given a \((p_0,\ldots)\in\varprojlim_n(AP_n,\sigma)\) and an \(x\in\mathbb{R}^d\) to act on the point. Let \(N\) be the level at which the substitution forces the border, and let \(n\) be the level at which \(\pi_n\) on every patch of radius \(x+\epsilon\) is contractible. Then, for all \(m\geq n\), patches of radius \(x+\epsilon\) also have their images under \(\pi_m\) contractible. For each such \(m\), let \(U\subseteq AP_{m+N}\) be the preimage of the set of \(d\)-cells incident to \(p_m\) that contains \(p_{m+N}\). By border-forcing, its image under \(\sigma^N\) is unbranched. By construction, this image contains the well-defined translation \(p_m+x\) (\cref{figure:inverse-limit-action}). Repeating this for all \(m\geq n\) picks a branch in each of the \(AP_m\)-complexes and a point in it compatible with \(\sigma\). Applying \(\sigma\) sufficiently many times yields the rest of the coordinates between \(0\) and \(n\).\\
\indent We denote this action by \((p_0,\ldots)+x\), as the coordinates themselves represent the location of the origin in the tiling relative to each of the supertiles.\\
\indent We have proven the following corollary to \cref{theorem:anderson-putnam} that is due to \cite{kellendonk95}.
\begin{corollary}[Kellendonk]
\label{corollary:robinson-map-conjugacy}
Under the action \(\mathbb{R}^d\ \rotatebox[origin=c]{-90}{\(\circlearrowright\)}\ \varprojlim_n(AP_n,\sigma)\), where
\[
x\cdot(p_0,\ldots)=(p_0,\ldots)+x,
\]
\noindent the Robinson map is a topological conjugacy.
\end{corollary}
\indent Summarized succinctly, border-forcing allows us to construct, starting at a sufficiently-high index, a well-defined (up to extension) sequence of \hyperlink{acu}{acu} sets compatible with substitution (up to extension), each of which containing the corresponding coordinate \(p_n\) and its translation \(p_n+x\).

\section{Groupoid \(C^\ast\)-algebras}
\label{section:groupoids}
\begin{figure}[t]
\centering
\begin{tikzpicture}
\fill[lightgray](.125,.125)--(.625,.125)--(.625,.375)--(.375,.375)--(.375,.625)--(.125,.625)--(.125,.125);
\draw(0,0)--(2,0)--(2,1)--(1,1)--(1,2)--(0,2)--(0,0);
\draw(1,0)--(1,0.5)--(0.5,0.5)--(0.5,1)--(0,1);
\draw(0.5,1)--(0.5,1.5)--(1,1.5);
\draw(1,0.5)--(1.5,0.5)--(1.5,1);
\draw[->](2.25,1)--(3.75,1)node[midway,above]{\(+B_\epsilon^{(d)}(x)\)};
\fill[lightgray](4.125,.125)--(4.625,.125)--(4.625,.375)--(4.375,.375)--(4.375,.625)--(4.125,.625)--(4.125,.125);
\draw(4,0)--(6,0)--(6,1)--(5,1)--(5,2)--(4,2)--(4,0);
\draw(5,0)--(5,.5)--(4.5,.5)--(4.5,1)--(4,1);
\draw(4.5,1)--(4.5,1.5)--(5,1.5);
\draw(5,.5)--(5.5,.5)--(5.5,1);
\draw(6.25,1)node{\(+\)};
\fill[lightgray](6.75,1)circle[radius=0.25];
\draw[->](6.5,.75)--(6.75,1);
\draw(7.25,1)node{\(=\)};
\fill[lightgray](7.625,.375)--(7.875,.125)--(8.375,.125)--(8.625,.375)--(8.625,.625)--(8.125,1.125)--(7.875,1.125)--(7.625,.875)--(7.625,.375);
\fill[lightgray](7.875,.375)circle[radius=0.25];
\fill[lightgray](8.375,.375)circle[radius=0.25];
\fill[lightgray](8.375,.625)circle[radius=0.25];
\fill[lightgray](8.125,.875)circle[radius=0.25];
\fill[lightgray](7.875,.875)circle[radius=0.25];
\draw(7.5,0)--(9.5,0)--(9.5,1)--(8.5,1)--(8.5,2)--(7.5,2)--(7.5,0);
\draw(8.5,0)--(8.5,.5)--(8.0,.5)--(8.0,1)--(7.5,1);
\draw(8.0,1)--(8.0,1.5)--(8.5,1.5);
\draw(8.5,.5)--(9.0,.5)--(9.0,1);
\end{tikzpicture}
\caption{A partial translation taking a cylinder set (choice of origin in light gray) to the cylinder set plus a ball.}
\label{figure:translation-cylinder-set}
\end{figure}
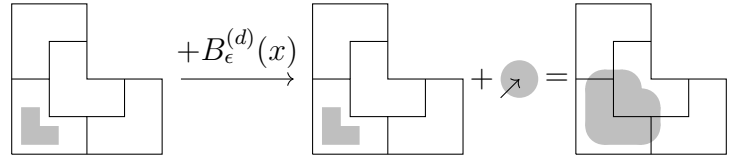%

\begin{figure}[t]
\centering
\begin{tikzpicture}
\draw[lightgray,line width=0.1cm](0.5,0.5)--(1,0);
\draw[->](0,1.0)--(.50,.50);
\draw(.50,.50)--(1.0,0);
\draw[->](0,-1.0)--(.50,-.50);
\draw(.50,-.50)--(1.0,0);
\draw[->](1.0,0)--(1.50,.50);
\draw(1.50,.50)--(2.0,1.0);
\draw[->](1.0,0)--(1.50,-.50);
\draw(1.50,-.50)--(2.0,-1.0);
\draw[->](2.25,0)--(3.75,0)node[midway,above]{\(+B_\epsilon^{(d)}(x)\)};
\draw[dotted,lightgray,line width=0.1cm](4.5,0.5)--(5,0);
\draw[lightgray,line width=0.1cm](5,0)--(6,1);
\draw[lightgray,line width=0.1cm](5,0)--(6,-1);
\draw[->](4.000,1.0)--(4.500,.50);
\draw(4.500,.50)--(5.000,0);
\draw[->](4.000,-1.0)--(4.500,-.50);
\draw(4.500,-.50)--(5.000,0);
\draw[->](5.000,0)--(5.500,.50);
\draw(5.500,.50)--(6.000,1.0);
\draw[->](5.000,0)--(5.500,-.50);
\draw(5.500,-.50)--(6.000,-1.0);
\end{tikzpicture}
\caption{A partial translation under the Robinson map. We have a cylinder set that does not intersect a branching point (left), translated by a sufficiently small ball so that the resulting set passes through a branching point (right). The dotted segment indicates that while origin placement is no longer allowed there, we still retain knowledge of its ``pattern'', and thus should include it as part of the set of partial translations. This set does not witness all of the branches, hence our usage of the term ``partial''.}
\label{figure:partial-translation}
\end{figure}
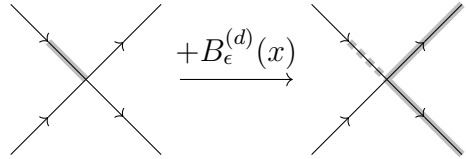%

\indent In this section, we define two main types of groupoids, one associated to the action, the other to the substitution rule, then discuss their associated \(C^\ast\)-algebras. While the reader with sufficient background may choose to skim this section, we will not follow the more standard perspective of defining the groupoids from punctured tiling spaces using the so-called doubly-pointed patterns introduced in \cite{kellendonk95}. Rather, we apply the Robinson map on the entire tiling space and obtain the analogue of the unstable and \(AF\)-groupoids, then restrict to the punctures to obtain a groupoid equivalence between the unstable groupoid on the punctured tiling space and the groupoid arising from a quotient of the \(AF\)-equivalence relation present on the associated Bratteli diagram constructed in \cite{bellissardjuliensavinien10} and \cite{juliensavinien10}. We hope that this approach provides more intuition to the relationship between the unstable and the \(AF\)-groupoids.

\subsection{Continuous groupoids}
\indent We start with a proposition on sets of the form \(C(P,U)+B_r(x)\) that we call \emph{partial translations} (\cref{figure:translation-cylinder-set}). \cref{figure:partial-translation} shows why sets of this form are given such a name. These are not the same as partial translations as defined in \cite{kellendonkputnam00}.
\begin{proposition}
\label{proposition:partial-translation}
The set of partial translations and the set of cylinder sets of the form \(C(P,U+B_r(x))\) where \(U+B_r(x)\subseteq\support P\) are mutual refinements.
\begin{proof}
Certainly cylinder sets of the latter form are partial translations. For the direction that the latter set is finer than the former, due to border-forcing and extension, the only case we have to check is that addition by sufficiently large \(x\) still returns unions of sets of the desired form.\\
\indent Let \(n\) be sufficiently large so that the border that is forced contains balls of radius \(R=\diameter U+x+r+\epsilon\) centered around the boundaries of each of the supertiles. Consider the set of all inclusions of \(P\) into patches of the form of unions of level-\(2n\) supertiles, so that the border that is forced is twice as large. For each such inclusion into \(\bigcup_i\varsigma^{2n}(t_i)\), for addition to make sense, let us consider the patch where we include the \(R\)-bordering tiles into the union. Then addition by \(B_r(x)\) never moves the resulting set outside of half of the forced border. However, two patches differing by the forced border return the same cylinder sets by extension, provided we pick the sets of origins to be outside of the forced border, thus while the patch of union of supertiles has its border grown by \(R\) to \(2R\), it is the same cylinder set as if its border stayed at \(R\). Repeating for each such inclusion gives that it can be written as a union of sets of the desired form.
\end{proof}
\end{proposition}
\noindent For notational purposes, by \(C(P,U)+B_r(x)\) we will assume that we work with sets of the latter form.\\
\indent We recall a few basic definitions on groupoids.
\begin{definition}
A \emph{groupoid} \(G\) is a set with a subset \(G^2\subseteq G\times G\) of \emph{composable elements of G}, a multiplication \(G^2\rightarrow G\) defined via \((g_1,g_2)\mapsto g_2g_1\), and an inversion \(G\rightarrow G\) defined via \(g\mapsto g^{-1}\) so that
\begin{itemize}
\item
For all \(g\in G\), \((g^{-1})^{-1}=g\),
\item
For all \((g_1,g_2),(g_2,g_3)\in G^2\), \((g_2g_1,g_3),(g_1,g_3g_2)\in G^2\), and \((g_3g_2)g_1=g_3(g_2g_1)\), and
\item
For all \(g\in G\), \((g,g^{-1})\in G^2\), and for all \((g_1,g_2)\in G^2\), \((g_1,g_2)g_1^{-1}=g_2\) and \(g_2^{-1}(g_1,g_2)=g_1\).
\end{itemize}%

\noindent \(G^0=\{g^{-1}g:g\in G\}\) is called the \emph{unit space}, and its elements are \emph{units}. The \emph{source} is the map \(s(g)=g^{-1}g\), and the \emph{range} is the map \(r(g)=gg^{-1}\).\\
\indent \(G\) is a \emph{topological groupoid} if it is a topological space with the source, the range, and the inversion maps continuous, and the multiplication map continuous with respect to the subspace topology of the product topology on \(G\times G\).\\
\indent \(G\) is \emph{\'etale} if it is a topological groupoid and the range map is a \emph{local homeomorphism}, or every \(g\in G\) has an open neighborhood \(U\) so that \(r(U)\) is open, and \(r:U\rightarrow r(U)\) is a homeomorphism.
\end{definition}
\noindent Note that if a groupoid is defined by an action on a topological space, then the unit space can be identified to the topological space itself.\\
\indent There is an important notion of equivalence in the sense of \cite[Example 2.7]{muhlyrenaultwilliams87} that relates one groupoid to another.
\begin{definition}
Let \(G\) be a locally compact Hausdorff groupoid. A closed subset \(N\subseteq G^0\) is an \emph{abstract transversal} if for all \(g_1\in G\), there exists \((g_1,g_2)\in G^2\) so that \(r(g_2g_1)\in N\). \(H=s^{-1}(N)\cap r^{-1}(N)\) is \emph{the restriction of \(G\) to \(N\)}, and \(G\) and \(H\) are \emph{equivalent (via \(s^{-1}(N)\) or just \(N\))}.
\end{definition}
\noindent It turns out that \(K\)-theory is invariant under equivalence of groupoids, so we will consider groupoids up to restriction to abstract transversals.\\
\indent We now define the groupoid on \(\Omega_T\) that we wish to study.
\begin{definition}
The \emph{unstable groupoid}, \(G_u\), on \(\Omega_T\) is the topological groupoid consisting of pairs of tilings that are translations of each other, i.e.
\[
G_u=\{(T',T'')\in\Omega_T^2:\exists x\in\mathbb{R}^d\textnormal{ such that }T''=T'-x\}.
\]
\noindent Its topology has the basis consisting of sets of the form
\[
(C(P,U)+B_{r_1}(x_1),C(P,U)+B_{r_2}(x_2))
\]
\noindent where \(r_1,r_2>0\) and \(x_1,x_2\in\mathbb{R}^d\).
\end{definition}
\indent We can replace this basis with the collection of partial homeomorphisms.
\begin{definition}
Given a homeomorphism of cylinder sets \(\phi:C(P,U_1)\rightarrow C(P,U_2)\), the \emph{source} is the map \(s(\phi)=C(P,U_1)\) (also denoted \(\phi^{-1}\phi\)), and the \emph{range} is the map \(r(\phi)=C(P,U_2)\) (also denoted \(\phi\phi^{-1}\)).\\
\indent A \emph{partial homeomorphism} on \(\Omega_T\) is an equivalence class of homeomorphisms of cylinder sets, where we quotient by the source and the range.
\end{definition}
\begin{figure}[t]
\centering
\begin{tikzpicture}
\fill[lightgray](0.5,0.5)circle[radius=0.25];
\draw(0,0)--(2,0)--(2,1)--(1,1)--(1,2)--(0,2)--(0,0);
\draw(1,0)--(1,0.5)--(0.5,0.5)--(0.5,1)--(0,1);
\draw(0.5,1)--(0.5,1.5)--(1,1.5);
\draw(1,0.5)--(1.5,0.5)--(1.5,1);
\draw[->](2.25,1)--(2.75,1)node[midway,above]{\(+B_1\)};
\fill[lightgray](3.5,0.5)circle[radius=0.25];
\draw(3,0)--(5,0)--(5,1)--(4,1)--(4,2)--(3,2)--(3,0);
\draw(4,0)--(4,.5)--(3.5,.5)--(3.5,1)--(3,1);
\draw(3.5,1)--(3.5,1.5)--(4,1.5);
\draw(4,.5)--(4.5,.5)--(4.5,1);
\draw(5.25,1)node{\(+\)};
\fill[lightgray](5.5,0.75)--(6,0.75)--(6,1)--(5.75,1)--(5.75,1.25)--(5.5,1.25)--(5.5,0.75);
\fill(5.625,0.875)circle[radius=0.05];
\draw(2.5,-1.25)node{\(=\)};
\fill[lightgray](3.125,-1.875)--(3.375,-2.125)--(3.875,-2.125)--(4.125,-1.875)--(4.125,-1.625)--(3.625,-1.125)--(3.375,-1.125)--(3.125,-1.375)--(3.125,-1.875);
\fill[lightgray](3.375,-1.875)circle[radius=0.25];
\fill[lightgray](3.875,-1.875)circle[radius=0.25];
\fill[lightgray](3.875,-1.625)circle[radius=0.25];
\fill[lightgray](3.625,-1.375)circle[radius=0.25];
\fill[lightgray](3.375,-1.375)circle[radius=0.25];
\draw(3,-2.25)--(5,-2.25)--(5,-1.25)--(4,-1.25)--(4,-.25)--(3,-.25)--(3,-2.25);
\draw(4,-2.25)--(4,-1.75)--(3.5,-1.75)--(3.5,-1.25)--(3,-1.25);
\draw(3.5,-1.25)--(3.5,-.75)--(4,-.75);
\draw(4,-1.75)--(4.5,-1.75)--(4.5,-1.25);
\draw(2.5,-3.5)node{\(=\)};
\fill[lightgray](3.375,-4.125)--(3.875,-4.125)--(3.875,-3.875)--(3.625,-3.875)--(3.625,-3.625)--(3.375,-3.625)--(3.375,-4.125);
\draw(3,-4.5)--(5,-4.5)--(5,-3.5)--(4,-3.5)--(4,-2.5)--(3,-2.5)--(3,-4.5);
\draw(4,-4.5)--(4,-4.0)--(3.5,-4.0)--(3.5,-3.5)--(3,-3.5);
\draw(3.5,-3.5)--(3.5,-3.0)--(4,-3.0);
\draw(4,-4.0)--(4.5,-4.0)--(4.5,-3.5);
\draw(5.25,-3.5)node{\(+\)};
\fill[lightgray](5.75,-3.5)circle[radius=0.25];
\fill(5.75,-3.5)circle[radius=0.05];
\draw[->](2.25,-5.75)--(2.75,-5.75)node[midway,above]{\(+B_2\)};
\fill[lightgray](3.375,-6.375)--(3.875,-6.375)--(3.875,-6.125)--(3.625,-6.125)--(3.625,-5.875)--(3.375,-5.875)--(3.375,-6.375);
\draw(3,-6.75)--(5,-6.75)--(5,-5.75)--(4,-5.75)--(4,-4.75)--(3,-4.75)--(3,-6.75);
\draw(4,-6.75)--(4,-6.25)--(3.5,-6.25)--(3.5,-5.75)--(3,-5.75);
\draw(3.5,-5.75)--(3.5,-5.25)--(4,-5.25);
\draw(4,-6.25)--(4.5,-6.25)--(4.5,-5.75);
\end{tikzpicture}
\caption{A construction of a basic open set in the unstable groupoid (modulo moving down sufficiently many levels) from a partial homeomorphism between the first and the last cylinder sets.}
\label{figure:partial-homeomorphism-to-translation}
\end{figure}
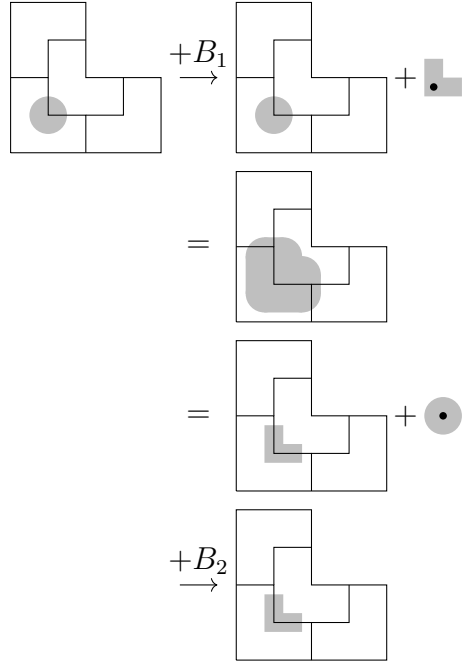%

\begin{proposition}
\label{proposition:unstable-groupoid-topology-partial-homeomorphism}
The basis of the unstable groupoid and the set of partial homeomorphisms (up to sources and ranges) are mutual refinements.
\begin{proof}
\(U+B_{r_1}(x_1)\) and \(U+B_{r_2}(x_2)\) are both homeomorphic to open balls, thus there exists a homeomorphism taking one to the other, giving us a partial homeomorphism between the cylinder sets.\\
\indent Conversely, given a partial homeomorphism \(\phi:C(P,U_1)\rightarrow C(P,U_2)\), let \(\{Q_i\}_i\) be an extension of \(P\) so that for each \(i\), \(\support Q_i\) contains balls of radius \(\diameter P\) centered around the boundary of \(\support P\). Pick any \(x_1\in U_1\) and \(x_2\in U_2\), then for each \(i\), it is not too hard to see (\cref{figure:partial-homeomorphism-to-translation}) that
\[
C(Q_i,U_1)+U_2-(x_1+x_2)/2=C(Q_i,U_2)+U_1-(x_1+x_2)/2
\]
\noindent and \(U_1+U_2-(x_1+x_2)/2\subseteq\support Q_i\), forming a basic open set in the unstable groupoid.
\end{proof}
\end{proposition}
\indent Using the topological conjugacy that is the Robinson map, we bring this over to \(\varprojlim_n(AP_n,\sigma)\) and obtain the unstable groupoid, also denoted \(G_u\), whose topology is generated by partial homeomorphisms. However, due to the possible existence of branches on each of the \(AP\)-complexes, our definitions look slightly different from what is expected if we applied the Robinson map directly. Applying the description given towards the end of the last section gives that the two are equivalent.
\begin{definition}
The \emph{unstable groupoid} on \(\varprojlim_n(AP_n,\sigma)\) is the set pairs of points \((p_0,\ldots),(q_0,\ldots)\) such that there exists a set \(\sigma^n(U_n)\times\cdots\times\sigma(U_n)\times U_n\times U_{n+1}\times\ldots\) such that for all \(m\geq n\), \(U_m\subseteq AP_m\) is \hyperlink{acu}{acu}, compatible with \(\sigma\), and \(p_m,q_m\in U_m\), i.e.
\[
G_u=\left\{((p_0,\ldots),(q_0,\ldots))\in\varprojlim_n(AP_n,\sigma)^2:\begin{array}{c}\exists U_n\times U_{n+1}\times\ldots\textnormal{ such that}\\\forall m\geq n\textnormal{, }U_m\subseteq AP_m\textnormal{ is acu,}\\\sigma(U_{m+1})=U_m\textnormal{, and }p_m,q_m\in U_m\end{array}\right\}.
\]
\indent A \emph{partial homeomorphism} on \(\varprojlim_n(AP_n,\sigma)\) is an equivalence class of homeomorphisms \(\phi:C(U_{1,n})\rightarrow C(U_{2,n})\) together with some \hyperlink{acu}{acu} \(U_n\subseteq AP_n\) such that \(U_{1,n},U_{2,n}\subseteq U_n\) that is its \emph{support}, denoted \(\support\phi\), modulo the source, the range, and the support.
\end{definition}
\noindent When we refer to sets of the form \(\sigma^n(U_n)\times\cdots\times\sigma(U_n)\times U_n\times U_{n+1}\times\ldots\), we will always assume \hyperlink{acu}{acu} starting at the level where \(\sigma\) is not used. An analogous proposition for \(\varprojlim_n(AP_n,\sigma)\) under the Robinson map shows that the set of partial homeomorphism forms a basis of this unstable groupoid.\\
\indent A concrete description of a partial homeomorphism can be given as follows. Let \(\phi:U_{1,n}\rightarrow U_{2,n}\) denote the underlying partial homeomorphism on \(AP_n\). Then for any \((p_0,\ldots)\in\varprojlim_n(AP_n,\sigma)\), \(\phi(p_0,\ldots)=(q_0,\ldots)\) has its coordinates \(\phi(p_n)=q_n\), and for all \(m\geq 0\), \(\sigma^{n-m}(q_n)=q_m\). This is possibly not well-defined for \(m>n\), but our choice of support has a single preimage under \(\sigma^{n-m}\) that contains \(p_m\). We pick \(q_m\) to be the unique point in this preimage. That is to say, if \(\phi:U_{1,n}\rightarrow U_{2,n}\) is a homeomorphism, then for all \(m>n\), each of the preimages of the support under \(\sigma^{n-m}\) determines a homeomorphism between the pair of corresponding preimages of \(U_{1,n}\) and \(U_{2,n}\) under \(\sigma^{n-m}\).\\
\indent Before proceeding further, since we will utilize it later, this description lets us extend partial homeomorphisms.
\begin{proposition}[Extension]
Suppose that \(\phi:C(U_{1,n})\rightarrow C(U_{2,n})\) is a partial homeomorphism.
\begin{itemize}
\item
If \(m\geq n\) and \(\psi:C(V_{1,m})\rightarrow C(V_{2,m})\) is a partial homeomorphism so that \(\sigma^{n-m}(\support\phi)\subseteq\support\psi\) and \(V_{i,m}=\sigma^{n-m}(U_{i,n})\) for \(i=1,2\), then \(\phi=\psi\).
\item
More generally, if \(m\geq n\) and \(\{\psi_i\}\) is a collection of partial homeomorphisms in \(AP_m\) so that \(\sigma^{n-m}(\support\phi)\subseteq\bigcup_i\support\psi_i\) with each \(\support\psi_i\) containing exactly one preimage in \(\sigma^{n-m}(\support\phi)\), and the sources and ranges are exactly \(\sigma^{n-m}(U_{i,n})\), respectively, then \(\phi=\bigcup_i\psi_i\).
\end{itemize}%

\begin{proof}
Extending the cylinder set \(C(\support\phi)\) to some level \(m\geq n\) then restricting the preimages of \(C(U_{1,n})\) and \(C(U_{2,n})\) to those that belong to \(\sigma^{n-m}(\support\phi)\) and are compatible under \(\sigma^{m-n}\) gives our desired collection of partial homeomorphisms.
\end{proof}
\end{proposition}
\begin{remark}
One may notice that this is a restricted version of a similar proposition that exists on the unstable groupoid on \(\Omega_T\). This proposition can be generalized to allow mixing of different levels of \(\support\psi_i\) in \(\varprojlim_n(AP_n,\sigma)\), but it is harder to state. We will also not require this generality.
\end{remark}
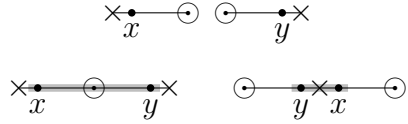
\begin{figure}[t]
\centering
\begin{tikzpicture}
\draw(-0.25,0)--(0.75,0);
\draw(-0.25,0)node{\(\times\)};
\draw(0.75,0)node{\(\odot\)};
\fill(0,0)circle[radius=0.05]node[below]{\(x\)};
\draw(1.25,0)--(2.25,0);
\draw(1.25,0)node{\(\odot\)};
\draw(2.25,0)node{\(\times\)};
\fill(2,0)circle[radius=0.05]node[below]{\(y\)};
\draw[lightgray,line width=0.1cm](-1.375,-1)--(0.375,-1);
\draw(-1.5,-1)--(0.5,-1);
\draw(-1.5,-1)node{\(\times\)};
\draw(0.5,-1)node{\(\times\)};
\draw(-0.5,-1)node{\(\odot\)};
\fill(-1.25,-1)circle[radius=0.05]node[below]{\(x\)};
\fill(0.25,-1)circle[radius=0.05]node[below]{\(y\)};
\draw[lightgray,line width=0.1cm](2.125,-1)--(2.875,-1);
\draw(1.5,-1)--(3.5,-1);
\draw(1.5,-1)node{\(\odot\)};
\draw(3.5,-1)node{\(\odot\)};
\draw(2.5,-1)node{\(\times\)};
\fill(2.25,-1)circle[radius=0.05]node[below]{\(y\)};
\fill(2.75,-1)circle[radius=0.05]node[below]{\(x\)};
\end{tikzpicture}
\caption{An \(AP_0\)-complex with two prototiles with the \(\times\)'s and \(\odot\)'s identified and two distinct points (top) whose translation is not well-defined until we choose a support that contains them (bottom).}
\label{figure:ap-complex-difference}
\end{figure}%

\begin{remark}
The support is vital in the definition of a partial homeomorphism. If one does not choose the support, then being a compact branched manifold without boundary, \(AP_n\) allows multiple ways to translate one point in \(U_{1,n}\) to another in \(U_{2,n}\) (\cref{figure:ap-complex-difference}). Choosing the support fixes the translation vector and lets the action to be well-defined. Alternatively, one can use the correspondence between \hyperlink{acu}{acu} sets in \(AP_n\) and patches in \(\Omega_T\).
\end{remark}
\indent Noting that if \(\support\phi\) resides in the interior of a \(d\)-cell, its preimages under \(\sigma\) necessarily also only reside in the interiors of \(d\)-cells, there is a dichotomy on the set of partial homeomorphisms, based on the support. To describe it, we first define a new topological groupoid.
\begin{definition}
The \emph{tail\textsubscript{\(n\)} groupoid}, denoted \(G_{\textnormal{tail},n}\), is the subgroupoid of the unstable groupoid with its topology generated by partial homeomorphisms whose supports reside within a single \(d\)-cell in \(AP_n\). This groupoid consists of pairs of points \((p_0,\ldots),(q_0,\ldots)\in\varprojlim_n(AP_n,\sigma)\) so that there exists a \(\sigma^n(U_n)\times\cdots\times\sigma(U_n)\times U_n\times U_{n+1}\times\cdots\) containing both points and for all \(m\geq n\), \(U_m\subseteq AP_m\) in the interior of a single \(d\)-cell, i.e.
\[
G_{\textnormal{tail},n}=\left\{((p_0,\ldots),(q_0,\ldots))\in G_u:\begin{array}{c}\exists U_n\times U_{n+1}\times\ldots\textnormal{ such that}\\\forall m\geq n\textnormal{, }p_m,q_m\in U_m\textnormal{ and}\\\exists d\textnormal{-cells }e_{d,m}\subseteq AP_m\textnormal{ such that }U_m\subseteq e_{d,m}\end{array}\right\}.
\]
\noindent \(U_n\times U_{n+1}\times\ldots\) is the \emph{tail}.
\end{definition}
\indent There exists a bonding map \(\sigma^\top:G_{\textnormal{tail},n}\rightarrow G_{\textnormal{tail},n+1}\). Given any two points that belong to \(G_{\textnormal{tail},n}\) via some \(\sigma^n(U_n)\times\cdots\times\sigma(U_n)\times U_n\times U_{n+1}\times\cdots\), since \(\sigma(U_{n+1})=U_n\), one can instead write it as \(\sigma^{n+1}(U_{n+1})\times\cdots\times\sigma^2(U_{n+1})\times\sigma(U_{n+1})\times U_{n+1}\times\cdots\), so the two points belong to \(G_{\textnormal{tail},n+1}\) as well, and \(\sigma^\top\) is the inclusion map. While this notation for the map may be misleading, as it is not the dual to the substitution map, on inducing it on the level of \(K\)-theory it becomes the dual.
\begin{definition}
We can form the \emph{tail groupoid}, denoted \(G_\textnormal{tail}\), by taking the inductive limit
\[
G_\textnormal{tail}=\varinjlim_n(G_{\textnormal{tail},n},\sigma^\top)=\bigcup_nG_{\textnormal{tail},n}.
\]
\indent The \emph{residual groupoid}, denoted \(G_\textnormal{res}\), is generated, as a groupoid, by \(G_u\backslash G_\textnormal{tail}\). More concretely, this groupoid is generated from pairs of points such that for all \(\sigma^n(U_n)\times\cdots\times\sigma(U_n)\times U_n\times U_{n+1}\times\cdots\) containing them, each set in each coordinate intersects the interiors of two \(d\)-cells, i.e.
\[
G_\textnormal{res}=\left\langle
\left\{((p_0,\ldots),(q_0,\ldots))\in G_u:\begin{array}{c}\forall U_n\times U_{n+1}\times\ldots\textnormal{ such that }\\\forall m\geq n\textnormal{, }p_m,q_m\in U_m\textnormal{,}\\U_m\cap AP_m^{(d-1)}\neq\emptyset
\end{array}\right\}
\right\rangle.
\]
\end{definition}
\begin{warning}
The tail\textsubscript{\(n\)} and tail groupoids on \(\varprojlim_n(AP_n,\sigma)\) are very much not finite- or approximately finite-dimensional! We use this term since on restricting to punctures, they become so, and the terminologies tail and \(AF\) coincide there.
\end{warning}
\indent This dichotomy allows us to put a filtration by dimension on the unstable groupoid. Consider the unstable groupoid on the \(d-1\)-skeleton, \(G_u^{(d-1)}\), that is the restriction of \(G_u\) to \(\varprojlim_n(AP_n^{(d-1)},\sigma)\). By unrestricting to a small neighborhood of the \(d-1\)-skeleton, we can view \(G_u^{(d-1)}\) as a subgroupoid of \(G_u\). Then the groupoid generated by \(G_\textnormal{tail}\) and \(G_u^{(d-1)}\) is exactly \(G_u\), and we can replace \(G_\textnormal{res}\) by \(G_u^{(d-1)}\) since they present the same topological information in \(G_u\).\\
\indent Proceeding in the same way \(G_\textnormal{tail}\) and \(G_\textnormal{res}\) are constructed, we obtain \(G_\textnormal{tail}^{(d-1)}\) and \(G_\textnormal{res}^{(d-1)}\). To obtain the filtration by dimension, we repeat this construction for each dimension \(0\leq k<d\) inductively in descending order. \(G_\textnormal{res}^{(0)}\) is empty, since \(G_u^{(0)}\) has its diagonal discrete points and there is no nontrivial action among them, so \(G_\textnormal{tail}^{(0)}\) sees the entire groupoid \(G_u^{(0)}\). We will use \(G_\textnormal{tail}^{(0)}\) instead of \(G_u^{(0)}\).
\begin{figure}[t]
\centering
\begin{tikzpicture}
\fill[lightgray](-0.125,0.375)--(0.375,-0.125)--(0.625,0.125)--(0.125,0.625)--(-0.125,0.375);
\fill[lightgray](0,0.5)circle[radius=0.1768];
\fill[lightgray](0.5,0)circle[radius=0.1768];
\draw(0,0)--(1,0)--(1,0.5)--(0.5,0.5)--(0.5,1)--(0,1)--(0,0);
\draw(0,1)--(-0.5,1)--(-0.5,0.5)--(-1,0.5)--(-1,0)--(0,0);
\draw(-1,0)--(-1,-0.5)--(-0.5,-0.5)--(-0.5,-1)--(0,-1)--(0,0);
\draw(0,-1)--(0.5,-1)--(0.5,-0.5)--(1,-0.5)--(1,0);
\fill(0.5,0)circle[radius=0.05];
\fill(0,0.5)circle[radius=0.05];
\draw[->](1.25,0)--(1.75,0)node[midway,above]{\(\textnormal{res}\)};
\draw[lightgray,line width=0.1cm](3,0.25)--(3,0.6768);
\draw[lightgray,line width=0.1cm](3.25,0)--(3.6768,0);
\draw(2,0)--(4,0);
\draw(3,-1)--(3,1);
\fill(3.5,0)circle[radius=0.05];
\fill(3,0.5)circle[radius=0.05];
\fill[lightgray](-0.1768,-3)--(0.1768,-3)--(0.1768,-2.5)--(-0.1768,-2.5)--(-0.1768,-3);
\fill[lightgray](0,-2.5)circle[radius=0.1768];
\fill[lightgray](0,-3.1768)--(0,-2.8232)--(0.5,-2.8232)--(0.5,-3.1768)--(0,-3.1768);
\fill[lightgray](0.5,-3)circle[radius=0.1768];
\fill[lightgray](0,-3)circle[radius=0.1768];
\draw(0,-3)--(1,-3)--(1,-2.5)--(.5,-2.5)--(.5,-2)--(0,-2)--(0,-3);
\draw(0,-2)--(-.5,-2)--(-.5,-2.5)--(-1,-2.5)--(-1,-3)--(0,-3);
\draw(-1,-3)--(-1,-3.5)--(-.5,-3.5)--(-.5,-4)--(0,-4)--(0,-3);
\draw(0,-4)--(.5,-4)--(.5,-3.5)--(1,-3.5)--(1,-3);
\fill(.5,-3)circle[radius=0.05];
\fill(0,-2.5)circle[radius=0.05];
\draw[->](1.25,-3)--(1.75,-3)node[midway,above]{\(\textnormal{res}\)};
\draw[lightgray,line width=0.1cm](3,-3.1768)--(3,-2.3232);
\draw[lightgray,line width=0.1cm](2.8232,-3)--(3.6768,-3);
\draw(2,-3)--(4,-3);
\draw(3,-4)--(3,-2);
\fill(3.5,-3)circle[radius=0.05];
\fill(3,-2.5)circle[radius=0.05];
\draw(4.25,-3)node{\(=\)};
\draw[lightgray,line width=0.1cm](5.5,-3.1768)--(5.5,-2.3232);
\draw(4.5,-3)--(6.5,-3);
\draw(5.5,-4)--(5.5,-2);
\fill(5.5,-3)circle[radius=0.05];
\fill(5.5,-2.5)circle[radius=0.05];
\draw(6.75,-3)node{\(\cdot\)};
\draw[lightgray,line width=0.1cm](7.8232,-3)--(8.6768,-3);
\draw(7.0,-3)--(9.0,-3);
\draw(8.0,-4)--(8.0,-2);
\fill(8.5,-3)circle[radius=0.05];
\fill(8.0,-3)circle[radius=0.05];
\end{tikzpicture}
\caption{If we do not impose contractibility on the supports of partial homeomorphisms, then a partial homeomorphism on \(G_u\) may not yield a partial homeomorphism on \(G_u^{(1)}\) from restriction to the \(1\)-skeleton (top). Assuming contractibility gives a sequence of partial homeomorphisms in \(G_u^{(1)}\) that composes to the desired partial homeomorphism (bottom).}
\label{figure:partial-homeomorphism-restriction}
\end{figure}
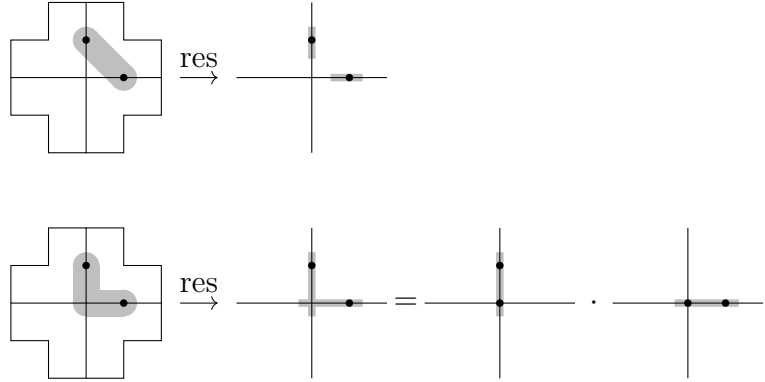%

\indent Notice that due to the contractibility condition, if we are given two points in the \(d-1\)-skeleton that belong in \(G_u\) via some \(\sigma^n(U_n)\times\cdots\times\sigma(U_n)\times U_n\times U_{n+1}\times\cdots\), restriction of this set to the skeleton gives a similar collection that is connected. Had we not imposed contractibility, the resulting restriction might be disconnected (\cref{figure:partial-homeomorphism-restriction}).\\
\indent Unfortunately, in the same vein as the warning in the previous section regarding lower-dimensional skeleta, this restriction presents an issue when we want to put a topology on \(G_u^{(k)}\) for \(0\leq k<d\). If we have a partial homeomorphism whose source and range cylinder sets restrict to interiors of different \(k\)-cells, thus are open balls in \(\mathbb{R}^k\), then its support crosses a \(k-1\)-cell. If the \(k\)-skeleton is branched at the \(k-1\)-cell, then the support, coming from restriction, necessarily sees all of the branches.\\
\indent As we want the partial homeomorphism to still arise from a \(k\)-dimensional translation action, one could fix the branching by deformation retracting the ``unused'' branches until we obtain an open ball in \(\mathbb{R}^k\). This is then an actual partial homeomorphism. However, this allows for sources or ranges, and therefore cylinder sets, to be centered about branch points. At the branching \(k-1\)-cell, intersecting two cylinder sets that share a single branch yields a set that is not open.\\
\indent Using a trick adapted from \cite{giordanomatuiputnamskau08} and \cite{giordanomatuiputnamskau09}, one easy way to resolve this issue is to, rather than restricting to skeleta, instead restrict to their \(\epsilon\)-neighborhoods that become \hyperlink{acu}{acu} upon border-forcing, where contractibility and unbranching are conditions per cell, not on the overall skeleton. Thus, while we morally want to work with groupoids restricted to skeleta, we will instead consider the subgroupoids of \(G_u\) with topology generated by partial homeomorphisms with supports in \(\epsilon\)-neighborhoods of the skeleta. Each \(k\)-skeleton forms an abstract transversal of \(G_u\) restricted to the \(\epsilon\)-neighborhood of the skeleton, and the new groupoids that we form by restricting to the \(\epsilon\)-neighborhood, still called \(G_u^{(k)}\), \(G_\textnormal{tail}^{(k)}\), and \(G_\textnormal{res}^{(k)}\), respectively, are therefore groupoid equivalent (in the sense of \cite{muhlyrenaultwilliams87}) to the ones defined on the \(k\)-skeleton.\\
\indent We now state these definitions precisely.
\begin{definition}
Let \(\epsilon>0\).\\
\indent Two open balls in the \(k\)-skeleton are \emph{partially homeomorphic} if there exists a partial homeomorphism in \(G_u\), supported on an \(\epsilon\)-neighborhood of the \(k\)-skeleton so that its source and range, restricted to the skeleton, are the two open balls.
\begin{itemize}
\item
The \emph{unstable groupoid of the \(k\)-skeleton}, \(G_u^{(k)}\), is the topological space with a basis of such partial homeomorphisms.
\item
The \emph{tail groupoid of the \(k\)-skeleton}, \(G_\textnormal{tail}^{(k)}\), is the subgroupoid of \(G_u^{(k)}\) generated by partial homeomorphisms (that live in \(G_u\)) whose support does not intersect the \(k-1\)-skeleton.
\item
The \emph{residual groupoid of the \(k\)-skeleton}, \(G_\textnormal{res}^{(k)}\), is generated by the difference of the two previous groupoids.
\end{itemize}
\end{definition}
\noindent One can also define the tail\textsubscript{\(n\)} groupoid of the \(k\)-skeleton, \(G_{\textnormal{tail},n}^{(k)}\), similarly, by restricting to supports that do not intersect \(AP_n^{(k)}\). The same process as before allows us to replace the residual groupoid \(G_\textnormal{res}^{(k)}\), for \(1\leq k\leq d\) by \(G_u^{(k-1)}\).\\
\indent It is crucial that we define each of the groupoids starting from \(G_u\), followed by restriction, since otherwise one can form compositions of partial homeomorphisms on \(G_u^{(k)}\) whose sources and ranges belong to different \(k\)-cells but avoid the \(k-1\)-skeleton. This is in fact the key reason our groupoids are stated as restrictions from \(G_u\). This is guaranteed to not occur since the supports of partial homeomorphisms in \(G_u\) are contractible.
\begin{remark}
The term ``partial homeomorphism'' is no longer accurate, since there exist partial homeomorphisms in \(G_u\) taking open balls intersecting \(k-1\)-cells to neighboring \(k\)-cells. If a \(k-1\)-cell is branched, the restriction of the source is branched, which is not homeomorphic to the restriction of the range, which is unbranched.\\
\indent If one forgoes restriction, then this is not an issue. However, the fact that partial homeomorphisms can ``bypass'' lower-dimensional skeleta makes the groupoid equivalence to the punctured tiling space harder to state.\\
\indent Lastly, while a sequence of partial homeomorphisms can still avoid the \(k-2\)-skeleton (or lower), our filtration on the unstable groupoid allows us to build skeleta dimension-by-dimension. More precisely, inductively, since the attaching maps are given in terms of consecutive-dimensional cells, this avoidance is a nonissue.
\end{remark}
\indent We now state an obvious but more useful version of extension for each \(G_u^{(k)}\) and \(G_\textnormal{tail}^{(k)}\).
\begin{proposition}[Extension]
Any partial homeomorphism in \(G_u^{(k)}\) can have its support extended to the entire interiors of (the \(\epsilon\)-neighborhood of) the \(k\)-cells it intersects. For partial homeomorphisms in \(G_\textnormal{tail}^{(k)}\), we may assume that the support is the interior of (the \(\epsilon\)-neighborhood of) a single \(k\)-cell.
\end{proposition}
\noindent From now on, we will assume that the supports of partial homeomorphisms in the two groupoids have these forms.

\subsection{Discrete groupoids}
\indent Let us now restrict to punctures to obtain the more familiar discrete groupoids. A considerable portion of the definitions and the results in this section are due to the seminal work of \cite{kellendonk95} and \cite{kellendonk97}, which provided the foundations of our understanding of the (discrete) unstable groupoid and the \(C^\ast\)-algebraic structure of aperiodic tilings.\\
\indent Since we will want to work with lower-dimensional skeleta, let us puncture each of the lower-dimensional boundaries of the tiles in the tiling as well, called \emph{\(k\)-punctures} for dimension \(k\), and pass them through the Robinson map (or \(\pi_n\) if we want to work with the \(AP_n\)-complex), so as to mimic our filtration by dimension. On \(AP_n\), there are as many punctures on each \(k\)-cell as there are \(k\)-cells on applying \(\sigma^n\), counting each occurrence of the same \(k\)-cells in \(AP_0\) separately. On the tiling space, this would be the number of \(k\)-dimensional boundaries of tiles that are in a \(k\)-dimensional boundary of a level-\(n\) supertile.\\
\indent The set of \(k\)-punctures forms an abstract transversal of the \(k\)-skeleton, thus restriction gives a groupoid equivalence. We will, however, give a slightly different presentation of the partial homeomorphisms that form a basis of the topology. To properly restrict, we will slightly amend the partial homeomorphisms that generate the topology on \(G_u\) to be those whose sources and ranges are sufficiently small, so that they contain one puncture each at the most. We also pick \(\epsilon\) (in the definitions of lower-dimensional groupoids) sufficiently small so that they only contain punctures from their respective skeleton. Then, on restriction, cylinder sets that contain punctures have their support turn into a collection of punctures and their open ball of origins turn into punctures themselves.
\begin{definition}
The \emph{discrete unstable groupoid of the \(k\)-skeleton}, \(\dot{G}_u^{(k)}\), is the restriction of \(G_u^{(k)}\) to the \(k\)-punctures. A \emph{partial homeomorphism} on the \(k\)-skeleton is a triple \([P,t',t]\), where \(P\) is a set of punctures inside an \hyperlink{acu}{acu} set, or a \emph{punctured acu set}, in some \(AP_n^{(k)}\), and \(t\) and \(t'\) are elements of \(P\). Composition of partial homeomorphisms is given by
\[
[P,t'',t']\cdot[P,t',t]=[P,t'',t].
\]
\noindent Let \(s([P,t',t])=t\) be the \emph{source}, and \(r([P,t',t])=t'\) be the \emph{range}. We omit the superscript if \(k=d\).
\end{definition}
\noindent It is well-known that \(\dot{G}_u\) is \'etale, although the unstable groupoids of the skeleta may not be.
\begin{remark}
A partial homeomorphism defined this way first appeared in \cite{kellendonk95}, and is also called a \emph{doubly-pointed pattern} in \cite{kellendonk97}. In order to maintain the set of doubly-pointed patterns as a left action on the set of punctures, we have elected to use the second puncture as the source and the first as the range, rather than the reverse as in \cite{kellendonkputnam00}.\\
\indent It is important to bear in mind that a punctured acu set comes from an underlying \hyperlink{acu}{acu} set, since there are numerous sets that restrict down to the same set of punctures, some of which may not be allowed, contractible, or unbranched. To prevent ambiguity, we will specify the support using the underlying acu set.
\end{remark}
\indent We can similarly define the restricted tail\textsubscript{\(n\)}, tail, and residual groupoids, but we switch to different terminology for the first two that is more suggestive of their defining features.
\begin{definition}
\leavevmode
\begin{itemize}
\item
The \emph{finite-dimensional groupoid at level-\(n\) of the \(k\)-skeleton}, or \(AF_n\)-groupoid, denoted \(\dot{G}_{AF,n}^{(k)}\), is the restriction of \(G_{\textnormal{tail},n}^{(k)}\) to the \(k\)-punctures.
\item
The \emph{approximately finite-dimensional groupoid of the \(k\)-skeleton}, or \(AF\)-groupoid, denoted \(\dot{G}_{AF}^{(k)}\), is the restriction of \(G_\textnormal{tail}^{(k)}\) to the \(k\)-punctures.
\item
The \emph{discrete residual groupoid of the \(k\)-skeleton}, denoted \(\dot{G}_\textnormal{res}^{(k)}\), is the restriction of \(G_\textnormal{res}^{(k)}\) to the \(k\)-punctures.
\end{itemize}
\end{definition}
\noindent Each is groupoid equivalent (in the sense of \cite{muhlyrenaultwilliams87}) to their unrestricted analogues. Whenever it is unclear, to distinguish from the discrete groupoid, we will refer to the unrestricted groupoids as \emph{continuous}.\\
\indent We now look at their associated (reduced) groupoid \(C^\ast\)-algebras. Recall that given a \'etale groupoid \(G\), we can form the continuous compactly-supported functions on \(G\) with values in \(\mathbb{C}\), denoted \(C_c(G)\), and give it the convolution product
\[
fg(a)=\sum_{bc=a}f(b)g(c)
\]
\noindent where the sum is finite due to local compactness, and involution
\[
f^\ast(a)=\overline{f(a^{-1})}.
\]
\noindent From \cite[Proposition 3.3.1]{sims17}, for each \(x\in G^0\), we have \emph{the regular representation associated to \(x\)}, \(\pi_x:C_c(G)\rightarrow\mathcal{B}(\ell^2(s^{-1}(x)))\), such that
\[
\pi_x(f)\delta_a=\sum_{s(b)=r(a)}f(b)\delta_{ba}.
\]
\noindent Forming the direct sum over all elements of \(G^0\) gives us the \emph{(reduced) groupoid \(C^\ast\)-algebra of \(G\)}, \(C_r^\ast(G)\), as the completion of
\[
\bigoplus_{x\in G^0}\pi_x(C_c(G))\subseteq\bigoplus_{x\in G^0}\mathcal{B}(\ell^2(s^{-1}(x)))
\]
\noindent under the supremum of the operator norm.\\
\indent Our goal will be to obtain \(K\)-theoretic information of \(C_r^\ast(\dot{G}_u)\) from that of \(C_r^\ast(\dot{G}_{AF}^{(k)})\). For most of the remainder of this subsection, we establish some elementary facts and consequences. To do so, let us denote the characteristic function on partial homeomorphisms by \(e[P,t',t]\). Note that these are elements of \(C_r^\ast(\dot{G}_u)\).\\
\indent There are three types of elements that are going to be crucial in our discussion,
\begin{itemize}
\item
\emph{Projections}, i.e. elements \(p\) so that \(p^2=p\),
\item
\emph{Partial isometries}, i.e. elements \(v\) so that the \emph{source}, \(v^\ast v\), and the \emph{range}, \(vv^\ast\), are projections, and
\item
\emph{Unitaries}, i.e. elements \(u\) so that \(u^\ast u=uu^\ast=1\).
\end{itemize}%

\noindent Let \(\mathcal{P}(\cdot)\) denote the set of projections of a \(C^\ast\)-algebra, and let \(\mathcal{U}(\cdot)\) denote the group of unitaries.\\
\indent We have several properties inherited from partial homeomorphisms that are due to \cite{kellendonk95} and we will not prove.
\begin{proposition}[Kellendonk]
\leavevmode
\begin{enumerate}
\item
\(e[P,t'',t']e[P,t',t]=e[P,t'',t]\).
\item
\(e[P,t,t]\) is a projection, denoted \(e[P,t]\), and \(e[P,t',t]\) is a partial isometry.
\item
\(e[P,t',t]^\ast=e[P,t,t']\).
\end{enumerate}%

\end{proposition}
\noindent Furthermore, the extension properties also hold here. We skip their statements since they are the same, but restricted to punctures.\\
\indent These allow us to show that \(C_r^\ast(\dot{G}_{AF}^{(k)})\) is a familiar \(C^\ast\)-algebra.
\begin{proposition}[Putnam, Kellendonk, Julien--Savinien]
\label{proposition:af-groupoid-algebra}
\(C_r^\ast(\dot{G}_{AF}^{(k)})\) is the closure of the direct limit of direct sums of finite-dimensional algebras, where the bonding maps are specified by the transpose of the induced substitution map on \(k\)-cells.
\begin{proof}
Given any partial homeomorphism in \(\dot{G}_{AF,n}^{(k)}\) where the support is some \(k\)-cell in \(AP_n^{(k)}\), suppose that we are given an order on the punctures, say \(\{t^i\}_{i=1}^m\). Consider the matrix algebra \(M_m(\mathbb{C})\) where the diagonal entries represent the ordered set of punctures. Then the elementary matrix \(e_{ji}\) is exactly the element \(e[P,t^j,t^i]\). Repeating this argument for each \(k\)-cell and noting that there are no partial homeomorphisms between different \(k\)-cells shows that \(C_r^\ast(\dot{G}_{AF,n}^{(k)})\) is the direct sum of matrix algebras consisting of the punctures inside each \(k\)-cell, with one summand for each \(k\)-cell. As we will see, the punctures cannot be arbitrarily ordered.\\
\indent There are as many punctures in a single \(k\)-cell in \(AP_{n+1}^{(k)}\) as there are punctures in each of the \(k\)-cells that it subdivides to in \(AP_n^{(k)}\). Given any partial homeomorphism in \(\dot{G}_{AF,n}^{(k)}\) with support \(P\) an entire \(k\)-cell, by extension, we can write it as a sum \(\sum_ie[P_i,t_i',t_i]\), with each \(P_i\) \(k\)-cells in \(AP_{n+1}^{(k)}\), possibly duplicated. Applying this to any pair of punctures in \(P\), we see that each \(P_i\) contains copies of the matrix algebra coming from \(P\), specified by \(\sigma^\top\). The order on the punctures is thus contingent on the order on the previous level, and we are only allowed a choice of order between sets of punctures of consecutive levels. By stationarity of the substitution rule, we maintain the same order between any two consecutive levels, and the order on the punctures on \(AP_1^{(k)}\) completely determines the rest.
\end{proof}
\end{proposition}
\noindent Closures of direct limits of direct sums of finite-dimensional algebras are called \emph{\(AF\)-algebras}.\\
\indent The above construction, due to \cite{putnam89} for \(k=d=1\), \cite{kellendonk95} for arbitrary \(k=d\), and \cite{juliensavinien10}\footnote{Their construction is slightly different due to wanting \emph{simplicity}.} in general, motivates a type of graphs associated to \(AF\)-algebras.
\begin{definition}
A \emph{Bratteli diagram} is a graph \((V,E)\) consisting of a vertex set \(V=\{V_n\}_{n=0}^\infty\), each of the \(V_n\) nonempty and consisting of finitely many vertices, and an edge set \(E=\{E_n\}_{n=0}^\infty\), each of the \(E_n\) nonempty and consisting of finitely many edges between \(V_n\) and \(V_{n+1}\), so that for all \(n\geq 0\), each vertex in each \(V_n\) is incident to an edge in \(E_n\), and for all \(n\geq 1\), each vertex in each \(V_n\) is incident to an edge in \(E_{n-1}\).\\
\indent The \emph{path space} \(X_E\) is the subset of \(\prod_{n=0}^\infty E_n\) consisting of paths in the graph. The \emph{\(AF_n\)-equivalence relation} on the path space is the set of pairs of paths in the path space whose edges coincide after level-\(n\), i.e.
\[
\{((p_0,\ldots),(q_0,\ldots))\in X_E^2:\forall m\geq n\textnormal{, }p_m=q_m\}.
\]
\noindent The \emph{\(AF\)-equivalence relation} is the union of all \(AF_n\)-equivalence relations.
\end{definition}
\indent For each \(AF\)-algebra, we can form an associated Bratteli diagram by taking each vertex at level-\(n\) to be each direct summand in the same level of the direct system, with each edge representing a single inclusion of a direct summand in one level into a direct summand in the next.\\
\indent Fixing a \(0\leq k\leq d\), let us form a Bratteli diagram for the direct system for \(C_r^\ast(\dot{G}_{AF}^{(k)})\), where at each level there are exactly as many vertices as there are \(k\)-cells, and the incidence matrix between consecutive levels is the transpose of the substitution matrix \(\sigma^\top\). It follows from the proof of \cref{proposition:af-groupoid-algebra} that each finite path in \(\prod_{n=0}^{N-1}E_n\) between \(V_0\) and \(V_N\) can be identified to a \(k\)-puncture at level-\(N\), and each vertex in the Bratteli diagram at level-\(N\) is a matrix algebra whose elementary matrices represent the characteristic functions on partial homeomorphisms among each pair of \(k\)-punctures. That is, the \(AF\)-equivalence relation on the Bratteli diagram formed from \(C_r^\ast(\dot{G}_{AF}^{(k)})\) is exactly the \(AF\)-groupoid of the \(k\)-skeleton.
\begin{remark}
For \(d=1\), there is an order one should choose on the \(d\)-punctures in \(AP_1\) that is the order of the punctures under \(\pi^{-1}\). This gives a well-defined (up to a sign) \emph{ordered} Bratteli diagram where the order of the punctures yields a \emph{Bratteli--Vershik map}.
\end{remark}
\indent As an immediate consequence, together with the description of the unstable and residual groupoids, we have a weaker form (groupoid equivalence rather than homeomorphism) of \cite[Theorem 3.22]{bellissardjuliensavinien10} and \cite[Theorem 4.9]{juliensavinien10}.
\begin{theorem}[Bellissard--Julien--Savinien, Julien--Savinien]
\leavevmode
\begin{itemize}
\item
The unstable groupoid on the punctured tiling space is groupoid equivalent to a quotient of the \(AF\)-equivalence relation on the top-dimensional Bratteli diagram.
\item
The unstable groupoid on the punctured tiling space is groupoid equivalent to an inductive quotient of the \(AF\)-equivalence relation on the \(k\)-dimensional Bratteli diagram by a quotient of the \(AF\)-equivalence relation on the \(k-1\)-dimensional Bratteli diagram.
\end{itemize}%

\begin{proof}
We start from the continuous groupoid on \(\varprojlim_n(AP_n,\sigma)\), \(G_u^{(k)}=G_\textnormal{tail}^{(k)}\cup G_\textnormal{res}^{(k)}\). Let us consider \(G_{u,\epsilon}^{(k)}\leq G_u^{(k)}\) whose diagonal contains exactly the \(k\)-cells neighboring the \(k-1\)-skeleton. Restricting this groupoid to the \(k\)- and \(k-1\)-punctures, we see that the set of \(k-1\)-punctures in the \(k-1\)-skeleton forms an abstract transversal for \(\dot{G}_{u,\epsilon}^{(k)}\), thus it is groupoid equivalent to \(\dot{G}_u^{(k-1)}\). Furthermore, any element in \(\dot{G}_\textnormal{res}^{(k)}\) can be written as a composition of elements from \(\dot{G}_{AF}^{(k)}\) and \(\dot{G}_{u,\epsilon}^{(k)}\). Since \(\dot{G}_{AF}^{(k)}\) is the \(AF\)-equivalence relation on the \(k\)-dimensional Bratteli diagram, \(\dot{G}_u^{(k)}\) is obtained by having certain \(AF\)-equivalence classes identified by \(\dot{G}_{u,\epsilon}^{(k)}\), and hence by \(\dot{G}_u^{(k-1)}\) by groupoid equivalence.\\
\indent To obtain the overall groupoid equivalence in the theorem, we restrict the unstable groupoid on \(\Omega_T\) to the \(d\)-punctures, and compose this with the Robinson map applied to the above explanation.
\end{proof}
\end{theorem}
\indent In the interest of giving an interpretation that will be useful in the next section, we will provide an explicit description purely in terms of doubly-pointed patterns. Given a partial homeomorphism \([P,t',t]\in\dot{G}_u^{(k)}\backslash\dot{G}_{AF,n}^{(k)}\) for \(k\geq 1\), we can decompose it into a sequence of partial homeomorphisms \(\{[P,t_i',t_i]\}_{i=1}^m\) so that \(t_i\) and \(t_i'\) are \(k\)-punctures neighboring a \(k-1\)-cell on projecting to \(AP_0\), \(t_{i+1}=t_i'\) for \(1\leq i<m\), and \(t_1=t\) and \(t_m'=t\), so that
\[
[P,t_m',t_m]\cdots[P,t_2',t_2]\cdot[P,t_1',t_1]=[P,t',t].
\]
Since the partial homeomorphism \([P,t',t]\) does not belong to the \(AF_n\)-groupoid of the \(k\)-skeleton, \(t\) and \(t'\) belong to different \(k\)-cells in \(AP_n\), and therefore there exists at least one index \(1\leq i\leq m\) so that \(t_i\) and \(t_i'\) belong to different \(k\)-cells. In particular, they neighbor a \(k-1\)-puncture in a \(k-1\)-cell in \(AP_n\). We can then think of such a \([P,t_i',t_i]\) as flowing from \(t_i\) to this \(k-1\)-puncture, then to \(t_i'\). Thus this sequence of partial homeomorphisms forms a path of neighboring \(k\)- and \(k-1\)-punctures, with segments that stay within a single \(k\)-cell, and to move between different \(k\)-cells, we move first to a \(k-1\)-puncture, then to the other \(k\)-cell. In other words, \(k-1\)-punctures in \(AP_n\) can be thought of as the nontrivial partial homeomorphisms in \(\dot{G}_u^{(k)}\) relative to \(\dot{G}_{AF,n}^{(k)}\).
\begin{remark}
An arbitrary sequence is unlikely to be efficient in terms of crossing as few \(k-1\)-boundaries as possible. However, one can ``subtract'' by loops until it is.
\end{remark}
\indent Unfortunately, this interpretation highlights the insufficiency of this theorem towards a na\"ive calculation of the \(K\)-theory of \(\dot{G}_u\) using the dimension filtration. That is, to compute the \(K\)-theory of \(\dot{G}_u\), one needs to inductively compute the \(K\)-theory of \(\dot{G}_u^{(k)}\) using \(\dot{G}_{AF}^{(k)}\) and \(\dot{G}_u^{(k-1)}\), where one thinks of \(\dot{G}_u^{(k-1)}\) as ``gluing'' elements of \(\dot{G}_{AF}^{(k)}\) together. If \(AP_n^{(k)}\) is branched at some \(k-1\)-cell with \(m\)-branches, then that corresponding puncture simultaneously represents, up to orientation, \(m-1\) partial homeomorphisms (all others can be generated by combining). Furthermore, this \(m\) is not consistent between different \(k-1\)-cells, which presents difficulties when there are cylinder sets in the \(k\)-skeleton that witness different \(k-1\)-cells when extending.\\
\indent One way to remedy this is to use a more general \(K\)-theoretic object called \emph{relative \(K\)-theory}. We will then see a more precise topological explanation for why the na\"ive computation does not work.
\begin{remark}
There is a way around where one employs a process similar to \emph{state-splitting} that \emph{does not preserve the \(K\)-theory of the \(AF\)-groupoid} (but does preserve the \(K\)-theory of the unstable groupoid), to the substitution rule to ensure the existence of a subgroupoid of \(\dot{G}_u^{(k)}\) closed under a \(k\)-parameter subgroup of \(\mathbb{R}^d\), e.g. \cite{juliensavinien16} where one takes the chair tiling and splits each supertile into three squares. It is unclear if this process always returns a stationary substitution rule.
\end{remark}

\section{\(K\)-theory}
\label{section:k-theory}
\indent In this section we recall ordinary \(K\)-theory, then briefly outline the key facts from relative \(K\)-theory and excision, leaving the details for the reader to fill in from \cite{haslehurst21}, \cite{putnam97}, and \cite{putnam21}. We will only sketch the relative \(K_0\)-group, and leave the relative \(K_1\)-group to be deduced from topology and exactness. We will also only state excision as it appears in \cite{putnam98} since we do not need its full power.

\subsection{Ordinary \(K\)-theory}
\indent Given a unital \(C^\ast\)-algebra \(A\), let \(M_n(A)\) denote the algebra of \(n\times n\) matrices with entries in \(A\). Let
\[
M_\infty(A)=\bigcup_{n=1}^\infty M_n(A)
\]
\noindent under the embedding \(a\mapsto\diagonal(a,0)\). Let the binary operation on \(\mathcal{P}(M_\infty(A))\) be direct sum, i.e.
\[
p\oplus q=\left(\begin{array}{cc}p&0\\0&q\end{array}\right).
\]
\noindent Two projections are \emph{Murray--von Neumann (MvN) equivalent} if there exists a partial isometry whose source and range are the two projections. The \emph{\(K_0\)-group of \(A\)}, \(K_0(A)\), is the Grothendieck group of the abelian semigroup formed from \(\mathcal{P}(M_\infty(A))\) with direct sum, under Murray--von Neumann equivalence in \(M_\infty(A)\).\\
\indent We now consider
\[
M_\infty(A)=\bigcup_{n=1}^\infty M_n(A)
\]
\noindent under the embedding \(a\mapsto\diagonal(a,1)\). Let the binary operation on \(\mathcal{U}(M_\infty(A))\) still be direct sum. Two unitaries \(u,v\in\mathcal{U}(M_\infty(A))\) are \emph{stably homotopic} if there exists a sufficiently large \(n\) so that \(\diagonal(u,1_k),\diagonal(v,1_\ell)\in\mathcal{U}(M_n(A))\) are homotopic within \(\mathcal{U}(M_n(A))\). The \emph{\(K_1\)-group of \(A\)}, \(K_1(A)\), is the group of equivalence classes of stably homotopic unitaries in \(\mathcal{U}(M_\infty(A))\) with direct sum.\\
\indent We have the following standard properties. See, for example, \cite{rordamlarsenlaustsen00}.
\begin{proposition}
For \(i=0,1\), \(K_i\) is a continuous functor. That is, if \(\varinjlim_n(A_n,\phi_n)\) is the closure of a direct limit of \(\ast\)-homomorphisms of \(C^\ast\)-algebras \(\phi_n:A_n\rightarrow A_{n+1}\), then \(K_i(\varinjlim_n(A_n,\phi_n))\cong\varinjlim_n(K_i(A_n),{\phi_n}_\ast)\).
\end{proposition}
\begin{proposition}
\(K_0(\mathbb{C})=K_0(M_n(\mathbb{C}))=\mathbb{Z}\), and \(K_1(\mathbb{C})=K_1(M_n(\mathbb{C}))=0\).
\end{proposition}
\begin{corollary}
If an \(AF\)-algebra \(A\) arises from a Bratteli diagram \((V,E)\) with incidence matrices \(\{h_{E_n}\}\), then \(K_0(A)=\varinjlim_n(\mathbb{Z}V_n,h_{E_n})\) and \(K_1(A)=0\).
\end{corollary}
\noindent The following is \cite[Theorem 2.8]{muhlyrenaultwilliams87} adapted for abstract transversals.
\begin{theorem}[Muhly--Renault--Williams]
If \(H\) is a subgroupoid of an \'etale groupoid \(G\), and the two are equivalent via an abstract transversal \(N\), then for \(i=0,1\), \(K_i(C_r^\ast(G))\cong K_i(C_r^\ast(H))\).
\end{theorem}

\subsection{Relative \(K\)-theory}
\indent We now move on to relative \(K\)-theory. Given a map of unital \(C^\ast\)-algebras \(\phi:A\rightarrow B\), consider triples of the form \((p,q,v)\) where \(p,q\in M_\infty(A)\) are projections and \(v\in M_\infty(B)\) is a partial isometry so that \(s(v)=\phi(p)\) and \(r(v)=\phi(q)\). That is, \(p\) and \(q\) are Murray--von Neumann equivalent under \(\phi\). Two triples \((p_1,q_1,v_1)\) and \((p_2,q_2,v_2)\) are \emph{isomorphic} if there exist partial isometries \(w_1,w_2\in M_\infty(A)\) so that \(s(w_1)=p_1\), \(r(w_1)=p_2\), \(s(w_2)=q_1\), \(r(w_2)=q_2\), and \(\phi(w_2)v_1=v_2\phi(w_1)\). A triple \((p,q,v)\) is \emph{elementary} if \(p=q\) and there is a path of partial isometries \(\{v_t\}_{0\leq t\leq 1}\subseteq B\) so that \(v_0=\phi(p)\), \(v_1=v\), and \(s(v_t)=r(v_t)=\phi(p)\). Letting the binary operation be component-wise direct sum, we define the \emph{relative \(K_0\)-group of \(\phi\)}, denoted \(K_0(A;B)\), to be the equivalence classes of triples under isomorphism and stabilization with respect to taking direct sums with elementary triples. Note that the map \(\phi\) affects the group, but we leave it implicit in the notation.\\
\indent We have several key properties that we will use as black boxes.
\begin{proposition}
\leavevmode
\begin{enumerate}
\item
\([p_1,q_1,v_1]+[p_2,q_2,v_2]=[p_1,q_2,v_2v_1]\) if \(q_1=p_2\).
\item
\(-[p,q,v]=[q,p,v^\ast]\).
\item
Relative \(K\)-theory is functorial. That is, if the diagram
\[
\begin{tikzcd}
A_1\arrow[r,"\phi_1"]\arrow[d,"\psi_1"]&B_1\arrow[d,"\psi_2"]\\
A_2\arrow[r,"\phi_2"]&B_2
\end{tikzcd}
\]
\noindent of \(\ast\)-homomorphisms between \(C^\ast\)-algebras commutes, then there exists a \(\psi_\ast:K_i(A_1;B_1)\rightarrow K_i(A_2;B_2)\) for \(i=0,1\) such that \(\psi_\ast([p,q,v])=[\psi_1(p),\psi_1(q),\psi_2(v)]\).
\end{enumerate}%

\end{proposition}
\noindent The main tool is the \emph{six-term sequence in relative \(K\)-theory} (\cite[Theorem 2.1]{haslehurst21}).
\begin{theorem}[Haslehurst]
There is an exact sequence
\[
\begin{tikzcd}
K_0(A;B)\arrow[r,"\evaluation"]&K_0(A)\arrow[r,"\phi_\ast"]&K_0(B)\arrow[d]\\
K_1(B)\arrow[u]&K_1(A)\arrow[l]&K_1(A;B)\arrow[l]
\end{tikzcd}
\]
where \(\evaluation([p,q,v])=[p]-[q]\) is the \emph{evaluation map}.
\end{theorem}

\subsection{Excision}
\indent We now proceed with excision in the context of open subgroupoids in (reduced) groupoid \(C^\ast\)-algebras in \cite{putnam98}. Let \(G\) be a groupoid that is an equivalence relation on a topological space, and let \(L\subseteq G\) as a subspace (\emph{not as a subgroupoid!}) be such that
\begin{itemize}
\item
\(L\) is closed,
\item
\(r(L)\cap s(L)=\emptyset\),
\item
\(G'=G\backslash(L\cup L^{-1})\) is such that \(G'G'\leq G\), and
\item
\(LG',G'L\subseteq L\).
\end{itemize}%

\noindent Let the topology on \(s(L)=L^{-1}L\) be given by the convergence of sequences, \(\{x_n\}_n\) converges to \(x\) if and only if there exist sequences \(\{y_n\}_n\) and \(\{z_n\}_n\) in \(L\) converging in \(G\) so that \(x_n=y_n^{-1}z_n\), and similarly for \(r(L)=LL^{-1}\). Let \(H'=s(L)\cup r(L)\) and \(H=H'\cup L\cup L^{-1}\) be given the disjoint union topology. Then we have the following \emph{excision} property from \cite{putnam21}.
\begin{theorem}[Putnam]
For \(i=0,1\),
\[
K_i(C_r^\ast(G');C_r^\ast(G))\cong K_i(C_r^\ast(H');C_r^\ast(H)).
\]
\end{theorem}
\indent The following corollary from \cite{putnam98} will be very useful. We will provide a sketch of the proof as we will need the map explicitly.
\begin{theorem}[Putnam]
\label{theorem:excision-isomorphism}
For \(i=0,1\),
\[
K_i(C_r^\ast(G');C_r^\ast(G))\cong K_i(C_r^\ast(H)).
\]
\begin{proof}
\(H'\) has two components, and under the new topology, each forms an abstract transversal to \(H\). Thus \(K_i(C_r^\ast(H'))\cong K_i(C_r^\ast(H))\oplus K_i(C_r^\ast(H))\). Each of the rows of the six-term sequence in relative \(K\)-theory then reads
\[
\begin{tikzcd}
K_i(C_r^\ast(H');C_r^\ast(H))\arrow[r,hook,"\Delta^\top"]&K_i(C_r^\ast(H))\oplus K_i(C_r^\ast(H))\arrow[r,two heads,"+"]&K_i(C_r^\ast(H))
\end{tikzcd}
\]
\noindent where first map \(\Delta^\top=(1,-1)\) is the skew-diagonal map and the second is summation. Thus \(K_i(C_r^\ast(H');C_r^\ast(H))\cong K_i(C_r^\ast(H))\), and excision gives the statement.
\end{proof}
\end{theorem}
\noindent The way we will use this theorem is that there will be an obvious abstract transversal to \(H\) that comes from our choice of \(L\) (rather, \(L\cup L^{-1}\)) that we then restrict to, where the restriction is continuous with respect to the new topology on \(H\), resulting in a groupoid equivalence.

\section{\v Cech cohomology and \(K\)-theory}
\label{section:k-theory-cohomology}
\indent Since \(\dot{G}_{AF}\leq\dot{G}_u\) is an inclusion of an open subgroupoid, a continuous compactly-supported function on \(\dot{G}_{AF}\) is also one on \(\dot{G}_u\), giving us the inclusion \(C_c(\dot{G}_{AF})\hookrightarrow C_c(\dot{G}_u)\) and therefore the inclusion \(\iota:C_r^\ast(\dot{G}_{AF})\hookrightarrow C_r^\ast(\dot{G}_u)\). Applying the six-term sequence in relative \(K\)-theory to \(\iota\) gives us
\[
\begin{tikzcd}
K_0(C_r^\ast(\dot{G}_{AF});C_r^\ast(\dot{G}_u))\arrow[r,"\evaluation"]&K_0(C_r^\ast(\dot{G}_{AF}))\arrow[r,"\iota_\ast"]&K_0(C_r^\ast(\dot{G}_u))\arrow[d,two heads]\\
K_1(C_r^\ast(\dot{G}_u))\arrow[u,hook]&0\arrow[l]&K_1(C_r^\ast(\dot{G}_{AF});C_r^\ast(\dot{G}_u))\arrow[l].
\end{tikzcd}
\]
\noindent The goal of this section is to construct an isomorphic exact sequence in \v Cech cohomology for \(d=1,2\), then use it to explicitly illustrate the roles of each \(K_0(C_r^\ast(\dot{G}_{AF}^{(k)}))\). We endeavor to describe as many things in generality as possible before focusing on \(d=1,2\), and begin by establishing some notation and convention. We will occasionally make use of the boundary hyperplane condition, and will state explicitly where it is used. Although it factors in to some of the proofs, the only place it is required is the reconstruction of \(K_1(C_r^\ast(\dot{G}_u))\) from \(K_0(C_r^\ast(\dot{G}_{AF}^{(1)}))\).\\
\indent If \(P_1\) and \(P_2\) are punctured \hyperlink{acu}{acu} sets in \(AP_n\), denote \(P=P_1\cup P_2\) the union of the two punctured acu sets. \emph{This is not the union of cylinder sets!} There may be many such \(P\), and we will let context dictate the relevant one. Furthermore, noting that in \(AP_0\), there is exactly one \(k\)-puncture for each \(k\)-cell, allowing us to identify the two, we will make use of the fact that we have a substitution rule and write \(\varsigma^n(t)\), \(t\) a \(k\)-cell in \(AP_0\), to denote the punctured acu set in \(AP_n^{(k)}\) corresponding to subdividing \(t\) \(n\)-times. By extension, we may replace any punctured acu set with \(\bigcup_i\varsigma^n(t_i)\), for some appropriate \(\{t_i\}_i\).

\subsection{Thickening cochains}
\indent Cohomology, and therefore the cellular structure on the \(AP\)-complexes, will be integral in this section. Let us assume that our \(k\)-cells are formed from the intersection of multiple \(d\)-cells. There is a type of \(k\)-cell that behaves particularly nicely with respect to substitutions.
\begin{definition}
A \(k\)-cell satisfies the \emph{boundary hyperplane condition} if it belongs to a hyperplane. \(T\) satisfies the \emph{boundary hyperplane condition} if, for all \(0\leq k<d\), all \(k\)-cells satisfy the boundary hyperplane condition.
\end{definition}
\noindent Note that \(d\)-cells automatically satisfy the boundary hyperplane condition. With this definition, we will assign orientations to \(k\)-cells, for \(k<d\), one of two ways.
\begin{itemize}
\item
If a \(k\)-cell does not satisfy the boundary hyperplane condition, then its orientation is assigned arbitrarily.
\item
If a \(k\)-cell satisfies the boundary hyperplane condition, then its orientation is assigned so that if it and another \(k\)-cell have preimages under some \(\sigma^n\) belonging in the same \(k\)-cell, then they, and each of their preimages, carry the same orientation. For example, assigning the same orientation for parallel \(k\)-cells guarantees this. In other words, we want orientations to stay the same within ``primitive components''.
\end{itemize}%

\noindent For the sake of signs working out, we will pick the right-handed orientation for the \(d\)-cells.\\
\indent Given the cochain complexes for the \(AP\)-complexes, one easily verifies that, from the conditions we placed on the orientations of the cells, the substitution map induces a map between consecutive levels that can be identified with \(\sigma^\top\), letting us form a cochain complex of direct limits
\[
\begin{tikzcd}
0\arrow[r]&C_0^0\arrow[r,"\delta_0^0"]\arrow[d,"\sigma^\top"]&C_0^1\arrow[r,"\delta_0^1"]\arrow[d,"\sigma^\top"]&\cdots\arrow[r,"\delta_0^{d-1}"]&C_0^d\arrow[r]\arrow[d,"\sigma^\top"]&0\\
0\arrow[r]&C_1^0\arrow[r,"\delta_1^0"]\arrow[d,"\sigma^\top"]&C_1^1\arrow[r,"\delta_1^1"]\arrow[d,"\sigma^\top"]&\cdots\arrow[r,"\delta_1^{d-1}"]&C_1^d\arrow[r]\arrow[d,"\sigma^\top"]&0\\
&\vdots\arrow[d,"\sigma^\top"]&\vdots\arrow[d,"\sigma^\top"]&&\vdots\arrow[d,"\sigma^\top"]&\\
0\arrow[r]&C_n^0\arrow[r,"\delta_n^0"]\arrow[d,"\sigma^\top"]&C_n^1\arrow[r,"\delta_n^1"]\arrow[d,"\sigma^\top"]&\cdots\arrow[r,"\delta_n^{d-1}"]&C_n^d\arrow[r]\arrow[d,"\sigma^\top"]&0\\
&\vdots\arrow[d]&\vdots\arrow[d]&&\vdots\arrow[d]&\\
0\arrow[r]&\varinjlim_n(C_n^0,\sigma^\top)\arrow[r,"\delta^0"]&\varinjlim_n(C_n^1,\sigma^\top)\arrow[r,"\delta^1"]&\cdots\arrow[r,"\delta^{d-1}"]&\varinjlim_n(C_n^d,\sigma^\top)\arrow[r]&0.
\end{tikzcd}
\]
\noindent To simplify notation, let us call \(C^k=\varinjlim_n(C_n^k,\sigma^\top)\) for each \(0\leq k\leq d\).\\
\indent We have a key observation that is the motivation behind most things we do.
\begin{lemma}[The Cochain--Finite-dimensional Algebra Isomorphism (\(\CF\))]
\label{lemma:af-cochain}
\(K_0(C_r^\ast(\dot{G}_{AF}))\cong C^d\). For \(k<d\), \(\varinjlim_n(K_0(C_r^\ast(\dot{G}_{AF,n}^{(k)})),\sigma^\top)\cong C^k\), and if the boundary hyperplane condition is satisfied, then \(K_0(C_r^\ast(\dot{G}_{AF}^{(k)}))\cong C^k\).
\begin{proof}
We picked the bonding maps \(K_0(C_r^\ast(\dot{G}_{AF,n}^{(k)}))\rightarrow K_0(C_r^\ast(\dot{G}_{AF,n+1}^{(k)}))\) to be the one given by \(\sigma^\top:C_n^k\rightarrow C_{n+1}^k\), so we only need the identification \(K_0(C_r^\ast(\dot{G}_{AF,n}^{(k)}))\cong C_n^k\).\\
\indent Given a \(k\)-cochain corresponding to a \(k\)-cell \(AP_n\), we pick a projection in \(C_r^\ast(\dot{G}_{AF,n}^{(k)})\) whose support is the entire \(k\)-cell, which is a \(k\)-punctured \hyperlink{acu}{acu} set, together with a puncture in it. \(K_0\) asserts that all punctures in this set are the same, making this identification well-defined.\\
\indent If the boundary hyperplane condition is satisfied, then each of the entries in \(\sigma^\top:C_n^k\rightarrow C_{n+1}^k\) is nonnegative, and the map \(K_0(C_r^\ast(\dot{G}_{AF,n}^{(k)}))\rightarrow K_0(C_r^\ast(\dot{G}_{AF,n+1}^{(k)}))\) is induced by the algebra map \(C_r^\ast(\dot{G}_{AF,n}^{(k)})\rightarrow C_r^\ast(\dot{G}_{AF,n+1}^{(k)})\).\\
\indent Finally, we chose each of the \(d\)-cells to have the same orientations, therefore \(\sigma^\top\) always has nonnegative entries.
\end{proof}
\end{lemma}
\indent In other words, if we wish to show that \(K_i(C_r^\ast(\dot{G}_u))\) can be constructed out of \(K_0(C_r^\ast(\dot{G}_{AF}^{(k)}))\) for \(0\leq k\leq d\) using direct sums, quotients, and subgroups, it suffices to demonstrate this for \v Cech cohomology and the cochains groups, then pass back to \(K\)-theory using an isomorphism between \(K\)-theory and \v Cech cohomology. In order to make use of this observation, we need an identification between \(C^{d-1}\) and \(K_0(C_r^\ast(\dot{G}_{AF});C_r^\ast(\dot{G}_u))\) that arises from extension.\\
\indent Let us provide a simple picture on certain types of elements of \(K_0(C_r^\ast(\dot{G}_{AF});C_r^\ast(\dot{G}_u))\). Observe that a partial isometry \(e[P,t',t]\in C_r^\ast(\dot{G}_u)\) is a triple \((p,q,v)\) by taking \(p=e[P,t]=e[P,t,t']e[P,t',t]\), \(q=e[P,t']=e[P,t',t]e[P,t,t']\), and \(v\) to be the partial isometry itself. Partial isometries with support single punctured \(d\)-cells are elementary triples, thus are trivial in \(K_0(C_r^\ast(\dot{G}_{AF});C_r^\ast(\dot{G}_u))\). However, it is not evident that if \(e[P,t',t]\notin C_r^\ast(\dot{G}_{AF})\), we still have that \(e[P,t],e[P,t']\in C_r^\ast(\dot{G}_{AF})\), given that the support \(P\) may intersect multiple \(d\)-cells. This is resolved in the next few definitions and propositions.
\begin{definition}
Given a \(d-1\)-cell in the \(AP_n\)-complex, let \(\{\varsigma^n(t_i)\}_i\) be the set of punctured \(d\)-cells, counting multiplicity, incident to and with orientations agreeing with it, and let \(\{\varsigma^n(t_j)\}_j\) be the set of punctured \(d\)-cells, counting multiplicity and with orientations disagreeing. We call such a pair a \emph{full translation}. The first collection satisfies the \emph{right-hand rule}, and the second satisfies the \emph{left-hand rule}.
\end{definition}
\begin{figure}[t]
\centering
\begin{tikzpicture}
\draw[->](0.5176,1.9319)--(0.2588,0.9659);
\draw(0.2588,0.9659)--(0,0);
\path(0,0)--(0.5176,1.9319)node[midway,above,rotate=75]{\((a)b(b)\)};
\draw[->](0,0)--(0.9659,0.2588);
\draw(0.9659,0.2588)--(1.9319,0.5176);
\path(0,0)--(1.9319,0.5176)node[midway,below,rotate=15]{\((b)b(a)\)};
\draw[->](1.9319,0.5176)--(1.2247,1.2247);
\draw(1.2247,1.2247)--(0.5176,1.9319);
\path(0.5176,1.9319)--(1.9319,0.5176)node[midway,above,rotate=-45]{\((b)a(b)\)};
\draw[->](0,0)arc(45:225:0.3183);
\draw(-0.4502,-0.4502)arc(-135:45:0.3183);
\path(-0.9003,0)--(0,-0.9003)node[midway,below,rotate=-45]{\((b)b(b)\)};
\fill(0,0)circle[radius=0.05];
\end{tikzpicture}
\caption{The \(AP\)-complex of the collared Silver Mean substitution, \(a\mapsto b\), \(b\mapsto bab\). The basic open ball around the marked point that is a \(0\)-cell is formed from the three \(1\)-cells \((a)b(b)\), \((b)b(a)\), and \((b)b(b)\). The preimages in the tiling space under the Robinson map are \((a)b.b(a)\), \((a)b.b(b)\), \((b)b.b(a)\), so \((b)b(b)\) is counted twice. The \(0\)-cell partitions this set into \(\{(a)b(b),(b)b(b)\}\) (source) and \(\{(b)b(b),(b)b(a)\}\) (range), which is a full translation.}
\label{figure:silver-mean-loop}
\end{figure}
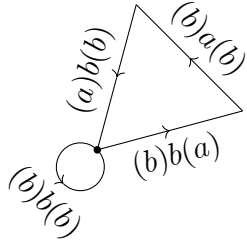%

\noindent A single \(d\)-cell may have multiple parts of its boundary attached to a single \(d-1\)-cell (\cref{figure:silver-mean-loop}). By \emph{counting multiplicity}, we mean that we count each occurrence separately. Taking union of a full translation yields an open ball in \(AP_n\). In other words, given the set of \(d\)-cells incident to a \(d-1\)-cell forming an open ball, counting multiplicity, we are forming a partition based on handedness.
\begin{remark}
Rather than counting multiplicities, one can instead apply a version of state-splitting. This amounts to decorating each \(d\)-cell in \(AP_n\) with its image \(d\)-cells in \(AP_{n-m}\), where \(m\) is sufficiently large so that none of the decorations is incident to the \(d-1\)-cell in more than one way. That is, for this particular \(m\), we pass to the induced substitution applied to sets of finite coordinates of length \(n-m\) in \(\varprojlim_n(AP_n,\sigma)\) that do not necessarily start at level-\(0\).
\end{remark}
\begin{proposition}[Extension]
Given a partial isometry \(e[P,t',t]\in C_r^\ast(\dot{G}_u)\), if \(t\) and \(t'\) belong to the same punctured \(d\)-cell for some \(AP_n\), then \(e[P,t',t]\in C_r^\ast(\dot{G}_{AF})\), up to shifting by a level given by border-forcing. More generally, if \(t\) and \(t'\) belong to different punctured \(d\)-cells, it can be written as a sum of partial isometries
\[
\sum_{i=1}^{m-1}e[\varsigma^n(t_i)\cup\varsigma^n(t_{i+1}),t_i'',t_i']
\]
where for \(1\leq i\leq m\), \(t_i'\in\varsigma^n(t_i)\) and \(t_i''\in\varsigma^n(t_{i+1})\), and for \(1\leq i<m\), each pair \(\varsigma^n(t_i)\) and \(\varsigma^n(t_{i+1})\) intersects at a \(d-1\)-cell and \(t_{i+1}'=t_i''\).
\begin{proof}
We will prove the second statement, since it implies the first.\\
\indent Let \(N\) be the level at which the substitution forces the border. By extension, we may assume that we are given a partial isometry of the form \(e[\bigcup_i\varsigma^n(t_i),t',t]\). Form a sequence \(\{i_j\}_{j=1}^m\subseteq\{i\}_i\) so that the corresponding punctured \(d\)-cells pairwise neighbor a \(d-1\)-cell, with \(t\in\varsigma^n(t_{i_1})\) and \(t'\in\varsigma^n(t_{i_m})\). Consider a single preimage \(U\) of \(\sigma^{-N}(\bigcup_i\varsigma^n(t_i))\). By border-forcing, \(U\) is an open ball in the subspace topology of \(AP_{n+N}\). Furthermore, since \(\bigcup_i\varsigma^n(t_i)\) contains each of the \(d-1\)-cells that are pairwise intersections of our sequence of punctured \(d\)-cells, \(\sigma^{-N}(\varsigma^n(t_{i_j})\cup\varsigma^n(t_{i_{j+1}}))\cap U\) are also open balls for each \(1\leq j<m\). By extension, we may reduce the support from \(U\) to \(\bigcup_{j=1}^{m-1}\sigma^{-N}(\varsigma^n(t_{i_j})\cup\varsigma^n(t_{i_{j+1}}))\cap U\), then for each partial isometry between adjacent \(d\)-cells, we may further reduce support to the corresponding \(\sigma^{-N}(\varsigma^n(t_{i_j})\cup\varsigma^n(t_{i_{j+1}}))\cap U\).\\
\indent These are full translations around \(d-1\)-cells. To obtain partial isometries, we convert each full translation into punctured \hyperlink{acu}{acu} sets around each of the \(d-1\)-cells.
\end{proof}
\end{proposition}
\begin{proposition}
\label{proposition:full-translation-mvn}
Given any \(d-1\)-cell in the \(AP_n\)-complex, let \(\{\varsigma^n(t_i)\}_i\) and \(\{\varsigma^n(t_j)\}_j\) be a full translation. Then \(\sum_ie[\varsigma^n(t_i),t_i']\) is Murray--von Neumann equivalent to \(\sum_je[\varsigma^n(t_j),t_j'']\) in \(C_r^\ast(\dot{G}_u)\), for any punctures \(t_i'\) in \(\varsigma^n(t_i)\) and \(t_j''\) in \(\varsigma^n(t_j)\).
\begin{proof}
Let \(I\) denote the pairs of indices \((i,j)\) such that \(\varsigma^n(t_i)\cup\varsigma^n(t_j)\) forms a punctured \hyperlink{acu}{acu} set, and let \(I^\top\) be the set with index order swapped. Noting that a sum in \(I\) and one in \(I^\top\) differ in the way individual projections and partial isometries are grouped, namely using the right-handed and the left-handed cells, we have that
\begin{align*}
\sum_ie[\varsigma^n(t_i),t_i']&=\sum_Ie[\varsigma^n(t_i)\cup\varsigma^n(t_j),t_i']\tag*{(Extension)}\\
&=\sum_Ie[\varsigma^n(t_i)\cup\varsigma^n(t_j),t_i',t_j'']e[\varsigma^n(t_i)\cup\varsigma^n(t_j),t_j'',t_i']\\
&\sim_\textnormal{MvN}\sum_{I^\top}e[\varsigma^n(t_i)\cup\varsigma^n(t_j),t_j'',t_i']e[\varsigma^n(t_i)\cup\varsigma^n(t_j),t_i',t_j'']\tag*{(Reordering summands)}\\
&=\sum_{I^\top}e[\varsigma^n(t_i)\cup\varsigma^n(t_j),t_j'']\\
&=\sum_je[\varsigma^n(t_j),t_j''].\tag*{(Extension)}
\end{align*}
\end{proof}
\end{proposition}
\noindent Together, these say that partial isometries given by doubly-pointed patterns that are nontrivial in \(K_0(C_r^\ast(\dot{G}_{AF});C_r^\ast(\dot{G}_u))\) are exactly those that arise from full translations, whose source and range projections land back in \(K_0(C_r^\ast(\dot{G}_{AF}))\).
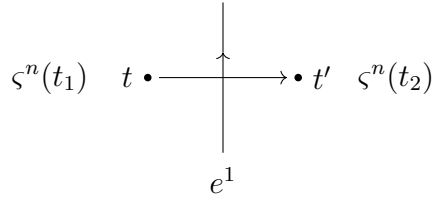
\begin{figure}[t]
\centering
\begin{tikzpicture}
\draw[->](0,-1)--(0,0.3333);
\draw(0,0.3333)--(0,1);
\fill(-1,0)circle[radius=0.05];
\fill(1,0)circle[radius=0.05];
\draw[->](-0.85,0)--(0.85,0);
\draw(0,-1.05)node[below]{\(e^1\)};
\draw(-1.05,0)node[left]{\(t\)};
\draw(1.05,0)node[right]{\(t'\)};
\draw(-2.3,0)node{\(\varsigma^n(t_1)\)};
\draw(2.3,0)node{\(\varsigma^n(t_2)\)};
\end{tikzpicture}
\caption{Assuming no branches occur at \(e^1\), \(\thickening_n\) maps it to relative \(K_0\)-equivalence class of the partial isometry \(e[\varsigma^n(t_1)\cup\varsigma^n(t_2),t',t]\). In other words, \(\varsigma^n(t_1)\) satisfies the right-hand rule with the direction of \(e^1\), and \(t\) is defined to be the source, and \(\varsigma^n(t_2)\) satisfies the left-hand rule, and \(t'\) is defined to be the range. Inverting the partial isometry results in a sign on \(e^1\).}
\label{figure:thickening}
\end{figure}%

\indent We next define the key map between the \(d-1\)-cochain group and the relative \(K_0\)-group that puts these propositions into topological context. We begin by defining it at a finite level, then extend it using functoriality of relative \(K\)-theory.
\begin{definition}
The \emph{thickening map}\footnote{This terminology is borrowed from \cite{juliensavinien16}, where one thickens a codimension-\(1\) cell into an open ball in the \(AP\)-complex.}, \(\thickening_n:C_n^{d-1}\rightarrow K_0(C_r^\ast(\dot{G}_{AF,n});C_r^\ast(\dot{G}_u))\), sends a \(d-1\)-cell to the relative \(K_0\)-equivalence class of partial isometries formed from the full translation surrounding it, with sources in each of the right-handed cells and ranges in the left-handed ones (\cref{figure:thickening}), i.e.
\[
\thickening_n(e^{d-1})=\sum_Ie[\varsigma^n(t_i)\cup\varsigma^n(t_j),t_j'',t_i']
\]
\noindent where the right sum is under the relative \(K_0\)-equivalence, \(\{\varsigma^n(t_i)\}_i\) and \(\{\varsigma^n(t_j)\}_j\) is a full translation around the \(d-1\)-cochain \(e^{d-1}\in C_n^{d-1}\), \(I\) is the set of pairs of indices \((i,j)\) such that \(\varsigma^n(t_i)\cup\varsigma^n(t_j)\) forms a punctured \hyperlink{acu}{acu} set, and \(t_i'\) and \(t_j''\) are punctures in \(\varsigma^n(t_i)\) and \(\varsigma^n(t_j)\) for each \(i\) and \(j\), respectively.
\end{definition}
\begin{remark}
If one reverses the order of the sources and ranges, then the map \(\thickening_n\) will have the opposite sign. More generally, the map will carry signs on the \(d-1\)-cells whose orientations are reversed.
\end{remark}
\begin{figure}[t]
\centering
\begin{tikzpicture}
\draw(-1,1.7321)--(0,0);
\draw(-1,-1.7321)--(0,0);
\draw(0,0)--(2,0);
\draw(2,0)--(3,1.7321);
\draw(2,0)--(3,-1.7321);
\draw(-3.5,0)node{\(\varsigma^n(t_1)\)};
\draw(1,-1.7321)node{\(\varsigma^n(t_2)\)};
\draw(1,1.7321)node{\(\varsigma^n(t_4)\)};
\draw(5.5,0)node{\(\varsigma^n(t_3)\)};
\fill(-1,0)circle[radius=0.05];
\fill(1,-0.8660)circle[radius=0.05];
\fill(3,0)circle[radius=0.05];
\fill(1,0.8660)circle[radius=0.05];
\draw[->](-0.8624,-0.0596)--(0.8624,-0.8064);
\draw[->](1.1376,-0.8064)--(2.8624,-0.0596);
\draw[->](2.8624,0.0596)--(1.1376,0.8064);
\draw[->](0.8624,0.8064)--(-0.8624,0.0596);
\draw(-1.05,0)node[left]{\(t_1'=t_4''\)};
\draw(1,-0.9160)node[below]{\(t_2'=t_1''\)};
\draw(3.05,0)node[right]{\(t_3'=t_2''\)};
\draw(1,0.9160)node[above]{\(t_4'=t_3''\)};
\draw[->,dashed](0.925,-0.7361)--(0.925,0.7361);
\draw[->,dashed](1.075,0.7361)--(1.075,-0.7361);
\end{tikzpicture}
\caption{A relation in \(K_0(C_r^\ast(\dot{G}_{AF});C_r^\ast(\dot{G}_u))\) (solid arrows forming a quadrilateral) that surrounds two \(0\)-cells decomposing into a sum of two relations (solid arrows and dashed arrows forming two triangles), each of which surrounds only one \(0\)-cell. The vertical partial isometries (dashed arrows) are negatives of each other, which cancel when we add the two smaller relations together to obtain the larger relation.}
\label{figure:relation-decompose}
\end{figure}
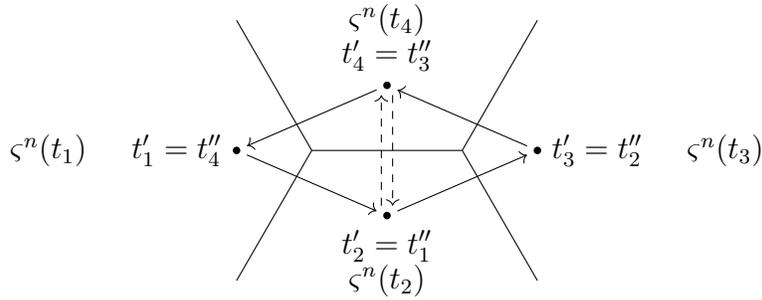%

\begin{proposition}
\(K_0(C_r^\ast(\dot{G}_{AF});C_r^\ast(\dot{G}_u))=\varinjlim_n(K_0(C_r^\ast(\dot{G}_{AF,n});C_r^\ast(\dot{G}_u)),\sigma^\top)\).
\begin{proof}
Since the diagram
\[
\begin{tikzcd}
C_r^\ast(\dot{G}_{AF,n})\arrow[r]\arrow[d,hook]&C_r^\ast(\dot{G}_u)\arrow[d,equal]\\
C_r^\ast(\dot{G}_{AF,n+1})\arrow[r]\arrow[d,hook]&C_r^\ast(\dot{G}_u)\arrow[d,equal]\\
\vdots\arrow[d,hook]&\vdots\arrow[d,equal]\\
C_r^\ast(\dot{G}_{AF})\arrow[r]&C_r^\ast(\dot{G}_u)
\end{tikzcd}
\]
commutes, functoriality of relative \(K\)-theory gives the existence of a map \(\sigma^\top:K_0(C_r^\ast(\dot{G}_{AF,n});C_r^\ast(\dot{G}_u))\rightarrow K_0(C_r^\ast(\dot{G}_{AF,n+1});C_r^\ast(\dot{G}_u))\), and together with the Universal Property of Direct Limits and the Five Lemma, we get the result.
\end{proof}
\end{proposition}
\indent Before stating the thickening map as extended to direct limits, we notice that there are relations among full translations that one can fashion by forming a ``loop'' so that the source of the first full translation is the range of the last (\cref{figure:relation-decompose}). Any such loop can always be decomposed to elementary ones that surround exactly one \(d-2\)-cell, thus we will only work with relations of this type. The thickening map identifies such loops with \(\image\delta_n^{d-2}\). We package this with compatibility with \(\sigma^\top\) into the following proposition.
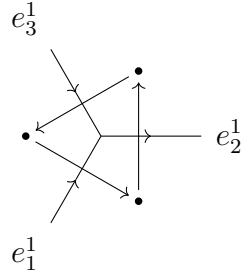
\begin{figure}[t]
\centering
\begin{tikzpicture}
\draw[->](-0.6667,1.1547)--(-0.3333,0.5774);
\draw(-0.3333,0.5774)--(0,0);
\draw[->](-0.6667,-1.1547)--(-0.3333,-0.5774);
\draw(-0.3333,-0.5774)--(0,0);
\draw[->](0,0)--(0.6667,0);
\draw(0.6667,0)--(1.3333,0);
\fill(-1,0)circle[radius=0.05];
\fill(0.5,-0.8660)circle[radius=0.05];
\fill(0.5,0.8660)circle[radius=0.05];
\draw[->](-0.8701,-0.075)--(0.3701,-0.7910);
\draw[->](0.5,-0.7160)--(0.5,0.7160);
\draw[->](0.3701,0.7910)--(-0.8701,0.075);
\draw(-0.6917,-1.1980)node[below left]{\(e_1^1\)};
\draw(1.3833,0)node[right]{\(e_2^1\)};
\draw(-0.6917,1.1980)node[above left]{\(e_3^1\)};
\end{tikzpicture}
\caption{The relation given by the sum of the partial isometries identifies with \(e_1^1-e_2^1+e_3^1=0\) under \(\thickening_n\), which coincides with \(\delta_n^0\) applied to the center \(0\)-cell. Here we arbitrarily prescribed some orientations to the \(1\)-cells.}
\label{figure:full-translation-image-delta0}
\end{figure}%

\begin{figure}[t]
\centering
\begin{tikzpicture}
\draw[->](0,0)--(0.5,0);
\draw(0.5,0)--(1,0);
\draw(0.5,0)node[below]{\(e^1\)};
\draw[|->](1.25,0)--(1.75,0);
\draw[->](2,0)--(2.5,0);
\draw(2.5,0)--(3,0);
\draw(3,0)--(3.5,0);
\draw[->](4,0)--(3.5,0);
\draw(3,-0.1)--(3,0.1);
\draw(2.5,0)node[below]{\(e_1^1\)};
\draw(3.5,0)node[below]{\(-e_2^1\)};
\draw(0,-1.5)--(1,-1.5);
\fill(0.5,-2)circle[radius=0.05];
\fill(0.5,-1)circle[radius=0.05];
\draw[->](0.5,-1.15)--(0.5,-1.85);
\draw[|->](1.25,-1.5)--(1.75,-1.5);
\draw(2,-1.5)--(4,-1.5);
\draw(3,-1.4)--(3,-1.6);
\fill(2.5,-2)circle[radius=0.05];
\fill(3.5,-2)circle[radius=0.05];
\fill(2.5,-1)circle[radius=0.05];
\fill(3.5,-1)circle[radius=0.05];
\draw[->](2.5,-1.15)--(2.5,-1.85);
\draw[->](3.5,-1.85)--(3.5,-1.15);
\draw(2.5,-2)node[below]{\(+1\)};
\draw(3.5,-2)node[below]{\(-1\)};
\end{tikzpicture}
\caption{The substitution map is \(e^1\mapsto e_1^1-e_2^1\) with the orientation of \(e_2^1\) reversed in comparison to \(e^1\) (top). Under thickening, the partial isometry \(\textnormal{th}(e_2^1)\) carries a sign and therefore has its direction reversed to agree with the substitution map applied to the partial isometry \(\textnormal{th}(e^1)\) (bottom).}
\label{figure:substitution-reverse}
\end{figure}
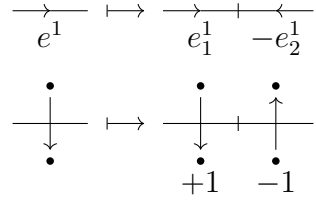%

\begin{proposition}
\label{proposition:thickening-compatible-substitution}
We have the commutative diagram
\[
\begin{tikzcd}
C_n^{d-1}\arrow[r,"\sigma^\top"]\arrow[rd,two heads]\arrow[rddd,swap,"\thickening_n"]&C_{n+1}^{d-1}\arrow[rd,two heads]\arrow[rddd,swap,"\thickening_{n+1}"]&\\
&C_n^{d-1}/\image\delta_n^{d-2}\arrow[r,crossing over,"\sigma^\top"]\arrow[dd,"\thickening_n"]&C_{n+1}^{d-1}/\image\delta_{n+1}^{d-2}\arrow[dd,"\thickening_{n+1}"]\\
&&\\
&K_0(C_r^\ast(\dot{G}_{AF,n});C_r^\ast(\dot{G}_u))\arrow[r]&K_0(C_r^\ast(\dot{G}_{AF,n+1});C_r^\ast(\dot{G}_u))
\end{tikzcd}
\]
\noindent where the quotients exist for \(d\geq 2\).
\begin{proof}
We first show, for \(d\geq 2\), the existence of the map \(\thickening_n:C_n^{d-1}/\image\delta_n^{d-2}\rightarrow K_0(C_r^\ast(\dot{G}_{AF,n});C_r^\ast(\dot{G}_u))\) by showing that \(\thickening_n(\image\delta_n^{d-2})=0\). Given a \(d-2\)-cell whose image under \(\delta_n^{d-2}\) is a signed sum \(\sum_i(-1)^{m_i}e_i^{d-1}\), we can apply the thickening map summand-wise to each \(e_i^{d-1}\) to obtain an associated full translation \(\thickening_n(e_i^{d-1})\). We then obtain the sum
\[
\sum_i(-1)^{m_i}\thickening_n(e_i^{d-1})=\sum_i\thickening_n(e_i^{d-1})^{\ast m_i}
\]
\noindent with the exponent on the right denoting that we take adjoint \(m_i\)-times. Noticing that taking the union of all such full translations forms an open ball around the \(d-2\)-cell, we see that if two full translations intersect nontrivially at either their sources or ranges (possibly source-to-source or range-to-range), then they must share the entire source or range. Consider any two full translations with associated indices \(i\) and \(i+1\) so that they share either the source or the range. We wish to show that with \(\ast m_i\), the sharing must always be source-to-range (\cref{figure:full-translation-image-delta0}). Suppose that \(\delta_n^{d-2}\) assigns the same sign to both \(d-1\)-cells. Since we chose the orientations of each of the \(d\)-cells to be the same, the attaching maps to the two \(d-1\)-cells assign opposite signs to the shared source or range of the two full translations, thus one must be the source and the other the range. If the \(d-1\)-cells have opposite signs, then the shared source or range are either both the source or both the range. However, on applying \(\ast m\), the one with the sign has the source and range flipped, giving that the sharing is source-to-range. Procedurally applying adjoints as \(m_i\) requires gives us a loop, which, by principality of our groupoid, is trivial. Thus, by the Universal Property of Quotients, thickening descends to the quotient.\\
\indent For \(d\geq 2\), the top square commutes because while the map \(\sigma:AP_{n+1}\rightarrow AP_n\) may not be cellular, one can easily make it cellular by introducing \(k\)-cells in \(AP_{n+1}\) as they exist under \(\sigma\).\\
\indent We next show that the bottom slanted square also commutes, that is that thickening is compatible with \(\sigma^\top:C_n^{d-1}\rightarrow C_{n+1}^{d-1}\). Traversing the bottom composition, by continuity, \(\sigma^\top\) sends the underlying open ball of the full translation that is the thickened \(d-1\)-cell to a collection of open balls. Such an open ball either belongs to the interior of a single \(d\)-cell, or contains a \(d-1\)-cell that is the image of the \(d-1\)-cell under \(\sigma^\top\). In the former case, it becomes a partial isometry whose source and range belongs in a single \(d\)-cell, and is thus trivial. In the latter, the partitioning of the full translation as inherited from its preimage agrees with the top composition since if \(\sigma^\top\) reverses orientation on a \(d-1\)-cell, the cell will carry a sign, which undoes the reversed orientation of the corresponding partial isometry (\cref{figure:substitution-reverse}).\\
\indent Finally, a diagram chase shows the front square commutes.
\end{proof}
\end{proposition}
\noindent Thus by functoriality and the Universal Property of Direct Limits, we obtain a map \(\thickening:C^{d-1}/\image\delta^{d-2}\rightarrow K_0(C_r^\ast(\dot{G}_{AF});C_r^\ast(\dot{G}_u))\), allowing us to state the theorem that relates the top two \v Cech cohomology groups to \(K\)-theory.
\begin{theorem}
\label{theorem:k-theory-cohomology}
We have the commutative diagram
\[
\begin{tikzcd}
0\arrow[r]&\check{H}^{d-1}(\Omega_T)\arrow[r,hook]\arrow[d]&C^{d-1}/\image\delta^{d-2}\arrow[r,"\delta^{d-1}"]\arrow[d,"\thickening"]&C^d\arrow[r,two heads]\arrow[d,"\CF"]&\check{H}^d(\Omega_T)\arrow[r]\arrow[d]&0&\\
0\arrow[r]&K_{1,u}\arrow[r,hook]&K_{0,AF;u}\arrow[r,"\evaluation"]&K_{0,AF}\arrow[r,"\iota_\ast"]&K_{0,u}\arrow[r,two heads]&K_{1,AF;u}\arrow[r]&0
\end{tikzcd}
\]
\noindent where \(K_{i,u}=K_i(C_r^\ast(\dot{G}_u))\), \(K_{i,AF}=K_i(C_r^\ast(\dot{G}_{AF}))\), and \(K_{i,AF;u}=K_i(C_r^\ast(\dot{G}_{AF});C_r^\ast(\dot{G}_u))\).
\begin{proof}
Since \(\delta^{d-1}\circ\delta^{d-2}=0\), we obtain the existence of the map \(\delta^{d-1}:C^{d-1}/\image\delta^{d-2}\rightarrow C^d\). Then the top exact sequence is none other than the exact sequence
\[
\begin{tikzcd}
0\arrow[r]&\kernel\phi\arrow[r,hook]&A\arrow[r,"\phi"]&B\arrow[r,two heads]&\cokernel\phi\arrow[r]&0
\end{tikzcd}
\]
\noindent applied to this map. It suffices to check that the middle square commutes at a finite level, which it does since thickening counts multiplicities, so \(\evaluation\circ\thickening_n([e^{d-1}])=[\sum_ic_ie_i^d]-[\sum_jc_je_j^d]\) where \(c_i\) is the number of occurrences of \(e^{d-1}\) in \(e_i^d\) with the orientations of the attaching maps agreeing, and \(c_j\) is disagreeing, thus this difference is \(\delta_n^{d-1}([e^{d-1}])\).
\end{proof}
\end{theorem}

\subsection{Isomorphism for \(d\leq 2\)}
\indent Over \(\mathbb{Z}\), it is known that for \(d\leq 3\) and for \(i=0,1\), \(K_i(C_r^\ast(\dot{G}_u))\cong\bigoplus_{k=0}^{\lfloor (d-i)/2\rfloor}\check{H}^{d-i-2k}(\Omega_T)\). For example, by \cite{sadunwilliams03}, tiling spaces are fiber bundles over \(\mathbb{Z}^d\), thus are decorations of \(\mathbb{Z}^d\), and one can then apply \cite{vanelst94}\footnote{Modulo mistakes for \(d=3\).} and observe that the generators and relations given are in terms of cohomology, or one can use the Chern and Connes--Thom isomorphisms. For an incomplete list, see \cite{putnam89}, \cite{kellendonk97}, \cite{andersonputnam98}, and \cite{forresthunton99}. However, in the interest of explicitly computing the map between generators and relations of the two, and of promoting relative \(K\)-theory, we will reprove this fact for \(d\leq 2\).\\
\indent The argument follows the same basic structure demonstrated in \cite{juliensavinien16} for the chair tiling, and is similar in spirit to \cite{vanelst94}, where one inductively ``peels off'' dimensions by appealing to the techniques demonstrated in the following theorem whose conclusion has been known since \cite{putnam89}.
\begin{theorem}[\cref{theorem:isomorphism-d-1}]
For \(d=1\), we have an isomorphism of exact sequences
\[
\begin{tikzcd}
0\arrow[r]&\check{H}^0(\Omega_T)\arrow[r,hook]\arrow[d]&C^0\arrow[r,"\delta^0"]\arrow[d,"\thickening"]&C^1\arrow[r,two heads]\arrow[d,"\CF"]&\check{H}^1(\Omega_T)\arrow[r]\arrow[d]&0\\
0\arrow[r]&K_{1,u}\arrow[r,hook]&K_{0,AF;u}\arrow[r,"\evaluation"]&K_{0,AF}\arrow[r,"\iota_\ast"]&K_{0,u}\arrow[r]&0
\end{tikzcd}
\]
\noindent where \(K_{i,u}=K_i(C_r^\ast(\dot{G}_u))\), \(K_{i,AF}=K_i(C_r^\ast(\dot{G}_{AF}))\), and \(K_{i,AF;u}=K_i(C_r^\ast(\dot{G}_{AF});C_r^\ast(\dot{G}_u))\). Thus
\begin{align*}
K_0(C_r^\ast(\dot{G}_u))&\cong K_0(C_r^\ast(\dot{G}_{AF}))/\image\evaluation\\
K_1(C_r^\ast(\dot{G}_u))&\cong\mathbb{Z}.
\end{align*}
\begin{proof}
Noting that the nontrivial elements in \(K_0(C_r^\ast(\dot{G}_{AF});C_r^\ast(\dot{G}_u))\) are partial isometries that intersect \(\varprojlim_n(AP_n^{(0)},\sigma)\), and thus \(\dot{G}_{AF}^{(0)}\) forms an abstract transversal, we apply \cref{theorem:excision-isomorphism} followed by restriction to \(\varprojlim_n(AP_n^{(0)},\sigma)\) that is a groupoid equivalence. The thickening map is the inverse map.\\
\indent More precisely, let \(L\) be the set of elements in \(G_u\) whose source belongs to the left and range belongs to the right of \(\varprojlim_n(AP_n^{(0)},\sigma)\). Then, forming the corresponding subgroupoids \(H\) and \(H'\), \(H\leq\dot{G}_u\) is the subgroupoid supported on points in \(\varprojlim_n(AP_n,\sigma)\) that are translations of \(\varprojlim_n(AP_n^{(0)},\sigma)\), and \(H'\leq\dot{G}_{AF}\) is similar, with appropriate topologies that is measured with respect to \(\varprojlim_n(AP_n^{(0)},\sigma)\) so that \(\varprojlim_n(AP_n^{(0)},\sigma)\) forms an abstract transversal to \(H\) and gives a groupoid equivalence between \(H\) and \(\dot{G}_{AF}^{(0)}\).\\
\indent By \cref{theorem:excision-isomorphism}, we have the sequence of isomorphisms
\[
K_i(C_r^\ast(\dot{G}_{AF});C_r^\ast(\dot{G}_u))\cong K_i(C_r^\ast(H');C_r^\ast(H))\cong K_i(C_r^\ast(H))\cong K_i(C_r^\ast(\dot{G}_{AF}^{(0)})).
\]
\noindent Applying \cref{lemma:af-cochain} and reversing the sequence of isomorphisms gives the thickening map.
\end{proof}
\end{theorem}
\begin{remark}
This theorem allows us to prove, for \(d=1\), facts regarding proper substitutions and substitutions that do not force the border.
\begin{itemize}
\item
\emph{If a substitution is \emph{proper}\footnote{The resulting Bratteli--Vershik system is proper.}, i.e. every supertile begins with the same prototile and ends with the same prototile, then \(K_0(C_r^\ast(\dot{G}_{AF}))\cong K_0(C_r^\ast(\dot{G}_u))\).} Properness holds for all higher-level supertiles as well, thus there is only one nontrivial partial isometry. \(K_1(C_r^\ast(\dot{G}_u))\cong\mathbb{Z}\leq K_0(C_r^\ast(\dot{G}_{AF});C_r^\ast(\dot{G}_u))\), and \(K_0(C_r^\ast(\dot{G}_{AF}))\) is torsion-free, therefore \(\evaluation=0\).
\item
\emph{If a substitution does not force the border, then \(K_0(C_r^\ast(\dot{G}_u))\cong K_0(C_r^\ast(\dot{G}_{AF}))/\image\evaluation\) is a nontrivial quotient.} There are multiple partial isometries that persist in the limit, thus \(K_1(C_r^\ast(\dot{G}_u))\lneq K_0(C_r^\ast(\dot{G}_{AF});C_r^\ast(\dot{G}_u))\), and \(\evaluation\neq 0\).
\end{itemize}%

\end{remark}
\indent For \(d=2\), we now pass to \(\mathbb{Z}^2\)-decorations using \cite{sadunwilliams03}. Summarizing, the set of vectors forming the boundaries of prototiles in \(T\) forms loops, thus they correspond to an element in the kernel of a linear system of equations in \(\mathbb{R}^2\) solved over \(\mathbb{R}\). There exists an element in the kernel over \(\mathbb{Q}\) that is arbitrarily close to the original element, and deforming\footnote{Such deformations can be viewed as elements of \(\check{H}^1(\Omega_T;\mathbb{R}^2)\).} the prototiles to this new element gives a homeomorphism between \(\Omega_T\) and the new tiling space formed from a tiling whose \(0\)-cells of tiles belong to \(\mathbb{Q}^2\). One can then uniformly scale the coordinates of the \(0\)-cells until they become integral, there being finitely many, up to translation, due to finite local complexity, then replace each \(1\)-cell with a collection of piecewise-constant \(1\)-cells where at least one of the coordinates of each piece is integral. Finally, one splits each of the new prototiles into unit squares decorated by the prototile they reside in. Each step is at least a homeomorphism of tiling spaces, and so the resulting tiling space is homeomorphic to the original. Let us denote this new tiling by \(T_\square\), and the homeomorphism by \(\square:\Omega_T\rightarrow\Omega_{T_\square}\).\\
\indent Rather than attempting to bring the substitution structure from \(T\) through \(\square\), we construct a new inverse system that is a hybrid of the Anderson--Putnam and G\"ahler (\cite{gahler02}) inverse systems. Let \(\Gamma_n\) be the \(CW\)-complex where the \(2\)-cells are patches of squares tiles residing in \([i\cdot 2^n,(i+1)\cdot 2^n]\times[j\cdot 2^n,(j+1)\cdot 2^n]\) collared by their neighboring \(2\)-cells of the same size, for \(i,j\in\mathbb{Z}\), up to translation, and we quotient by identifications as they exist in \(T_\square\). The substitution map \(\sigma_\square:\Gamma_n\rightarrow\Gamma_{n-1}\) is subdividing each \(2\)-cell into \(\Gamma_n\) into their respective \(2\)-cells in \(\Gamma_{n-1}\). One easily checks, following the arguments given in \cite[Theorem 4.3]{andersonputnam98} and \cref{corollary:robinson-map-conjugacy}, that, under the translation action, \(\varprojlim_n(\Gamma_n,\sigma_\square)\) is topologically conjugate to \(\Omega_{T_\square}\). Note that \(T_\square\) satisfies the boundary hyperplane condition, thus the cells in each \(\Gamma_n\) are assigned orientations as aforedescribed.\\
\indent Let us denote the unstable and tail groupoids of each dimension from this inverse limit \(\varprojlim_n(\Gamma_n,\sigma_\square)\) by the same notation as before, and introduce two new types of objects that are refinements of these groupoids by observing that in \(\Gamma_n\), the \(1\)-cells are either horizontal or vertical. Let us denote \(\Gamma_{n,h}^{(1)}\) and \(\Gamma_{n,v}^{(1)}\) the horizontal and vertical parts of the \(1\)-skeleton of \(\Gamma_n\), respectively. Let \(\bullet\) be either \(h\) or \(v\), and let \(h^\top=v\) and \(v^\top=h\).
\begin{itemize}
\item
Let \(G_{u,\bullet}^{(1)}\) be the subgroupoid of \(G_u^{(1)}\) restricted to \(\varprojlim_n(\Gamma_{n,\bullet}^{(1)},\sigma_\square)\), i.e., with the source and range maps taken to be on \(G_u^{(1)}\),
\[
G_{u,\bullet}^{(1)}=s^{-1}(\varprojlim_n(\Gamma_{n,\bullet}^{(1)},\sigma_\square))\cap r^{-1}(\varprojlim_n(\Gamma_{n,\bullet}^{(1)},\sigma_\square)),
\]
\noindent where we have identified the unit space with the underlying topological space.
\item
Let \(G_{\textnormal{tail},\bullet}^{(1)}\) be the tail subgroupoid of \(G_{u,\bullet}^{(1)}\) induced by \(\sigma\), i.e.
\[
G_{\textnormal{tail},\bullet}^{(1)}=\left\{((p_0,\ldots),(q_0,\ldots))\in G_{u,\bullet}^{(1)}:\begin{array}{c}\exists U_n\times U_{n+1}\times\ldots\textnormal{ such that}\\\forall m\geq n\textnormal{, }p_m,q_m\in U_m\textnormal{ and}\\\exists 1\textnormal{-cells }e_{1,m}\subseteq\Gamma_{m,\bullet}^{(1)}\textnormal{ such that }U_m\subseteq e_{1,m}\end{array}\right\}.
\]
\item
Let \(G_{u,\bullet}\) be the subgroupoid of \(G_u\) with the elements crossing \(\varprojlim_n(\Gamma_{n,\bullet}^{(1)},\sigma_\square)\) removed, e.g. for \(\bullet=h\),
\[
G_{u,h}=\left\{((p_0,\ldots),(q_0,\ldots))\in G_u:\begin{array}{c}\exists y_1,y_2\textnormal{ such that }\sign(y_1)=\sign(y_2)\textnormal{ and}\\(p_0,\ldots)-(0,y_1),(q_0,\ldots)-(0,y_2)\in\varprojlim_n(\Gamma_{n,h}^{(1)},\sigma_\square)\end{array}\right\}.
\]
\item
Let \(C_\bullet^1\) be the subgroup of \(C^1\) restricted to \(\varprojlim_n(\Gamma_{n,\bullet}^{(1)},\sigma_\square)\), and \(\delta_\bullet^0:C^0\rightarrow C_\bullet^1\) and \(\delta_\bullet^1:C_\bullet^1\rightarrow C^2\) be the induced coboundary maps on the corresponding cochains. Note that \(\delta_\bullet^1\circ\delta_\bullet^0\) is not necessarily zero.
\end{itemize}%

\noindent As before, we introduce punctures to each of the cells, the natural choice being the exact middle of each of the cells, and denote the restricted groupoids with a dot.\\
\indent Let us consider the three inclusions of groupoids
\begin{gather}
\dot{G}_{u,\bullet}\hookrightarrow\dot{G}_u\\
\dot{G}_{AF}\hookrightarrow\dot{G}_{u,\bullet^\top}\\
\dot{G}_{AF,\bullet}^{(1)}\hookrightarrow\dot{G}_{u,\bullet}^{(1)}
\end{gather}
\noindent and their associated six-term sequences in relative \(K\)-theory.\\
\indent We have the following proposition that simplifies the relative \(K\)-groups.
\begin{proposition}
For \(i=0,1\),
\begin{enumerate}
\item
\(K_i(C_r^\ast(\dot{G}_{u,\bullet});C_r^\ast(\dot{G}_u))\cong K_i(C_r^\ast(\dot{G}_{u,\bullet}^{(1)}))\),
\item
\(K_i(C_r^\ast(\dot{G}_{AF});C_r^\ast(\dot{G}_{u,\bullet^\top}))\cong K_i(C_r^\ast(\dot{G}_{AF,\bullet}^{(1)}))\), and
\item
\(K_i(C_r^\ast(\dot{G}_{AF,\bullet}^{(1)});C_r^\ast(\dot{G}_{u,\bullet}^{(1)}))\cong K_i(C_r^\ast(\dot{G}_{AF}^{(0)}))\).
\end{enumerate}%

\begin{proof}
These all follow from \cref{theorem:excision-isomorphism} with appropriate choices of \(L\), so that \(L\cup L^{-1}\) is
\begin{enumerate}
\item
The subset of the groupoid \(\dot{G}_u\) that cross \(\varprojlim_n(\Gamma_{n,\bullet}^{(1)},\sigma_\square)\),
\item
The subset of the groupoid \(\dot{G}_{u,\bullet^\top}\) that cross \(\varprojlim_n(\Gamma_{n,\bullet}^{(1)},\sigma_\square)\), and
\item
The subset of the groupoid \(\dot{G}_{u,\bullet}^{(1)}\) that cross \(\varprojlim_n(\Gamma_n^{(0)},\sigma_\square)\).
\end{enumerate}%

\indent Let us explain the first one with \(\bullet=h\). The others are similar. Let \(L\) be the subset of \(\dot{G}_u\) consisting of elements of the form \(((p_0,\ldots)+(0,y_1),(q_0,\ldots)+(0,y_2))\), with \((p_0,\ldots),(q_0,\ldots)\in\varprojlim_n(\Gamma_{n,h}^{(1)},\sigma_\square)\) and \(y_1<0\) and \(y_2>0\). In the notation of \cref{theorem:excision-isomorphism}, \(G=\dot{G}_u\), \(G'=\dot{G}_{u,h}\), and \(H=\dot{G}_{u,h}^{(1)}\), and the isomorphism is as desired.
\end{proof}
\end{proposition}
\indent Composing these isomorphisms (in reverse) with the evaluation maps (and the counterparts from the relative \(K_1\)-groups) for each of the three six-term sequences, we obtain the three maps
\begin{align*}
\delta_i&:K_i(C_r^\ast(\dot{G}_{u,\bullet}^{(1)}))\rightarrow K_i(C_r^\ast(\dot{G}_{u,\bullet}))\\
\delta_\bullet^1&:K_i(C_r^\ast(\dot{G}_{AF,\bullet}^{(1)}))\rightarrow K_i(C_r^\ast(\dot{G}_{AF}))\\
\delta_\bullet^0&:K_i(C_r^\ast(\dot{G}_{AF}^{(0)}))\rightarrow K_i(C_r^\ast(\dot{G}_{AF,\bullet}^{(1)}))
\end{align*}
\noindent where the latter two names are as given because they are the coboundary maps restricted to either the horizontal or the vertical \(1\)-cells.\\
\indent The three six-term sequences in relative \(K\)-theory then read
\begin{equation}
\label{equation:six-term-d-2}
\begin{tikzcd}
K_0(C_r^\ast(\dot{G}_{u,h}^{(1)}))\arrow[r,"\delta_0"]&K_0(C_r^\ast(\dot{G}_{u,h}))\arrow[r,"\iota_\ast"]&K_0(C_r^\ast(\dot{G}_u))\arrow[d]\\
K_1(C_r^\ast(\dot{G}_u))\arrow[u]&K_1(C_r^\ast(\dot{G}_{u,h}))\arrow[l,swap,"\iota_\ast"]&K_1(C_r^\ast(\dot{G}_{u,h}^{(1)}))\arrow[l,swap,"\delta_1"]
\end{tikzcd}
\end{equation}
\begin{equation}
\label{equation:six-term-d-1}
\begin{tikzcd}
K_0(C_r^\ast(\dot{G}_{AF,v}^{(1)}))\arrow[r,"\delta_v^1"]&K_0(C_r^\ast(\dot{G}_{AF}))\arrow[r,"\iota_\ast"]&K_0(C_r^\ast(\dot{G}_{u,h}))\arrow[d,two heads]\\
K_1(C_r^\ast(\dot{G}_{u,h}))\arrow[u,hook]&0\arrow[l]&0\arrow[l]
\end{tikzcd}
\end{equation}
\begin{equation}
\label{equation:six-term-d-0}
\begin{tikzcd}
K_0(C_r^\ast(\dot{G}_{AF}^{(0)}))\arrow[r,"\delta_h^0"]&K_0(C_r^\ast(\dot{G}_{AF,h}^{(1)}))\arrow[r,"\iota_\ast"]&K_0(C_r^\ast(\dot{G}_{u,h}^{(1)}))\arrow[d,two heads]\\
K_1(C_r^\ast(\dot{G}_{u,h}^{(1)}))\arrow[u,hook]&0\arrow[l]&0\arrow[l].
\end{tikzcd}
\end{equation}
\noindent From \cref{equation:six-term-d-0},
\begin{align*}
K_0(C_r^\ast(\dot{G}_{u,h}^{(1)}))&\cong\cokernel\delta_h^0\\
K_1(C_r^\ast(\dot{G}_{u,h}^{(1)}))&\cong\kernel\delta_h^0,
\end{align*}
\noindent and from \cref{equation:six-term-d-1},
\begin{align*}
K_0(C_r^\ast(\dot{G}_{u,h}))&\cong\cokernel\delta_v^1\\
K_1(C_r^\ast(\dot{G}_{u,h}))&\cong\kernel\delta_v^1.
\end{align*}
\noindent Using the identification from \cref{lemma:af-cochain}, let us show that \(\delta_0\) can be identified with \(\widetilde{\delta}_h^1\) in the diagram
\[
\begin{tikzcd}
C_h^1\arrow[rr,"\delta_h^1"]\arrow[d,two heads]&&C^2\arrow[d,two heads]\\
C_h^1/\image\delta_h^0\arrow[rr,dashed,"\widetilde{\delta}_h^1"]\arrow[d,"\CF"]&&C^2/\image\delta_v^1\arrow[d,"\CF"]\\
K_0(C_r^\ast(\dot{G}_{AF,h}^{(1)}))/\image\delta_h^0\arrow[rr]\arrow[d,"\sim"]&&K_0(C_r^\ast(\dot{G}_{AF}))/\image\delta_v^1\arrow[d,"\sim"]\\
K_0(C_r^\ast(\dot{G}_{u,h}^{(1)}))\arrow[r,"\sim"]\arrow[rr,bend right=15,swap,"\delta_0"]&K_0(C_r^\ast(\dot{G}_{u,h});C_r^\ast(\dot{G}_u))\arrow[r]&K_0(C_r^\ast(\dot{G}_{u,h})).
\end{tikzcd}
\]
\noindent The map \(\delta_0\) is thickening a projection supported on a horizontal \(1\)-cell, then evaluating by taking formal differences of projections supported on neighboring \(2\)-cells, while noting that this is modulo horizontal translations. Interpreted topologically, this is given by taking a horizontal \(1\)-cell \(c_h^1\), taking its coboundary under \(\delta_h^1\), then taking equivalence classes modulo \(\image\delta_v^1\). Formally, this is given by the diagram
\[
\begin{tikzcd}
C^1\arrow[rr,"\delta^1"]\arrow[dd,two heads]&[-2em]&C^2\arrow[rr,two heads]&[-2em]&C^2/\image\delta_v^1&[-2em]\\
[-2em]&(c_h^1,0)\arrow[rr,mapsto]&&\delta^1((c_h^1,0))&&{[\delta_h^1(c_h^1)]}\\
C_h^1\arrow[rr,two heads]&&C_h^1/\image\delta_h^0\arrow[rruu,dashed,bend right=15,swap,"\widetilde{\delta}_h^1"]&&&\\
[-2em]&c_h^1\arrow[rr,mapsto]\arrow[from=uu,mapsto,crossing over]&&{[c_h^1]}\arrow[rruu,dashed,mapsto,bend right=15]\arrow[from=uu,to=rruu,mapsto,crossing over]&&
\end{tikzcd}
\]
\begin{figure}[t]
\centering
\begin{tikzpicture}
\fill(0,0)circle[radius=0.05];
\draw(-1,0)--(1,0);
\draw(0,-1)--(0,1);
\draw(-0.5,-0.05)node[below]{\(e_{h,1}^1\)};
\draw(0.5,-0.05)node[below]{\(e_{h,2}^1\)};
\draw[<->](-0.5,0.05)to[bend left=30](0.5,0.05);
\draw[|->](1.25,0)--(1.75,0)node[midway,above]{\(\delta_h^1\)};
\fill(3.0,0)circle[radius=0.05];
\draw(2.0,0)--(4.0,0);
\draw(3.0,-1)--(3.0,1);
\draw(3.0,-1.05)node[below]{\(e_{v,1}^1\)};
\draw(3.0,1.05)node[above]{\(e_{v,2}^1\)};
\draw(2.45,-.5)node[left]{\(e_1^2\)};
\draw(2.45,.5)node[left]{\(e_3^2\)};
\draw(3.55,-.5)node[right]{\(e_2^2\)};
\draw(3.55,.5)node[right]{\(e_4^2\)};
\draw[<->](2.5,-.5)to(3.5,-.5);
\draw[<->](2.5,.5)to(3.5,.5);
\end{tikzpicture}
\caption{The identification of two \(1\)-cells in \(C_h^1\) across \(\image\delta_h^0\) (left) is mapped, under \(\delta_h^1\), to the identification of the corresponding \(2\)-cells in \(C^2\) across \(\image\delta_v^1\) (right).}
\label{figure:delta-h-1}
\end{figure}
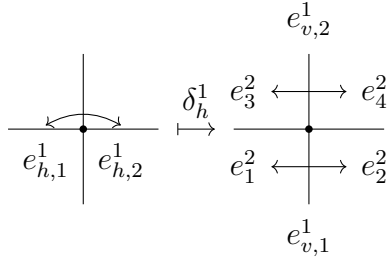%

\noindent where this is independent of the representative chosen in \([c_h^1]\) (\cref{figure:delta-h-1}). In words, given \([c_h^1]\in C_h^1/\image\delta_h^0\), we lift it to an element \(c_h^1\in C_h^1\). This is \((c_h^1,0)\in C^1\). Then \(\delta^1((c_h^1,0))=\delta_h^1(c_h^1)\). Finally, we take its quotient \([\delta_h^1(c_h^1)]\in C^2/\image\delta_v^1\). Thus the map is \(\widetilde{\delta}_h^1([c_h^1])=[\delta_h^1(c_h^1)]\).\\
\indent Similarly, again using \cref{lemma:af-cochain}, \(\delta_1\) can be identified with \(\widetilde{\delta}_v^0\) from the diagram
\[
\begin{tikzcd}
\kernel\delta_h^0\arrow[rr,dashed,"\widetilde{\delta}_v^0"]\arrow[d,"\sim"]&&\kernel\delta_v^1\arrow[d,"\sim"]\\
K_1(C_r^\ast(\dot{G}_{u,h}^{(1)}))\arrow[r,"\sim"]\arrow[rr,bend left=15,"\delta_1"]\arrow[d,hook]&K_1(C_r^\ast(\dot{G}_{u,h});C_r^\ast(\dot{G}_u))\arrow[r]\arrow[d,hook]&K_1(C_r^\ast(\dot{G}_{u,h}))\arrow[d,hook]\\
K_0(C_r^\ast(\dot{G}_{AF}^{(0)}))\arrow[r,"\sim"]\arrow[d,"\CF^{-1}"]&K_0(C_r^\ast(\dot{G}_{AF,v}^{(1)});C_r^\ast(\dot{G}_{u,v}^{(1)}))\arrow[r]&K_0(C_r^\ast(\dot{G}_{AF,v}^{(1)}))\arrow[d,"\CF^{-1}"]\\
C^0\arrow[rr,"\delta_v^0"]&&C_v^1
\end{tikzcd}
\]
\noindent with \(\widetilde{\delta}_v^0\) given by the diagram
\[
\begin{tikzcd}
\kernel\delta_h^0\arrow[rr,hook]\arrow[rrdd,dashed,bend right=30,swap,"\widetilde{\delta}_v^0"]&[-2em]&C^0\arrow[rr,"\delta^0"]&[-2em]&C^1\arrow[rr,"\delta^1"]\arrow[dd,two heads]&[-2em]&C^2&[-2em]\\
[-2em]&c^0\arrow[rr,mapsto]\arrow[rrdd,dashed,mapsto,bend right=30,crossing over]&&c^0\arrow[rr,mapsto,crossing over]&&(0,c_v^1)\arrow[rr,mapsto]\arrow[dd,mapsto]&&0\\
&&\kernel\delta_v^1\arrow[rr,hook]&&C_v^1&&&\\
[-2em]&&&c_v^1\arrow[rr,mapsto]&&c_v^1&&.
\end{tikzcd}
\]
\noindent In words, given \(c^0\in\kernel\delta_h^0\), let \(\delta^0(c^0)=(0,c_v^1)\). \(\delta^1\circ\delta^0(c^0)=\delta^1((0,c_v^1))=\delta_v^1(c_v^1)=0\), so \(c_v^1\in\kernel\delta_v^1\). The map is \(\widetilde{\delta}_v^0(c^0)=c_v^1\).
\begin{remark}
\label{remark:coboundary-induced}
The maps \(\widetilde{\delta}_h^1\) and \(\widetilde{\delta}_v^0\) essentially arise from this diagram that \emph{does not commute}
\[
\begin{tikzcd}
\kernel\delta_\bullet^0\arrow[r,"\widetilde{\delta}_{\bullet^\top}^0"]\arrow[d,hook]&\kernel\delta_{\bullet^\top}^1\arrow[d,hook]\\
C^0\arrow[r,"\delta_{\bullet^\top}^0"]\arrow[d,"\delta_\bullet^0"]\arrow[rd,phantom,"\xcancel{\rotatebox[origin=c]{-90}{\(\circlearrowleft\)}}"]&C_{\bullet^\top}^1\arrow[d,"\delta_{\bullet^\top}^1"]\\
C_\bullet^1\arrow[r,"\delta_\bullet^1"]\arrow[d,two heads]&C^2\arrow[d,two heads]\\
C_\bullet^1/\image\delta_\bullet^0\arrow[r,"\widetilde{\delta}_\bullet^1"]&C^2/\image\delta_{\bullet^\top}^1.
\end{tikzcd}
\]
\end{remark}
\indent Finally, phrased purely in terms of topological objects, the six-term sequence in \cref{equation:six-term-d-2} becomes
\[
\begin{tikzcd}
\cokernel\delta_h^0\arrow[r,"\widetilde{\delta}_h^1"]&\cokernel\delta_v^1\arrow[r]&K_0(C_r^\ast(\dot{G}_u))\arrow[d]\\
K_1(C_r^\ast(\dot{G}_u))\arrow[u]&\kernel\delta_v^1\arrow[l,swap]&\kernel\delta_h^0\arrow[l,swap,"\widetilde{\delta}_v^0"].
\end{tikzcd}
\]
\noindent From the definition of the maps,
\begin{itemize}
\item
\(\cokernel\widetilde{\delta}_h^1\cong C^2/(\image\widetilde{\delta}_h^1+\image\delta_v^1)=C^2/(\image\delta_h^1+\image\delta_v^1)\), and
\item
\(\kernel\widetilde{\delta}_v^0\cong\kernel\delta_h^0\cap\kernel\delta_v^0=\kernel\delta^0\).
\end{itemize}%

\noindent This allows us to break the six-term sequence into two short exact sequences
\begin{gather}
\label{equation:short-exact-k-0}
\begin{tikzcd}[ampersand replacement=\&]
0\arrow[r]\&C^2/(\image\delta_h^1+\image\delta_v^1)\arrow[r,hook]\&K_0(C_r^\ast(\dot{G}_u))\arrow[r,two heads]\&\kernel\delta^0\arrow[r]\&0
\end{tikzcd}\\
\label{equation:short-exact-k-1}
\begin{tikzcd}[ampersand replacement=\&]
0\arrow[r]\&\kernel\delta_v^1/\image\widetilde{\delta}_v^0\arrow[r,hook]\&K_1(C_r^\ast(\dot{G}_u))\arrow[r,two heads]\&\kernel\widetilde{\delta}_h^1\arrow[r]\&0.
\end{tikzcd}
\end{gather}
\noindent \cref{equation:short-exact-k-0} splits because \(\kernel\delta^0\cong\mathbb{Z}\), giving us
\[
K_0(C_r^\ast(\dot{G}_u))\cong C^2/(\image\delta_h^1+\image\delta_v^1)\oplus\kernel\delta^0.
\]
\indent Towards \cref{equation:short-exact-k-1}, we observe that by reversing the roles of \(h\) and \(v\), we obtain the analogous induced maps \(\widetilde{\delta}_h^0\) and \(\widetilde{\delta}_v^1\), giving us the following proposition.
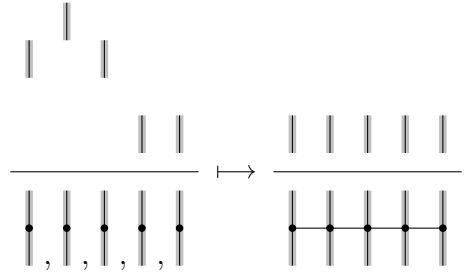
\begin{figure}[t]
\centering
\begin{tikzpicture}
\draw[lightgray,line width=0.1cm](0,1.25)--(0,1.75);
\draw(0,1.25)--(0,1.75);
\draw[lightgray,line width=0.1cm](0.5,1.75)--(0.5,2.25);
\draw(0.5,1.75)--(0.5,2.25);
\draw[lightgray,line width=0.1cm](1,1.25)--(1,1.75);
\draw(1,1.25)--(1,1.75);
\draw[lightgray,line width=0.1cm](1.5,0.25)--(1.5,0.75);
\draw(1.5,0.25)--(1.5,0.75);
\draw[lightgray,line width=0.1cm](2,0.25)--(2,0.75);
\draw(2,0.25)--(2,0.75);
\draw(-0.25,0)--(2.25,0);
\draw[lightgray,line width=0.1cm](0,-1.25)--(0,-0.25);
\draw(0,-1.25)--(0,-0.25);
\draw[lightgray,line width=0.1cm](0.5,-1.25)--(0.5,-0.25);
\draw(0.5,-1.25)--(0.5,-0.25);
\draw[lightgray,line width=0.1cm](1,-1.25)--(1,-0.25);
\draw(1,-1.25)--(1,-0.25);
\draw[lightgray,line width=0.1cm](1.5,-1.25)--(1.5,-0.25);
\draw(1.5,-1.25)--(1.5,-0.25);
\draw[lightgray,line width=0.1cm](2,-1.25)--(2,-0.25);
\draw(2,-1.25)--(2,-0.25);
\fill(0,-0.75)circle[radius=0.05];
\fill(0.5,-0.75)circle[radius=0.05];
\fill(1,-0.75)circle[radius=0.05];
\fill(1.5,-0.75)circle[radius=0.05];
\fill(2,-0.75)circle[radius=0.05];
\draw(0.25,-1.25)node{\(,\)};
\draw(0.75,-1.25)node{\(,\)};
\draw(1.25,-1.25)node{\(,\)};
\draw(1.75,-1.25)node{\(,\)};
\draw[|->](2.5,0)--(3,0);
\draw[lightgray,line width=0.1cm](3.5,0.25)--(3.5,0.75);
\draw(3.5,0.25)--(3.5,0.75);
\draw[lightgray,line width=0.1cm](4,0.25)--(4,0.75);
\draw(4,0.25)--(4,0.75);
\draw[lightgray,line width=0.1cm](4.5,0.25)--(4.5,0.75);
\draw(4.5,0.25)--(4.5,0.75);
\draw[lightgray,line width=0.1cm](5,0.25)--(5,0.75);
\draw(5,0.25)--(5,0.75);
\draw[lightgray,line width=0.1cm](5.5,0.25)--(5.5,0.75);
\draw(5.5,0.25)--(5.5,0.75);
\draw(3.25,0)--(5.75,0);
\draw[lightgray,line width=0.1cm](3.5,-1.25)--(3.5,-0.25);
\draw(3.5,-1.25)--(3.5,-0.25);
\draw[lightgray,line width=0.1cm](4,-1.25)--(4,-0.25);
\draw(4,-1.25)--(4,-0.25);
\draw[lightgray,line width=0.1cm](4.5,-1.25)--(4.5,-0.25);
\draw(4.5,-1.25)--(4.5,-0.25);
\draw[lightgray,line width=0.1cm](5,-1.25)--(5,-0.25);
\draw(5,-1.25)--(5,-0.25);
\draw[lightgray,line width=0.1cm](5.5,-1.25)--(5.5,-0.25);
\draw(5.5,-1.25)--(5.5,-0.25);
\draw(3.5,-0.75)--(5.5,-0.75);
\fill(3.5,-0.75)circle[radius=0.05];
\fill(4,-0.75)circle[radius=0.05];
\fill(4.5,-0.75)circle[radius=0.05];
\fill(5,-0.75)circle[radius=0.05];
\fill(5.5,-0.75)circle[radius=0.05];
\end{tikzpicture}
\caption{An isomorphism \(\kernel(\pi_h^2\circ\delta_v^1)/\image\delta_v^0\xrightarrow{\sim}\kernel\delta_v^1/\image\widetilde{\delta}_v^0\). The \(1\)-cochains in the numerators and denominators are highlighted in light gray. On the left, we are allowed identification across \(\image\delta_v^0\) for each individual \(1\)-cochain. On the right, we identify only if they are joined by \(\delta_v^0(\kernel\delta_h^0)\) (unhighlighted horizontal line).}
\label{figure:cochain-isomorphism}
\end{figure}%

\begin{proposition}
Denote \(\pi_\bullet^2:C^2\twoheadrightarrow C^2/\image\delta_\bullet^1\). Then \(\kernel\widetilde{\delta}_\bullet^1=\kernel(\pi_{\bullet^\top}^2\circ\delta_\bullet^1)/\image\delta_\bullet^0\cong\kernel\delta_\bullet^1/\image\widetilde{\delta}_\bullet^0\).
\begin{proof}
We will prove this for \(\bullet=v\). If \(c_v^1\in\kernel\delta_v^1\), then \(c_v^1\) is a collection of \(1\)-cochains in \(C_v^1\) that are vertical boundaries of \(2\)-cochains closed under horizontal translations.\\
\indent Suppose that \(c_v^1\in\kernel(\pi_h^2\circ\delta_v^1)\). The difference from \(\kernel\delta_v^1\) is that its image in \(C^2\) under \(\delta_v^1\) is allowed identification across \(\image\delta_h^1\). In other words, the \(1\)-cochains in \(c_v^1\) are allowed to be shifted across \(\image\delta_v^0\). \(\kernel\widetilde{\delta}_v^1\) then quotients by identifying across \(\image\delta_v^0\).\\
\indent On the other hand, \(\image\widetilde{\delta}_v^0=\delta_v^0(\kernel\delta_h^0)\), and \(\kernel\delta_h^0\) is the collection of \(0\)-cochains in \(C^0\) that are boundaries of horizontal \(1\)-cochains closed under horizontal translations. That is, \(\image\widetilde{\delta}_v^0\) shifts each \(1\)-cochain in \(c_v^1\in\kernel\delta_v^1\) vertically by the same amount.\\
\indent The isomorphism between the two quotients is as follows: given an element in \(\kernel(\pi_v^2\circ\delta_h^1)/\image\delta_h^0\) and a \(1\)-cochain in it, form a new sum consisting of all vertical \(1\)-cochains that are horizontal translations of it, then take its equivalence class under \(\image\widetilde{\delta}_v^0\) (\cref{figure:cochain-isomorphism}). This is independent of representatives by visual inspection.
\end{proof}
\end{proposition}
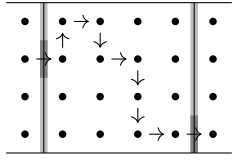
\begin{figure}[t]
\centering
\begin{tikzpicture}
\draw[lightgray,line width=0.1cm](0,0)--(0,2);
\draw[lightgray,line width=0.1cm](2,0)--(2,2);
\draw[gray,line width=0.1cm](0,1)--(0,1.5);
\draw[gray,line width=0.1cm](2,0)--(2,0.5);
\draw(-0.5,0)--(2.5,0);
\draw(-0.5,2)--(2.5,2);
\draw(0,0)--(0,2);
\draw(2,0)--(2,2);
\fill(-0.25,0.25)circle[radius=0.05];
\fill(0.25,0.25)circle[radius=0.05];
\fill(0.75,0.25)circle[radius=0.05];
\fill(1.25,0.25)circle[radius=0.05];
\fill(1.75,0.25)circle[radius=0.05];
\fill(2.25,0.25)circle[radius=0.05];
\fill(-0.25,0.75)circle[radius=0.05];
\fill(0.25,0.75)circle[radius=0.05];
\fill(0.75,0.75)circle[radius=0.05];
\fill(1.25,0.75)circle[radius=0.05];
\fill(1.75,0.75)circle[radius=0.05];
\fill(2.25,0.75)circle[radius=0.05];
\fill(-0.25,1.25)circle[radius=0.05];
\fill(0.25,1.25)circle[radius=0.05];
\fill(0.75,1.25)circle[radius=0.05];
\fill(1.25,1.25)circle[radius=0.05];
\fill(1.75,1.25)circle[radius=0.05];
\fill(2.25,1.25)circle[radius=0.05];
\fill(-0.25,1.75)circle[radius=0.05];
\fill(0.25,1.75)circle[radius=0.05];
\fill(0.75,1.75)circle[radius=0.05];
\fill(1.25,1.75)circle[radius=0.05];
\fill(1.75,1.75)circle[radius=0.05];
\fill(2.25,1.75)circle[radius=0.05];
\draw[->](-0.1,1.25)--(0.1,1.25);
\draw[->](0.25,1.4)--(0.25,1.6);
\draw[->](0.4,1.75)--(0.6,1.75);
\draw[->](0.75,1.6)--(0.75,1.4);
\draw[->](0.9,1.25)--(1.1,1.25);
\draw[->](1.25,1.1)--(1.25,0.9);
\draw[->](1.25,0.6)--(1.25,0.4);
\draw[->](1.4,0.25)--(1.6,0.25);
\draw[->](1.9,0.25)--(2.1,0.25);
\end{tikzpicture}
\caption{A unitary in \(K_1(C_r^\ast(\dot{G}_{u,h}))\) (sequence of arrows on punctures) formed from parallel \(1\)-cochains on \(\Gamma_{2,v}^{(1)}\) (light gray). We can freely move the ``exit points'' (arrows crossing the \(1\)-cochains), and therefore the corresponding \(1\)-cochains on \(\Gamma_{0,v}^{(1)}\) (gray), up and down, provided we do not move to a different \(1\)-cochain on \(\Gamma_{2,v}^{(1)}\). This is the quotient in \(\kernel(\pi_h^2\circ\delta_v^1)/\image\delta_v^0\).}
\label{figure:unitary-cochain}
\end{figure}%

\indent Recall that we have an inductive limit structure on \(K_1(C_r^\ast(\dot{G}_{u,h}))\) using \(\sigma_\square\), by picking a cochain on \(\Gamma_{n,v}^{(1)}\), considering the set of all cochains in it that are horizontal translations of each other, forming the (equivalence class of) unitaries in \(K_1(C_r^\ast(\dot{G}_{u,h}))\) that cross each of these cochains (\cref{figure:unitary-cochain}), then taking the union over all \(n\). In other words, the more appropriate topological object to use in place of the domain of the map \(K_1(C_r^\ast(\dot{G}_{u,h}))/\image\delta_1\rightarrow K_1(C_r^\ast(\dot{G}_u))\) induced by the map \(K_1(C_r^\ast(\dot{G}_{u,h}))\rightarrow K_1(C_r^\ast(\dot{G}_u))\) in \cref{equation:six-term-d-2} is \(\kernel\widetilde{\delta}_v^1\). \cref{equation:short-exact-k-1} then becomes
\[
\begin{tikzcd}
0\arrow[r]&\kernel\widetilde{\delta}_v^1\arrow[r,hook]&K_1(C_r^\ast(\dot{G}_u))\arrow[r,two heads]&\kernel\widetilde{\delta}_h^1\arrow[r]&0.
\end{tikzcd}
\]
\noindent Reversing the roles of \(h\) and \(v\) then gives us the splitting of this short exact sequence, and
\[
K_1(C_r^\ast(\dot{G}_u))\cong\kernel\widetilde{\delta}_v^1\oplus\kernel\widetilde{\delta}_h^1.
\]
\indent This gives us the following theorem.
\begin{figure}[t]
\centering
\begin{tikzpicture}
\draw[lightgray,line width=0.1cm](0,0)--(0,0.5);
\draw[lightgray,line width=0.1cm](0.5,0.5)--(0.5,1);
\draw[lightgray,line width=0.1cm](1,0)--(1,0.5);
\draw(0,0)--(1,0);
\draw(0,0.5)--(1,0.5);
\draw(0,1)--(1,1);
\draw[->](0,0)--(0,0.25);
\draw(0,0.25)--(0,1);
\draw[->](0.5,0)--(0.5,0.75);
\draw(0.5,0.75)--(0.5,1);
\draw[->](1,0)--(1,0.25);
\draw(1,0.25)--(1,1);
\draw(0.25,0.25)node{\(-1\)};
\draw(0.25,0.75)node{\(1\)};
\draw(0.75,0.25)node{\(1\)};
\draw(0.75,0.75)node{\(-1\)};
\draw(1.25,0.5)node{\(+\)};
\draw[lightgray,line width=0.1cm](1.5,0.5)--(2,0.5);
\draw[lightgray,line width=0.1cm](2,0.5)--(2.5,0.5);
\draw(1.5,0)--(2.5,0)--(2.5,1)--(1.5,1)--(1.5,0);
\draw(2,0)--(2,1);
\draw(1.5,0.5)--(1.75,0.5);
\draw[->](2,0.5)--(1.75,0.5);
\draw[->](2,0.5)--(2.25,0.5);
\draw(2.25,0.5)--(2.5,0.5);
\draw(1.75,0.25)node{\(1\)};
\draw(1.75,0.75)node{\(-1\)};
\draw(2.25,0.25)node{\(-1\)};
\draw(2.25,0.75)node{\(1\)};
\draw(2.75,0.5)node{\(=\)};
\draw(3,0)--(4,0)--(4,1)--(3,1)--(3,0);
\draw(3.5,0)--(3.5,1);
\draw(3,0.5)--(4,0.5);
\draw(3.25,0.25)node{\(0\)};
\draw(3.25,0.75)node{\(0\)};
\draw(3.75,0.25)node{\(0\)};
\draw(3.75,0.75)node{\(0\)};
\end{tikzpicture}
\caption{\(\delta_v^2\) applied to an element of \(\kernel(\pi_h^2\circ\delta_v^1)\) that does not belong to \(\kernel\delta_v^1\) (first summand) whose triviality under \(\pi_h^2\circ\delta_v^1\) is by adding by elements of \(\image\delta_h^1\) (second summand) that turn out to belong to \(\image\delta_h^0\). The \(1\)-cochains are highlighted in light gray, and the signs are indicated by the arrows.}
\label{figure:pi-h-2-trivial}
\end{figure}%

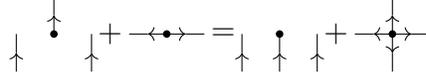
\begin{figure}[t]
\centering
\begin{tikzpicture}

\draw[->](0,0)--(0,.25);
\draw(0,.25)--(0,.5);
\draw[->](.5,.5)--(.5,.75);
\draw(.5,.75)--(.5,1);
\draw[->](1,0)--(1,.25);
\draw(1,.25)--(1,.5);
\fill(.5,.5)circle[radius=0.05];
\draw(1.25,.5)node{\(+\)};
\draw(1.5,.5)--(1.75,.5);
\draw[->](2,.5)--(1.75,.5);
\draw[->](2,.5)--(2.25,.5);
\draw(2.25,.5)--(2.5,.5);
\fill(2,.5)circle[radius=0.05];
\draw(2.75,0.5)node{\(=\)};
\draw[->](3,0)--(3,.25);
\draw(3,.25)--(3,.5);
\draw[->](3.5,0)--(3.5,.25);
\draw(3.5,.25)--(3.5,.5);
\draw[->](4,0)--(4,.25);
\draw(4,.25)--(4,.5);
\fill(3.5,.5)circle[radius=0.05];
\draw(4.25,.5)node{\(+\)};
\draw(4.5,.5)--(4.75,.5);
\draw[->](5,.5)--(4.75,.5);
\draw[->](5,.5)--(5.25,.5);
\draw(5.25,.5)--(5.5,.5);
\draw(5,0)--(5,.25);
\draw[->](5,.5)--(5,.25);
\draw[->](5,.5)--(5,.75);
\draw(5,.75)--(5,1);
\fill(5,.5)circle[radius=0.05];
\end{tikzpicture}
\caption{\(\kernel(\pi_h^2\circ\delta_v^1)/\image\delta_v^0\) (left, first summand) added by \(\image\delta_h^0\) (left, second summand) gives a sum of an element of \(\kernel\delta_v^1\) (right, first summand) and an element of \(\image\delta^0\) (right, second summand), and vice versa (up to elements in \(\image\delta_v^0\) and \(\image\delta_h^0\)). The signs on the \(1\)-cochains are indicated by the arrows.}
\label{figure:delta-v-1-trivial}
\end{figure}%

\begin{theorem}
For \(d=2\), \(K_0(C_r^\ast(\dot{G}_u))\cong\check{H}^2(\Omega_{T_\square})\oplus\check{H}^0(\Omega_{T_\square})\) and \(K_1(C_r^\ast(\dot{G}_u))\cong\check{H}^1(\Omega_{T_\square})\).
\begin{proof}
It remains to show that \(\check{H}^1(\Omega_{T_\square})=\kernel\widetilde{\delta}_v^1\oplus\kernel\widetilde{\delta}_h^1\). We will show this for \(\bullet=v\).\\
\indent Let us begin with the observation that there is an interaction between elements that go to zero only at \(\pi_h^2\) in \(\kernel(\pi_h^2\circ\delta_v^1)\) and elements of \(\image\delta_h^0\). Suppose that \(\sum_ic_{v,i}^1\in\kernel\delta_v^1\) and for some \(j\) and some \(0\)-cochain \(c^0\), \(\delta_v^0(c^0)=c_{v,j}^1-{c_{v,j}^1}'\). Then \(\sum_{i\neq j}c_v^1+{c_{v,j}^1}'\in\kernel(\pi_h^2\circ\delta_v^1)\) but \(\sum_{i\neq j}c_v^1+{c_{v,j}^1}'\notin\kernel\delta_v^1\). In particular, the quotient from \(\pi_h^2\) is by two pairs of neighboring \(2\)-cochains that arise from some \(\delta_h^1(c_{h,j_1}^1)\) and \(\delta_h^1(c_{h,j_2}^1)\), and the difference of the two horizontal \(1\)-cochains is \(\delta_h^0(c^0)\). In other words, changing representatives of elements in \(\kernel(\pi_h^2\circ\delta_v^1)/\image\delta_v^0\) implicitly includes elements of \(\image\delta_h^0\) (\cref{figure:pi-h-2-trivial}).\\
\indent Furthermore, we notice that
\[
\delta^0(c^0)=c_{v,j}^1-{c_{v,j}^1}'+c_{h,j_1}^1-c_{h,j_2}^1,
\]
\noindent i.e., this changing of representatives can be further realized as adding elements of \(\image\delta^0\) (\cref{figure:delta-v-1-trivial}).\\
\indent From these observations, \(\kernel\widetilde{\delta}_v^1\oplus\kernel\widetilde{\delta}_h^1\leq\check{H}^1(\Omega_{T_\square})\).\\
\indent For the converse, let us consider the inductive limit structure on \(\check{H}^1(\Omega_{T_\square})\) afforded by \(\sigma_\square\), and recall that elements of \(\kernel\delta^1\) correspond to collections of closed curves intersecting transversely with the \(1\)-cochains, with the oriented intersection numbers coinciding with the values assigned by the \(1\)-cochains (\cite[pp.~188--189]{hatcher02}). Given a closed curve in \(\Gamma_n\) as described, let \(c_n^1=(c_{n,v}^1,c_{n,h}^1)\) denote the collection of \(1\)-cochains it intersects. Being a closed curve, \(c_{n,v}^1\) is formed from the vertical boundaries of a collection of \(2\)-cells that are neighboring, up to vertical translations. That is, \(c_{n,v}^1\in\kernel\widetilde{\delta}_{n,v}^1\). Same holds for \(c_{n,h}^1\). Thus \(\check{H}^1(\Omega_{T_\square})\leq\kernel\widetilde{\delta}_v^1\oplus\kernel\widetilde{\delta}_h^1\).
\end{proof}
\end{theorem}
\noindent We can then compose these isomorphisms with the map \(\square^{-1}\) induced on cohomology to get the generators and relations in terms of the cochains and coboundary maps in \(\Omega_T\), from which we see that the isomorphism \(K_1(C_r^\ast(\dot{G}_u))\cong\check{H}^1(\Omega_T)\) is induced by thickening.\\
\indent We now revert to groupoids on \(\Omega_T\). Using this theorem and \cref{theorem:k-theory-cohomology}, we can concisely package the relevant topological objects in the context of relative \(K\)-theory into the following corollary that is a generalization of \cref{theorem:isomorphism-d-1} to \(d=2\).
\begin{corollary}[\cref{theorem:isomorphism-d-2}]
For \(d=2\), we have an isomorphism of exact sequences
\[
\begin{tikzcd}
0\arrow[r]&\check{H}^1(\Omega_T)\arrow[r,hook]\arrow[ddd]&C^1/\image\delta^0\arrow[r,"\delta^1"]\arrow[ddd,"\thickening"]&C^2\arrow[r,two heads]&\check{H}^2(\Omega_T)\arrow[r]&0&\\[-28pt]
&&&\oplus&\oplus&\oplus&\\[-28pt]
&&\vphantom{C^1/\image\delta^0}&0\arrow[r]\arrow[d,"\CF"]&\check{H}^0(\Omega_T)\arrow[r,hook,two heads]\arrow[d]&\check{H}^0(\Omega_T)\arrow[r]\arrow[d]&0\\
0\arrow[r]&K_{1,u}\arrow[r,hook]&K_{0,AF;u}\arrow[r,"\evaluation"]&K_{0,AF}\arrow[r,"\iota_\ast"]&K_{0,u}\arrow[r,two heads]&K_{1,AF;u}\arrow[r]&0,
\end{tikzcd}
\]
\noindent where \(K_{i,u}=K_i(C_r^\ast(\dot{G}_u))\), \(K_{i,AF}=K_i(C_r^\ast(\dot{G}_{AF}))\), and \(K_{i,AF;u}=K_i(C_r^\ast(\dot{G}_{AF});C_r^\ast(\dot{G}_u))\). Thus, if \(T\) satisfies the boundary hyperplane condition,
\begin{align*}
K_0(C_r^\ast(\dot{G}_u))&\cong K_0(C_r^\ast(\dot{G}_{AF}))/\image\evaluation\oplus\mathbb{Z}\\
K_1(C_r^\ast(\dot{G}_u))&\cong\kernel[K_0(C_r^\ast(\dot{G}_{AF}^{(1)}))/\thickening(K_0(C_r^\ast(\dot{G}_{AF}^{(0)})))\xrightarrow{\evaluation} K_0(C_r^\ast(\dot{G}_{AF}))]
\end{align*}
\noindent where \(\thickening:K_0(C_r^\ast(\dot{G}_{AF}^{(0)}))\rightarrow K_0(C_r^\ast(\dot{G}_{AF}^{(1)}))\) returns an alternating sum of elements of \(K_0(C_r^\ast(\dot{G}_{AF}^{(1)}))\) that coincides with \(\delta^0\).
\end{corollary}
\noindent Therefore for \(d=2\), the \(K\)-theory of the unstable groupoid can be reconstructed from the \(K\)-theory of the \(AF\)-groupoids if \(T\) satisfies the boundary hyperplane condition. This includes all triangular substitutions where the tiles are aligned edge-to-edge.
\begin{remark}
The right-hand side of the isomorphism of exact sequences recovers the well-known fact from \cite{kellendonk97} that the contribution of \(K_0(C_r^\ast(\dot{G}_{AF}))\) towards \(K_0(C_r^\ast(\dot{G}_u))\) is as the top-dimensional cohomology, also called the \emph{integer group of coinvariants}.
\end{remark}
\begin{warning}
If \(T\) does not satisfy the boundary hyperplane condition, \(\sigma^\top:C_n^1\rightarrow C_{n+1}^1\) may carry signs that \emph{cannot be made nonnegative} through reorienting \(1\)-cells!\\
\indent For example, in the chair tiling, consider the two \(1\)-cells (light gray with neighboring tiles drawn)
\[
\begin{tikzpicture}
\draw[lightgray,line width=0.1cm](0,0.5)--(0,0)--(0.5,0);
\draw(0,.5)--(0,1.0);
\draw(.5,0)--(1.0,0);
\draw(1.0,0)--(1.0,.5);
\draw[->](0,.5)--(0,.25);
\draw(0,.25)--(0,0);
\draw(.5,.5)--(1.0,.5);
\draw(-.5,0)--(-.5,-.5);
\draw[->](0,0)--(.25,0);
\draw(.5,0)--(.25,0);
\draw(0,1.0)--(.5,1.0);
\draw(0,-.5)--(-.5,-.5);
\draw(0,-.5)--(.5,-.5);
\draw(.5,0)--(.5,-.5);
\draw(-.5,0)--(-.5,.5);
\draw(0,.5)--(-.5,.5);
\draw(.5,.5)--(.5,1.0);

\draw[lightgray,line width=0.1cm](2.5,0)--(3.0,0)--(3.0,.5);
\draw(3.0,.5)--(3.0,1.0);
\draw[->](3.0,0)--(3.0,.25);
\draw(3.0,.5)--(3.0,.25);
\draw(2.5,.5)--(2.5,1.0);
\draw(3.5,0)--(3.5,.5);
\draw(2.5,0)--(2.5,-.5);
\draw(3.0,1.0)--(2.5,1.0);
\draw(3.0,.5)--(3.5,.5);
\draw[->](2.5,0)--(2.75,0);
\draw(2.75,0)--(3.0,0);
\draw(2.5,.5)--(2.0,.5);
\draw(2.0,0)--(2.5,0);
\draw(3.5,0)--(3.5,-.5);
\draw(3.0,-.5)--(2.5,-.5);
\draw(3.0,-.5)--(3.5,-.5);
\draw(2.0,0)--(2.0,.5);
\end{tikzpicture}
\]
with the prescribed orientations. Note that one cannot further separate each of the two \(1\)-cells into smaller \(1\)-cells, since there is only a single partial isometry associated to each\footnote{Rather, one can, but doing so amounts to passing to a different cellular decomposition of the tiling, therefore a different substitution rule.}.\\
\indent On applying sufficiently many substitutions (\(3\) or more for the chair tiling), we observe the two \(1\)-cells
\[
\begin{tikzpicture}
\draw[lightgray,line width=0.1cm](1.5,0)--(1.5,1.0);
\draw(1.5,1.0)--(1.5,.5);
\draw(1.0,1.0)--(1.0,.5);
\draw(1.5,0)--(1.5,.5);
\draw(1.5,0)--(2.0,0);
\draw(.5,.5)--(1.0,.5);
\draw(1.0,0)--(.5,0);
\draw(2.5,0)--(2.0,0);
\draw(2.0,1.0)--(1.5,1.0);
\draw(1.0,1.0)--(1.5,1.0);
\draw(1.5,0)--(1.0,0);
\draw(2.5,.5)--(2.0,.5);
\draw(2.5,0)--(2.5,.5);
\draw(.5,0)--(.5,.5);
\draw(2.0,1.0)--(2.0,.5);

\draw[lightgray,line width=0.1cm](3.5,.5)--(4.5,.5);
\draw(4.5,0)--(4.0,0);
\draw(3.5,1.0)--(3.5,1.5);
\draw(4.0,.5)--(4.5,.5);
\draw(4.5,.5)--(4.5,1.0);
\draw(3.5,1.0)--(3.5,.5);
\draw(4.0,1.0)--(4.5,1.0);
\draw(4.0,0)--(4.0,-.5);
\draw(3.5,0)--(3.5,-.5);
\draw(4.5,0)--(4.5,.5);
\draw(3.5,.5)--(4.0,.5);
\draw(3.5,1.5)--(4.0,1.5);
\draw(3.5,-.5)--(4.0,-.5);
\draw(3.5,0)--(3.5,.5);
\draw(4.0,1.0)--(4.0,1.5);
\end{tikzpicture}
\]
in the approximate positions
\[
\begin{tikzpicture}
\draw[lightgray,line width=0.1cm](0,0.5)--(0,0)--(0.5,0);
\draw(0,.5)--(0,1.0);
\draw(.5,0)--(1.0,0);
\draw(1.0,0)--(1.0,.5);
\draw[->](0,.5)--(0,.25);
\draw(0,.25)--(0,0);
\draw(.5,.5)--(1.0,.5);
\draw(-.5,0)--(-.5,-.5);
\draw[->](0,0)--(.25,0);
\draw(.5,0)--(.25,0);
\draw(0,1.0)--(.5,1.0);
\draw(0,-.5)--(-.5,-.5);
\draw(0,-.5)--(.5,-.5);
\draw(.5,0)--(.5,-.5);
\draw(-.5,0)--(-.5,.5);
\draw(0,.5)--(-.5,.5);
\draw(.5,.5)--(.5,1.0);
\draw[|->](1.25,0)--(1.75,0)node[midway,above]{\(\varsigma^n\)};
\draw[lightgray,line width=0.1cm](3,2.5)--(3,0)--(5.5,0);
\draw(3,1.5)--(3,0)--(4.5,0);
\draw[->](3.0,2.5)--(3.0,2.0);
\draw(2.5,2.5)--(2.5,2.0);
\draw(3.0,1.5)--(3.0,2.0);
\draw(3.0,1.5)--(3.5,1.5);
\draw(2.0,2.0)--(2.5,2.0);
\draw(2.5,1.5)--(2.0,1.5);
\draw(4.0,1.5)--(3.5,1.5);
\draw(3.5,2.5)--(3.0,2.5);
\draw(2.5,2.5)--(3.0,2.5);
\draw(3.0,1.5)--(2.5,1.5);
\draw(4.0,2.0)--(3.5,2.0);
\draw(4.0,1.5)--(4.0,2.0);
\draw(2.0,1.5)--(2.0,2.0);
\draw(3.5,2.5)--(3.5,2.0);
\draw(5.5,-.5)--(5.0,-.5);
\draw(4.5,.5)--(4.5,1.0);
\draw(5.0,0)--(5.5,0);
\draw(5.5,0)--(5.5,.5);
\draw(4.5,.5)--(4.5,0);
\draw(5.0,.5)--(5.5,.5);
\draw(5.0,-.5)--(5.0,-1.0);
\draw(4.5,-.5)--(4.5,-1.0);
\draw(5.5,-.5)--(5.5,0);
\draw[->](4.5,0)--(5.0,0);
\draw(4.5,1.0)--(5.0,1.0);
\draw(4.5,-1.0)--(5.0,-1.0);
\draw(4.5,-.5)--(4.5,0);
\draw(5.0,.5)--(5.0,1.0);

\draw[lightgray,line width=0.1cm](7,0)--(7.5,0)--(7.5,0.5);
\draw(7.5,.5)--(7.5,1.0);
\draw[->](7.5,0)--(7.5,.25);
\draw(7.5,.5)--(7.5,.25);
\draw(7.0,.5)--(7.0,1.0);
\draw(8.0,0)--(8.0,.5);
\draw(7.0,0)--(7.0,-.5);
\draw(7.5,1.0)--(7.0,1.0);
\draw(7.5,.5)--(8.0,.5);
\draw[->](7.0,0)--(7.25,0);
\draw(7.25,0)--(7.5,0);
\draw(7.0,.5)--(6.5,.5);
\draw(6.5,0)--(7.0,0);
\draw(8.0,0)--(8.0,-.5);
\draw(7.5,-.5)--(7.0,-.5);
\draw(7.5,-.5)--(8.0,-.5);
\draw(6.5,0)--(6.5,.5);
\draw[|->](8.25,0)--(8.75,0)node[midway,above]{\(\varsigma^n\)};
\draw[lightgray,line width=0.1cm](9,0)--(11.5,0)--(11.5,2.5);
\draw(10,0)--(11.5,0)--(11.5,1.5);
\draw(11.5,2.0)--(11.5,2.5);
\draw(11.0,2.5)--(11.0,2.0);
\draw[->](11.5,1.5)--(11.5,2.0);
\draw(11.5,1.5)--(12.0,1.5);
\draw(10.5,2.0)--(11.0,2.0);
\draw(11.0,1.5)--(10.5,1.5);
\draw(12.5,1.5)--(12.0,1.5);
\draw(12.0,2.5)--(11.5,2.5);
\draw(11.0,2.5)--(11.5,2.5);
\draw(11.5,1.5)--(11.0,1.5);
\draw(12.5,2.0)--(12.0,2.0);
\draw(12.5,1.5)--(12.5,2.0);
\draw(10.5,1.5)--(10.5,2.0);
\draw(12.0,2.5)--(12.0,2.0);
\draw(10.0,-.5)--(9.5,-.5);
\draw(9.0,.5)--(9.0,1.0);
\draw(9.5,0)--(10.0,0);
\draw(10.0,0)--(10.0,.5);
\draw(9.0,.5)--(9.0,0);
\draw(9.5,.5)--(10.0,.5);
\draw(9.5,-.5)--(9.5,-1.0);
\draw(9.0,-.5)--(9.0,-1.0);
\draw(10.0,-.5)--(10.0,0);
\draw[->](9.0,0)--(9.5,0);
\draw(9.0,1.0)--(9.5,1.0);
\draw(9.0,-1.0)--(9.5,-1.0);
\draw(9.0,-.5)--(9.0,0);
\draw(9.5,.5)--(9.5,1.0);
\end{tikzpicture}
\]
in the substituted \(1\)-cells.\\
\indent Regardless of the orientations assigned to the horizontal and vertical \(1\)-cells, one of them will always have its orientation reversed with respect to the orientation of one of the substituted \(1\)-cells, resulting in a sign on the induced substitution matrix. E.g. having the horizontal and vertical \(1\)-cells oriented left-to-right and down-to-up results in the substitution on the left \(1\)-cell carrying a \(-1\) on the vertical \(1\)-cell.\\
\indent This is the reason we cared little for the orientations on the \(1\)-cells not satisfying the boundary hyperplane condition. Thus we may no longer have a map \(K_0(C_r^\ast(\dot{G}_{AF,n}^{(1)}))\rightarrow K_0(C_r^\ast(\dot{G}_{AF,n+1}^{(1)}))\) that is both compatible with the map on cochains and induced by a \(C^\ast\)-algebra map \(C_r^\ast(\dot{G}_{AF,n}^{(1)})\hookrightarrow C_r^\ast(\dot{G}_{AF,n+1}^{(1)})\).\\
\indent In brief, maps on cochains are always allowed signs since they are maps on free abelian groups, but that is not the case for \(C^\ast\)-algebra maps.
\end{warning}
\indent Despite this, \(K_0(C_r^\ast(\dot{G}_{AF}^{(1)}))\) and \(C^1\) may still be isomorphic\footnote{This isomorphism may not arise from shift equivalence.}! Let us call a direct limit under a matrix all of whose eigenvalues are integral \emph{completely split} if it is isomorphic to \(\bigoplus_\lambda\mathbb{Z}[1/\lambda]\) with \(\lambda\) the eigenvalues of the matrix.
\begin{proposition}
\label{proposition:af-cochain-nontrivial}
Let \(H\) be the set of \(1\)-cells satisfying the boundary hyperplane condition, and let \(\sigma^\top|_H\) be the induced substitution matrix restricted to \(H\). If \(\sigma^\top|_H\) has all of its eigenvalues integral and its direct limit splits completely, then \(K_0(C_r^\ast(\dot{G}_{AF}^{(1)}))\cong C^1\).
\begin{proof}
Let \(H^c\) be the complement of \(H\) in the set of \(1\)-cells. Elements of \(H\) cannot subdivide to elements of \(H^c\), so \(\sigma^\top\), written as a \(2\times 2\) block matrix, has the form
\[
\left(
\begin{array}{cc}
\sigma^\top|_{H^c}&A\\
0&\sigma^\top|_H
\end{array}
\right).
\]
\indent The substitution map, restricted to \(H^c\), is necessarily nonexpanding. Therefore, after removing the eventual kernel, \(\sigma^\top|_{H^c}\) is a permutation matrix. By applying a sufficiently-high power (and removing the eventual kernel), we will assume \(\sigma^\top|_{H^c}=I\). Letting \(B=\sigma^\top|_H\), our block matrix has the form
\[
\left(
\begin{array}{cc}
I&A\\
0&B
\end{array}
\right).
\]
\indent Noting that eigenvectors of unit eigenvalues contribute copies of \(\mathbb{Z}\) in the limit, suppose that all eigenvalues of \(B\) are nonzero and nonunit. Let \(v\) be an eigenvector for an eigenvalue \(\lambda\). Setting \(w=Av/(\lambda-1)\),
\begin{align*}
\left(
\begin{array}{cc}
I&A\\
0&B
\end{array}
\right)
\left(
\begin{array}{c}
w\\v
\end{array}
\right)&=
\left(
\begin{array}{c}
w+Av\\
\lambda v
\end{array}
\right)\\
&=
\left(
\begin{array}{c}
w+(\lambda-1)w\\
\lambda v
\end{array}
\right)\\
&=\lambda
\left(
\begin{array}{c}
w\\
v
\end{array}
\right),
\end{align*}
\noindent and we see that each such eigenvector \(v\) of \(B\) can be used to construct an eigenvector for the block matrix.\\
\indent To compute the direct limit under this block matrix, we are looking for the set of vectors \(c\in\mathbb{Z}^{|H^c|+|H|}\) so that
\[
\sum_{i=1}^{|H^c|}c_ie_i+\sum_{j=|H^c|+1}^{|H^c|+|H|}c_j\lambda_j^k\left(\begin{array}{c}w_j\\v_j\end{array}\right)\in\mathbb{Z}^{|H^c|+|H|}.
\]
\noindent By assumption, since the direct limit under \(B\) splits completely, each \(c_j\lambda_j^kw_j\) is integral. Thus, each of the \(c_i\) is integral (up to a fraction), and the direct limit of the block matrix also splits completely.\\
\indent Repeating the same argument with \(A\) replaced to represent the matrix induced by \(C_r^\ast(\dot{G}_{AF,n}^{(1)})\hookrightarrow C_r^\ast(\dot{G}_{AF,n+1}^{(1)})\), denoted \(A'\), yields our desired isomorphism, given by the map
\[
\left(
\begin{array}{c}
w_j\\
v_j
\end{array}
\right)
\mapsto
\left(
\begin{array}{c}
w_j'\\
v_j
\end{array}
\right)
\]
\noindent where \(w'=A'v/(\lambda-1)\).
\end{proof}
\end{proposition}
\begin{figure}[t]
\centering
\begin{tikzpicture}
\draw[lightgray,line width=0.1cm](-0.3833,0.6640)--(0,0);
\draw[lightgray,line width=0.1cm](0,0)--(0.7667,0);
\draw[->](-0.6667,1.1547)--(-0.5,0.8660);
\draw(-0.5,0.8660)--(0,0);
\draw[->](-0.6667,-1.1547)--(-0.5,-0.8660);
\draw(-0.5,-0.8660)--(0,0);
\draw[->](0,0)--(1,0);
\draw(1,0)--(1.3333,0);
\fill(-0.3333,0.5774)circle[radius=0.05];
\fill(0.6667,0)circle[radius=0.05];
\draw[lightgray,line width=0.1cm](2.6167,-0.6640)--(3,0);
\draw[lightgray,line width=0.1cm](2.6167,0.6640)--(3,0);
\draw[lightgray,line width=0.1cm](3,0)--(3.7667,0);
\draw[->](2.3333,1.1547)--(2.5,.8660);
\draw(2.5,.8660)--(3,0);
\draw[->](2.3333,-1.1547)--(2.5,-.8660);
\draw(2.5,-.8660)--(3,0);
\draw[->](3,0)--(4,0);
\draw(4,0)--(4.3333,0);
\fill(2.6667,-.5774)circle[radius=0.05];
\fill(2.6667,.5774)circle[radius=0.05];
\fill(3.6667,0)circle[radius=0.05];
\end{tikzpicture}
\caption{A partial homeomorphism yielding a partial isometry, as given by restriction (left), compared to \(\image\delta^0\) (right). We expect to quotient the set of all projections on the skeleton by the former, but the latter is what we actually quotient by. This partial isometry cannot be realized as an element of \(\image\delta^0\).}
\label{figure:k-theory-skeleton}
\end{figure}
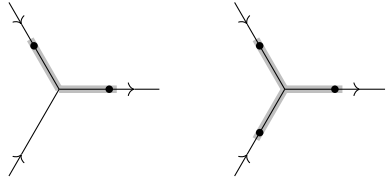%

\begin{warning}
In general, \emph{it does not appear that \(K_0(C_r^\ast(\dot{G}_{AF});C_r^\ast(\dot{G}_u))\) arises as the \(K_0\)-group of the \(d-1\)-skeleton!} This is true for \(d=1\), as we saw in the proof of \cref{theorem:isomorphism-d-1}, where it arises from an \(AF\)-algebra. However, for \(d=2\), we possibly have \(AF\)-equivalence classes that arise as the intersections of three or more infinite \(1\)-boundaries. If \(K_0(C_r^\ast(\dot{G}_{AF});C_r^\ast(\dot{G}_u))\cong C^1/\image\delta^0\) arises as the \(K_0\)-group of \(\varprojlim_n(AP_n^{(1)},\sigma)\), then we expect it to come from restriction.\\
\indent On the one hand, there is a partial isometry containing only two of the \(1\)-boundaries. On the other, \(\delta^0\) applied to this intersection gives a single signed sum of involving all three infinite \(1\)-boundaries. In particular, \emph{\(\image\delta^0\) must assign nonzero values to all incident \(1\)-cells,} whereas (generating) \emph{partial isometries select two \(1\)-cells to assign nonzero values.} Under Murray--von Neumann equivalence, we expect the latter to be the relation we quotient by in the \(K_0\)-group of the \(1\)-skeleton, but the relation we have is the former, and we cannot obtain all of the relations of the latter form from the former (\cref{figure:k-theory-skeleton}).\\
\indent However, in certain situations (e.g. \cite{juliensavinien16}), there exist substitution tiling spaces with stationary substitution rules within the same MLD-equivalence class that have a subset of its \(d-1\)-skeleton that intersects all other infinite boundaries and is closed under the induced group action by translation. This subset is a \(d-1\)-dimensional tiling space, and so the relative \(K\)-theory relative to this subgroupoid is isomorphic to the \(K\)-theory of this new groupoid.
\end{warning}
\begin{remark}
Since the Chern and Connes--Thom isomorphisms hold over \(\mathbb{Z}\) for \(d\leq 3\), we expect that the analogous theorems hold for \(d=3\) as well, at least abstractly,
\[
\begin{tikzcd}
0\arrow[r]&\check{H}^2(\Omega_T)\arrow[r,hook]&C^2/\image\delta^1\arrow[r,"\delta^2"]&C^3\arrow[r,two heads]&\check{H}^3(\Omega_T)\arrow[r]&0&\\[-28pt]
&\oplus&\oplus&\oplus&\oplus&\oplus&\\[-28pt]
0\arrow[r]&\check{H}^0(\Omega_T)\arrow[r,hook,two heads]\arrow[d]&\check{H}^0(\Omega_T)\arrow[r]\arrow[d,"\thickening\oplus\phi"]&0\arrow[r]\arrow[d,"\CF"]&\check{H}^1(\Omega_T)\arrow[r,hook,two heads]\arrow[d]&\check{H}^1(\Omega_T)\arrow[r]\arrow[d]&0\\
0\arrow[r]&K_{1,u}\arrow[r,hook]&K_{0,AF;u}\arrow[r,"\evaluation"]&K_{0,AF}\arrow[r,"\iota_\ast"]&K_{0,u}\arrow[r,two heads]&K_{1,AF;u}\arrow[r]&0,
\end{tikzcd}
\]
\noindent where assuming the boundary hyperplane condition and using \cref{lemma:af-cochain} gives the contributions of each of the \(AF\)-groupoids towards the \(K\)-theory of the unstable groupoid after undoing the isomorphisms.\\
\indent Unfortunately, as it is for \(d=2\), the component of \(\check{H}^0(\Omega_T)\) in \(K_1(C_r^\ast(\dot{G}_u))\), denoted by \(\phi\), is harder to identify. By writing \(C_r^\ast(\dot{G}_u)\) as a crossed product via \cite[Lemma 11]{kellendonk97}, peeling off one dimension, then applying the Connes--Thom isomorphism to this single dimension gives the \(K_0\)-group of a \(\mathbb{Z}^2\)-crossed product. The reasoning provided in \cite{kellendonkputnam00} results in a rank-\(1\) element that cannot be homotopic to the trivial rank-\(1\) projection. This should be the element seen by \(\phi\).
\end{remark}

\subsection{Interpretation}
\indent \cref{theorem:isomorphism-d-1} and \cref{theorem:isomorphism-d-2} allow us to give a very concise description of the six-term sequence in relative \(K\)-theory of \(\iota:C_r^\ast(\dot{G}_{AF})\hookrightarrow C_r^\ast(\dot{G}_u)\) for \(d\leq 2\). For motivation, let us first recall an alternate way to record a substitution rule with an example.
\begin{example}[Fibonacci]
The square of the collared Fibonacci substitution is
\begin{align*}
A&\mapsto (AD)AB(C)\\
B&\mapsto (B)CAB(C)\\
C&\mapsto (AB)CAD(A)\\
D&\mapsto (B)CAD(A)
\end{align*}
\noindent with the \(AP\)-complex
\[
\begin{tikzpicture}
\draw[->](0,0)--(0.75,0);
\draw[-](0.75,0)--(1.5,0);
\draw(0,0)--(1.5,0)node[midway,above]{\(A\)};
\draw[->](1.5,0)arc(-30:-90:0.8660);
\draw[-](0,0)arc(-150:-90:0.8660);
\draw(1.5,0)arc(-30:-150:0.8660)node[midway,below]{\(D\)};
\draw[->](1.5,0)--(1.125,0.75);
\draw[-](1.125,0.75)--(0.75,1.5);
\draw(1.5,0)--(0.75,1.5)node[midway,above,rotate=-63.4349]{\(B\)};
\draw[->](0.75,1.5)--(0.375,0.75);
\draw[-](0.375,0.75)--(0,0);
\draw(0,0)--(0.75,1.5)node[midway,above,rotate=63.4349]{\(C\)};
\end{tikzpicture}.
\]
\noindent We can condense these two pieces of information into a single picture
\[
\begin{tikzpicture}
\draw[->](0,0)--(1.50,0);
\draw[-](1.50,0)--(3.0,0);
\draw(0,0)--(1.50,0)node[midway,below,yshift=2pt]{\fontsize{7}{7}\(A\)};
\draw(1.50,0)--(3.0,0)node[midway,below,yshift=2pt]{\fontsize{7}{7}\(B\)};
\draw(0,0)--(3,0)node[midway,above]{\(A\)};
\draw[->](3.0,0)arc(-30:-90:1.7321);
\draw[-](0,0)arc(-150:-90:1.7321);
\draw(3.0,0)arc(-30:-70:1.7321)node[midway,above,rotate=45,yshift=-2pt]{\fontsize{7}{7}\(C\)};
\draw(2.0924,-.7616)arc(-70:-110:1.7321)node[midway,above,yshift=-2pt]{\fontsize{7}{7}\(A\)};
\draw(.9076,-.7616)arc(-110:-150:1.7321)node[midway,above,rotate=-45,yshift=-2pt]{\fontsize{7}{7}\(D\)};
\draw(3.0,0)arc(-30:-150:1.7321)node[midway,below]{\(D\)};
\draw[->](3.0,0)--(2.250,1.50);
\draw[-](2.250,1.50)--(1.50,3.0);
\draw(3.0,0)--(2.50,1.0)node[midway,below,rotate=-63.4349,yshift=2pt]{\fontsize{7}{7}\(C\)};
\draw(2.50,1.0)--(2,2)node[midway,below,rotate=-63.4349,yshift=2pt]{\fontsize{7}{7}\(A\)};
\draw(2,2)--(1.50,3.0)node[midway,below,rotate=-63.4349,yshift=2pt]{\fontsize{7}{7}\(B\)};
\draw(3,0)--(1.5,3)node[midway,above,rotate=-63.4349]{\(B\)};
\draw[->](1.50,3.0)--(.750,1.50);
\draw[-](.750,1.50)--(0,0);
\draw(0,0)--(.50,1.0)node[midway,below,rotate=63.4349,yshift=2pt]{\fontsize{7}{7}\(D\)};
\draw(.50,1.0)--(1.0,2)node[midway,below,rotate=63.4349,yshift=2pt]{\fontsize{7}{7}\(A\)};
\draw(1.0,2)--(1.50,3.0)node[midway,below,rotate=63.4349,yshift=2pt]{\fontsize{7}{7}\(C\)};
\draw(1.5,3)--(0,0)node[midway,above,rotate=63.4349]{\(C\)};
\end{tikzpicture}
\]
\noindent where the smaller letters are labels we further attach to each of the top-dimensional cells to indicate where they are sent under \(\sigma\). E.g.
\[
\begin{tikzpicture}
\draw[->,dashed](0,0)--(1.50,0);
\draw[-,dashed](1.50,0)--(3.0,0);
\draw[->,dashed](3.0,0)arc(-30:-90:1.7321);
\draw[-,dashed](0,0)arc(-150:-90:1.7321);
\draw[->,dashed](3.0,0)--(2.250,1.50);
\draw[-,dashed](2.250,1.50)--(1.50,3.0);
\draw[->](1.50,3.0)--(.750,1.50);
\draw[-](.750,1.50)--(0,0);
\draw(0,0)--(.50,1.0)node[midway,below,rotate=63.4349,yshift=2pt]{\fontsize{7}{7}\(D\)};
\draw(.50,1.0)--(1.0,2)node[midway,below,rotate=63.4349,yshift=2pt]{\fontsize{7}{7}\(A\)};
\draw(1.0,2)--(1.50,3.0)node[midway,below,rotate=63.4349,yshift=2pt]{\fontsize{7}{7}\(C\)};
\draw(1.5,3)--(0,0)node[midway,above,rotate=63.4349]{\(C\)};
\draw[|->](3.25,0)--(4.25,0)node[midway,above]{\(\sigma^2\)};
\draw[->](4.5,0)--(5.25,0);
\draw[-](5.25,0)--(6.0,0);
\draw(4.5,0)--(6.0,0)node[midway,above]{\(A\)};
\draw[->](6.0,0)arc(-30:-90:0.8660);
\draw[-](4.5,0)arc(-150:-90:0.8660);
\draw(6.0,0)arc(-30:-150:0.8660)node[midway,below]{\(D\)};
\draw[->,dashed](6.0,0)--(5.625,.75);
\draw[-,dashed](5.625,.75)--(5.25,1.5);
\draw[->](5.25,1.5)--(4.875,.75);
\draw[-](4.875,.75)--(4.5,0);
\draw(4.5,0)--(5.25,1.5)node[midway,above,rotate=63.4349]{\(C\)};
\end{tikzpicture}
\]
\noindent where the first third of the level-\(2\) supertile \(\varsigma^2(C)\) is subdivided to \(C\), the second third to \(A\), and the last third to \(D\), so that the path traced out in \(AP_0\) has orientations agreeing with the underlying orientations of the \(1\)-cells.
\end{example}
\indent When we add these labels on top of the supertiles, the substitution map turns into an evaluation map where one evaluates the labels. \emph{The thickening map is exactly labelling the \(d-1\)-cells with their images under \(\delta^{d-1}\),} and performing this process in low dimensions turns the exact sequence
\[
\begin{tikzcd}
0\arrow[r]&\check{H}^{d-1}(\Omega_T)\arrow[r,hook]&C^{d-1}/\image\delta^{d-2}\arrow[r,"\delta^{d-1}"]&C^d\arrow[r,two heads]&\check{H}^d(\Omega_T)\arrow[r]&0
\end{tikzcd}
\]
\noindent into the nontrivial part of the six-term sequence in relative \(K\)-theory that is witnessed by \(\evaluation:K_0(C_r^\ast(\dot{G}_{AF});C_r^\ast(\dot{G}_u))\rightarrow K_0(C_r^\ast(\dot{G}_{AF}))\)!

\section{Computations}
\label{section:computations}
\indent In this section, we first perform the computations by hand for \(d=1\) to illustrate the details, then resort to a script in Sage written by the author that automatically collars and computes all of the induced substitution matrices for \(d=2\), whose results are compiled in \cref{table:computation-d-2}. Note that the computations are performed over \emph{exact rings}, and therefore carry no floating point errors (aside from the eigenspace computations for the induced substitution matrices in Danzer \(7\)-fold and the GKM examples). Furthermore, it appears that the script works for certain \emph{pseudo}substitutions, and we perform the computations for \cite[Figure 4]{clarksadun06}.

\subsection{\(d=1\)}
\begin{example}[Fibonacci]
The Fibonacci substitution has collared prototiles
\[
\{A=(b)a(b),D=(a)b(a),B=(a)b(b),C=(b)b(a)\},
\]
\noindent with the substitution rule
\begin{align*}
A&\mapsto(B)C(A)\\
D&\mapsto AB\\
B&\mapsto AD(A)\\
C&\mapsto (D)AB
\end{align*}
\noindent whose square forces the border
\begin{align*}
A&\mapsto (AD)AB(C)\\
D&\mapsto (B)CAD(A)\\
B&\mapsto (B)CAB(C)\\
C&\mapsto (AB)CAD(A).
\end{align*}
\noindent The transpose of the induced substitution on them is given by the stationary ordered Bratteli diagram
\[
\begin{tikzcd}[nodes={inner sep=1pt},column sep=small]
\cdot\arrow[rd,dash]\arrow[rrd,dash]\arrow[rrrd,dash]&\cdot\arrow[rd,dash]&\cdot\arrow[ld,dash]\arrow[rd,dash]&\cdot\arrow[llld,dash]\\
\cdot\arrow[rd,dash]\arrow[rrd,dash]\arrow[rrrd,dash]&\cdot\arrow[rd,dash]&\cdot\arrow[ld,dash]\arrow[rd,dash]&\cdot\arrow[llld,dash]\\
\cdot&\cdot&\cdot&\cdot
\end{tikzcd}
\]
\noindent where the order of the vertices is the same as the one in the set and the corresponding Bratteli--Vershik system has the order of the edges from left to right as presented in the diagram.\\
\indent The corresponding matrix is
\[
\sigma^\top=\left(\begin{array}{cccc}
0&0&0&1\\
1&0&1&0\\
1&1&0&0\\
1&0&1&0
\end{array}\right)
\]
\noindent which has eigenspaces
\begin{align*}
\Lambda_{AF,\phi}&=\spn\{(\phi^{-1},1,1,1)\}\\
\Lambda_{AF,-1}&=\spn\{(-1,1,0,1)\}\\
\Lambda_{AF,-\phi^{-1}}&=\spn\{(-\phi,1,1,1)\}\\
\Lambda_{AF,0}&=\spn\{-1,1,1,0\}
\end{align*}
\noindent where \(\phi=(1+\sqrt{5})/2\). The product of the nontrivial eigenvalues is \(1\), so the induced map \(\sigma^\top:K_0(C_r^\ast(\dot{G}_{AF,n}))\rightarrow K_0(C_r^\ast(\dot{G}_{AF,n+1}))\) is invertible after removing the kernel, and
\[
K_0(C_r^\ast(\dot{G}_{AF}))\cong\mathbb{Z}^3=\Lambda_{AF,\phi}\oplus\Lambda_{AF,-1}\oplus\Lambda_{AF,-\phi^{-1}}.
\]
\indent We now look for potentially nontrivial partial isometries that are formed from neighboring prototiles, which are
\[
\{(A,D),(A,B),(D,A),(B,C),(C,A)\}
\]
\noindent giving us that the set of full translations is
\[
\{(A,D+B),(D+C,A),(B,C)\}.
\]
\noindent Since we only deal with basic projections, these can be represented by pairs of sets of matrices where we only keep track of the diagonal. Under \((\sigma^\top)^2\), reading off from the Bratteli diagram, the full translations in level-\(0\) are mapped by
\begin{align*}
(1,0,0,0),(0,1,1,0)\mapsto{}&((1,0),(0,1,0),(0,1,0),(0,1,0)),((0,1),(0,0,1),(0,0,1),(0,0,1))\\
={}&((0,0),(0,0,0),(0,0,0),(0,0,0)),((0,0),(0,0,0),(0,0,0),(0,0,0))\\
(0,1,0,1),(1,0,0,0)\mapsto{}&((0,0),(1,0,1),(1,0,0),(0,0,1)),((1,0),(0,1,0),(0,1,0),(0,1,0))\\
={}&((0,0),(0,0,1),(0,0,0),(0,0,1)),((1,0),(0,0,0),(0,0,0),(0,0,0))\\
(0,0,1,0),(0,0,0,1)\mapsto{}&((0,1),(0,0,0),(0,0,1),(0,0,0)),((0,0),(1,0,0),(1,0,0),(1,0,0))\\
={}&((0,1),(0,0,0),(0,0,0),(0,0,0)),((0,0),(1,0,0),(1,0,0),(0,0,0))\\
&+((0,0),(0,0,0),(0,0,1),(0,0,0)),((0,0),(0,0,0),(0,0,0),(1,0,0))
\end{align*}
\noindent where we eliminated the pairs of entries whose sources and ranges belong to \(C_r^\ast(\dot{G}_{AF,2})\), thus are trivial in \(K_0(C_r^\ast(\dot{G}_{AF,2});C_r^\ast(\dot{G}_u))\).\\
\indent Let us explicitly write out the images in terms of tiles in supertiles. For the respective level-\(2\) supertiles, let us denote the position corresponding to the basic projections by underscores and the absence of projections by zeroes. Then, for the second matrix algebra pair, for the second summand of the source, we have the projection
\[
\sigma^2(D)=(B)CA\underline{D}(A)
\]
\noindent which has range
\[
\sigma^2(A)=(AD)\underline{A}B(C)
\]
\noindent which is the first summand of the range. However, \((D,A)\) is not a full translation, and we need to include
\[
\sigma^2(C)=(AB)CA\underline{D}(A)
\]
\noindent in the source as well, which is the fourth summand, giving the triple
\[
[0+(B)CA\underline{D}(A)+0+(AB)CA\underline{D}(A),(AD)\underline{A}B(C)+0+0+0,1]
\]
\noindent as the image, where the last entry indicates that the partial isometry moves left to right by \(1\) unit. We can then denote this by
\[
[\sigma^2(D)+\sigma^2(C),\sigma^2(A),1].
\]
\noindent For the first term in the third matrix algebra pair, we have, for the first summand of the source,
\[
\sigma^2(A)=(AD)A\underline{B}(C),
\]
\noindent which has
\begin{align*}
\sigma^2(D)&=(B)\underline{C}AD(A)\\
\sigma^2(B)&=(B)\underline{C}AB(C),
\end{align*}
\noindent the second and third summands, both as its range. This is a full translation, and the triple is
\[
[(AD)A\underline{B}(C)+0+0+0,0+(B)\underline{C}AD(A)+(B)\underline{C}AB(C)+0,1]
\]
\noindent which we can denote as
\[
[\sigma^2(A),\sigma^2(D)+\sigma^2(B),1].
\]
\noindent For the second term in the third pair, we have, for the third summand,
\[
\sigma^2(B)=(B)CA\underline{B}(C)
\]
\noindent which has range
\[
\sigma^2(C)=(AB)\underline{C}AD(A).
\]
\noindent This is a full translation, so we obtain the triple
\[
[0+0+(B)CA\underline{B}(C)+0,0+0+0+(AB)\underline{C}AD(A),1]
\]
\noindent which we denote as
\[
[\sigma^2(B),\sigma^2(C),1].
\]
\noindent In conclusion, the map is
\begin{align*}
[A,D+B,1]\mapsto{}&0\\
[D+C,A,1]\mapsto{}&[\sigma^2(D)+\sigma^2(C),\sigma^2(A),1]\\
[B,C,1]\mapsto{}&[\sigma^2(A),\sigma^2(D)+\sigma^2(B),1]+[\sigma^2(B),\sigma^2(C),1].
\end{align*}
\indent Since the codomain is a pair of direct sums of matrix algebras, by realizing this as a direct sum of twice the size and repeating the computation for the other two elements, we obtain a map of matrix algebras, which we can represent by a Bratteli diagram
\[
\begin{tikzcd}[nodes={inner sep=1pt},column sep=small]
\cdot&\cdot\arrow[d,dash]&\cdot\arrow[lld,dash]\arrow[d,dash]\\
\cdot&\cdot&\cdot
\end{tikzcd}
\]
\noindent where the order of vertices is the same order we computed. Thus
\[
(\sigma^\top)^2=\left(\begin{array}{ccc}0&0&1\\0&1&0\\0&0&1\end{array}\right),
\]
\noindent and one computes that the eigenspaces are
\begin{align*}
\Lambda_{AF;u,1}&=\spn\{(0,1,0),(1,0,1)\}\\
\Lambda_{AF;u,0}&=\spn\{(1,0,0)\}.
\end{align*}
\noindent Evidently there are no nontrivial relations among the two eigenvectors of \(\Lambda_{AF;u,1}\), so
\[
K_0(C_r^\ast(\dot{G}_{AF});C_r^\ast(\dot{G}_u))\cong\mathbb{Z}^2=\Lambda_{AF;u,1}.
\]
\indent For each of the two eigenvectors, we have that
\begin{align*}
\evaluation((0,1,0))&=(-1,1,0,1)\\
\evaluation((1,0,1))&=(1,-1,0,-1)
\end{align*}
\noindent where the domain and codomain are both expressed in standard basis in full translations and prototiles, respectively. Thus, expressed in eigenbasis,
\[
\evaluation=\left(\begin{array}{cc}0&0\\1&-1\\0&0\end{array}\right).
\]
\indent Finally,
\begin{align*}
K_0(C_r^\ast(\dot{G}_u))&\cong\mathbb{Z}^2=\Lambda_{AF,\phi}\oplus\Lambda_{AF,-\phi^{-1}}\\
K_1(C_r^\ast(\dot{G}_u))&\cong\mathbb{Z}=\spn\{(1,1,1)\}=\Delta(\Lambda_{AF;u,1}).
\end{align*}
\end{example}
\begin{example}[Silver Mean]
The Silver Mean substitution has collared prototiles
\[
\{A=(b)b(a),B=(b)a(b),C=(a)b(b),D=(b)b(b)\}.
\]
\noindent The transpose of the induced substitution on them is given by the stationary ordered Bratteli diagram
\[
\begin{tikzcd}[nodes={inner sep=1pt},column sep=small]
\cdot\arrow[d,dash]\arrow[rrd,dash]\arrow[rrrd,dash]&\cdot\arrow[ld,dash]\arrow[rd,dash]\arrow[rrd,dash]&\cdot\arrow[lld,dash]\arrow[d,dash]\arrow[rd,dash]&\cdot\arrow[lld,dash]\\
\cdot\arrow[d,dash]\arrow[rrd,dash]\arrow[rrrd,dash]&\cdot\arrow[ld,dash]\arrow[rd,dash]\arrow[rrd,dash]&\cdot\arrow[lld,dash]\arrow[d,dash]\arrow[rd,dash]&\cdot\arrow[lld,dash]\\
\cdot&\cdot&\cdot&\cdot
\end{tikzcd}
\]
\noindent where the order of the vertices is the same as the one in the set, the order of the edges is from left to right, and we provided two levels because the substitution rule forces the border at level-\(2\). The corresponding matrix is
\[
\sigma^\top=\left(\begin{array}{cccc}
2&2&2&1\\
1&1&1&0\\
2&2&2&1\\
2&2&2&1
\end{array}\right).
\]
It has eigenspaces
\begin{align*}
\Lambda_{AF,\delta_S^2}&=\spn\{(1,\delta_S^{-1},1,1)\}\\
\Lambda_{AF,\delta_S^{-2}}&=\spn\{(1,-\delta_S,1,1)\}\\
\Lambda_{AF,0}&=\spn\{(-1,0,1,0),(-1,1,0,0)\}
\end{align*}
where \(\delta_S=1+\sqrt{2}\). The product of the nontrivial eigenvalues is again \(1\), so
\[
K_0(C_r^\ast(\dot{G}_{AF}))\cong\mathbb{Z}^2=\Lambda_{AF,\delta_S^2}\oplus\Lambda_{AF,\delta_S^{-2}}.
\]
\indent We look for partial isometries formed from neighboring prototiles, which are
\[
\{(A,B),(B,C),(C,A),(C,D),(D,A)\}
\]
\noindent giving us that the set of triples in \(K_0(C_r^\ast(\dot{G}_{AF});C_r^\ast(\dot{G}_u))\) that we compute the direct limit with is
\[
\{[A,B,1],[B,C,1],[C+D,A+D,1]\}.
\]
\noindent Representing by matrices where we only keep track of the diagonal, under \((\sigma^\top)^2\),
\begin{align*}
(1,0,0,0),(0,1,0,0)\mapsto{}&\begin{array}{c}((1,0,0,0,1,0,0),(0,0,0),(1,0,0,0,1,0,0),(1,0,0,0,1,0,0)),\\((0,1,0,0,0,1,0),(0,0,0),(0,1,0,0,0,1,0),(0,1,0,0,0,1,0))\end{array}\\
={}&\begin{array}{c}((0,0,0,0,0,0,0),(0,0,0),(0,0,0,0,0,0,0),(0,0,0,0,0,0,0)),\\((0,0,0,0,0,0,0),(0,0,0),(0,0,0,0,0,0,0),(0,0,0,0,0,0,0))\end{array}\\
(0,1,0,0),(0,0,1,0)\mapsto{}&\begin{array}{c}((0,1,0,0,0,1,0),(0,0,0),(0,1,0,0,0,1,0),(0,1,0,0,0,1,0)),\\((0,0,1,0,0,0,1),(0,0,0),(0,0,1,0,0,0,1),(0,0,1,0,0,0,1))\end{array}\\
={}&\begin{array}{c}((0,0,0,0,0,0,0),(0,0,0),(0,0,0,0,0,0,0),(0,0,0,0,0,0,0)),\\((0,0,0,0,0,0,0),(0,0,0),(0,0,0,0,0,0,0),(0,0,0,0,0,0,0))\end{array}\\
(0,0,1,1),(1,0,0,1)\mapsto{}&\begin{array}{c}((0,0,1,1,0,0,1),(0,0,1),(0,0,1,1,0,0,1),(0,0,1,1,0,0,1)),\\((1,0,0,1,1,0,0),(1,0,0),(1,0,0,1,1,0,0),(1,0,0,1,1,0,0))\end{array}\\
={}&\begin{array}{c}((0,0,0,0,0,0,1),(0,0,1),(0,0,0,0,0,0,1),(0,0,0,0,0,0,1)),\\((1,0,0,0,0,0,0),(1,0,0),(1,0,0,0,0,0,0),(1,0,0,0,0,0,0))\end{array}\\
={}&\begin{array}{c}((0,0,0,0,0,0,1),(0,0,0),(0,0,0,0,0,0,0),(0,0,0,0,0,0,0)),\\((0,0,0,0,0,0,0),(1,0,0),(0,0,0,0,0,0,0),(0,0,0,0,0,0,0))\end{array}\\
&+\begin{array}{c}((0,0,0,0,0,0,0),(0,0,1),(0,0,0,0,0,0,0),(0,0,0,0,0,0,0)),\\((0,0,0,0,0,0,0),(0,0,0),(1,0,0,0,0,0,0),(0,0,0,0,0,0,0))\end{array}\\
&+\begin{array}{c}((0,0,0,0,0,0,0),(0,0,0),(0,0,0,0,0,0,1),(0,0,0,0,0,0,1)),\\((1,0,0,0,0,0,0),(0,0,0),(0,0,0,0,0,0,0),(1,0,0,0,0,0,0))\end{array}
\end{align*}
\noindent where we eliminated the pairs of entries that are in \(C_r^\ast(\dot{G}_{AF,2})\). For clarity, let us illustrate the process again for computing the images in terms of full translations.\\
\indent Following the same convention as before, the projection from the first summand of the first term of the third matrix algebra pair is
\[
\sigma^2(A)=(BCDABC)ABCDAB\underline{C}(AB)
\]
\noindent which, under the partial isometry, is sent to
\[
\sigma^2(B)=(BCDABC)\underline{A}BC(ABCDAB).
\]
\noindent Since \((A,B)\) is a full translation, the triple is
\[
[(BCDABC)ABCDAB\underline{C}(AB)+0+0+0,0+(BCDABC)\underline{A}BC(ABCDAB)+0+0,1],
\]
\noindent denoted
\[
[\sigma^2(A),\sigma^2(B),1].
\]
\noindent For the second summand of the second term, we have
\[
\sigma^2(B)=(BCDABC)AB\underline{C}(ABCDAB),
\]
\noindent and
\[
\sigma^2(C)=(BC)\underline{A}BCDABC(ABCDAB)
\]
\noindent is the range, and since \((B,C)\) is a full translation, the triple is
\[
[0+(BCDABC)AB\underline{C}(ABCDAB)+0+0,0+0+(BC)\underline{A}BCDABC(ABCDAB)+0,1],
\]
\noindent denoted
\[
[\sigma^2(B),\sigma^2(C),1].
\]
\noindent Finally, the third summand of the third term is
\[
\sigma^2(C)=(BC)ABCDAB\underline{C}(ABCDAB)
\]
\noindent and has both
\begin{align*}
\sigma^2(A)&=(BCDABC)\underline{A}BCDABC(AB)\\
\sigma^2(D)&=(BCDABC)\underline{A}BCDABC(ABCDAB)
\end{align*}
\noindent as the range. However, \((C,A+D)\) is not a full translation. We need to include \(D\), therefore
\[
\sigma^2(D)=(BCDABC)ABCDAB\underline{C}(ABCDAB),
\]
\noindent in the source as well. Thus the corresponding triple is
\begin{align*}
&[0+0+(BC)\underline{A}BCDABC(ABCDAB)+(BCDABC)\underline{A}BCDABC(ABCDAB),\\
&(BCDABC)\underline{A}BCDABC(AB)+0+0+(BCDABC)\underline{A}BCDABC(ABCDAB),1],
\end{align*}
\noindent denoted
\[
[\sigma^2(C)+\sigma^2(D),\sigma^2(A)+\sigma^2(D),1].
\]
\indent Thus the map is
\begin{align*}
[A,B,1]\mapsto{}&0\\
[B,C,1]\mapsto{}&0\\
[C+D,A+D,1]\mapsto{}&[\sigma^2(A),\sigma^2(B),1]+[\sigma^2(B),\sigma^2(C),1]\\
&+[\sigma^2(C)+\sigma^2(D),\sigma^2(A)+\sigma^2(D),1].
\end{align*}
\indent Since the codomain is a pair of direct sums of matrix algebras, by realizing this as a direct sum of twice the size, we obtain a map of matrix algebras, which we can represent by a Bratteli diagram
\[
\begin{tikzcd}[nodes={inner sep=1pt},column sep=small]
\cdot&\cdot&\cdot\arrow[d,dash]\arrow[ld,dash]\arrow[lld,dash]\\
\cdot&\cdot&\cdot
\end{tikzcd}
\]
\noindent where the order of vertices is the same order we computed. Thus
\[
(\sigma^\top)^2=\left(\begin{array}{ccc}0&0&1\\0&0&1\\0&0&1\end{array}\right),
\]
\noindent and the eigenspaces are
\begin{align*}
\Lambda_{AF;u,1}&=\spn\{(1,1,1)\}\\
\Lambda_{AF;u,0}&=\spn\{(1,0,0),(0,1,0)\}.
\end{align*}
\noindent There are no nontrivial relations among the eigenvectors, so
\[
K_0(C_r^\ast(\dot{G}_{AF});C_r^\ast(\dot{G}_u))\cong\mathbb{Z}=\Lambda_{AF;u,1}.
\]
\indent For the single vector, we have that
\begin{align*}
\evaluation((1,1,1))&=(0,0,0,0)
\end{align*}
\noindent where the domain and codomain are both expressed in standard basis in full translations and prototiles, respectively.\\
\indent Finally,
\begin{align*}
K_0(C_r^\ast(\dot{G}_u))&\cong\mathbb{Z}^2=\Lambda_{AF,\delta_S^2}\oplus\Lambda_{AF,\delta_S^{-2}}\\
K_1(C_r^\ast(\dot{G}_u))&\cong\mathbb{Z}=\Lambda_{AF;u,1}.
\end{align*}
\end{example}

\subsection{\(d=2\)}
\begin{example}[Dyadic solenoid]
The only ``example'' we can do by hand\footnote{Rather, the only two-dimensional computation the author is willing to do by hand, particularly since the matrix representations of the translation action that identifies neighboring tiles in a supertile is not as readily transparent as it is for \(d=1\). See \cite{kellendonk97} for the computation for the Robinson triangle tiling.} is the two-dimensional dyadic solenoid, given by the substitution
\[
\begin{tikzpicture}
\draw(0,0)--(0.5,0)--(0.5,0.5)--(0,0.5)--(0,0);
\draw[|->](0.75,0.25)--(1.75,0.25)node[midway,above]{\(\sigma\)};
\draw(2,-0.25)--(3,-0.25)--(3,0.75)--(2,0.75)--(2,-0.25);
\draw(2,0.25)--(3,0.25);
\draw(2.5,-0.25)--(2.5,0.75);
\end{tikzpicture}
\]
where we call the puncture in the single square \(t\). The substitution forces the border, and it is not too hard to see that the \(AP\)-complex is a torus, giving us the two generating partial isometries that are simultaneously full translations
\[
\begin{tikzpicture}
\draw(0,0)--(1,0)--(1,0.5)--(0,0.5)--(0,0);
\draw(0.5,0)--(0.5,0.5);
\fill(0.25,0.25)circle[radius=0.05];
\fill(0.75,0.25)circle[radius=0.05];
\draw[->](0.4,0.25)--(0.6,0.25);
\draw(2.00,-.25)--(2.50,-.25)--(2.50,.75)--(2.00,.75)--(2.00,-.25);
\draw(2.00,.25)--(2.50,.25);
\fill(2.25,0)circle[radius=0.05];
\fill(2.25,.5)circle[radius=0.05];
\draw[->](2.25,.15)--(2.25,.35);
\end{tikzpicture}
\]
which we denote \([t,t,(1,0)]\) and \([t,t,(0,1)]\). Notice that as a substitution, there are multiple preimages that subdivide to each, for example for \([t,t,(1,0)]\),
\[
\begin{tikzpicture}
\draw(0,0)--(2,0)--(2,1)--(0,1)--(0,0);
\draw(1,0)--(1,1);
\draw[dashed](0,0.5)--(2,0.5);
\draw[dashed](0.5,0)--(0.5,1);
\draw[dashed](1.5,0)--(1.5,1);
\fill(0.25,0.25)circle[radius=0.05];
\fill(0.75,0.25)circle[radius=0.05];
\fill(0.25,0.75)circle[radius=0.05];
\fill(0.75,0.75)circle[radius=0.05];
\fill(1.25,0.25)circle[radius=0.05];
\fill(1.75,0.25)circle[radius=0.05];
\fill(1.25,0.75)circle[radius=0.05];
\fill(1.75,0.75)circle[radius=0.05];
\draw[->](0.4,0.25)--(0.6,0.25);
\draw[->](0.4,0.75)--(0.6,0.75);
\draw[->](0.9,0.25)--(1.1,0.25);
\draw[->](0.9,0.75)--(1.1,0.75);
\draw[->](1.4,0.25)--(1.6,0.25);
\draw[->](1.4,0.75)--(1.6,0.75);
\draw[|->](2.25,0.5)--(3.25,0.5)node[midway,above]{\(\sigma\)};
\draw(3.5,0.25)--(4.5,0.25)--(4.5,0.75)--(3.5,0.75)--(3.5,0.25);
\draw(4,0.25)--(4,0.75);
\fill(3.75,0.5)circle[radius=0.05];
\fill(4.25,0.5)circle[radius=0.05];
\draw[->](3.9,0.5)--(4.1,0.5);
\end{tikzpicture}
\]
but if the image is a generator in \(K_0(C_r^\ast(\dot{G}_{AF,n});C_r^\ast(\dot{G}_u))\), four of the partial isometries in the domain have sources and ranges belonging to \(C_r^\ast(\dot{G}_{AF,n+1})\), thus are trivial in \(K_0(C_r^\ast(\dot{G}_{AF,n+1});C_r^\ast(\dot{G}_u))\), and the transpose map becomes
\[
\begin{tikzpicture}
\draw(0,0.25)--(1,0.25)--(1,0.75)--(0,0.75)--(0,0.25);
\draw(0.5,0.25)--(0.5,0.75);
\fill(0.25,0.5)circle[radius=0.05];
\fill(0.75,0.5)circle[radius=0.05];
\draw[->](0.4,0.5)--(0.6,0.5);
\draw[|->](1.25,0.5)--(2.25,0.5)node[midway,above]{\(\sigma^\top\)};
\draw(2.5,0)--(4.5,0)--(4.5,1)--(2.5,1)--(2.5,0);
\draw(3.5,0)--(3.5,1);
\draw[dashed](2.5,0.5)--(4.5,0.5);
\draw[dashed](3,0)--(3,1);
\draw[dashed](4,0)--(4,1);
\fill(3.25,0.25)circle[radius=0.05];
\fill(3.25,0.75)circle[radius=0.05];
\fill(3.75,0.25)circle[radius=0.05];
\fill(3.75,0.75)circle[radius=0.05];
\draw[->](3.4,0.25)--(3.6,0.25);
\draw[->](3.4,0.75)--(3.6,0.75);
\end{tikzpicture}.
\]
There is a single generating relation, given by
\[
\begin{tikzpicture}
\draw(0,0)--(1,0)--(1,1)--(0,1)--(0,0);
\draw(0,0.5)--(1,0.5);
\draw(0.5,0)--(0.5,1);
\fill(0.25,0.25)circle[radius=0.05];
\fill(0.75,0.25)circle[radius=0.05];
\fill(0.25,0.75)circle[radius=0.05];
\fill(0.75,0.75)circle[radius=0.05];
\draw[->](0.4,0.25)--(0.6,0.25);
\draw[->](0.75,0.4)--(0.75,0.6);
\draw[->](0.6,0.75)--(0.4,0.75);
\draw[->](0.25,0.6)--(0.25,0.4);
\end{tikzpicture}
\]
where the top full translation is \([t,t,-(1,0)]\) and the left is \([t,t,-(0,1)]\), giving us
\begin{align*}
[t,t,(1,0)]&+[t,t,(0,1)]+[t,t,-(1,0)]+[t,t,-(0,1)]\\
&=[t,t,(1,0)]+[t,t,(0,1)]-[t,t,(1,0)]-[t,t,(0,1)]\\
&=0.
\end{align*}
This sum of full translations is always trivial since addition in \(K\)-theory is commutative, thus we do not need to find the induced substitution map. The evaluation map is trivial, since \([t,t,(1,0)]\mapsto[t]-[t]=0\) and \([t,t,(0,1)]\mapsto[t]-[t]=0\). Thus, taking direct limits, the six-term sequence is
\[
\begin{tikzcd}
\mathbb{Z}[1/2]^2\arrow[r,"0"]&\mathbb{Z}[1/4]\arrow[r,"\iota_\ast"]&K_0(C_r^\ast(\dot{G}_u))\arrow[d,two heads]\\
K_1(C_r^\ast(\dot{G}_u))\arrow[u,hook]&0\arrow[l]&\mathbb{Z}\arrow[l].
\end{tikzcd}
\]
Therefore \(K_0(C_r^\ast(\dot{G}_u))\cong\mathbb{Z}[1/4]\oplus\mathbb{Z}\) and \(K_1(C_r^\ast(\dot{G}_u))\cong\mathbb{Z}[1/2]^2\), and the roles \(K_0(C_r^\ast(\dot{G}_{AF}))\) and \(K_0(C_r^\ast(\dot{G}_{AF}^{(1)}))\) are as
\[
K_0(C_r^\ast(\dot{G}_u))\cong K_0(C_r^\ast(\dot{G}_{AF}))\oplus\mathbb{Z}
\]
and
\[
K_1(C_r^\ast(\dot{G}_u))\cong K_0(C_r^\ast(\dot{G}_{AF}^{(1)})).
\]
\end{example}
\begin{example}[Half-hex]
Even though the half-hex substitution forces the border, the script collars everything, yielding \(24\) collared prototiles, \(60\) partial isometries, and \(38\) ``loops'' of partial isometries. Out of the partial isometries, there are only \(42\) full translations, meaning that there are at least \(60-42=18\) partial isometries that do not come from subdividing full translations. Similarly, out of the loops, there are only \(20\) that persist. Note that \emph{this does not mean that they do not belong to the kernel of \(\sigma^\top\)!}\\
\indent In \(C^2\), there is one eigenvector of eigenvalue \(4\), one of eigenvalue \(2\), and four of \(1\). In \(C^1\), there are three each of eigenvalues \(2\) and \(1\). \(\image\delta^0\cong\mathbb{Z}\) as one of the eigenvectors of eigenvalue \(1\) in \(C^1\). One easily checks that the direct limits of the induced substitution matrices split completely, giving us the six-term sequence
\[
\begin{tikzcd}
\mathbb{Z}[1/2]^3\oplus\mathbb{Z}^2\arrow[r,"\evaluation"]&\mathbb{Z}[1/4]\oplus\mathbb{Z}[1/2]\oplus\mathbb{Z}^4\arrow[r,"\iota_\ast"]&K_0(C_r^\ast(\dot{G}_u))\arrow[d,two heads]\\
K_1(C_r^\ast(\dot{G}_u))\arrow[u,hook]&0\arrow[l]&\mathbb{Z}\arrow[l]
\end{tikzcd}
\]
where the evaluation map maps an eigenvector of eigenvalue \(2\) and both of the eigenvectors of eigenvalue \(1\) isomorphically onto their images, with the rest mapped to \(0\), giving us \(K_0(C_r^\ast(\dot{G}_u))\cong\mathbb{Z}[1/4]\oplus\mathbb{Z}^3\) and \(K_1(C_r^\ast(\dot{G}_u))\cong\mathbb{Z}[1/2]^2\). The component of \(K_0(C_r^\ast(\dot{G}_{AF}))\) in \(K_0(C_r^\ast(\dot{G}_u))\) is as
\begin{align*}
K_0(C_r^\ast(\dot{G}_u))&\cong(\mathbb{Z}[1/4]\oplus\mathbb{Z}[1/2]\oplus\mathbb{Z}^4)/(\mathbb{Z}[1/2]\oplus\mathbb{Z}^2)\oplus\mathbb{Z}\\
&\cong K_0(C_r^\ast(\dot{G}_{AF}))/(\mathbb{Z}[1/2]\oplus\mathbb{Z}^2)\oplus\mathbb{Z}.
\end{align*}
\noindent Since the substitution rule satisfies the boundary hyperplane condition, we can read off the role of \(K_0(C_r^\ast(\dot{G}_{AF}^{(1)}))\) in \(K_1(C_r^\ast(\dot{G}_u))\) easily as
\begin{align*}
K_1(C_r^\ast(\dot{G}_u))&\cong(\mathbb{Z}[1/2]^2\oplus\mathbb{Z})/\mathbb{Z}\\
&\leq(\mathbb{Z}[1/2]^3\oplus\mathbb{Z}^3)/\mathbb{Z}\\
&\cong K_0(C_r^\ast(\dot{G}_{AF}^{(1)}))/\mathbb{Z}.
\end{align*}
\end{example}
\begin{figure}[t]
\centering
\begin{tikzpicture}
\draw(1.0,.866)--(2.0,.866);
\draw(1.0,.866)--(0,.866);
\draw(.5,1.7321)--(0,.866);
\draw(2.0,.866)--(1.5,1.7321);
\draw(.5,1.7321)--(1.5,1.7321);
\draw[|->](2.25,1.299)--(2.75,1.299);
\draw(7.0,0)--(6.0,0);
\draw(6.5,2.5981)--(5.5,2.5981);
\draw(8.5,.866)--(9.0,0);
\draw(5.0,1.7321)--(4.5,.866);
\draw(3.5,.866)--(3.0,0);
\draw(7.0,0)--(6.5,.866);
\draw(7.5,.866)--(7.0,1.7321);
\draw(7.0,1.7321)--(6.5,2.5981);
\draw(7.5,2.5981)--(6.5,2.5981);
\draw(7.5,.866)--(8.0,1.7321);
\draw(5.0,0)--(5.5,.866);
\draw(4.0,1.7321)--(4.5,2.5981);
\draw(6.5,.866)--(5.5,.866);
\draw(5.0,1.7321)--(5.5,2.5981);
\draw(7.5,.866)--(8.0,0);
\draw(8.0,1.7321)--(8.5,.866);
\draw(7.5,2.5981)--(8.0,1.7321);
\draw(6.0,1.7321)--(7.0,1.7321);
\draw(8.0,0)--(9.0,0);
\draw(5.5,.866)--(5.0,1.7321);
\draw(4.0,0)--(3.0,0);
\draw(7.0,0)--(8.0,0);
\draw(4.0,1.7321)--(4.5,.866);
\draw(5.0,0)--(6.0,0);
\draw(4.0,0)--(4.5,.866);
\draw(3.5,.866)--(4.0,1.7321);
\draw(7.0,1.7321)--(6.5,.866);
\draw(4.0,0)--(5.0,0);
\draw(6.0,1.7321)--(5.0,1.7321);
\draw(5.5,2.5981)--(4.5,2.5981);
\end{tikzpicture}
\caption{Substitution rule for a variant of the half-hex tiling with expansion factor \(3\). The substitution rule is rotationally symmetric.}
\label{figure:halfhex3}
\end{figure}
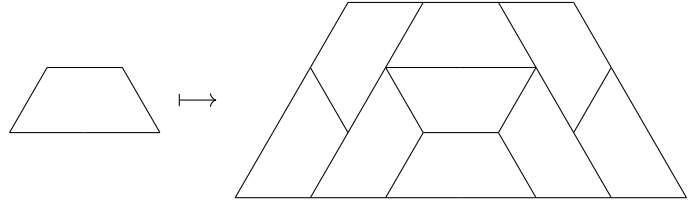%

\begin{example}[Half-hex \(3\times 3\)]
This is the substitution rule in \cite[Figure 9]{langford40} and is shown in \cref{figure:halfhex3}. A supertile is shown in \cref{figure:halfhex3-patch4}. On a finite level of the direct limit, \(\evaluation:K_0(C_r^\ast(\dot{G}_{AF,n});C_r^\ast(\dot{G}_u))\rightarrow K_0(C_r^\ast(\dot{G}_{AF,n}))\) contains a vector that is multiplied by \(2\), hence \(\image\iota_\ast\) has a torsion term. This torsion term has the eigenvalue of \(1\) under \(\sigma^\top\), which is invertible.\\
\indent Due to the sizes of the induced matrices, we assume each of the direct limits splits completely. The induced substitution matrix on \(K_0(C_r^\ast(\dot{G}_{AF}))\) has \(12\) nonintegral eigenvalues multiplying to \(15625\), all of which are greater than \(1\). Restricting the direct limit to the associated eigenvectors of the \(12\) eigenvalues that multiply to \(15625\) and computing it over \(\mathbb{Q}\) gives that \(K_0(C_r^\ast(\dot{G}_{AF}))\cong\mathbb{Z}[1/9]\oplus\mathbb{Z}[1/3]^7\oplus\mathbb{Z}[1/2]^{24}\oplus\mathbb{Z}^{306}\oplus\mathbb{Q}^{12}\). \(6\) of the eigenvectors of the \(12\) eigenvalues are mapped isomorphically onto its image in \(\image\iota_\ast\), the other \(6\) pulled back isomorphically onto its preimage in \(\coimage\evaluation\).\\
\indent The six-term sequence is
\[
\begin{tikzcd}
\begin{tabular}{@{}c@{}}\(\mathbb{Z}[1/3]^6\oplus\mathbb{Z}[1/2]^{15}\)\\\(\oplus\mathbb{Z}^{259}\oplus\mathbb{Q}^6\)\end{tabular}\arrow[r,"\evaluation"]&\begin{tabular}{@{}c@{}}\(\mathbb{Z}[1/9]\oplus\mathbb{Z}[1/3]^7\)\\\(\oplus\mathbb{Z}[1/2]^{24}\oplus\mathbb{Z}^{306}\oplus\mathbb{Q}^{12}\)\end{tabular}\arrow[r,"\iota_\ast"]&\begin{tabular}{@{}c@{}}\(\mathbb{Z}[1/9]\oplus\mathbb{Z}[1/3]^3\oplus\mathbb{Z}[1/2]^9\)\\\(\oplus\mathbb{Z}^{48}\oplus\mathbb{Z}_2\oplus\mathbb{Q}^6\)\end{tabular}\arrow[d,two heads]\\
\mathbb{Z}[1/3]^2\arrow[u,hook]&0\arrow[l]&\mathbb{Z}\arrow[l]
\end{tikzcd}
\]
\noindent where the evaluation map sends \(\mathbb{Z}[1/3]^4\oplus\mathbb{Z}[1/2]^{15}\oplus\mathbb{Z}^{258}\oplus\mathbb{Q}^6\leq K_0(C_r^\ast(\dot{G}_{AF});C_r^\ast(\dot{G}_u))\) isomorphically onto its image, multiplies \(\mathbb{Z}\leq K_0(C_r^\ast(\dot{G}_{AF});C_r^\ast(\dot{G}_u))\) by \(2\), and sends the rest to \(0\). We have that
\begin{align*}
K_0(C_r^\ast(\dot{G}_u))&\cong\left.\left(\begin{tabular}{@{}c@{}}\(\mathbb{Z}[1/9]\oplus\mathbb{Z}[1/3]^7\)\\\(\oplus\mathbb{Z}[1/2]^{24}\oplus\mathbb{Z}^{306}\oplus\mathbb{Q}^{12}\)\end{tabular}\right)\middle/(\mathbb{Z}[1/3]^4\oplus\mathbb{Z}[1/2]^{15}\oplus\mathbb{Z}^{258}\oplus 2\mathbb{Z}\oplus\mathbb{Q}^6)\right.\oplus\mathbb{Z}\\
&\cong K_0(C_r^\ast(\dot{G}_{AF}))/(\mathbb{Z}[1/3]^4\oplus\mathbb{Z}[1/2]^{15}\oplus\mathbb{Z}^{258}\oplus 2\mathbb{Z}\oplus\mathbb{Q}^6)\oplus\mathbb{Z}
\end{align*}
\noindent and, since the substitution rule satisfies the boundary hyperplane condition,
\begin{align*}
K_1(C_r^\ast(\dot{G}_u))&\cong(\mathbb{Z}[1/3]^2\oplus\mathbb{Z}^{143})/\mathbb{Z}^{143}\\
&\leq(\mathbb{Z}[1/3]^6\oplus\mathbb{Z}[1/2]^{15}\oplus\mathbb{Z}^{402}\oplus\mathbb{Q}^6)/\mathbb{Z}^{143}\\
&\cong K_0(C_r^\ast(\dot{G}_{AF}^{(1)}))/\mathbb{Z}^{143}.
\end{align*}
\indent The \(12\) nonintegral eigenvalues in \(K_0(C_r^\ast(\dot{G}_{AF}))\) are \(-2\pm i\), each with multiplicity \(6\). Three of each contribute to \(\check{H}^2(\Omega_T)\cong K_0(C_r^\ast(\dot{G}_{AF}))/\image\evaluation\leq K_0(C_r^\ast(\dot{G}_u))\).
\end{example}
\begin{example}[Chair]
Its \(K\)-theory has already been computed in \cite{juliensavinien16} using an MLD-equivalent substitution with squares, but here we compute directly with the original ``chair'' prototiles. Each of the induced substitution matrices have integral eigenvalues. One easily checks that the direct limits of the induced substitution matrices split completely. It turns out that \(\image\delta^0=0\), giving us \(K_0(C_r^\ast(\dot{G}_{AF});C_r^\ast(\dot{G}_u))\cong K_0(C_r^\ast(\dot{G}_{AF}^{(1)}))\).\\
\indent The six-term sequence is
\[
\begin{tikzcd}
\mathbb{Z}[1/2]^2\oplus\mathbb{Z}^4\arrow[r,"\evaluation"]&\mathbb{Z}[1/4]\oplus\mathbb{Z}[1/2]^2\oplus\mathbb{Z}^4\arrow[r,"\iota_\ast"]&\mathbb{Z}[1/4]\oplus\mathbb{Z}[1/2]^2\oplus\mathbb{Z}\arrow[d,two heads]\\
\mathbb{Z}[1/2]^2\arrow[u,hook]&0\arrow[l]&\mathbb{Z}\arrow[l]
\end{tikzcd}
\]
\noindent where the evaluation map sends \(\mathbb{Z}^4\leq K_0(C_r^\ast(\dot{G}_{AF});C_r^\ast(\dot{G}_u))\) isomorphically onto its image and the rest to \(0\). We have that
\begin{align*}
K_0(C_r^\ast(\dot{G}_u))&\cong(\mathbb{Z}[1/4]\oplus\mathbb{Z}[1/2]^2\oplus\mathbb{Z}^4)/\mathbb{Z}^4\oplus\mathbb{Z}\\
&\cong K_0(C_r^\ast(\dot{G}_{AF}))/\mathbb{Z}^4\oplus\mathbb{Z}
\end{align*}
\noindent and
\begin{align*}
K_1(C_r^\ast(\dot{G}_u))&\cong\mathbb{Z}[1/2]^2\\
&\leq\mathbb{Z}[1/2]^2\oplus\mathbb{Z}^4\\
&\cong C^1\\
&\cong K_0(C_r^\ast(\dot{G}_{AF}^{(1)}))\tag*{(\cref{proposition:af-cochain-nontrivial})}
\end{align*}
\noindent where we have to apply a slightly nontrivial isomorphism in the last step since this substitution does not satisfy the boundary hyperplane condition. The proposition applies since the component in \(C^1\) contributed by the \(1\)-cells satisfying the boundary hyperplane condition splits completely, being \(\mathbb{Z}[1/2]^2\).
\end{example}
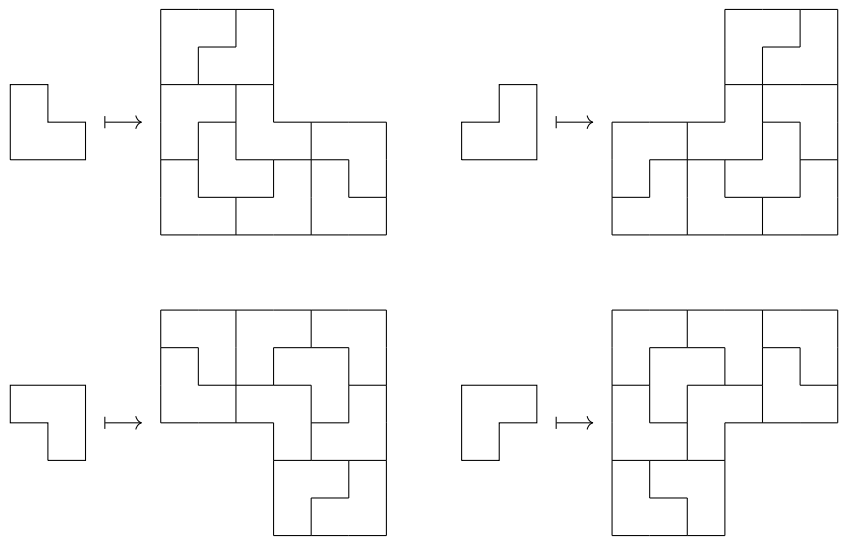
\begin{figure}[t]
\centering
\begin{tikzpicture}
\draw(0,1)--(1,1)--(1,1.5)--(.5,1.5)--(.5,2)--(0,2)--(0,1);
\draw[|->](1.25,1.5)--(1.75,1.5);
\draw(3.5,1.0)--(3.0,1.0);
\draw(4.0,0)--(3.5,0);
\draw(2.5,2.0)--(2.5,2.5);
\draw(4.5,1.5)--(4.0,1.5);
\draw(2.5,3.0)--(2.0,3.0);
\draw(4.5,1.5)--(5.0,1.5);
\draw(3.0,1.5)--(3.0,1.0);
\draw(3.0,3.0)--(3.5,3.0);
\draw(5.0,.5)--(4.5,.5);
\draw(5.0,1.0)--(5.0,1.5);
\draw(4.0,1.0)--(4.0,1.5);
\draw(5.0,.5)--(5.0,1.0);
\draw(2.0,2.0)--(2.0,2.5);
\draw(3.5,1.5)--(4.0,1.5);
\draw(3.5,1.5)--(3.5,2.0);
\draw(2.0,1.0)--(2.0,1.5);
\draw(3.5,1.0)--(4.0,1.0);
\draw(2.5,1.0)--(2.5,1.5);
\draw(4.0,1.0)--(4.5,1.0);
\draw(2.0,.5)--(2.0,1.0);
\draw(3.5,.5)--(3.5,1.0);
\draw(3.0,0)--(3.0,.5);
\draw(3.5,2.0)--(3.5,2.5);
\draw(4.5,.5)--(4.5,1.0);
\draw(4.5,0)--(5.0,0);
\draw(5.0,.5)--(5.0,0);
\draw(2.5,0)--(2.0,0);
\draw(2.0,1.0)--(2.5,1.0);
\draw(3.0,2.0)--(2.5,2.0);
\draw(4.0,0)--(4.5,0);
\draw(2.0,1.5)--(2.0,2.0);
\draw(2.5,3.0)--(3.0,3.0);
\draw(2.5,.5)--(2.5,1.0);
\draw(3.0,0)--(3.5,0);
\draw(2.0,2.5)--(2.0,3.0);
\draw(2.5,0)--(3.0,0);
\draw(4.0,0)--(4.0,.5);
\draw(2.0,.5)--(2.0,0);
\draw(2.0,2.0)--(2.5,2.0);
\draw(2.5,.5)--(3.0,.5);
\draw(4.0,.5)--(4.0,1.0);
\draw(3.5,2.5)--(3.5,3.0);
\draw(3.0,2.0)--(3.5,2.0);
\draw(3.0,2.5)--(3.0,3.0);
\draw(3.5,.5)--(3.0,.5);
\draw(3.0,2.5)--(2.5,2.5);
\draw(3.0,1.5)--(3.0,2.0);
\draw(3.0,1.5)--(2.5,1.5);
\draw(6.0,1)--(7.0,1)--(7.0,2)--(6.5,2)--(6.5,1.5)--(6.0,1.5)--(6.0,1);
\draw[|->](7.25,1.5)--(7.75,1.5);
\draw(8.5,1.5)--(9.0,1.5);
\draw(8.0,.5)--(8.0,1.0);
\draw(9.0,0)--(9.5,0);
\draw(9.0,0)--(9.0,.5);
\draw(11.0,1.0)--(10.5,1.0);
\draw(10.0,0)--(9.5,0);
\draw(10.0,0)--(10.0,.5);
\draw(9.0,1.0)--(9.0,1.5);
\draw(8.5,0)--(9.0,0);
\draw(10.5,1.0)--(10.5,1.5);
\draw(10.0,1.0)--(10.0,1.5);
\draw(9.5,2.5)--(9.5,3.0);
\draw(10.5,3.0)--(10.5,2.5);
\draw(9.5,1.0)--(10.0,1.0);
\draw(11.0,2.0)--(11.0,2.5);
\draw(8.0,1.5)--(8.0,1.0);
\draw(8.0,0)--(8.0,.5);
\draw(10.5,0)--(11.0,0);
\draw(11.0,1.0)--(11.0,1.5);
\draw(10.0,.5)--(9.5,.5);
\draw(9.5,1.5)--(9.0,1.5);
\draw(8.5,.5)--(8.0,.5);
\draw(11.0,.5)--(11.0,1.0);
\draw(9.5,1.0)--(9.0,1.0);
\draw(10.0,1.5)--(10.0,2.0);
\draw(8.5,0)--(8.0,0);
\draw(10.5,2.5)--(10.0,2.5);
\draw(10.5,3.0)--(10.0,3.0);
\draw(10.0,1.5)--(10.5,1.5);
\draw(10.0,.5)--(10.5,.5);
\draw(10.5,2.0)--(10.0,2.0);
\draw(10.0,0)--(10.5,0);
\draw(11.0,1.5)--(11.0,2.0);
\draw(10.0,3.0)--(9.5,3.0);
\draw(9.5,2.0)--(9.5,2.5);
\draw(9.5,2.0)--(10.0,2.0);
\draw(10.0,2.0)--(10.0,2.5);
\draw(8.0,1.5)--(8.5,1.5);
\draw(11.0,2.5)--(11.0,3.0);
\draw(9.0,.5)--(9.0,1.0);
\draw(10.5,.5)--(10.5,1.0);
\draw(8.5,1.0)--(9.0,1.0);
\draw(11.0,.5)--(11.0,0);
\draw(9.5,1.5)--(9.5,2.0);
\draw(10.5,3.0)--(11.0,3.0);
\draw(11.0,2.0)--(10.5,2.0);
\draw(9.5,1.0)--(9.5,.5);
\draw(8.5,.5)--(8.5,1.0);
\draw(.5,-3)--(1,-3)--(1,-2)--(0,-2)--(0,-2.5)--(.5,-2.5)--(.5,-3);
\draw[|->](1.25,-2.5)--(1.75,-2.5);
\draw(4.0,-3.0)--(4.0,-2.5);
\draw(5.0,-3.5)--(5.0,-4.0);
\draw(5.0,-2.5)--(5.0,-2.0);
\draw(3.0,-1.0)--(3.5,-1.0);
\draw(3.5,-2.5)--(3.0,-2.5);
\draw(4.0,-1.0)--(3.5,-1.0);
\draw(4.0,-3.0)--(3.5,-3.0);
\draw(3.5,-2.5)--(3.5,-3.0);
\draw(3.0,-2.0)--(2.5,-2.0);
\draw(3.5,-4.0)--(3.5,-3.5);
\draw(2.0,-2.5)--(2.0,-2.0);
\draw(2.0,-2.5)--(2.5,-2.5);
\draw(4.5,-3.5)--(4.5,-3.0);
\draw(2.5,-1.0)--(3.0,-1.0);
\draw(3.0,-1.5)--(3.0,-1.0);
\draw(3.5,-4.0)--(4.0,-4.0);
\draw(4.5,-1.0)--(5.0,-1.0);
\draw(5.0,-3.5)--(5.0,-3.0);
\draw(5.0,-2.0)--(4.5,-2.0);
\draw(4.5,-2.5)--(4.5,-2.0);
\draw(3.5,-2.0)--(4.0,-2.0);
\draw(4.0,-2.5)--(4.0,-2.0);
\draw(3.5,-1.5)--(4.0,-1.5);
\draw(3.0,-2.0)--(3.5,-2.0);
\draw(3.0,-2.0)--(3.0,-2.5);
\draw(5.0,-4.0)--(4.5,-4.0);
\draw(2.5,-1.0)--(2.0,-1.0);
\draw(3.5,-1.5)--(3.5,-2.0);
\draw(4.0,-1.5)--(4.5,-1.5);
\draw(4.5,-2.5)--(4.0,-2.5);
\draw(4.5,-1.5)--(4.5,-2.0);
\draw(4.5,-1.0)--(4.0,-1.0);
\draw(5.0,-2.0)--(5.0,-1.5);
\draw(4.5,-3.5)--(4.0,-3.5);
\draw(4.5,-4.0)--(4.0,-4.0);
\draw(2.5,-2.0)--(2.5,-1.5);
\draw(3.0,-2.0)--(3.0,-1.5);
\draw(3.5,-3.5)--(3.5,-3.0);
\draw(5.0,-3.0)--(4.5,-3.0);
\draw(5.0,-3.0)--(5.0,-2.5);
\draw(4.0,-3.5)--(4.0,-4.0);
\draw(4.0,-1.0)--(4.0,-1.5);
\draw(2.0,-1.5)--(2.0,-1.0);
\draw(2.0,-1.5)--(2.5,-1.5);
\draw(4.5,-3.0)--(4.0,-3.0);
\draw(5.0,-1.5)--(5.0,-1.0);
\draw(2.0,-1.5)--(2.0,-2.0);
\draw(2.5,-2.5)--(3.0,-2.5);
\draw(6.0,-3)--(6.5,-3)--(6.5,-2.5)--(7.0,-2.5)--(7.0,-2)--(6.0,-2)--(6.0,-3);
\draw[|->](7.25,-2.5)--(7.75,-2.5);
\draw(9.0,-1.5)--(8.5,-1.5);
\draw(10.0,-1.0)--(9.5,-1.0);
\draw(11.0,-2.5)--(11.0,-2.0);
\draw(8.0,-3.5)--(8.0,-4.0);
\draw(8.0,-2.5)--(8.0,-2.0);
\draw(9.0,-2.5)--(8.5,-2.5);
\draw(8.5,-4.0)--(8.0,-4.0);
\draw(9.0,-2.5)--(9.0,-2.0);
\draw(8.5,-4.0)--(9.0,-4.0);
\draw(10.0,-1.5)--(10.5,-1.5);
\draw(8.0,-3.5)--(8.0,-3.0);
\draw(9.5,-3.0)--(9.5,-2.5);
\draw(10.5,-2.5)--(10.0,-2.5);
\draw(10.0,-2.0)--(10.0,-2.5);
\draw(9.5,-3.5)--(9.5,-3.0);
\draw(8.5,-3.5)--(8.5,-3.0);
\draw(10.5,-2.5)--(11.0,-2.5);
\draw(9.5,-2.5)--(10.0,-2.5);
\draw(9.0,-2.0)--(9.5,-2.0);
\draw(9.0,-4.0)--(9.0,-3.5);
\draw(9.0,-1.5)--(9.5,-1.5);
\draw(8.5,-2.0)--(8.5,-1.5);
\draw(9.5,-2.0)--(9.5,-1.5);
\draw(10.5,-2.0)--(10.5,-1.5);
\draw(10.0,-1.5)--(10.0,-2.0);
\draw(10.5,-1.0)--(11.0,-1.0);
\draw(10.0,-1.0)--(10.5,-1.0);
\draw(11.0,-1.5)--(11.0,-2.0);
\draw(8.5,-3.5)--(9.0,-3.5);
\draw(9.0,-2.5)--(9.0,-3.0);
\draw(8.0,-1.5)--(8.0,-2.0);
\draw(8.5,-1.0)--(8.0,-1.0);
\draw(8.5,-2.5)--(8.5,-2.0);
\draw(9.5,-3.5)--(9.5,-4.0);
\draw(8.0,-3.0)--(8.0,-2.5);
\draw(9.0,-1.0)--(9.5,-1.0);
\draw(9.0,-4.0)--(9.5,-4.0);
\draw(8.5,-1.0)--(9.0,-1.0);
\draw(8.0,-3.0)--(8.5,-3.0);
\draw(9.5,-3.0)--(9.0,-3.0);
\draw(8.5,-3.0)--(9.0,-3.0);
\draw(10.5,-2.0)--(11.0,-2.0);
\draw(9.0,-1.5)--(9.0,-1.0);
\draw(8.0,-2.0)--(8.5,-2.0);
\draw(10.0,-1.5)--(10.0,-1.0);
\draw(8.0,-1.5)--(8.0,-1.0);
\draw(10.0,-2.0)--(9.5,-2.0);
\draw(11.0,-1.5)--(11.0,-1.0);
\end{tikzpicture}
\caption{Substitution rule for a variant of the chair tiling with expansion factor \(3\).}
\label{figure:chair3}
\end{figure}%

\begin{example}[Chair \(3\times 3\)\footnote{Supplied by Rodrigo Trevi\~no.}]
The substitution rule is a variant of the chair tiling with expansion factor \(3\), and is shown in \cref{figure:chair3}. A supertile is shown in \cref{figure:chair3-patch4}. Unlike the standard chair tiling, rotational symmetry has been eliminated. Each of the induced substitution matrices have integral eigenvalues. Due to the sizes of the induced matrices, we assume each of the direct limits splits completely.\\
\indent The six-term sequence is
\[
\begin{tikzcd}
\mathbb{Z}[1/3]^4\oplus\mathbb{Z}[1/2]^2\oplus\mathbb{Z}^{76}\arrow[r,"\evaluation"]&\begin{tabular}{@{}c@{}}\(\mathbb{Z}[1/9]\oplus\mathbb{Z}[1/3]^6\)\\\(\oplus\mathbb{Z}[1/2]^4\oplus\mathbb{Z}^{101}\)\end{tabular}\arrow[r,"\iota_\ast"]&\begin{tabular}{@{}c@{}}\(\mathbb{Z}[1/9]\oplus\mathbb{Z}[1/3]^4\)\\\(\oplus\mathbb{Z}[1/2]^2\oplus\mathbb{Z}^{26}\)\end{tabular}\arrow[d,two heads]\\
\mathbb{Z}[1/3]^2\arrow[u,hook]&0\arrow[l]&\mathbb{Z}\arrow[l]
\end{tikzcd}
\]
\noindent where the evaluation map sends \(\mathbb{Z}[1/3]^2\oplus\mathbb{Z}[1/2]^2\oplus\mathbb{Z}^{76}\leq K_0(C_r^\ast(\dot{G}_{AF});C_r^\ast(\dot{G}_u))\) isomorphically onto its image and the rest to \(0\). We have that
\begin{align*}
K_0(C_r^\ast(\dot{G}_u))&\cong(\mathbb{Z}[1/9]\oplus\mathbb{Z}[1/3]^6\oplus\mathbb{Z}[1/2]^4\oplus\mathbb{Z}^{101})/(\mathbb{Z}[1/3]^2\oplus\mathbb{Z}[1/2]^2\oplus\mathbb{Z}^{76})\oplus\mathbb{Z}\\
&\cong K_0(C_r^\ast(\dot{G}_{AF}))/(\mathbb{Z}[1/3]^2\oplus\mathbb{Z}[1/2]^2\oplus\mathbb{Z}^{76})\oplus\mathbb{Z}
\end{align*}
\noindent and
\begin{align*}
K_1(C_r^\ast(\dot{G}_u))&\cong(\mathbb{Z}[1/3]^2\oplus\mathbb{Z}^{34})/\mathbb{Z}^{34}\\
&\leq(\mathbb{Z}[1/3]^4\oplus\mathbb{Z}[1/2]^2\oplus\mathbb{Z}^{110})/\mathbb{Z}^{34}\\
&\cong C^1/\mathbb{Z}^{34}\\
&\cong K_0(C_r^\ast(\dot{G}_{AF}^{(1)}))/\mathbb{Z}^{34}\tag*{(\cref{proposition:af-cochain-nontrivial})}
\end{align*}
\noindent where we have to apply a slightly nontrivial isomorphism in the last step since this substitution does not satisfy the boundary hyperplane condition. The proposition applies since we assumed that \(C^1\) splits completely.
\end{example}
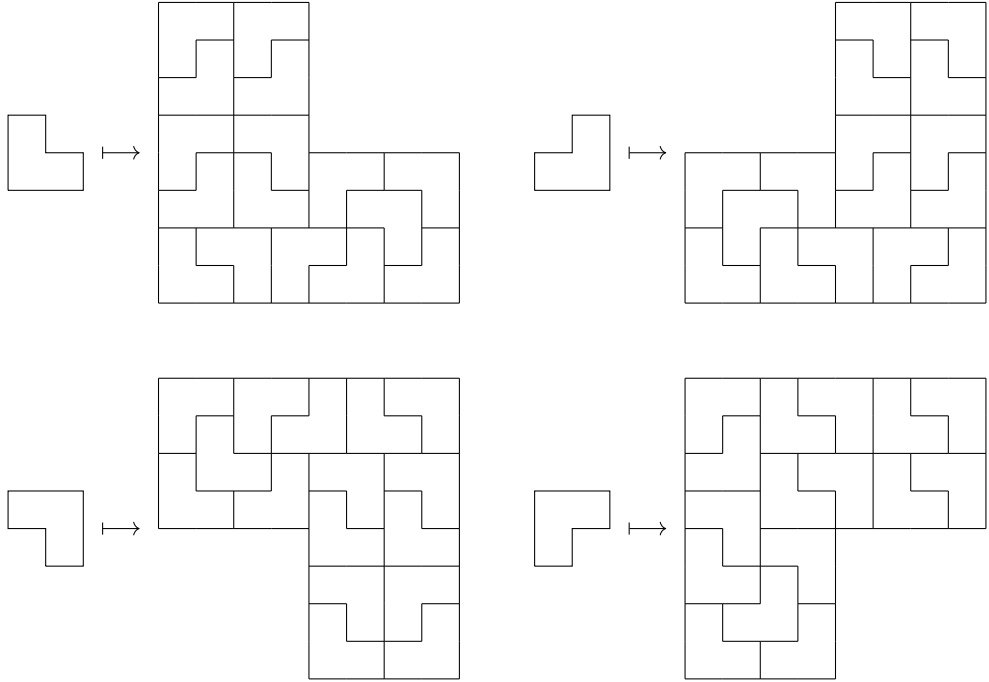
\begin{figure}[t]
\centering
\begin{tikzpicture}
\draw(0,1.5)--(1,1.5)--(1,2.0)--(.5,2.0)--(.5,2.5)--(0,2.5)--(0,1.5);
\draw[|->](1.25,2.0)--(1.75,2.0);
\draw(2.5,3.0)--(2.0,3.0);
\draw(3.0,1.5)--(3.0,1.0);
\draw(4.0,1.0)--(4.0,1.5);
\draw(5.0,.5)--(5.0,1.0);
\draw(4.5,2.0)--(5.0,2.0);
\draw(2.0,1.0)--(2.0,1.5);
\draw(3.0,3.5)--(3.0,4.0);
\draw(2.5,4.0)--(3.0,4.0);
\draw(2.5,0)--(2.0,0);
\draw(3.0,2.0)--(2.5,2.0);
\draw(4.0,0)--(4.5,0);
\draw(6.0,1.0)--(5.5,1.0);
\draw(2.0,1.5)--(2.0,2.0);
\draw(4.0,2.5)--(3.5,2.5);
\draw(2.5,.5)--(2.5,1.0);
\draw(4.0,0)--(4.0,.5);
\draw(2.0,.5)--(2.0,0);
\draw(4.0,3.5)--(4.0,4.0);
\draw(2.5,3.0)--(2.5,3.5);
\draw(2.0,3.5)--(2.0,3.0);
\draw(3.0,2.5)--(3.0,3.0);
\draw(3.0,1.5)--(3.0,2.0);
\draw(2.5,1.5)--(2.5,2.0);
\draw(4.0,3.0)--(4.0,3.5);
\draw(2.5,1.0)--(3.0,1.0);
\draw(5.5,2.0)--(5.0,2.0);
\draw(5.0,.5)--(5.5,.5);
\draw(5.0,1.5)--(5.0,2.0);
\draw(6.0,1.5)--(6.0,2.0);
\draw(5.0,1.5)--(5.5,1.5);
\draw(4.0,1.0)--(4.5,1.0);
\draw(3.0,2.0)--(3.0,2.5);
\draw(3.5,.5)--(3.5,1.0);
\draw(2.0,2.5)--(2.5,2.5);
\draw(2.0,1.0)--(2.5,1.0);
\draw(2.0,4.0)--(2.5,4.0);
\draw(2.0,2.5)--(2.0,3.0);
\draw(5.5,0)--(6.0,0);
\draw(3.5,1.5)--(4.0,1.5);
\draw(6.0,0)--(6.0,.5);
\draw(5.5,1.0)--(5.5,1.5);
\draw(4.0,0)--(3.5,0);
\draw(3.5,1.0)--(3.0,1.0);
\draw(3.5,4.0)--(3.0,4.0);
\draw(3.5,.5)--(3.5,0);
\draw(6.0,1.0)--(6.0,1.5);
\draw(4.5,1.5)--(5.0,1.5);
\draw(3.0,3.0)--(3.5,3.0);
\draw(2.0,2.0)--(2.0,2.5);
\draw(3.5,1.5)--(3.5,2.0);
\draw(2.0,1.5)--(2.5,1.5);
\draw(4.5,1.5)--(4.5,1.0);
\draw(2.0,3.5)--(2.0,4.0);
\draw(3.0,0)--(3.0,.5);
\draw(4.5,.5)--(4.5,1.0);
\draw(5.0,.5)--(5.0,0);
\draw(6.0,1.0)--(6.0,.5);
\draw(3.5,3.5)--(4.0,3.5);
\draw(3.0,2.5)--(2.5,2.5);
\draw(3.5,4.0)--(4.0,4.0);
\draw(3.5,1.0)--(4.0,1.0);
\draw(4.0,2.5)--(4.0,3.0);
\draw(3.5,3.5)--(3.5,3.0);
\draw(5.5,0)--(5.0,0);
\draw(3.0,3.0)--(3.0,3.5);
\draw(5.5,.5)--(5.5,1.0);
\draw(2.0,.5)--(2.0,1.0);
\draw(2.5,3.5)--(3.0,3.5);
\draw(4.0,2.0)--(4.0,2.5);
\draw(4.0,2.0)--(4.0,1.5);
\draw(4.5,0)--(5.0,0);
\draw(5.5,2.0)--(6.0,2.0);
\draw(5.0,1.0)--(4.5,1.0);
\draw(3.0,2.5)--(3.5,2.5);
\draw(4.0,2.0)--(4.5,2.0);
\draw(3.0,0)--(3.5,0);
\draw(2.5,0)--(3.0,0);
\draw(2.5,.5)--(3.0,.5);
\draw(4.0,.5)--(4.5,.5);
\draw(3.0,2.0)--(3.5,2.0);
\draw(7.0,1.5)--(8.0,1.5)--(8.0,2.5)--(7.5,2.5)--(7.5,2.0)--(7.0,2.0)--(7.0,1.5);
\draw[|->](8.25,2.0)--(8.75,2.0);
\draw(11.0,0)--(11.0,.5);
\draw(12.0,0)--(11.5,0);
\draw(12.0,0)--(12.0,.5);
\draw(9.5,.5)--(9.5,1.0);
\draw(11.5,0)--(11.5,.5);
\draw(12.0,3.0)--(12.0,2.5);
\draw(10.0,1.5)--(9.5,1.5);
\draw(13.0,1.0)--(13.0,1.5);
\draw(11.5,1.5)--(11.0,1.5);
\draw(12.5,4.0)--(12.0,4.0);
\draw(12.0,3.5)--(12.5,3.5);
\draw(12.0,1.5)--(12.5,1.5);
\draw(12.0,1.0)--(12.5,1.0);
\draw(12.0,2.5)--(12.5,2.5);
\draw(12.0,.5)--(12.5,.5);
\draw(12.0,0)--(12.5,0);
\draw(13.0,1.5)--(13.0,2.0);
\draw(11.5,4.0)--(11.0,4.0);
\draw(11.0,2.0)--(10.5,2.0);
\draw(12.5,.5)--(12.5,1.0);
\draw(13.0,.5)--(13.0,0);
\draw(13.0,4.0)--(12.5,4.0);
\draw(11.0,2.5)--(11.5,2.5);
\draw(10.5,1.0)--(10.5,1.5);
\draw(11.5,3.0)--(11.5,3.5);
\draw(13.0,3.5)--(13.0,3.0);
\draw(10.0,.5)--(10.0,1.0);
\draw(13.0,1.0)--(12.5,1.0);
\draw(9.5,1.0)--(9.5,1.5);
\draw(11.5,4.0)--(12.0,4.0);
\draw(11.0,3.0)--(11.0,3.5);
\draw(9.0,1.0)--(9.0,1.5);
\draw(11.0,2.5)--(11.0,3.0);
\draw(11.0,2.0)--(11.0,1.5);
\draw(9.5,.5)--(10.0,.5);
\draw(12.0,3.5)--(12.0,4.0);
\draw(10.5,0)--(10.0,0);
\draw(12.0,3.0)--(11.5,3.0);
\draw(11.5,2.0)--(12.0,2.0);
\draw(13.0,2.5)--(13.0,3.0);
\draw(10.5,1.0)--(11.0,1.0);
\draw(10.5,1.0)--(10.0,1.0);
\draw(11.5,1.0)--(11.5,.5);
\draw(10.5,.5)--(10.5,1.0);
\draw(9.0,1.0)--(9.5,1.0);
\draw(10.5,.5)--(11.0,.5);
\draw(10.0,2.0)--(9.5,2.0);
\draw(11.0,3.5)--(11.0,4.0);
\draw(11.0,1.0)--(11.0,1.5);
\draw(9.0,.5)--(9.0,1.0);
\draw(9.5,0)--(10.0,0);
\draw(10.5,0)--(11.0,0);
\draw(9.0,0)--(9.0,.5);
\draw(13.0,2.0)--(13.0,2.5);
\draw(13.0,3.5)--(13.0,4.0);
\draw(12.5,3.0)--(12.5,3.5);
\draw(12.0,3.0)--(12.0,3.5);
\draw(12.5,3.0)--(13.0,3.0);
\draw(13.0,2.0)--(12.5,2.0);
\draw(10.0,2.0)--(10.5,2.0);
\draw(11.0,0)--(11.5,0);
\draw(9.0,1.5)--(9.0,2.0);
\draw(9.0,0)--(9.5,0);
\draw(12.0,1.0)--(12.0,1.5);
\draw(10.0,0)--(10.0,.5);
\draw(11.5,1.0)--(12.0,1.0);
\draw(12.5,0)--(13.0,0);
\draw(9.0,2.0)--(9.5,2.0);
\draw(11.0,2.0)--(11.0,2.5);
\draw(11.0,3.5)--(11.5,3.5);
\draw(11.5,2.5)--(12.0,2.5);
\draw(13.0,.5)--(13.0,1.0);
\draw(11.5,1.0)--(11.0,1.0);
\draw(12.5,1.5)--(12.5,2.0);
\draw(10.0,1.5)--(10.0,2.0);
\draw(12.0,1.5)--(12.0,2.0);
\draw(12.0,2.0)--(12.0,2.5);
\draw(10.0,1.5)--(10.5,1.5);
\draw(11.5,1.5)--(11.5,2.0);
\draw(13.0,2.5)--(12.5,2.5);
\draw(.5,-3.5)--(1,-3.5)--(1,-2.5)--(0,-2.5)--(0,-3.0)--(.5,-3.0)--(.5,-3.5);
\draw[|->](1.25,-3.0)--(1.75,-3.0);
\draw(2.5,-3.0)--(3.0,-3.0);
\draw(4.5,-4.5)--(5.0,-4.5);
\draw(5.0,-1.0)--(4.5,-1.0);
\draw(2.0,-3.0)--(2.5,-3.0);
\draw(4.0,-1.5)--(3.5,-1.5);
\draw(4.0,-2.0)--(3.5,-2.0);
\draw(2.0,-1.0)--(2.0,-1.5);
\draw(6.0,-4.5)--(6.0,-4.0);
\draw(2.0,-2.5)--(2.0,-2.0);
\draw(4.0,-2.0)--(4.0,-2.5);
\draw(5.0,-3.5)--(5.5,-3.5);
\draw(5.0,-2.5)--(5.5,-2.5);
\draw(3.0,-2.0)--(3.0,-1.5);
\draw(2.5,-2.0)--(2.5,-1.5);
\draw(5.5,-1.0)--(5.0,-1.0);
\draw(5.5,-1.5)--(5.0,-1.5);
\draw(5.5,-2.0)--(5.5,-1.5);
\draw(6.0,-2.0)--(6.0,-1.5);
\draw(5.0,-2.0)--(5.5,-2.0);
\draw(6.0,-5.0)--(5.5,-5.0);
\draw(6.0,-3.0)--(5.5,-3.0);
\draw(4.0,-4.5)--(4.0,-5.0);
\draw(3.0,-1.5)--(2.5,-1.5);
\draw(6.0,-1.5)--(6.0,-1.0);
\draw(6.0,-5.0)--(6.0,-4.5);
\draw(6.0,-2.5)--(6.0,-2.0);
\draw(5.0,-5.0)--(5.0,-4.5);
\draw(5.0,-3.5)--(5.0,-3.0);
\draw(5.0,-3.5)--(4.5,-3.5);
\draw(3.5,-3.0)--(4.0,-3.0);
\draw(5.0,-2.5)--(5.0,-2.0);
\draw(4.5,-5.0)--(4.0,-5.0);
\draw(3.5,-1.0)--(3.0,-1.0);
\draw(5.0,-3.5)--(5.0,-4.0);
\draw(5.0,-1.0)--(5.0,-1.5);
\draw(3.0,-2.0)--(3.5,-2.0);
\draw(3.0,-1.5)--(3.0,-1.0);
\draw(4.5,-1.5)--(4.5,-1.0);
\draw(2.0,-3.0)--(2.0,-2.5);
\draw(5.5,-3.0)--(5.5,-2.5);
\draw(5.0,-3.0)--(5.0,-2.5);
\draw(4.0,-4.5)--(4.0,-4.0);
\draw(4.5,-2.5)--(4.0,-2.5);
\draw(3.0,-2.5)--(3.5,-2.5);
\draw(2.5,-1.0)--(3.0,-1.0);
\draw(3.5,-1.0)--(4.0,-1.0);
\draw(4.0,-3.5)--(4.0,-3.0);
\draw(4.0,-1.5)--(4.0,-1.0);
\draw(3.5,-2.5)--(3.5,-2.0);
\draw(3.0,-2.5)--(2.5,-2.5);
\draw(6.0,-3.5)--(6.0,-3.0);
\draw(4.5,-2.0)--(5.0,-2.0);
\draw(4.0,-3.0)--(4.0,-2.5);
\draw(2.5,-2.5)--(2.5,-2.0);
\draw(3.0,-2.5)--(3.0,-3.0);
\draw(2.0,-2.0)--(2.5,-2.0);
\draw(6.0,-3.5)--(6.0,-4.0);
\draw(4.0,-4.0)--(4.5,-4.0);
\draw(4.0,-1.0)--(4.5,-1.0);
\draw(5.0,-4.5)--(5.0,-4.0);
\draw(5.0,-3.0)--(4.5,-3.0);
\draw(5.5,-4.5)--(5.5,-4.0);
\draw(4.5,-2.5)--(4.5,-3.0);
\draw(4.0,-4.0)--(4.0,-3.5);
\draw(2.0,-1.0)--(2.5,-1.0);
\draw(4.0,-3.5)--(4.5,-3.5);
\draw(5.5,-1.0)--(6.0,-1.0);
\draw(6.0,-2.0)--(5.5,-2.0);
\draw(4.0,-2.0)--(4.5,-2.0);
\draw(6.0,-4.0)--(5.5,-4.0);
\draw(4.5,-4.5)--(4.5,-4.0);
\draw(6.0,-3.5)--(5.5,-3.5);
\draw(4.5,-5.0)--(5.0,-5.0);
\draw(6.0,-3.0)--(6.0,-2.5);
\draw(3.5,-3.0)--(3.0,-3.0);
\draw(2.0,-1.5)--(2.0,-2.0);
\draw(5.0,-4.5)--(5.5,-4.5);
\draw(5.5,-5.0)--(5.0,-5.0);
\draw(3.5,-1.5)--(3.5,-2.0);
\draw(4.5,-1.5)--(4.5,-2.0);
\draw(7.0,-3.5)--(7.5,-3.5)--(7.5,-3.0)--(8.0,-3.0)--(8.0,-2.5)--(7.0,-2.5)--(7.0,-3.5);
\draw[|->](8.25,-3.0)--(8.75,-3.0);
\draw(12.0,-2.5)--(12.0,-2.0);
\draw(10.5,-3.0)--(11.0,-3.0);
\draw(11.0,-4.5)--(11.0,-5.0);
\draw(9.0,-4.5)--(9.0,-4.0);
\draw(13.0,-1.5)--(13.0,-1.0);
\draw(9.5,-3.5)--(9.5,-3.0);
\draw(10.0,-2.0)--(10.5,-2.0);
\draw(10.5,-2.5)--(10.5,-2.0);
\draw(11.0,-1.5)--(11.0,-2.0);
\draw(11.5,-2.0)--(11.5,-1.5);
\draw(10.0,-1.5)--(10.0,-2.0);
\draw(12.5,-1.5)--(12.5,-2.0);
\draw(11.0,-1.0)--(11.5,-1.0);
\draw(9.0,-1.5)--(9.0,-2.0);
\draw(13.0,-2.0)--(13.0,-1.5);
\draw(9.5,-2.0)--(9.5,-1.5);
\draw(9.5,-1.0)--(9.0,-1.0);
\draw(9.0,-3.0)--(9.5,-3.0);
\draw(10.0,-1.5)--(10.0,-1.0);
\draw(11.5,-2.0)--(12.0,-2.0);
\draw(9.0,-1.5)--(9.0,-1.0);
\draw(12.0,-2.5)--(12.5,-2.5);
\draw(11.0,-2.0)--(10.5,-2.0);
\draw(11.0,-1.5)--(10.5,-1.5);
\draw(12.0,-3.0)--(11.5,-3.0);
\draw(9.0,-5.0)--(9.0,-4.5);
\draw(10.5,-1.5)--(10.5,-1.0);
\draw(9.0,-2.5)--(9.0,-2.0);
\draw(10.0,-5.0)--(10.0,-4.5);
\draw(9.5,-4.0)--(9.0,-4.0);
\draw(10.0,-2.5)--(9.5,-2.5);
\draw(11.0,-4.0)--(10.5,-4.0);
\draw(11.5,-2.5)--(11.5,-3.0);
\draw(12.0,-3.0)--(12.5,-3.0);
\draw(13.0,-2.5)--(13.0,-2.0);
\draw(10.5,-2.5)--(11.0,-2.5);
\draw(12.5,-2.0)--(12.0,-2.0);
\draw(12.5,-1.5)--(12.0,-1.5);
\draw(11.0,-3.0)--(11.0,-2.5);
\draw(12.5,-1.0)--(13.0,-1.0);
\draw(10.0,-4.5)--(9.5,-4.5);
\draw(12.0,-1.5)--(12.0,-1.0);
\draw(10.0,-5.0)--(9.5,-5.0);
\draw(11.0,-1.0)--(10.5,-1.0);
\draw(10.5,-5.0)--(11.0,-5.0);
\draw(9.0,-5.0)--(9.5,-5.0);
\draw(12.5,-2.0)--(13.0,-2.0);
\draw(10.0,-2.5)--(10.0,-2.0);
\draw(11.5,-1.5)--(11.5,-1.0);
\draw(11.0,-2.0)--(11.5,-2.0);
\draw(9.0,-3.5)--(9.0,-3.0);
\draw(13.0,-3.0)--(12.5,-3.0);
\draw(11.0,-4.5)--(11.0,-4.0);
\draw(10.0,-4.0)--(10.0,-3.5);
\draw(10.5,-3.5)--(10.0,-3.5);
\draw(10.5,-4.5)--(10.5,-4.0);
\draw(9.0,-2.5)--(9.5,-2.5);
\draw(12.5,-3.0)--(12.5,-2.5);
\draw(9.0,-2.0)--(9.5,-2.0);
\draw(10.0,-5.0)--(10.5,-5.0);
\draw(10.0,-1.5)--(9.5,-1.5);
\draw(9.0,-3.5)--(9.0,-4.0);
\draw(10.0,-4.5)--(10.5,-4.5);
\draw(9.5,-4.0)--(10.0,-4.0);
\draw(12.5,-1.0)--(12.0,-1.0);
\draw(11.0,-3.5)--(11.0,-3.0);
\draw(11.5,-1.0)--(12.0,-1.0);
\draw(10.0,-2.5)--(10.0,-3.0);
\draw(9.5,-4.0)--(9.5,-4.5);
\draw(9.5,-3.5)--(10.0,-3.5);
\draw(10.5,-3.5)--(10.5,-4.0);
\draw(11.5,-2.5)--(11.5,-2.0);
\draw(9.0,-3.0)--(9.0,-2.5);
\draw(10.0,-1.0)--(10.5,-1.0);
\draw(9.5,-1.0)--(10.0,-1.0);
\draw(13.0,-3.0)--(13.0,-2.5);
\draw(10.5,-3.0)--(10.0,-3.0);
\draw(11.0,-3.0)--(11.5,-3.0);
\draw(10.0,-3.5)--(10.0,-3.0);
\draw(11.0,-4.0)--(11.0,-3.5);
\end{tikzpicture}
\caption{Substitution rule for a variant of the chair tiling with expansion factor \(4\).}
\label{figure:chair4}
\end{figure}%

\begin{example}[Chair \(4\times 4\)\footnote{Supplied by Rodrigo Trevi\~no.}]
The substitution rule is a variant of the chair tiling with expansion factor \(4\), and is shown in \cref{figure:chair4}. A supertile is shown in \cref{figure:chair4-patch3}. Unlike the standard chair tiling, rotational symmetry has been eliminated. Each of the induced substitution matrices have integral eigenvalues. Due to the sizes of the induced matrices, we assume each of the direct limits splits completely. On a finite level of the direct limit, \(\evaluation:K_0(C_r^\ast(\dot{G}_{AF,n});C_r^\ast(\dot{G}_u))\rightarrow K_0(C_r^\ast(\dot{G}_{AF,n}))\) contains a vector that is multiplied by \(2\), hence \(\image\iota_\ast\) has a torsion term. (Un)fortunately, this torsion term belongs to the kernel of \(\sigma^\top\). It also turns out that \(\image\delta^0=0\), giving us \(K_0(C_r^\ast(\dot{G}_{AF});C_r^\ast(\dot{G}_u))\cong K_0(C_r^\ast(\dot{G}_{AF}^{(1)}))\).\\
\indent The six-term sequence is
\[
\begin{tikzcd}
\mathbb{Z}[1/4]^2\oplus\mathbb{Z}[1/2]\oplus\mathbb{Z}^{19}\arrow[r,"\evaluation"]&\mathbb{Z}[1/16]\oplus\mathbb{Z}[1/2]^4\oplus\mathbb{Z}^{30}\arrow[r,"\iota_\ast"]&\mathbb{Z}[1/16]\oplus\mathbb{Z}[1/2]^3\oplus\mathbb{Z}^{12}\arrow[d,two heads]\\
\mathbb{Z}[1/4]^2\arrow[u,hook]&0\arrow[l]&\mathbb{Z}\arrow[l]
\end{tikzcd}
\]
\noindent where the evaluation map sends \(\mathbb{Z}[1/2]\oplus\mathbb{Z}^{19}\leq K_0(C_r^\ast(\dot{G}_{AF});C_r^\ast(\dot{G}_u))\) isomorphically onto its image and the rest to \(0\). We have that
\begin{align*}
K_0(C_r^\ast(\dot{G}_u))&\cong(\mathbb{Z}[1/16]\oplus\mathbb{Z}[1/2]^4\oplus\mathbb{Z}^{30})/(\mathbb{Z}[1/2]\oplus\mathbb{Z}^{19})\oplus\mathbb{Z}\\
&\cong K_0(C_r^\ast(\dot{G}_{AF}))/(\mathbb{Z}[1/2]\oplus\mathbb{Z}^{19})\oplus\mathbb{Z}
\end{align*}
\noindent and
\begin{align*}
K_1(C_r^\ast(\dot{G}_u))&\cong\mathbb{Z}[1/4]^2\\
&\leq\mathbb{Z}[1/4]^2\oplus\mathbb{Z}[1/2]^2\oplus\mathbb{Z}^{19}\\
&\cong C^1\\
&\cong K_0(C_r^\ast(\dot{G}_{AF}^{(1)}))\tag*{(\cref{proposition:af-cochain-nontrivial})}
\end{align*}
\noindent where we have to apply a slightly nontrivial isomorphism in the last step since this substitution does not satisfy the boundary hyperplane condition. The proposition applies since we assumed that \(C^1\) splits completely.
\end{example}
\begin{example}[Table]
This is also known as the domino tiling. On a finite level of the direct limit, \(\evaluation:K_0(C_r^\ast(\dot{G}_{AF,n});C_r^\ast(\dot{G}_u))\rightarrow K_0(C_r^\ast(\dot{G}_{AF,n}))\) contains a vector that is multiplied by \(2\), hence \(\image\iota_\ast\) has a torsion term. This torsion term has the eigenvalue of \(1\) under \(\sigma^\top\), which is invertible. Due to the sizes of the induced matrices, we assume each of the direct limits splits completely.\\
\indent The six-term sequence is
\[
\begin{tikzcd}
\mathbb{Z}[1/2]^6\oplus\mathbb{Z}^{29}\arrow[r,"\evaluation"]&\mathbb{Z}[1/4]\oplus\mathbb{Z}[1/2]^8\oplus\mathbb{Z}^{32}\arrow[r,"\iota_\ast"]&\mathbb{Z}[1/4]\oplus\mathbb{Z}[1/2]^4\oplus\mathbb{Z}^4\oplus\mathbb{Z}_2\arrow[d,two heads]\\
\mathbb{Z}[1/2]^2\arrow[u,hook]&0\arrow[l]&\mathbb{Z}\arrow[l]
\end{tikzcd}
\]
\noindent where the evaluation map sends \(\mathbb{Z}[1/2]^4\oplus\mathbb{Z}^{28}\leq K_0(C_r^\ast(\dot{G}_{AF});C_r^\ast(\dot{G}_u))\) isomorphically onto its image, multiplies \(\mathbb{Z}\leq K_0(C_r^\ast(\dot{G}_{AF});C_r^\ast(\dot{G}_u))\) by \(2\), and sends the rest to \(0\). We have that
\begin{align*}
K_0(C_r^\ast(\dot{G}_u))&\cong(\mathbb{Z}[1/4]\oplus\mathbb{Z}[1/2]^8\oplus\mathbb{Z}^{32})/(\mathbb{Z}[1/2]^4\oplus\mathbb{Z}^{28}\oplus 2\mathbb{Z})\oplus\mathbb{Z}\\
&\cong K_0(C_r^\ast(\dot{G}_{AF}))/(\mathbb{Z}[1/2]^4\oplus\mathbb{Z}^{28}\oplus 2\mathbb{Z})\oplus\mathbb{Z}
\end{align*}
\noindent and, since the substitution rule satisfies the boundary hyperplane condition,
\begin{align*}
K_1(C_r^\ast(\dot{G}_u))&\cong(\mathbb{Z}[1/2]^2\oplus\mathbb{Z}^{21})/\mathbb{Z}^{21}\\
&\leq(\mathbb{Z}[1/2]^6\oplus\mathbb{Z}^{50})/\mathbb{Z}^{21}\\
&\cong K_0(C_r^\ast(\dot{G}_{AF}^{(1)}))/\mathbb{Z}^{21}.
\end{align*}
\end{example}
\begin{example}[Robinson triangle]
This is MLD-equivalent to the Penrose tiling. This substitution rule is known to be border-forcing, but we still collar everything, resulting in a large number of collared prototiles. The induced substitution matrices, after removing the eventual kernels, are all of determinant \(1\).\\
\indent The six-term sequence is
\[
\begin{tikzcd}
\mathbb{Z}^{37}\arrow[r,"\evaluation"]&\mathbb{Z}^{40}\arrow[r,"\iota_\ast"]&\mathbb{Z}^9\arrow[d,two heads]\\
\mathbb{Z}^5\arrow[u,hook]&0\arrow[l]&\mathbb{Z}\arrow[l]
\end{tikzcd}
\]
\noindent where the evaluation map sends \(\mathbb{Z}^{32}\leq K_0(C_r^\ast(\dot{G}_{AF});C_r^\ast(\dot{G}_u))\) isomorphically onto its image and the rest to \(0\). We have that
\begin{align*}
K_0(C_r^\ast(\dot{G}_u))&\cong\mathbb{Z}^{40}/\mathbb{Z}^{32}\oplus\mathbb{Z}\\
&\cong K_0(C_r^\ast(\dot{G}_{AF}))/\mathbb{Z}^{32}\oplus\mathbb{Z}
\end{align*}
\noindent and, since the substitution rule satisfies the boundary hyperplane condition,
\begin{align*}
K_1(C_r^\ast(\dot{G}_u))&\cong\mathbb{Z}^8/\mathbb{Z}^3\\
&\leq\mathbb{Z}^{40}/\mathbb{Z}^3\\
&\cong K_0(C_r^\ast(\dot{G}_{AF}^{(1)}))/\mathbb{Z}^3.
\end{align*}
\end{example}
\begin{example}[T\"ubingen triangle]
A supertile is shown in \cref{figure:tubingen-patch9}.\\
\indent On a finite level of the direct limit, \(\evaluation:K_0(C_r^\ast(\dot{G}_{AF,n});C_r^\ast(\dot{G}_u))\rightarrow K_0(C_r^\ast(\dot{G}_{AF,n}))\) contains \(2\) vectors that are each multiplied by \(5\), hence \(\image\iota_\ast\) has \(2\) torsion terms. These torsion terms both have the same eigenvalue of \(4\) under \(\sigma^\top\), which is invertible. Aside from this term, the induced substitution matrices, after removing the eventual kernels, are all of determinant \(1\).\\
\indent The six-term sequence is
\[
\begin{tikzcd}
\mathbb{Z}^{101}\arrow[r,"\evaluation"]&\mathbb{Z}^{120}\arrow[r,"\iota_\ast"]&\mathbb{Z}^{25}\oplus\mathbb{Z}_5^2\arrow[d,two heads]\\
\mathbb{Z}^5\arrow[u,hook]&0\arrow[l]&\mathbb{Z}\arrow[l]
\end{tikzcd}
\]
\noindent where the evaluation map sends \(\mathbb{Z}^{94}\leq K_0(C_r^\ast(\dot{G}_{AF});C_r^\ast(\dot{G}_u))\) isomorphically onto its image, multiplies \(\mathbb{Z}^2\leq K_0(C_r^\ast(\dot{G}_{AF});C_r^\ast(\dot{G}_u))\) each by \(5\), and sends the rest to \(0\). We have that
\begin{align*}
K_0(C_r^\ast(\dot{G}_u))&\cong\mathbb{Z}^{120}/(\mathbb{Z}^{94}\oplus(5\mathbb{Z})^2)\oplus\mathbb{Z}\\
&\cong K_0(C_r^\ast(\dot{G}_{AF}))/(\mathbb{Z}^{94}\oplus(5\mathbb{Z}^2))\oplus\mathbb{Z}
\end{align*}
\noindent and, since the substitution rule satisfies the boundary hyperplane condition,
\begin{align*}
K_1(C_r^\ast(\dot{G}_u))&\cong\mathbb{Z}^{14}/\mathbb{Z}^9\\
&\leq\mathbb{Z}^{110}/\mathbb{Z}^9\\
&\cong K_0(C_r^\ast(\dot{G}_{AF}^{(1)}))/\mathbb{Z}^9.
\end{align*}
\end{example}
\begin{example}[Danzer \(7\)-fold]
We are using the original version of the Danzer \(7\)-fold substitution for simplicity, as it is a substitution of triangles, rather than of (degenerate) quadrilaterals. A supertile is shown in \cref{figure:danzer7-patch4}. The expansion factor is \(\lambda=1+2\cos\pi/7\), which is \emph{not a PV number}, thus is weak mixing. Due to the sizes of the induced substitution matrices on \(K_0(C_r^\ast(\dot{G}_{AF,n}))\) and \(K_0(C_r^\ast(\dot{G}_{AF,n});C_r^\ast(\dot{G}_u))\), their eigenvalues had to be computed numerically. We assume that any eigenvalue with order \(\leq 10^{-7}\) is nonzero due to floating point errors.\\
\indent The determinant of the induced substitution matrix on \(K_0(C_r^\ast(\dot{G}_{AF}))\), upon removing the eventual kernel, is \(32761\). Fewer than \(18\) of the nonintegral eigenvalues multiply to it, with the remaining multiplying to \(1\). Restricting the direct limit to the associated eigenvectors of the \(18\) eigenvalues that multiply to \(32761\) and computing it over \(\mathbb{Q}\) gives that \(K_0(C_r^\ast(\dot{G}_{AF}))\cong\mathbb{Z}^{\geq 2213}\oplus\mathbb{Q}^{\leq 18}\). These eigenvectors of the \(18\) eigenvalues are mapped isomorphically onto its image in \(\image\iota_\ast\). The induced substitution matrix on \(K_0(C_r^\ast(\dot{G}_{AF});C_r^\ast(\dot{G}_u))\), after removing the eventual kernel, has determinant \(1\).\\
\indent The six-term sequence is
\[
\begin{tikzcd}
\mathbb{Z}^{2017}\arrow[r,"\evaluation"]&\mathbb{Z}^{\geq 2222}\oplus\mathbb{Q}^{\leq 18}\arrow[r,"\iota_\ast"]&\mathbb{Z}^{\geq 212}\oplus\mathbb{Q}^{\leq 18}\arrow[d,two heads]\\
\mathbb{Z}^6\arrow[u,hook]&0\arrow[l]&\mathbb{Z}\arrow[l]
\end{tikzcd}
\]
\noindent where the evaluation map sends \(\mathbb{Z}^{2011}\leq K_0(C_r^\ast(\dot{G}_{AF});C_r^\ast(\dot{G}_u))\) isomorphically onto its image and the rest to \(0\). We have that
\begin{align*}
K_0(C_r^\ast(\dot{G}_u))&\cong(\mathbb{Z}^{\geq 2222}\oplus\mathbb{Q}^{\leq 18})/\mathbb{Z}^{2011}\oplus\mathbb{Z}\\
&\cong K_0(C_r^\ast(\dot{G}_{AF}))/\mathbb{Z}^{2011}\oplus\mathbb{Z}
\end{align*}
\noindent and, since the substitution rule satisfies the boundary hyperplane condition,
\begin{align*}
K_1(C_r^\ast(\dot{G}_u))&\cong\mathbb{Z}^{215}/\mathbb{Z}^{209}\\
&\leq\mathbb{Z}^{2226}/\mathbb{Z}^{209}\\
&\cong K_0(C_r^\ast(\dot{G}_{AF}^{(1)}))/\mathbb{Z}^{209}.
\end{align*}
\indent Among the \(229\) nontrivial eigenvalues contributing to \(\check{H}^2(\Omega_T)\cong K_0(C_r^\ast(\dot{G}_{AF}))/\image\evaluation\leq K_0(C_r^\ast(\dot{G}_u))\), all but \(21\) belong to \(\mathbb{Z}[\lambda]\). The contribution to the determinant of \(32761\) is purely from them. There are three that multiply to the integral value of \(1\) which we include into the integral part of \(K_0(C_r^\ast(\dot{G}_u))\),
\begin{align*}
0.2864\cdots\textnormal{, }{}&2.9122\cdots\textnormal{, }\\
{}&\textnormal{roots of }x^3-2x^2-3x+1,\\
1.1986\cdots\textnormal{, }{}&\\
{}&\textnormal{root of }x^3+2x^2-3x-1,
\end{align*}
\noindent with the rest, each with multiplicity \(2\),
\begin{align*}
0.4826\cdots\textnormal{, }{}&1.1310\cdots\textnormal{, }2.0456\cdots\textnormal{, }2.4534\cdots\textnormal{, }\\
{}&\textnormal{roots of }x^9+6x^8-4x^7-73x^6-49x^5+277x^4+271x^3-364x^2-289x+181,\\
1.7331\cdots\textnormal{, }{}&1.8534\cdots\textnormal{, }1.8872\cdots\textnormal{, }2.9696\cdots\textnormal{, }3.6693\cdots\textnormal{, }\\
{}&\textnormal{roots of }x^9-6x^8-4x^7+73x^6-49x^5-277x^4+271x^3+364x^2-289x-181.
\end{align*}
\noindent There are four eigenvalues, \(\lambda^2-\lambda-1\) with multiplicity \(4\), \(3.6694\cdots\) with multiplicity \(2\), \(2.9696\cdots\) with multiplicity \(2\), and \(2.9122\cdots\), that are strictly between \(\lambda\) and \(\lambda^2\).\\
\indent Among the six nontrivial eigenvalues contributing to \(K_1(C_r^\ast(\dot{G}_u))\cong\check{H}^1(\Omega_T)\), they multiply to \(1\) with four expanding and two contracting.
\end{example}
\begin{example}[Ammann A\(2\)]
Per \cite[Section 10.3]{andersonputnam98}, we use the square of the substitution rule and the \(8\) prototiles it is primitive on. This substitution rule is known to be border-forcing, but we still collar everything, resulting in a large number of collared prototiles. The induced substitution matrices, after removing the eventual kernels, are all of determinant \(1\).\\
\indent The six-term sequence is
\[
\begin{tikzcd}
\mathbb{Z}^{6}\arrow[r,"\evaluation"]&\mathbb{Z}^8\arrow[r,"\iota_\ast"]&\mathbb{Z}^{7}\arrow[d,two heads]\\
\mathbb{Z}^4\arrow[u,hook]&0\arrow[l]&\mathbb{Z}\arrow[l]
\end{tikzcd}
\]
\noindent where the evaluation map sends \(\mathbb{Z}^2\leq K_0(C_r^\ast(\dot{G}_{AF});C_r^\ast(\dot{G}_u))\) isomorphically onto its image and the rest to \(0\). We have that
\begin{align*}
K_0(C_r^\ast(\dot{G}_u))&\cong\mathbb{Z}^8/\mathbb{Z}^2\oplus\mathbb{Z}\\
&\cong K_0(C_r^\ast(\dot{G}_{AF}))/\mathbb{Z}^2\oplus\mathbb{Z}
\end{align*}
\noindent and
\begin{align*}
K_1(C_r^\ast(\dot{G}_u))&\cong\mathbb{Z}^6/\mathbb{Z}^2\\
&\leq\mathbb{Z}^8/\mathbb{Z}^2\\
&\cong C^1/\mathbb{Z}^2\\
&\cong K_0(C_r^\ast(\dot{G}_{AF}^{(1)}))/\mathbb{Z}^2
\end{align*}
\noindent where we have to apply a slightly nontrivial isomorphism in the last step since this substitution does not satisfy the boundary hyperplane condition. \cref{proposition:af-cochain-nontrivial} \emph{does not apply}, since the eigenvalues are not integral. However, the induced substitution matrix on \(C^1\) has determinant \(1\) (up to the eventual kernel), thus as does the induced substitution matrix restricted to the \(1\)-cells that satisfy the boundary hyperplane condition. A proof similar to the one in the proposition gives an isomorphism \(K_0(C_r^\ast(\dot{G}_{AF}^{(1)}))\cong C^1\).
\end{example}
\begin{example}[Ammann A\(5\) decorated]
We are using the rhombus-triangle substitution. This is MLD-equivalent to the Ammann--Beenker tiling. The induced substitution matrices, after removing the eventual kernels, are all of determinant \(1\).\\
\indent The six-term sequence is
\[
\begin{tikzcd}
\mathbb{Z}^{65}\arrow[r,"\evaluation"]&\mathbb{Z}^{80}\arrow[r,"\iota_\ast"]&\mathbb{Z}^{24}\arrow[d,two heads]\\
\mathbb{Z}^8\arrow[u,hook]&0\arrow[l]&\mathbb{Z}\arrow[l]
\end{tikzcd}
\]
\noindent where the evaluation map sends \(\mathbb{Z}^{57}\leq K_0(C_r^\ast(\dot{G}_{AF});C_r^\ast(\dot{G}_u))\) isomorphically onto its image and the rest to \(0\). We have that
\begin{align*}
K_0(C_r^\ast(\dot{G}_u))&\cong\mathbb{Z}^{80}/\mathbb{Z}^{57}\oplus\mathbb{Z}\\
&\cong K_0(C_r^\ast(\dot{G}_{AF}))/\mathbb{Z}^{57}\oplus\mathbb{Z}
\end{align*}
\noindent and, since the substitution rule satisfies the boundary hyperplane condition,
\begin{align*}
K_1(C_r^\ast(\dot{G}_u))&\cong\mathbb{Z}^{15}/\mathbb{Z}^7\\
&\leq\mathbb{Z}^{72}/\mathbb{Z}^7\\
&\cong K_0(C_r^\ast(\dot{G}_{AF}^{(1)}))/\mathbb{Z}^7.
\end{align*}
\end{example}
\begin{example}[Ammann A\(5\) undecorated]
We are using the rhombus-triangle substitution. This is MLD-equivalent to the Ammann--Beenker tiling \emph{with the tile labels removed}. The induced substitution matrices, after removing the eventual kernels, are all of determinant \(1\). It turns out that \(\image\delta^0=0\), giving us \(K_0(C_r^\ast(\dot{G}_{AF});C_r^\ast(\dot{G}_u))\cong K_0(C_r^\ast(\dot{G}_{AF}^{(1)}))\).\\
\indent The six-term sequence is
\[
\begin{tikzcd}
\mathbb{Z}^{16}\arrow[r,"\evaluation"]&\mathbb{Z}^{20}\arrow[r,"\iota_\ast"]&\mathbb{Z}^{10}\arrow[d,two heads]\\
\mathbb{Z}^5\arrow[u,hook]&0\arrow[l]&\mathbb{Z}\arrow[l]
\end{tikzcd}
\]
\noindent where the evaluation map sends \(\mathbb{Z}^{11}\leq K_0(C_r^\ast(\dot{G}_{AF});C_r^\ast(\dot{G}_u))\) isomorphically onto its image and the rest to \(0\). We have that
\begin{align*}
K_0(C_r^\ast(\dot{G}_u))&\cong\mathbb{Z}^{20}/\mathbb{Z}^{11}\oplus\mathbb{Z}\\
&\cong K_0(C_r^\ast(\dot{G}_{AF}))/\mathbb{Z}^{11}\oplus\mathbb{Z}
\end{align*}
\noindent and, since the substitution rule satisfies the boundary hyperplane condition,
\begin{align*}
K_1(C_r^\ast(\dot{G}_u))&\cong\mathbb{Z}^8\\
&\leq\mathbb{Z}^{16}\\
&\cong K_0(C_r^\ast(\dot{G}_{AF}^{(1)})).
\end{align*}
\end{example}
\begin{example}[GKM \(9.1.1.1\)]
This is the substitution rule in \cite[Figure 9]{gahlerkwanmaloney15}, where the enumeration is, from left to right, the first triangle, followed by the first triangle connected by an arc, followed by the first triangle connected by an arc. A supertile is shown in \cref{figure:gkm9.1.1.1-patch4}. The expansion factor is \(\lambda=1+2\cos\pi/7\), which is \emph{not a PV number}, thus is weak mixing. Due to the sizes of the induced substitution matrices on \(K_0(C_r^\ast(\dot{G}_{AF,n}))\) and \(K_0(C_r^\ast(\dot{G}_{AF,n});C_r^\ast(\dot{G}_u))\), their eigenvalues had to be computed numerically. The induced substitution matrices, after removing the eventual kernels, are all of determinant \(1\).\\
\indent The six-term sequence is
\[
\begin{tikzcd}
\mathbb{Z}^{602}\arrow[r,"\evaluation"]&\mathbb{Z}^{672}\arrow[r,"\iota_\ast"]&\mathbb{Z}^{84}\arrow[d,two heads]\\
\mathbb{Z}^{12}\arrow[u,hook]&0\arrow[l]&\mathbb{Z}^2\arrow[l]
\end{tikzcd}
\]
\noindent where the evaluation map sends \(\mathbb{Z}^{590}\leq K_0(C_r^\ast(\dot{G}_{AF});C_r^\ast(\dot{G}_u))\) isomorphically onto its image and the rest to \(0\). Note that the substitution \emph{appears to have two components, and is thus not primitive!} We have that
\begin{align*}
K_0(C_r^\ast(\dot{G}_u))&\cong\mathbb{Z}^{672}/\mathbb{Z}^{590}\oplus\mathbb{Z}^2\\
&\cong K_0(C_r^\ast(\dot{G}_{AF}))/\mathbb{Z}^{590}\oplus\mathbb{Z}^2
\end{align*}
\noindent and, since the substitution rule satisfies the boundary hyperplane condition,
\begin{align*}
K_1(C_r^\ast(\dot{G}_u))&\cong\mathbb{Z}^{96}/\mathbb{Z}^{84}\\
&\leq\mathbb{Z}^{686}/\mathbb{Z}^{84}\\
&\cong K_0(C_r^\ast(\dot{G}_{AF}^{(1)}))/\mathbb{Z}^{84}.
\end{align*}
\indent Among the \(84\) nontrivial eigenvalues contributing to \(\check{H}^2(\Omega_T)\cong K_0(C_r^\ast(\dot{G}_{AF}))/\image\evaluation\leq K_0(C_r^\ast(\dot{G}_u))\), they all belong to \(\mathbb{Z}[\lambda]\), and there is an eigenvalue \(\lambda^2-\lambda-1\) with multiplicity \(8\) that is strictly between \(\lambda\) and \(\lambda^2\).\\
\indent Among the twelve nontrivial eigenvalues contributing to \(K_1(C_r^\ast(\dot{G}_u))\cong\check{H}^1(\Omega_T)\), they multiply to \(1\) with eight expanding and four contracting.
\end{example}
\begin{example}[GKM \(9.2.1.1\)]
This is the substitution rule in \cite[Figure 9]{gahlerkwanmaloney15}, where the enumeration is, from left to right, the second triangle, followed by the first triangle connected by an arc, followed by the first triangle connected by an arc. A supertile is shown in \cref{figure:gkm9.2.1.1-patch4}. The expansion factor is \(\lambda=1+2\cos\pi/7\), which is \emph{not a PV number}, thus is weak mixing. Due to the sizes of the induced substitution matrices on \(K_0(C_r^\ast(\dot{G}_{AF,n}))\) and \(K_0(C_r^\ast(\dot{G}_{AF,n});C_r^\ast(\dot{G}_u))\), their eigenvalues had to be computed numerically. The induced substitution matrices, after removing the eventual kernels, are all of determinant \(1\).\\
\indent The six-term sequence is
\[
\begin{tikzcd}
\mathbb{Z}^{1022}\arrow[r,"\evaluation"]&\mathbb{Z}^{1148}\arrow[r,"\iota_\ast"]&\mathbb{Z}^{140}\arrow[d,two heads]\\
\mathbb{Z}^{12}\arrow[u,hook]&0\arrow[l]&\mathbb{Z}^2\arrow[l]
\end{tikzcd}
\]
\noindent where the evaluation map sends \(\mathbb{Z}^{1010}\leq K_0(C_r^\ast(\dot{G}_{AF});C_r^\ast(\dot{G}_u))\) isomorphically onto its image and the rest to \(0\). Note that the substitution \emph{appears to have two components, and is thus not primitive!} We have that
\begin{align*}
K_0(C_r^\ast(\dot{G}_u))&\cong\mathbb{Z}^{1148}/\mathbb{Z}^{1010}\oplus\mathbb{Z}^2\\
&\cong K_0(C_r^\ast(\dot{G}_{AF}))/\mathbb{Z}^{1010}\oplus\mathbb{Z}^2
\end{align*}
\noindent and, since the substitution rule satisfies the boundary hyperplane condition,
\begin{align*}
K_1(C_r^\ast(\dot{G}_u))&\cong\mathbb{Z}^{124}/\mathbb{Z}^{112}\\
&\leq\mathbb{Z}^{1134}/\mathbb{Z}^{112}\\
&\cong K_0(C_r^\ast(\dot{G}_{AF}^{(1)}))/\mathbb{Z}^{112}.
\end{align*}
\indent Among the \(84\) nontrivial eigenvalues contributing to \(\check{H}^2(\Omega_T)\cong K_0(C_r^\ast(\dot{G}_{AF}))/\image\evaluation\leq K_0(C_r^\ast(\dot{G}_u))\), they all belong to \(\mathbb{Z}[\lambda]\), and there is an eigenvalue \(\lambda^2-\lambda-1\) with multiplicity \(8\) that is strictly between \(\lambda\) and \(\lambda^2\).\\
\indent Among the twelve nontrivial eigenvalues contributing to \(K_1(C_r^\ast(\dot{G}_u))\cong\check{H}^1(\Omega_T)\), they multiply to \(1\) with eight expanding and four contracting.
\end{example}
\begin{example}[GKM \(10.1.1.1\)]
This is the substitution rule in \cite[Figure 10]{gahlerkwanmaloney15}, where the enumeration is, from left to right, the first triangle, followed by the first triangle connected by an arc, followed by the first triangle connected by an arc. The expansion factor is \(\lambda=1+2\cos\pi/7\), which is \emph{not a PV number}, thus is weak mixing. Due to the sizes of the induced substitution matrices on \(K_0(C_r^\ast(\dot{G}_{AF,n}))\) and \(K_0(C_r^\ast(\dot{G}_{AF,n});C_r^\ast(\dot{G}_u))\), their eigenvalues had to be computed numerically.\\
\indent The determinant of the induced substitution matrix on \(K_0(C_r^\ast(\dot{G}_{AF}))\), upon removing the eventual kernel, is \(5547\). Fewer than \(27\) of the nonintegral eigenvalues multiply to it, with the remaining multiplying to \(1\). Restricting the direct limit to the associated eigenvectors of the \(27\) eigenvalues that multiply to \(5547\) and computing it over \(\mathbb{Q}\) gives that \(K_0(C_r^\ast(\dot{G}_{AF}))\cong\mathbb{Z}^{\geq 1457}\oplus\mathbb{Q}^{\leq 27}\). These eigenvectors of the \(27\) eigenvalues are mapped isomorphically onto its image in \(\image\iota_\ast\). The induced substitution matrix on \(K_0(C_r^\ast(\dot{G}_{AF});C_r^\ast(\dot{G}_u))\), after removing the eventual kernel, has determinant \(1\).\\
\indent The six-term sequence is
\[
\begin{tikzcd}
\mathbb{Z}^{1331}\arrow[r,"\evaluation"]&\mathbb{Z}^{\geq 1457}\oplus\mathbb{Q}^{\leq 27}\arrow[r,"\iota_\ast"]&\mathbb{Z}^{\geq 133}\oplus\mathbb{Q}^{\leq 27}\arrow[d,two heads]\\
\mathbb{Z}^6\arrow[u,hook]&0\arrow[l]&\mathbb{Z}\arrow[l]
\end{tikzcd}
\]
\noindent where the evaluation map sends \(\mathbb{Z}^{1325}\leq K_0(C_r^\ast(\dot{G}_{AF});C_r^\ast(\dot{G}_u))\) isomorphically onto its image and the rest to \(0\). We have that
\begin{align*}
K_0(C_r^\ast(\dot{G}_u))&\cong(\mathbb{Z}^{\geq 1457}\oplus\mathbb{Q}^{\leq 27})/\mathbb{Z}^{1325}\oplus\mathbb{Z}\\
&\cong K_0(C_r^\ast(\dot{G}_{AF}))/\mathbb{Z}^{1325}\oplus\mathbb{Z}
\end{align*}
\noindent and, since the substitution rule satisfies the boundary hyperplane condition,
\begin{align*}
K_1(C_r^\ast(\dot{G}_u))&\cong\mathbb{Z}^{145}/\mathbb{Z}^{139}\\
&\leq\mathbb{Z}^{1470}/\mathbb{Z}^{139}\\
&\cong K_0(C_r^\ast(\dot{G}_{AF}^{(1)}))/\mathbb{Z}^{139}.
\end{align*}
\indent Among the \(159\) nontrivial eigenvalues contributing to \(\check{H}^2(\Omega_T)\cong K_0(C_r^\ast(\dot{G}_{AF}))/\image\evaluation\leq K_0(C_r^\ast(\dot{G}_u))\), all but \(27\) belong to \(\mathbb{Z}[\lambda]\). The contribution to the determinant of \(5547\) is purely from them. There are three that multiply to the integral value of \(3\),
\begin{align*}
0.7608\cdots\textnormal{, }{}&2.6996\cdots\textnormal{, }\\
{}&\textnormal{roots of }x^3-2x^2-3x+3,\\
1.4605\cdots\textnormal{, }{}&\\
{}&\textnormal{root of }x^3+2x^2-3x-3,
\end{align*}
\noindent with the rest, each with multiplicity \(2\),
\begin{align*}
0.3856\cdots\textnormal{, }{}&0.5451\cdots\textnormal{, }1.1032\cdots\textnormal{, }1.5015\cdots\textnormal{, }2.6495\cdots\textnormal{, }2.7518\cdots\textnormal{, }2.7699\cdots\textnormal{, }\\
{}&\textnormal{roots of }x^{12}-3x^{11}-19x^{10}+64x^9+90x^8-387x^7-66x^6+797x^5-161x^4\\
{}&\phantom{\textnormal{roots of }}-573x^3+166x^2+126x-43,\\
0.6932\cdots\textnormal{, }{}&0.7677\cdots\textnormal{, }1.4062\cdots\textnormal{, }2.3209\cdots\textnormal{, }3.5188\cdots\textnormal{, }\\
{}&\textnormal{roots of }x^{12}+3x^{11}-19x^{10}-64x^9+90x^8+387x^7-66x^6-797x^5-161x^4\\
{}&\phantom{\textnormal{roots of }}+573x^3+166x^2-126x-43.
\end{align*}
\noindent There are two eigenvalues, \(\lambda^2-\lambda-1\) with multiplicity \(4\) and \(3.5188\cdots\) with multiplicity \(2\), that are strictly between \(\lambda\) and \(\lambda^2\).\\
\indent Among the six nontrivial eigenvalues contributing to \(K_1(C_r^\ast(\dot{G}_u))\cong\check{H}^1(\Omega_T)\), they multiply to \(1\) with four expanding and two contracting.
\end{example}
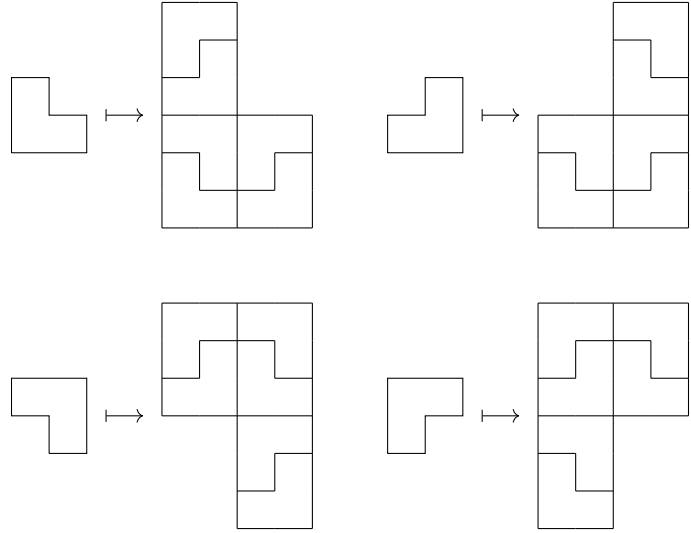
\begin{figure}[t]
\centering
\begin{tikzpicture}
\draw(0,1)--(1,1)--(1,1.5)--(.5,1.5)--(.5,2)--(0,2)--(0,1);
\draw[|->](1.25,1.5)--(1.75,1.5);
\draw(4.0,0)--(3.5,0);
\draw(2.5,2.0)--(2.5,2.5);
\draw(2.5,3.0)--(2.0,3.0);
\draw(3.0,1.5)--(3.0,1.0);
\draw(4.0,1.0)--(4.0,1.5);
\draw(3.5,1.5)--(4.0,1.5);
\draw(2.0,2.0)--(2.0,2.5);
\draw(2.0,1.5)--(2.5,1.5);
\draw(2.0,1.0)--(2.0,1.5);
\draw(3.0,.5)--(3.0,1.0);
\draw(3.0,2.0)--(3.0,2.5);
\draw(2.0,.5)--(2.0,1.0);
\draw(3.5,.5)--(3.5,1.0);
\draw(3.0,0)--(3.0,.5);
\draw(3.0,1.5)--(2.5,1.5);
\draw(2.5,0)--(2.0,0);
\draw(2.0,1.0)--(2.5,1.0);
\draw(2.0,1.5)--(2.0,2.0);
\draw(3.0,1.5)--(3.5,1.5);
\draw(2.5,3.0)--(3.0,3.0);
\draw(2.5,.5)--(2.5,1.0);
\draw(3.0,0)--(3.5,0);
\draw(2.0,2.5)--(2.0,3.0);
\draw(2.5,0)--(3.0,0);
\draw(4.0,0)--(4.0,.5);
\draw(2.0,2.0)--(2.5,2.0);
\draw(2.0,.5)--(2.0,0);
\draw(2.5,.5)--(3.0,.5);
\draw(4.0,.5)--(4.0,1.0);
\draw(3.0,2.5)--(3.0,3.0);
\draw(3.5,.5)--(3.0,.5);
\draw(3.0,2.5)--(2.5,2.5);
\draw(3.0,1.5)--(3.0,2.0);
\draw(3.5,1.0)--(4.0,1.0);
\draw(5.0,1)--(6.0,1)--(6.0,2)--(5.5,2)--(5.5,1.5)--(5.0,1.5)--(5.0,1);
\draw[|->](6.25,1.5)--(6.75,1.5);
\draw(7.0,0)--(7.5,0);
\draw(7.0,0)--(7.0,.5);
\draw(9.0,1.0)--(8.5,1.0);
\draw(8.0,0)--(7.5,0);
\draw(8.0,0)--(8.0,.5);
\draw(7.0,1.0)--(7.0,1.5);
\draw(8.0,1.0)--(8.0,1.5);
\draw(8.0,3.0)--(8.0,2.5);
\draw(9.0,2.0)--(9.0,2.5);
\draw(9.0,1.5)--(8.5,1.5);
\draw(8.5,0)--(9.0,0);
\draw(9.0,1.0)--(9.0,1.5);
\draw(8.0,.5)--(7.5,.5);
\draw(7.5,1.5)--(7.0,1.5);
\draw(9.0,.5)--(9.0,1.0);
\draw(7.5,1.5)--(8.0,1.5);
\draw(7.5,1.0)--(7.0,1.0);
\draw(8.0,1.5)--(8.0,2.0);
\draw(8.5,3.0)--(8.0,3.0);
\draw(8.5,2.5)--(8.0,2.5);
\draw(8.5,1.5)--(8.0,1.5);
\draw(8.0,.5)--(8.5,.5);
\draw(8.0,0)--(8.5,0);
\draw(9.0,1.5)--(9.0,2.0);
\draw(8.5,2.0)--(8.5,2.5);
\draw(8.0,2.0)--(8.0,2.5);
\draw(9.0,2.5)--(9.0,3.0);
\draw(7.0,.5)--(7.0,1.0);
\draw(8.0,.5)--(8.0,1.0);
\draw(8.5,.5)--(8.5,1.0);
\draw(9.0,.5)--(9.0,0);
\draw(8.5,3.0)--(9.0,3.0);
\draw(9.0,2.0)--(8.5,2.0);
\draw(7.5,1.0)--(7.5,.5);
\draw(.5,-3)--(1,-3)--(1,-2)--(0,-2)--(0,-2.5)--(.5,-2.5)--(.5,-3);
\draw[|->](1.25,-2.5)--(1.75,-2.5);
\draw(3.0,-3.0)--(3.0,-2.5);
\draw(4.0,-2.5)--(3.5,-2.5);
\draw(4.0,-3.5)--(4.0,-4.0);
\draw(4.0,-2.5)--(4.0,-2.0);
\draw(2.0,-1.0)--(2.5,-1.0);
\draw(2.5,-2.5)--(2.0,-2.5);
\draw(3.0,-1.0)--(2.5,-1.0);
\draw(3.0,-3.5)--(3.0,-3.0);
\draw(3.5,-3.5)--(3.5,-3.0);
\draw(2.0,-1.5)--(2.0,-1.0);
\draw(3.5,-1.0)--(4.0,-1.0);
\draw(4.0,-3.5)--(4.0,-3.0);
\draw(4.0,-2.0)--(3.5,-2.0);
\draw(3.0,-2.5)--(3.0,-2.0);
\draw(2.5,-1.5)--(3.0,-1.5);
\draw(2.0,-2.0)--(2.5,-2.0);
\draw(2.0,-2.0)--(2.0,-2.5);
\draw(4.0,-4.0)--(3.5,-4.0);
\draw(3.0,-3.5)--(3.5,-3.5);
\draw(3.0,-2.5)--(3.5,-2.5);
\draw(2.0,-1.5)--(2.0,-2.0);
\draw(2.5,-2.0)--(2.5,-1.5);
\draw(3.5,-1.0)--(3.0,-1.0);
\draw(3.5,-1.5)--(3.0,-1.5);
\draw(3.5,-2.0)--(3.5,-1.5);
\draw(4.0,-2.0)--(4.0,-1.5);
\draw(3.5,-4.0)--(3.0,-4.0);
\draw(3.0,-1.5)--(3.0,-2.0);
\draw(4.0,-3.0)--(3.5,-3.0);
\draw(4.0,-3.0)--(4.0,-2.5);
\draw(3.0,-3.5)--(3.0,-4.0);
\draw(3.0,-1.0)--(3.0,-1.5);
\draw(2.5,-2.5)--(3.0,-2.5);
\draw(4.0,-1.5)--(4.0,-1.0);
\draw(5.0,-3)--(5.5,-3)--(5.5,-2.5)--(6.0,-2.5)--(6.0,-2)--(5.0,-2)--(5.0,-3);
\draw[|->](6.25,-2.5)--(6.75,-2.5);
\draw(8.0,-1.5)--(7.5,-1.5);
\draw(9.0,-1.0)--(8.5,-1.0);
\draw(7.0,-3.5)--(7.0,-4.0);
\draw(7.0,-2.5)--(7.0,-2.0);
\draw(8.0,-2.5)--(7.5,-2.5);
\draw(7.5,-4.0)--(7.0,-4.0);
\draw(8.0,-2.5)--(8.0,-2.0);
\draw(7.5,-4.0)--(8.0,-4.0);
\draw(7.0,-3.5)--(7.0,-3.0);
\draw(9.0,-2.0)--(9.0,-2.5);
\draw(7.5,-3.5)--(7.5,-3.0);
\draw(8.5,-2.5)--(9.0,-2.5);
\draw(8.0,-1.5)--(8.5,-1.5);
\draw(8.0,-4.0)--(8.0,-3.5);
\draw(8.5,-1.5)--(8.5,-2.0);
\draw(9.0,-2.0)--(9.0,-1.5);
\draw(7.0,-1.5)--(7.0,-2.0);
\draw(7.5,-1.5)--(7.5,-2.0);
\draw(7.5,-1.0)--(7.0,-1.0);
\draw(8.0,-2.5)--(8.0,-3.0);
\draw(8.0,-2.0)--(8.0,-1.5);
\draw(7.5,-3.5)--(8.0,-3.5);
\draw(7.0,-3.0)--(7.0,-2.5);
\draw(8.0,-1.0)--(8.5,-1.0);
\draw(7.5,-1.0)--(8.0,-1.0);
\draw(7.0,-3.0)--(7.5,-3.0);
\draw(8.0,-2.5)--(8.5,-2.5);
\draw(7.0,-2.5)--(7.5,-2.5);
\draw(8.0,-3.5)--(8.0,-3.0);
\draw(8.0,-1.5)--(8.0,-1.0);
\draw(7.0,-2.0)--(7.5,-2.0);
\draw(9.0,-1.5)--(9.0,-1.0);
\draw(7.0,-1.5)--(7.0,-1.0);
\draw(9.0,-2.0)--(8.5,-2.0);
\end{tikzpicture}
\caption{Substitution rule for a variant of the chair tiling with affine expansion \(\diagonal(2,3)\).}
\label{figure:chair23-affine}
\end{figure}%

\begin{example}[Chair \(2\times 3\) (affine)\footnote{Supplied by Rodrigo Trevi\~no.}]
The substitution rule is a variant of the chair tiling with \emph{affine} expansion \(\diagonal(2,3)\), and is shown in \cref{figure:chair23-affine}. A supertile is shown in \cref{figure:chair23-patch4}. Note that by appropriately grouping pairs of prototiles into rectangles, one can reencode the substitution rule as an affine substitution on rectangles. Each of the induced substitution matrices have integral eigenvalues. One easily checks that the direct limits of the induced substitution matrices split completely.\\
\indent The six-term sequence is
\[
\begin{tikzcd}
\mathbb{Z}[1/3]\oplus\mathbb{Z}[1/2]\oplus\mathbb{Z}^3\arrow[r,"\evaluation"]&\mathbb{Z}[1/6]\oplus\mathbb{Z}[1/2]\oplus\mathbb{Z}^4\arrow[r,"\iota_\ast"]&\mathbb{Z}[1/6]\oplus\mathbb{Z}[1/2]\oplus\mathbb{Z}^2\arrow[d,two heads]\\
\mathbb{Z}[1/3]\oplus\mathbb{Z}[1/2]\arrow[u,hook]&0\arrow[l]&\mathbb{Z}\arrow[l]
\end{tikzcd}
\]
\noindent where the evaluation map sends \(\mathbb{Z}^3\leq K_0(C_r^\ast(\dot{G}_{AF});C_r^\ast(\dot{G}_u))\) isomorphically onto its image and the rest to \(0\). We have that
\begin{align*}
K_0(C_r^\ast(\dot{G}_u))&\cong(\mathbb{Z}[1/6]\oplus\mathbb{Z}[1/2]\oplus\mathbb{Z}^4)/\mathbb{Z}^3\oplus\mathbb{Z}\\
&\cong K_0(C_r^\ast(\dot{G}_{AF}))/\mathbb{Z}^3\oplus\mathbb{Z}
\end{align*}
\noindent and
\begin{align*}
K_1(C_r^\ast(\dot{G}_u))&\cong\mathbb{Z}[1/3]\oplus\mathbb{Z}[1/2]\\
&\leq\mathbb{Z}[1/3]\oplus\mathbb{Z}[1/2]\oplus\mathbb{Z}^3\\
&\cong C^1\\
&\cong K_0(C_r^\ast(\dot{G}_{AF}^{(1)}))\tag*{(\cref{proposition:af-cochain-nontrivial})}
\end{align*}
\noindent where we have to apply a slightly nontrivial isomorphism in the last step since this substitution does not satisfy the boundary hyperplane condition. The proposition applies since the component in \(C^1\) contributed by the \(1\)-cells satisfying the boundary hyperplane condition splits completely, being \(\mathbb{Z}[1/3]\oplus\mathbb{Z}[1/2]\).
\end{example}
\begin{figure}[t]
\centering
\begin{tikzpicture}
\draw(0,1)--(1,1)--(1,1.5)--(.5,1.5)--(.5,2)--(0,2)--(0,1);
\draw[|->](1.25,1.5)--(1.75,1.5);
\draw(3.5,.5)--(3.0,.5);
\draw(4.0,-.5)--(3.5,-.5);
\draw(2.5,1.5)--(2.5,2.0);
\draw(4.5,1.0)--(4.0,1.0);
\draw(3.5,3.5)--(3.0,3.5);
\draw(3.5,3.0)--(3.5,2.5);
\draw(2.5,2.5)--(2.0,2.5);
\draw(5.0,1.0)--(5.0,1.5);
\draw(3.0,1.0)--(3.0,.5);
\draw(3.0,2.5)--(3.5,2.5);
\draw(5.0,.5)--(5.0,1.0);
\draw(5.0,0)--(5.0,.5);
\draw(2.0,1.5)--(2.0,2.0);
\draw(3.5,.5)--(3.5,1.0);
\draw(4.5,1.5)--(5.0,1.5);
\draw(4.0,1.5)--(3.5,1.5);
\draw(2.0,.5)--(2.0,1.0);
\draw(3.0,3.0)--(3.0,3.5);
\draw(4.5,1.0)--(4.5,.5);
\draw(2.5,.5)--(2.5,1.0);
\draw(2.0,3.0)--(2.0,3.5);
\draw(3.5,3.0)--(3.5,3.5);
\draw(2.0,0)--(2.0,.5);
\draw(3.5,1.5)--(3.5,2.0);
\draw(3.0,-.5)--(3.0,0);
\draw(2.5,3.0)--(3.0,3.0);
\draw(3.5,0)--(3.5,.5);
\draw(4.5,0)--(4.5,.5);
\draw(3.0,1.0)--(2.5,1.0);
\draw(4.0,1.5)--(4.0,1.0);
\draw(5.0,0)--(5.0,-.5);
\draw(2.5,3.5)--(3.0,3.5);
\draw(4.5,-.5)--(5.0,-.5);
\draw(2.5,-.5)--(2.0,-.5);
\draw(2.0,.5)--(2.5,.5);
\draw(5.0,.5)--(4.5,.5);
\draw(3.0,1.5)--(2.5,1.5);
\draw(4.0,-.5)--(4.5,-.5);
\draw(2.0,1.0)--(2.0,1.5);
\draw(4.0,1.5)--(4.5,1.5);
\draw(2.5,2.5)--(3.0,2.5);
\draw(2.5,0)--(2.5,.5);
\draw(2.0,3.5)--(2.5,3.5);
\draw(3.0,-.5)--(3.5,-.5);
\draw(2.0,2.0)--(2.0,2.5);
\draw(2.5,-.5)--(3.0,-.5);
\draw(4.0,-.5)--(4.0,0);
\draw(3.5,0)--(3.0,0);
\draw(2.0,1.5)--(2.5,1.5);
\draw(2.0,0)--(2.0,-.5);
\draw(2.5,0)--(3.0,0);
\draw(4.0,0)--(4.0,.5);
\draw(2.5,2.5)--(2.5,3.0);
\draw(4.0,0)--(4.5,0);
\draw(3.5,2.0)--(3.5,2.5);
\draw(3.0,1.5)--(3.5,1.5);
\draw(2.0,3.0)--(2.0,2.5);
\draw(3.0,2.0)--(3.0,2.5);
\draw(3.5,1.0)--(4.0,1.0);
\draw(3.0,2.0)--(2.5,2.0);
\draw(3.0,1.0)--(3.0,1.5);
\draw(3.5,.5)--(4.0,.5);
\draw(6.0,1)--(7.0,1)--(7.0,2)--(6.5,2)--(6.5,1.5)--(6.0,1.5)--(6.0,1);
\draw[|->](7.25,1.5)--(7.75,1.5);
\draw(8.5,1.0)--(9.0,1.0);
\draw(8.0,0)--(8.0,.5);
\draw(8.0,1.5)--(8.5,1.5);
\draw(9.0,-.5)--(9.5,-.5);
\draw(9.0,-.5)--(9.0,0);
\draw(11.0,.5)--(10.5,.5);
\draw(10.0,-.5)--(9.5,-.5);
\draw(10.0,-.5)--(10.0,0);
\draw(9.5,3.5)--(10.0,3.5);
\draw(9.5,.5)--(9.5,1.0);
\draw(8.5,-.5)--(9.0,-.5);
\draw(10.5,.5)--(10.5,1.0);
\draw(10.0,.5)--(10.0,1.0);
\draw(9.5,2.0)--(9.5,2.5);
\draw(10.0,2.5)--(10.0,2.0);
\draw(11.0,1.5)--(11.0,2.0);
\draw(9.5,.5)--(10.0,.5);
\draw(8.0,1.0)--(8.0,.5);
\draw(8.0,-.5)--(8.0,0);
\draw(10.5,-.5)--(11.0,-.5);
\draw(11.0,.5)--(11.0,1.0);
\draw(10.0,0)--(9.5,0);
\draw(9.5,1.0)--(9.0,1.0);
\draw(9.0,1.5)--(9.0,1.0);
\draw(11.0,3.0)--(11.0,3.5);
\draw(11.0,0)--(11.0,.5);
\draw(8.5,.5)--(8.0,.5);
\draw(10.0,3.0)--(10.0,3.5);
\draw(9.5,.5)--(9.0,.5);
\draw(10.0,1.0)--(10.0,1.5);
\draw(8.0,1.0)--(8.0,1.5);
\draw(8.5,-.5)--(8.0,-.5);
\draw(10.0,2.0)--(10.5,2.0);
\draw(10.5,1.5)--(10.0,1.5);
\draw(10.0,2.5)--(10.5,2.5);
\draw(10.5,3.0)--(10.0,3.0);
\draw(10.5,1.0)--(10.0,1.0);
\draw(10.5,2.5)--(10.5,3.0);
\draw(10.5,3.5)--(10.0,3.5);
\draw(11.0,1.0)--(11.0,1.5);
\draw(10.5,0)--(10.0,0);
\draw(10.0,-.5)--(10.5,-.5);
\draw(9.0,1.5)--(9.5,1.5);
\draw(10.0,2.5)--(9.5,2.5);
\draw(9.5,1.5)--(9.5,2.0);
\draw(9.0,1.5)--(8.5,1.5);
\draw(9.5,1.5)--(10.0,1.5);
\draw(10.5,1.5)--(10.5,2.0);
\draw(11.0,2.0)--(11.0,2.5);
\draw(9.0,0)--(9.0,.5);
\draw(11.0,0)--(11.0,-.5);
\draw(10.5,0)--(10.5,.5);
\draw(11.0,3.5)--(10.5,3.5);
\draw(10.5,2.5)--(11.0,2.5);
\draw(11.0,1.5)--(10.5,1.5);
\draw(9.5,2.5)--(9.5,3.0);
\draw(9.5,3.5)--(9.5,3.0);
\draw(11.0,3.0)--(11.0,2.5);
\draw(8.5,.5)--(8.5,1.0);
\draw(9.5,.5)--(9.5,0);
\draw(8.5,0)--(8.5,.5);
\draw(8.5,0)--(9.0,0);
\draw(.5,-4)--(1,-4)--(1,-3)--(0,-3)--(0,-3.5)--(.5,-3.5)--(.5,-4);
\draw[|->](1.25,-3.5)--(1.75,-3.5);
\draw(4.0,-3.5)--(4.0,-3.0);
\draw(5.0,-5.5)--(5.0,-5.0);
\draw(3.5,-3.0)--(3.5,-2.5);
\draw(5.0,-4.0)--(5.0,-4.5);
\draw(5.0,-3.0)--(5.0,-2.5);
\draw(3.0,-1.5)--(3.5,-1.5);
\draw(3.5,-3.0)--(3.0,-3.0);
\draw(4.5,-5.0)--(4.5,-4.5);
\draw(4.0,-3.5)--(3.5,-3.5);
\draw(4.0,-1.5)--(3.5,-1.5);
\draw(3.0,-2.0)--(2.5,-2.0);
\draw(3.5,-4.5)--(3.5,-4.0);
\draw(2.0,-3.0)--(2.0,-2.5);
\draw(4.0,-5.5)--(4.0,-5.0);
\draw(4.5,-4.0)--(4.5,-3.5);
\draw(2.5,-1.5)--(3.0,-1.5);
\draw(3.0,-2.0)--(3.0,-1.5);
\draw(2.5,-3.0)--(2.5,-2.5);
\draw(2.5,-3.5)--(3.0,-3.5);
\draw(3.5,-4.5)--(4.0,-4.5);
\draw(4.5,-1.5)--(5.0,-1.5);
\draw(5.0,-5.0)--(5.0,-4.5);
\draw(5.0,-4.0)--(5.0,-3.5);
\draw(5.0,-2.5)--(4.5,-2.5);
\draw(4.0,-3.0)--(4.0,-2.5);
\draw(3.5,-2.5)--(4.0,-2.5);
\draw(3.0,-2.5)--(3.5,-2.5);
\draw(4.5,-3.0)--(4.5,-2.5);
\draw(3.5,-2.0)--(4.0,-2.0);
\draw(5.0,-4.5)--(4.5,-4.5);
\draw(3.5,-5.5)--(3.5,-5.0);
\draw(3.5,-5.0)--(3.5,-4.5);
\draw(2.5,-1.5)--(2.0,-1.5);
\draw(4.5,-1.5)--(4.0,-1.5);
\draw(4.5,-2.0)--(4.0,-2.0);
\draw(4.5,-2.5)--(4.5,-2.0);
\draw(5.0,-2.5)--(5.0,-2.0);
\draw(3.0,-2.5)--(3.0,-2.0);
\draw(2.5,-2.5)--(2.5,-2.0);
\draw(4.0,-5.0)--(4.5,-5.0);
\draw(4.5,-5.5)--(4.0,-5.5);
\draw(5.0,-5.5)--(4.5,-5.5);
\draw(3.5,-5.5)--(4.0,-5.5);
\draw(3.0,-3.5)--(3.5,-3.5);
\draw(5.0,-3.5)--(4.5,-3.5);
\draw(5.0,-3.5)--(5.0,-3.0);
\draw(3.5,-4.0)--(3.5,-3.5);
\draw(4.0,-4.0)--(4.0,-4.5);
\draw(3.0,-3.5)--(3.0,-3.0);
\draw(2.5,-3.5)--(2.0,-3.5);
\draw(4.0,-1.5)--(4.0,-2.0);
\draw(2.0,-2.0)--(2.0,-1.5);
\draw(2.0,-2.5)--(2.5,-2.5);
\draw(2.0,-3.0)--(2.0,-3.5);
\draw(4.0,-4.5)--(4.5,-4.5);
\draw(5.0,-2.0)--(5.0,-1.5);
\draw(4.0,-3.0)--(4.5,-3.0);
\draw(4.0,-3.5)--(4.5,-3.5);
\draw(4.0,-4.0)--(4.5,-4.0);
\draw(2.0,-2.0)--(2.0,-2.5);
\draw(2.5,-3.0)--(3.0,-3.0);
\draw(3.5,-2.5)--(3.5,-2.0);
\draw(6.0,-4)--(6.5,-4)--(6.5,-3.5)--(7.0,-3.5)--(7.0,-3)--(6.0,-3)--(6.0,-4);
\draw[|->](7.25,-3.5)--(7.75,-3.5);
\draw(9.0,-2.0)--(8.5,-2.0);
\draw(10.0,-1.5)--(9.5,-1.5);
\draw(8.0,-5.5)--(8.0,-5.0);
\draw(11.0,-3.0)--(11.0,-2.5);
\draw(9.0,-2.0)--(9.5,-2.0);
\draw(8.0,-3.0)--(8.0,-2.5);
\draw(8.0,-4.0)--(8.0,-4.5);
\draw(9.0,-5.5)--(9.0,-5.0);
\draw(8.0,-5.5)--(8.5,-5.5);
\draw(9.0,-3.0)--(8.5,-3.0);
\draw(8.5,-4.5)--(8.0,-4.5);
\draw(9.5,-3.5)--(10.0,-3.5);
\draw(9.0,-3.0)--(9.0,-2.5);
\draw(8.5,-4.5)--(9.0,-4.5);
\draw(10.0,-2.0)--(10.5,-2.0);
\draw(10.0,-2.5)--(9.5,-2.5);
\draw(8.0,-5.0)--(8.0,-4.5);
\draw(8.0,-4.0)--(8.0,-3.5);
\draw(10.5,-3.0)--(10.0,-3.0);
\draw(9.5,-4.0)--(9.5,-3.5);
\draw(8.5,-4.0)--(8.5,-3.5);
\draw(11.0,-3.5)--(10.5,-3.5);
\draw(9.5,-3.0)--(10.0,-3.0);
\draw(9.0,-2.5)--(9.5,-2.5);
\draw(9.0,-4.5)--(9.0,-4.0);
\draw(9.5,-3.0)--(9.5,-2.5);
\draw(8.0,-2.0)--(8.0,-2.5);
\draw(9.0,-3.0)--(9.0,-3.5);
\draw(8.5,-2.0)--(8.5,-2.5);
\draw(8.5,-1.5)--(8.0,-1.5);
\draw(10.5,-2.5)--(10.5,-2.0);
\draw(10.0,-2.0)--(10.0,-2.5);
\draw(9.5,-2.5)--(9.5,-2.0);
\draw(8.5,-4.5)--(8.5,-5.0);
\draw(8.5,-4.0)--(9.0,-4.0);
\draw(10.0,-1.5)--(10.5,-1.5);
\draw(8.5,-3.0)--(8.5,-2.5);
\draw(8.0,-3.5)--(8.0,-3.0);
\draw(10.5,-3.0)--(10.5,-2.5);
\draw(9.5,-4.0)--(9.5,-4.5);
\draw(10.5,-1.5)--(11.0,-1.5);
\draw(9.0,-1.5)--(9.5,-1.5);
\draw(9.0,-4.5)--(9.5,-4.5);
\draw(8.5,-1.5)--(9.0,-1.5);
\draw(9.5,-5.0)--(9.5,-4.5);
\draw(8.0,-3.5)--(8.5,-3.5);
\draw(9.5,-5.5)--(9.5,-5.0);
\draw(9.5,-3.5)--(9.0,-3.5);
\draw(10.0,-3.5)--(10.0,-3.0);
\draw(10.0,-3.5)--(10.5,-3.5);
\draw(8.5,-3.5)--(9.0,-3.5);
\draw(10.5,-2.5)--(11.0,-2.5);
\draw(11.0,-3.5)--(11.0,-3.0);
\draw(9.0,-2.0)--(9.0,-1.5);
\draw(8.0,-2.5)--(8.5,-2.5);
\draw(11.0,-2.0)--(11.0,-2.5);
\draw(8.0,-2.0)--(8.0,-1.5);
\draw(10.0,-2.0)--(10.0,-1.5);
\draw(11.0,-2.0)--(11.0,-1.5);
\draw(9.0,-5.0)--(8.5,-5.0);
\draw(9.0,-5.5)--(9.5,-5.5);
\draw(9.0,-5.5)--(8.5,-5.5);
\end{tikzpicture}
\caption{Substitution rule for a variant of the chair tiling with affine expansion \(\diagonal(3,4)\).}
\label{figure:chair34-affine}
\end{figure}
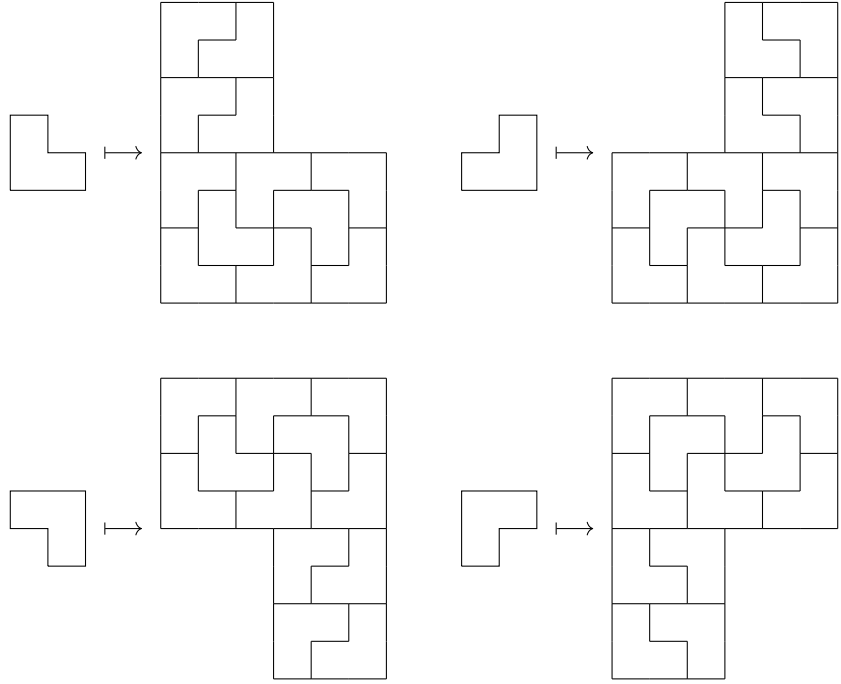%

\begin{example}[Chair \(3\times 4\) (affine)]
The substitution rule is a variant of the chair tiling with \emph{affine} expansion \(\diagonal(3,4)\), and is shown in \cref{figure:chair34-affine}. A supertile is shown in \cref{figure:chair34-patch3}. Each of the induced substitution matrices have integral eigenvalues. Due to the sizes of the induced matrices, we assume each of the direct limits splits completely.\\
\indent The six-term sequence is
\[
\begin{tikzcd}
\begin{tabular}{@{}c@{}}\(\mathbb{Z}[1/4]\oplus\mathbb{Z}[1/3]^2\)\\\(\oplus\mathbb{Z}[1/2]^3\oplus\mathbb{Z}^{16}\)\end{tabular}\arrow[r,"\evaluation"]&\begin{tabular}{@{}c@{}}\(\mathbb{Z}[1/12]\oplus\mathbb{Z}[1/4]\)\\\(\oplus\mathbb{Z}[1/3]^2\oplus\mathbb{Z}[1/2]^4\oplus\mathbb{Z}^{22}\)\end{tabular}\arrow[r,"\iota_\ast"]&\begin{tabular}{@{}c@{}}\(\mathbb{Z}[1/12]\oplus\mathbb{Z}[1/4]\)\\\(\oplus\mathbb{Z}[1/3]\oplus\mathbb{Z}[1/2]\oplus\mathbb{Z}^7\)\end{tabular}\arrow[d,two heads]\\
\mathbb{Z}[1/4]\oplus\mathbb{Z}[1/3]\arrow[u,hook]&0\arrow[l]&\mathbb{Z}\arrow[l]
\end{tikzcd}
\]
\noindent where the evaluation map sends \(\mathbb{Z}[1/3]\oplus\mathbb{Z}[1/2]^3\oplus\mathbb{Z}^{16}\leq K_0(C_r^\ast(\dot{G}_{AF});C_r^\ast(\dot{G}_u))\) isomorphically onto its image and the rest to \(0\). We have that
\begin{align*}
K_0(C_r^\ast(\dot{G}_u))&\cong\left.\left(\begin{tabular}{@{}c@{}}\(\mathbb{Z}[1/12]\oplus\mathbb{Z}[1/4]\)\\\(\oplus\mathbb{Z}[1/3]^2\oplus\mathbb{Z}[1/2]^4\oplus\mathbb{Z}^{22}\)\end{tabular}\right)\middle/(\mathbb{Z}[1/3]\oplus\mathbb{Z}[1/2]^3\oplus\mathbb{Z}^{16})\right.\oplus\mathbb{Z}\\
&\cong K_0(C_r^\ast(\dot{G}_{AF}))/(\mathbb{Z}[1/3]\oplus\mathbb{Z}[1/2]^3\oplus\mathbb{Z}^{16})\oplus\mathbb{Z}
\end{align*}
\noindent and
\begin{align*}
K_1(C_r^\ast(\dot{G}_u))&\cong(\mathbb{Z}[1/4]\oplus\mathbb{Z}[1/3]\oplus\mathbb{Z}^4)/\mathbb{Z}^4\\
&\leq(\mathbb{Z}[1/4]\oplus\mathbb{Z}[1/3]^2\oplus\mathbb{Z}[1/2]^3\oplus\mathbb{Z}^{20})/\mathbb{Z}^4\\
&\cong C^1/\mathbb{Z}^4\\
&\cong K_0(C_r^\ast(\dot{G}_{AF}^{(1)}))/\mathbb{Z}^4\tag*{(\cref{proposition:af-cochain-nontrivial})}
\end{align*}
\noindent where we have to apply a slightly nontrivial isomorphism in the last step since this substitution does not satisfy the boundary hyperplane condition. The proposition applies since we assumed that \(C^1\) splits completely.
\end{example}
\begin{figure}[t]
\centering
\begin{tikzpicture}
\draw(0,1)--(1,1)--(1,1.5)--(.5,1.5)--(.5,2)--(0,2)--(0,1);
\draw[|->](1.25,1.5)--(1.75,1.5);
\draw(3.5,0)--(3.0,0);
\draw(4.0,-1.0)--(3.5,-1.0);
\draw(2.5,1.0)--(2.5,1.5);
\draw(4.5,1.5)--(5.0,1.5);
\draw(2.5,.5)--(2.5,1.0);
\draw(2.5,-1.0)--(2.5,-.5);
\draw(2.0,3.5)--(2.0,4.0);
\draw(3.5,3.0)--(3.0,3.0);
\draw(2.5,0)--(3.0,0);
\draw(3.5,-.5)--(3.5,-1.0);
\draw(3.5,2.5)--(3.5,2.0);
\draw(2.5,2.0)--(2.0,2.0);
\draw(4.5,.5)--(5.0,.5);
\draw(5.0,.5)--(5.0,1.0);
\draw(3.0,.5)--(3.0,0);
\draw(3.0,2.0)--(3.5,2.0);
\draw(4.5,.5)--(4.5,1.0);
\draw(4.0,1.5)--(4.5,1.5);
\draw(5.0,0)--(5.0,.5);
\draw(3.5,3.5)--(3.5,4.0);
\draw(4.0,0)--(4.0,.5);
\draw(5.0,-.5)--(5.0,0);
\draw(3.0,2.0)--(3.0,2.5);
\draw(2.0,1.0)--(2.0,1.5);
\draw(3.5,3.0)--(3.5,3.5);
\draw(2.5,2.5)--(2.5,3.0);
\draw(3.5,.5)--(3.5,1.0);
\draw(3.0,4.0)--(2.5,4.0);
\draw(4.0,1.0)--(3.5,1.0);
\draw(2.0,3.0)--(2.0,3.5);
\draw(2.0,0)--(2.0,.5);
\draw(3.0,-.5)--(3.0,0);
\draw(4.0,0)--(4.5,0);
\draw(3.0,1.0)--(3.0,1.5);
\draw(2.0,2.5)--(2.0,3.0);
\draw(3.5,2.5)--(3.5,3.0);
\draw(2.0,-.5)--(2.0,0);
\draw(3.0,3.5)--(2.5,3.5);
\draw(3.5,-.5)--(3.5,0);
\draw(2.5,2.5)--(3.0,2.5);
\draw(4.5,-.5)--(4.5,0);
\draw(4.0,1.0)--(4.0,1.5);
\draw(3.0,.5)--(2.5,.5);
\draw(4.5,-1.0)--(5.0,-1.0);
\draw(5.0,-.5)--(5.0,-1.0);
\draw(2.5,3.0)--(3.0,3.0);
\draw(4.0,1.0)--(4.0,.5);
\draw(2.5,-1.0)--(2.0,-1.0);
\draw(2.0,0)--(2.5,0);
\draw(5.0,0)--(4.5,0);
\draw(4.0,-1.0)--(4.5,-1.0);
\draw(3.0,.5)--(3.5,.5);
\draw(3.0,4.0)--(3.5,4.0);
\draw(4.0,1.5)--(3.5,1.5);
\draw(2.0,.5)--(2.0,1.0);
\draw(2.0,3.0)--(2.5,3.0);
\draw(2.5,2.0)--(3.0,2.0);
\draw(3.0,-1.0)--(3.5,-1.0);
\draw(4.0,1.0)--(4.5,1.0);
\draw(2.0,1.5)--(2.0,2.0);
\draw(2.5,-1.0)--(3.0,-1.0);
\draw(4.0,-1.0)--(4.0,-.5);
\draw(2.0,-.5)--(2.0,-1.0);
\draw(2.0,1.0)--(2.5,1.0);
\draw(2.5,-.5)--(3.0,-.5);
\draw(3.0,3.5)--(3.0,3.0);
\draw(4.0,-.5)--(4.5,-.5);
\draw(5.0,1.0)--(5.0,1.5);
\draw(3.5,1.5)--(3.5,2.0);
\draw(3.0,1.0)--(3.5,1.0);
\draw(2.5,3.5)--(2.5,4.0);
\draw(2.0,2.5)--(2.0,2.0);
\draw(3.0,1.5)--(3.0,2.0);
\draw(2.0,4.0)--(2.5,4.0);
\draw(3.0,1.5)--(2.5,1.5);
\draw(3.5,0)--(4.0,0);
\draw(6.0,1)--(7.0,1)--(7.0,2)--(6.5,2)--(6.5,1.5)--(6.0,1.5)--(6.0,1);
\draw[|->](7.25,1.5)--(7.75,1.5);
\draw(8.0,-.5)--(8.0,0);
\draw(11.0,3.5)--(11.0,4.0);
\draw(9.0,-1.0)--(9.5,-1.0);
\draw(11.0,0)--(10.5,0);
\draw(9.0,-1.0)--(9.0,-.5);
\draw(10.0,-1.0)--(9.5,-1.0);
\draw(9.5,3.5)--(9.5,4.0);
\draw(10.5,-1.0)--(10.5,-.5);
\draw(9.0,0)--(9.0,.5);
\draw(9.5,3.0)--(10.0,3.0);
\draw(8.5,-.5)--(8.5,0);
\draw(9.5,-1.0)--(9.5,-.5);
\draw(8.5,-1.0)--(9.0,-1.0);
\draw(10.0,0)--(10.0,.5);
\draw(8.0,1.5)--(8.5,1.5);
\draw(9.5,0)--(10.0,0);
\draw(8.0,-1.0)--(8.0,-.5);
\draw(10.0,2.0)--(10.0,1.5);
\draw(9.5,1.5)--(9.5,2.0);
\draw(9.5,4.0)--(10.0,4.0);
\draw(8.0,.5)--(8.0,0);
\draw(11.0,1.0)--(11.0,1.5);
\draw(10.5,-1.0)--(11.0,-1.0);
\draw(11.0,3.0)--(11.0,3.5);
\draw(8.0,1.0)--(8.0,1.5);
\draw(11.0,0)--(11.0,.5);
\draw(9.0,1.0)--(9.0,1.5);
\draw(9.0,1.0)--(9.0,.5);
\draw(11.0,2.5)--(11.0,3.0);
\draw(11.0,-.5)--(11.0,0);
\draw(8.5,.5)--(8.5,1.0);
\draw(9.5,.5)--(10.0,.5);
\draw(10.5,2.5)--(10.5,3.0);
\draw(9.5,0)--(9.0,0);
\draw(10.5,.5)--(10.5,1.0);
\draw(8.0,.5)--(8.0,1.0);
\draw(10.0,3.0)--(10.0,3.5);
\draw(10.5,-.5)--(10.0,-.5);
\draw(10.0,-1.0)--(10.5,-1.0);
\draw(10.0,2.5)--(10.5,2.5);
\draw(10.5,3.0)--(10.0,3.0);
\draw(10.0,.5)--(10.5,.5);
\draw(10.0,2.0)--(10.0,2.5);
\draw(10.0,0)--(10.5,0);
\draw(11.0,.5)--(11.0,1.0);
\draw(10.5,3.5)--(10.0,3.5);
\draw(10.0,4.0)--(10.5,4.0);
\draw(9.0,1.0)--(9.5,1.0);
\draw(10.0,2.0)--(9.5,2.0);
\draw(10.5,3.5)--(10.5,4.0);
\draw(9.0,1.0)--(8.5,1.0);
\draw(9.5,3.0)--(9.5,3.5);
\draw(11.0,4.0)--(10.5,4.0);
\draw(9.5,1.0)--(10.0,1.0);
\draw(8.0,.5)--(8.5,.5);
\draw(10.5,1.0)--(10.5,1.5);
\draw(10.0,1.0)--(10.0,1.5);
\draw(11.0,1.5)--(11.0,2.0);
\draw(10.0,-.5)--(10.0,0);
\draw(8.5,0)--(9.0,0);
\draw(11.0,-.5)--(11.0,-1.0);
\draw(9.5,.5)--(9.5,1.0);
\draw(11.0,3.0)--(10.5,3.0);
\draw(10.5,2.0)--(11.0,2.0);
\draw(9.0,1.5)--(9.5,1.5);
\draw(11.0,1.0)--(10.5,1.0);
\draw(9.5,2.0)--(9.5,2.5);
\draw(9.5,3.0)--(9.5,2.5);
\draw(11.0,2.5)--(11.0,2.0);
\draw(8.5,0)--(8.0,0);
\draw(9.5,0)--(9.5,-.5);
\draw(10.0,2.0)--(10.5,2.0);
\draw(10.0,1.5)--(10.5,1.5);
\draw(9.0,1.5)--(8.5,1.5);
\draw(8.5,-.5)--(9.0,-.5);
\draw(8.5,-1.0)--(8.0,-1.0);
\draw(.5,-5)--(1,-5)--(1,-4)--(0,-4)--(0,-4.5)--(.5,-4.5)--(.5,-5);
\draw[|->](1.25,-4.5)--(1.75,-4.5);
\draw(4.5,-4.0)--(4.5,-3.5);
\draw(5.0,-6.0)--(5.0,-5.5);
\draw(5.0,-4.5)--(5.0,-5.0);
\draw(5.0,-3.5)--(5.0,-3.0);
\draw(2.5,-4.5)--(2.0,-4.5);
\draw(3.0,-2.0)--(3.5,-2.0);
\draw(4.0,-5.5)--(4.0,-5.0);
\draw(4.0,-4.0)--(3.5,-4.0);
\draw(4.0,-2.0)--(3.5,-2.0);
\draw(3.5,-6.5)--(3.5,-6.0);
\draw(2.5,-4.5)--(3.0,-4.5);
\draw(3.5,-3.5)--(3.5,-4.0);
\draw(3.0,-3.0)--(2.5,-3.0);
\draw(2.0,-3.5)--(2.0,-3.0);
\draw(3.0,-2.5)--(2.5,-2.5);
\draw(4.5,-6.0)--(4.5,-5.5);
\draw(3.5,-5.0)--(3.5,-4.5);
\draw(2.0,-3.5)--(2.5,-3.5);
\draw(4.5,-7.0)--(4.5,-6.5);
\draw(4.0,-4.5)--(4.0,-4.0);
\draw(4.5,-4.5)--(4.5,-4.0);
\draw(2.5,-2.0)--(3.0,-2.0);
\draw(3.0,-4.5)--(3.0,-4.0);
\draw(4.0,-7.0)--(3.5,-7.0);
\draw(3.0,-4.5)--(3.5,-4.5);
\draw(2.5,-4.0)--(3.0,-4.0);
\draw(3.0,-2.5)--(3.0,-2.0);
\draw(5.0,-7.0)--(4.5,-7.0);
\draw(3.5,-5.0)--(4.0,-5.0);
\draw(4.5,-2.0)--(5.0,-2.0);
\draw(5.0,-5.5)--(5.0,-5.0);
\draw(5.0,-4.5)--(5.0,-4.0);
\draw(5.0,-6.5)--(5.0,-6.0);
\draw(5.0,-3.0)--(4.5,-3.0);
\draw(4.0,-3.5)--(4.0,-3.0);
\draw(3.5,-3.0)--(4.0,-3.0);
\draw(3.0,-3.0)--(3.5,-3.0);
\draw(2.0,-4.5)--(2.0,-4.0);
\draw(3.0,-3.0)--(3.0,-3.5);
\draw(2.5,-4.0)--(2.5,-3.5);
\draw(5.0,-5.0)--(4.5,-5.0);
\draw(3.5,-6.0)--(3.5,-5.5);
\draw(3.5,-5.5)--(3.5,-5.0);
\draw(2.5,-2.0)--(2.0,-2.0);
\draw(4.5,-6.5)--(4.0,-6.5);
\draw(4.0,-6.0)--(4.5,-6.0);
\draw(4.5,-7.0)--(4.0,-7.0);
\draw(5.0,-6.0)--(4.5,-6.0);
\draw(4.0,-5.0)--(4.5,-5.0);
\draw(4.5,-5.5)--(4.0,-5.5);
\draw(4.0,-3.5)--(4.5,-3.5);
\draw(4.0,-4.5)--(4.5,-4.5);
\draw(5.0,-3.0)--(5.0,-2.5);
\draw(3.5,-6.0)--(4.0,-6.0);
\draw(3.0,-4.0)--(3.5,-4.0);
\draw(5.0,-4.0)--(4.5,-4.0);
\draw(5.0,-4.0)--(5.0,-3.5);
\draw(4.0,-4.5)--(4.0,-5.0);
\draw(5.0,-7.0)--(5.0,-6.5);
\draw(3.0,-4.0)--(3.0,-3.5);
\draw(4.0,-6.5)--(4.0,-6.0);
\draw(4.5,-2.0)--(4.5,-2.5);
\draw(2.0,-3.0)--(2.5,-3.0);
\draw(3.5,-3.5)--(4.0,-3.5);
\draw(2.0,-2.5)--(2.0,-2.0);
\draw(2.0,-3.5)--(2.0,-4.0);
\draw(4.0,-2.5)--(4.5,-2.5);
\draw(5.0,-2.5)--(5.0,-2.0);
\draw(4.0,-3.0)--(4.0,-2.5);
\draw(4.5,-3.0)--(4.0,-3.0);
\draw(2.0,-2.5)--(2.0,-3.0);
\draw(2.5,-2.5)--(2.5,-3.0);
\draw(3.5,-3.0)--(3.5,-2.5);
\draw(3.5,-2.5)--(3.5,-2.0);
\draw(4.5,-2.0)--(4.0,-2.0);
\draw(3.5,-6.5)--(3.5,-7.0);
\draw(6.0,-5)--(6.5,-5)--(6.5,-4.5)--(7.0,-4.5)--(7.0,-4)--(6.0,-4)--(6.0,-5);
\draw[|->](7.25,-4.5)--(7.75,-4.5);
\draw(9.0,-3.0)--(8.5,-3.0);
\draw(10.0,-2.0)--(9.5,-2.0);
\draw(9.0,-2.5)--(8.5,-2.5);
\draw(8.0,-7.0)--(8.5,-7.0);
\draw(9.5,-2.5)--(9.5,-2.0);
\draw(8.0,-6.0)--(8.0,-5.5);
\draw(8.0,-4.5)--(8.0,-5.0);
\draw(11.0,-3.5)--(11.0,-3.0);
\draw(8.0,-3.5)--(8.0,-3.0);
\draw(11.0,-4.5)--(11.0,-4.0);
\draw(11.0,-4.5)--(10.5,-4.5);
\draw(9.0,-3.5)--(8.5,-3.5);
\draw(8.0,-6.0)--(8.5,-6.0);
\draw(8.5,-4.0)--(8.5,-3.5);
\draw(8.5,-6.0)--(8.5,-5.5);
\draw(8.5,-5.0)--(8.0,-5.0);
\draw(9.5,-4.0)--(10.0,-4.0);
\draw(9.5,-7.0)--(9.5,-6.5);
\draw(9.0,-3.5)--(9.0,-3.0);
\draw(8.5,-5.0)--(9.0,-5.0);
\draw(9.5,-7.0)--(9.0,-7.0);
\draw(8.5,-7.0)--(9.0,-7.0);
\draw(9.0,-5.5)--(9.0,-5.0);
\draw(10.0,-3.0)--(10.5,-3.0);
\draw(10.0,-2.5)--(10.5,-2.5);
\draw(8.0,-5.5)--(8.0,-5.0);
\draw(8.0,-4.5)--(8.0,-4.0);
\draw(9.5,-4.0)--(9.5,-3.5);
\draw(10.5,-3.5)--(10.5,-4.0);
\draw(10.0,-3.0)--(10.0,-3.5);
\draw(8.0,-6.5)--(8.0,-6.0);
\draw(8.5,-4.5)--(8.5,-4.0);
\draw(10.5,-3.5)--(11.0,-3.5);
\draw(10.0,-4.5)--(10.0,-4.0);
\draw(8.5,-7.0)--(8.5,-6.5);
\draw(9.0,-3.0)--(9.5,-3.0);
\draw(9.0,-5.0)--(9.0,-4.5);
\draw(10.5,-2.0)--(11.0,-2.0);
\draw(10.5,-3.0)--(10.5,-2.5);
\draw(10.0,-4.5)--(10.5,-4.5);
\draw(9.5,-2.5)--(9.5,-3.0);
\draw(9.0,-3.0)--(9.0,-2.5);
\draw(8.5,-2.0)--(8.0,-2.0);
\draw(8.5,-2.0)--(8.5,-2.5);
\draw(8.0,-2.5)--(8.0,-3.0);
\draw(8.5,-4.5)--(9.0,-4.5);
\draw(10.0,-2.0)--(10.5,-2.0);
\draw(11.0,-2.5)--(11.0,-3.0);
\draw(9.5,-6.5)--(9.5,-6.0);
\draw(9.5,-4.5)--(9.5,-5.0);
\draw(8.0,-4.0)--(8.0,-3.5);
\draw(8.0,-7.0)--(8.0,-6.5);
\draw(9.0,-2.0)--(9.5,-2.0);
\draw(8.5,-6.5)--(9.0,-6.5);
\draw(9.0,-5.0)--(9.5,-5.0);
\draw(8.5,-2.0)--(9.0,-2.0);
\draw(9.5,-5.5)--(9.5,-5.0);
\draw(8.0,-4.0)--(8.5,-4.0);
\draw(9.0,-3.5)--(9.5,-3.5);
\draw(9.5,-4.5)--(10.0,-4.5);
\draw(9.5,-4.0)--(9.0,-4.0);
\draw(10.0,-4.0)--(10.5,-4.0);
\draw(10.0,-4.0)--(10.0,-3.5);
\draw(9.5,-6.0)--(9.5,-5.5);
\draw(9.0,-6.0)--(9.5,-6.0);
\draw(9.0,-4.5)--(9.0,-4.0);
\draw(10.5,-3.0)--(11.0,-3.0);
\draw(11.0,-4.0)--(11.0,-3.5);
\draw(9.0,-5.5)--(8.5,-5.5);
\draw(8.0,-3.0)--(8.5,-3.0);
\draw(9.0,-6.0)--(9.0,-6.5);
\draw(10.0,-2.5)--(10.0,-2.0);
\draw(8.0,-2.5)--(8.0,-2.0);
\draw(10.0,-3.0)--(9.5,-3.0);
\draw(11.0,-2.5)--(11.0,-2.0);
\draw(9.0,-6.0)--(8.5,-6.0);
\end{tikzpicture}
\caption{Substitution rule for a variant of the chair tiling with affine expansion \(\diagonal(3,5)\).}
\label{figure:chair35-variant-1-affine}
\end{figure}
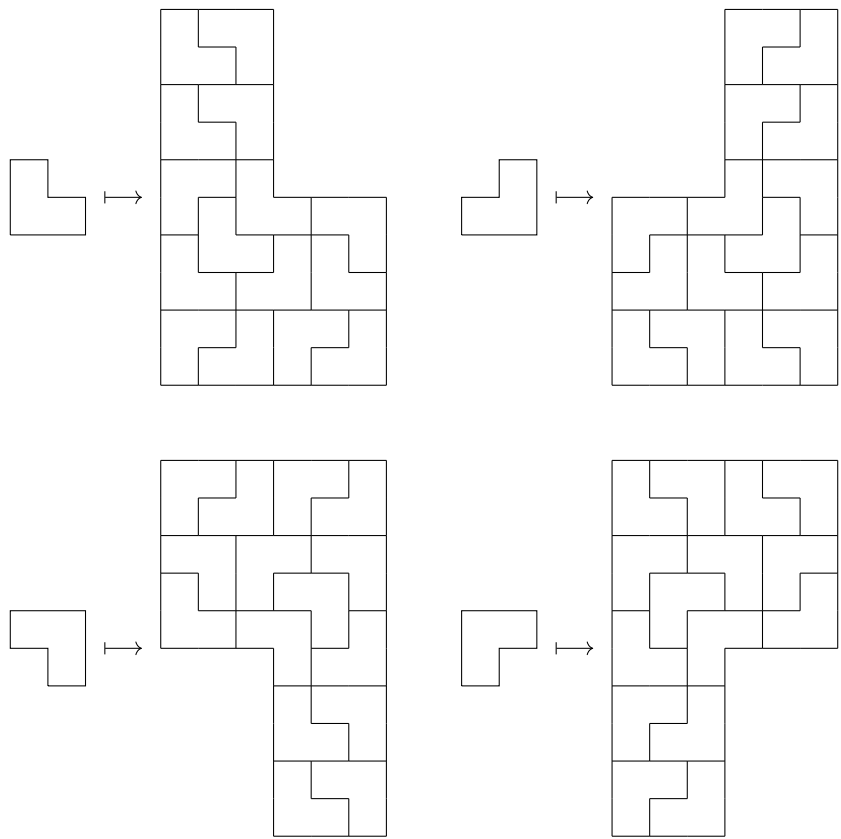%

\begin{example}[Chair \(3\times 5\) variant \(1\) (affine)]
The substitution rule is a variant of the chair tiling with \emph{affine} expansion \(\diagonal(3,5)\), and is shown in \cref{figure:chair35-variant-1-affine}. A supertile is shown in \cref{figure:chair35-1-patch3}. Each of the induced substitution matrices have integral eigenvalues. Due to the sizes of the induced matrices, we assume each of the direct limits splits completely.\\
\indent The six-term sequence is
\[
\begin{tikzcd}
\mathbb{Z}[1/5]^2\oplus\mathbb{Z}[1/3]^2\oplus\mathbb{Z}^{32}\arrow[r,"\evaluation"]&\begin{tabular}{@{}c@{}}\(\mathbb{Z}[1/15]\oplus\mathbb{Z}[1/5]^2\)\\\(\oplus\mathbb{Z}[1/3]^5\oplus\mathbb{Z}^{48}\)\end{tabular}\arrow[r,"\iota_\ast"]&\begin{tabular}{@{}c@{}}\(\mathbb{Z}[1/15]\oplus\mathbb{Z}[1/5]\)\\\(\oplus\mathbb{Z}[1/3]^4\oplus\mathbb{Z}^{17}\)\end{tabular}\arrow[d,two heads]\\
\mathbb{Z}[1/5]\oplus\mathbb{Z}[1/3]\arrow[u,hook]&0\arrow[l]&\mathbb{Z}\arrow[l]
\end{tikzcd}
\]
\noindent where the evaluation map sends \(\mathbb{Z}[1/5]\oplus\mathbb{Z}[1/3]\oplus\mathbb{Z}^{32}\leq K_0(C_r^\ast(\dot{G}_{AF});C_r^\ast(\dot{G}_u))\) isomorphically onto its image and the rest to \(0\). We have that
\begin{align*}
K_0(C_r^\ast(\dot{G}_u))&\cong(\mathbb{Z}[1/15]\oplus\mathbb{Z}[1/5]^2\oplus\mathbb{Z}[1/3]^5\oplus\mathbb{Z}^{48})/(\mathbb{Z}[1/5]\oplus\mathbb{Z}[1/3]\oplus\mathbb{Z}^{32})\oplus\mathbb{Z}\\
&\cong K_0(C_r^\ast(\dot{G}_{AF}))/(\mathbb{Z}[1/5]\oplus\mathbb{Z}[1/3]\oplus\mathbb{Z}^{32})\oplus\mathbb{Z}
\end{align*}
\noindent and
\begin{align*}
K_1(C_r^\ast(\dot{G}_u))&\cong(\mathbb{Z}[1/5]\oplus\mathbb{Z}[1/3]\oplus\mathbb{Z}^4)/\mathbb{Z}^4\\
&\leq(\mathbb{Z}[1/5]^2\oplus\mathbb{Z}[1/3]^2\oplus\mathbb{Z}^{36})/\mathbb{Z}^4\\
&\cong C^1/\mathbb{Z}^4\\
&\cong K_0(C_r^\ast(\dot{G}_{AF}^{(1)}))/\mathbb{Z}^4\tag*{(\cref{proposition:af-cochain-nontrivial})}
\end{align*}
\noindent where we have to apply a slightly nontrivial isomorphism in the last step since this substitution does not satisfy the boundary hyperplane condition. The proposition applies since we assumed that \(C^1\) splits completely.
\end{example}
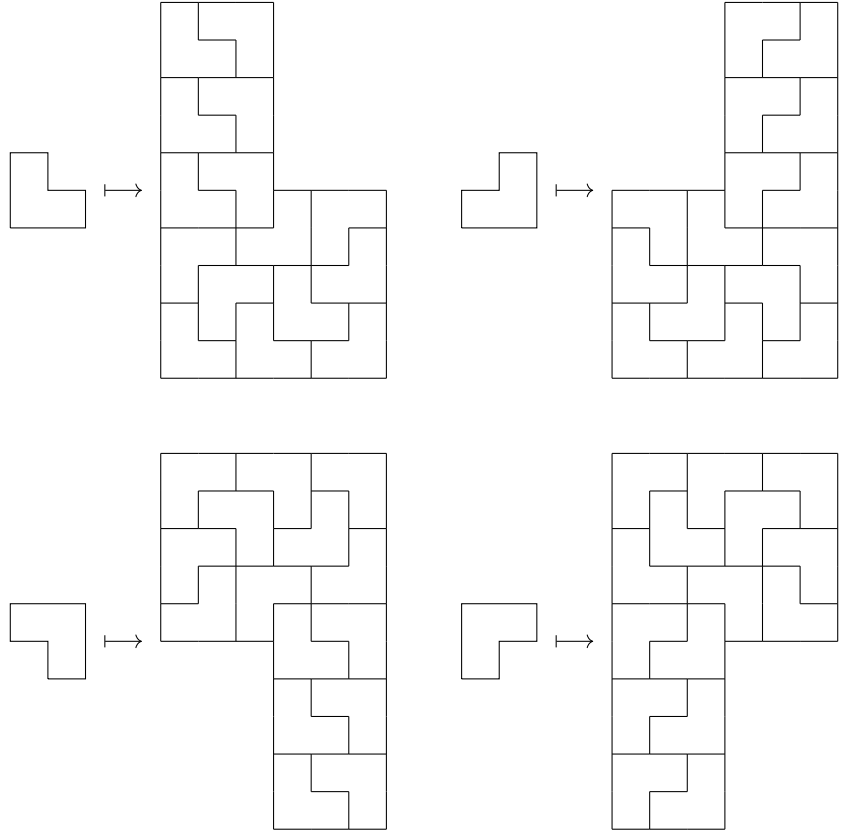
\begin{figure}[t]
\centering
\begin{tikzpicture}
\draw(0,1)--(1,1)--(1,1.5)--(.5,1.5)--(.5,2)--(0,2)--(0,1);
\draw[|->](1.25,1.5)--(1.75,1.5);
\draw(3.5,0)--(3.0,0);
\draw(4.0,-1.0)--(3.5,-1.0);
\draw(4.5,1.5)--(5.0,1.5);
\draw(4.5,.5)--(4.0,.5);
\draw(2.0,3.5)--(2.0,4.0);
\draw(3.5,3.0)--(3.0,3.0);
\draw(3.5,2.5)--(3.5,2.0);
\draw(5.0,.5)--(5.0,1.0);
\draw(2.5,2.0)--(2.0,2.0);
\draw(4.5,.5)--(4.5,1.0);
\draw(3.0,2.0)--(3.5,2.0);
\draw(4.0,1.5)--(4.5,1.5);
\draw(3.5,3.5)--(3.5,4.0);
\draw(4.0,0)--(4.0,.5);
\draw(5.0,-.5)--(5.0,0);
\draw(2.0,1.0)--(2.0,1.5);
\draw(3.0,2.0)--(3.0,2.5);
\draw(3.5,0)--(3.5,.5);
\draw(3.5,3.0)--(3.5,3.5);
\draw(2.5,2.5)--(2.5,3.0);
\draw(4.5,1.0)--(5.0,1.0);
\draw(3.0,4.0)--(2.5,4.0);
\draw(2.0,3.0)--(2.0,3.5);
\draw(2.0,0)--(2.0,.5);
\draw(3.0,-.5)--(3.0,0);
\draw(2.5,0)--(2.5,.5);
\draw(4.0,0)--(4.5,0);
\draw(3.0,1.0)--(3.0,1.5);
\draw(2.0,2.5)--(2.0,3.0);
\draw(2.5,2.0)--(2.5,1.5);
\draw(3.5,1.0)--(3.5,1.5);
\draw(3.5,-.5)--(3.5,0);
\draw(3.0,-1.0)--(3.0,-.5);
\draw(4.0,1.0)--(4.0,1.5);
\draw(4.5,-.5)--(4.5,0);
\draw(2.0,-.5)--(2.0,0);
\draw(3.0,3.5)--(2.5,3.5);
\draw(2.5,2.5)--(3.0,2.5);
\draw(3.5,2.5)--(3.5,3.0);
\draw(4.0,1.0)--(4.0,.5);
\draw(5.0,-.5)--(5.0,-1.0);
\draw(4.5,-1.0)--(5.0,-1.0);
\draw(2.5,3.0)--(3.0,3.0);
\draw(2.5,-1.0)--(2.0,-1.0);
\draw(2.0,0)--(2.5,0);
\draw(5.0,0)--(4.5,0);
\draw(4.0,-1.0)--(4.5,-1.0);
\draw(3.0,.5)--(3.5,.5);
\draw(4.0,1.5)--(3.5,1.5);
\draw(3.0,4.0)--(3.5,4.0);
\draw(3.5,-.5)--(4.0,-.5);
\draw(2.0,.5)--(2.0,1.0);
\draw(3.0,1.0)--(2.5,1.0);
\draw(2.5,2.0)--(3.0,2.0);
\draw(2.5,-.5)--(2.5,0);
\draw(2.0,3.0)--(2.5,3.0);
\draw(3.0,-1.0)--(3.5,-1.0);
\draw(2.0,1.5)--(2.0,2.0);
\draw(2.5,-1.0)--(3.0,-1.0);
\draw(4.0,-1.0)--(4.0,-.5);
\draw(2.0,-.5)--(2.0,-1.0);
\draw(2.0,1.0)--(2.5,1.0);
\draw(2.5,-.5)--(3.0,-.5);
\draw(3.0,3.5)--(3.0,3.0);
\draw(4.0,-.5)--(4.5,-.5);
\draw(5.0,1.0)--(5.0,1.5);
\draw(3.5,1.5)--(3.5,2.0);
\draw(3.0,1.0)--(3.5,1.0);
\draw(2.5,3.5)--(2.5,4.0);
\draw(2.0,2.5)--(2.0,2.0);
\draw(5.0,0)--(5.0,.5);
\draw(3.5,.5)--(4.0,.5);
\draw(2.0,4.0)--(2.5,4.0);
\draw(3.0,1.5)--(2.5,1.5);
\draw(3.0,.5)--(3.0,1.0);
\draw(3.0,.5)--(2.5,.5);
\draw(6.0,1)--(7.0,1)--(7.0,2)--(6.5,2)--(6.5,1.5)--(6.0,1.5)--(6.0,1);
\draw[|->](7.25,1.5)--(7.75,1.5);
\draw(8.5,.5)--(9.0,.5);
\draw(8.0,-.5)--(8.0,0);
\draw(8.0,1.0)--(8.5,1.0);
\draw(9.0,-1.0)--(9.5,-1.0);
\draw(11.0,3.5)--(11.0,4.0);
\draw(9.0,-1.0)--(9.0,-.5);
\draw(11.0,0)--(10.5,0);
\draw(10.0,-1.0)--(9.5,-1.0);
\draw(9.5,3.5)--(9.5,4.0);
\draw(10.0,-1.0)--(10.0,-.5);
\draw(9.0,-.5)--(9.5,-.5);
\draw(9.0,0)--(9.0,.5);
\draw(9.5,3.0)--(10.0,3.0);
\draw(9.5,0)--(9.5,.5);
\draw(8.5,-1.0)--(9.0,-1.0);
\draw(10.5,0)--(10.5,.5);
\draw(8.0,1.5)--(8.5,1.5);
\draw(9.5,0)--(10.0,0);
\draw(8.0,-1.0)--(8.0,-.5);
\draw(11.0,1.0)--(11.0,1.5);
\draw(10.5,2.0)--(10.5,1.5);
\draw(9.5,4.0)--(10.0,4.0);
\draw(8.0,.5)--(8.0,0);
\draw(9.5,1.5)--(9.5,2.0);
\draw(10.5,-1.0)--(11.0,-1.0);
\draw(11.0,3.0)--(11.0,3.5);
\draw(8.0,1.0)--(8.0,1.5);
\draw(11.0,0)--(11.0,.5);
\draw(9.0,1.0)--(9.0,1.5);
\draw(9.5,.5)--(9.0,.5);
\draw(9.0,1.0)--(9.0,.5);
\draw(11.0,2.5)--(11.0,3.0);
\draw(11.0,-.5)--(11.0,0);
\draw(9.5,.5)--(10.0,.5);
\draw(8.5,.5)--(8.5,1.0);
\draw(10.5,2.5)--(10.5,3.0);
\draw(8.5,-1.0)--(8.0,-1.0);
\draw(10.0,.5)--(10.0,1.0);
\draw(10.0,3.0)--(10.0,3.5);
\draw(10.5,3.5)--(10.0,3.5);
\draw(10.0,4.0)--(10.5,4.0);
\draw(10.5,2.5)--(10.0,2.5);
\draw(10.0,2.0)--(10.5,2.0);
\draw(10.0,.5)--(10.5,.5);
\draw(10.0,2.0)--(10.0,2.5);
\draw(10.0,-.5)--(10.5,-.5);
\draw(10.0,1.0)--(10.5,1.0);
\draw(11.0,.5)--(11.0,1.0);
\draw(8.0,.5)--(8.0,1.0);
\draw(9.5,1.0)--(9.5,1.5);
\draw(10.5,3.5)--(10.5,4.0);
\draw(10.0,2.0)--(9.5,2.0);
\draw(9.5,3.0)--(9.5,3.5);
\draw(10.5,3.0)--(10.0,3.0);
\draw(9.5,1.0)--(10.0,1.0);
\draw(10.0,1.5)--(10.5,1.5);
\draw(10.0,1.0)--(10.0,1.5);
\draw(11.0,4.0)--(10.5,4.0);
\draw(11.0,1.5)--(11.0,2.0);
\draw(10.0,-.5)--(10.0,0);
\draw(8.5,0)--(9.0,0);
\draw(10.5,-.5)--(10.5,0);
\draw(11.0,-.5)--(11.0,-1.0);
\draw(11.0,3.0)--(10.5,3.0);
\draw(10.5,2.0)--(11.0,2.0);
\draw(9.0,1.5)--(9.5,1.5);
\draw(11.0,1.0)--(10.5,1.0);
\draw(10.0,-1.0)--(10.5,-1.0);
\draw(9.5,3.0)--(9.5,2.5);
\draw(11.0,2.5)--(11.0,2.0);
\draw(8.5,0)--(8.0,0);
\draw(9.5,0)--(9.5,-.5);
\draw(8.5,-.5)--(8.5,0);
\draw(9.5,2.0)--(9.5,2.5);
\draw(9.0,1.5)--(8.5,1.5);
\draw(8.5,-.5)--(9.0,-.5);
\draw(.5,-5)--(1,-5)--(1,-4)--(0,-4)--(0,-4.5)--(.5,-4.5)--(.5,-5);
\draw[|->](1.25,-4.5)--(1.75,-4.5);
\draw(4.0,-4.0)--(4.0,-3.5);
\draw(5.0,-6.0)--(5.0,-5.5);
\draw(3.5,-3.5)--(3.5,-3.0);
\draw(2.5,-4.5)--(2.0,-4.5);
\draw(5.0,-3.5)--(5.0,-3.0);
\draw(5.0,-4.5)--(5.0,-5.0);
\draw(3.0,-2.0)--(3.5,-2.0);
\draw(3.5,-3.5)--(3.0,-3.5);
\draw(4.0,-2.0)--(3.5,-2.0);
\draw(4.0,-4.0)--(3.5,-4.0);
\draw(4.0,-5.5)--(4.0,-5.0);
\draw(3.5,-6.5)--(3.5,-6.0);
\draw(3.0,-3.0)--(2.5,-3.0);
\draw(3.5,-5.0)--(3.5,-4.5);
\draw(2.0,-3.5)--(2.0,-3.0);
\draw(3.0,-2.5)--(2.5,-2.5);
\draw(4.5,-6.0)--(4.5,-5.5);
\draw(4.5,-7.0)--(4.5,-6.5);
\draw(4.0,-4.5)--(4.0,-4.0);
\draw(2.5,-2.0)--(3.0,-2.0);
\draw(3.0,-4.5)--(3.0,-4.0);
\draw(3.0,-4.5)--(3.5,-4.5);
\draw(4.0,-7.0)--(3.5,-7.0);
\draw(3.0,-2.5)--(3.0,-2.0);
\draw(5.0,-7.0)--(4.5,-7.0);
\draw(3.5,-5.0)--(4.0,-5.0);
\draw(4.5,-2.0)--(5.0,-2.0);
\draw(5.0,-5.5)--(5.0,-5.0);
\draw(5.0,-4.5)--(5.0,-4.0);
\draw(5.0,-6.5)--(5.0,-6.0);
\draw(5.0,-3.0)--(4.5,-3.0);
\draw(4.5,-3.5)--(4.5,-3.0);
\draw(3.5,-3.0)--(4.0,-3.0);
\draw(3.0,-2.5)--(3.5,-2.5);
\draw(2.0,-4.5)--(2.0,-4.0);
\draw(3.0,-3.0)--(3.0,-3.5);
\draw(2.5,-4.0)--(2.5,-3.5);
\draw(5.0,-5.0)--(4.5,-5.0);
\draw(3.5,-6.0)--(3.5,-5.5);
\draw(2.5,-2.0)--(2.0,-2.0);
\draw(3.5,-5.5)--(3.5,-5.0);
\draw(4.5,-5.0)--(4.0,-5.0);
\draw(4.5,-4.5)--(4.0,-4.5);
\draw(4.0,-6.0)--(4.5,-6.0);
\draw(5.0,-6.0)--(4.5,-6.0);
\draw(4.5,-6.5)--(4.0,-6.5);
\draw(4.5,-2.5)--(4.5,-3.0);
\draw(4.0,-2.5)--(4.5,-2.5);
\draw(3.5,-4.5)--(3.5,-4.0);
\draw(4.0,-3.0)--(4.0,-2.5);
\draw(4.5,-3.5)--(4.0,-3.5);
\draw(3.5,-6.0)--(4.0,-6.0);
\draw(5.0,-4.0)--(4.5,-4.0);
\draw(5.0,-4.0)--(5.0,-3.5);
\draw(4.5,-4.5)--(4.5,-5.0);
\draw(5.0,-7.0)--(5.0,-6.5);
\draw(3.0,-4.0)--(3.0,-3.5);
\draw(3.5,-6.5)--(3.5,-7.0);
\draw(4.0,-6.5)--(4.0,-6.0);
\draw(2.5,-4.0)--(2.0,-4.0);
\draw(4.0,-2.0)--(4.0,-2.5);
\draw(2.0,-3.0)--(2.5,-3.0);
\draw(3.5,-3.5)--(4.0,-3.5);
\draw(2.0,-2.5)--(2.0,-2.0);
\draw(2.0,-3.5)--(2.0,-4.0);
\draw(3.5,-2.5)--(3.5,-3.0);
\draw(2.0,-2.5)--(2.0,-3.0);
\draw(2.5,-2.5)--(2.5,-3.0);
\draw(4.0,-4.0)--(4.5,-4.0);
\draw(4.5,-7.0)--(4.0,-7.0);
\draw(4.5,-5.5)--(4.0,-5.5);
\draw(2.5,-3.5)--(3.0,-3.5);
\draw(4.5,-2.0)--(4.0,-2.0);
\draw(5.0,-3.0)--(5.0,-2.5);
\draw(2.5,-4.5)--(3.0,-4.5);
\draw(5.0,-2.5)--(5.0,-2.0);
\draw(6.0,-5)--(6.5,-5)--(6.5,-4.5)--(7.0,-4.5)--(7.0,-4)--(6.0,-4)--(6.0,-5);
\draw[|->](7.25,-4.5)--(7.75,-4.5);
\draw(9.0,-2.5)--(8.5,-2.5);
\draw(10.0,-2.0)--(9.5,-2.0);
\draw(8.0,-7.0)--(8.5,-7.0);
\draw(8.0,-6.0)--(8.0,-5.5);
\draw(11.0,-3.5)--(11.0,-3.0);
\draw(8.0,-4.5)--(8.0,-5.0);
\draw(8.0,-3.5)--(8.0,-3.0);
\draw(8.5,-6.0)--(8.5,-5.5);
\draw(8.0,-6.0)--(8.5,-6.0);
\draw(9.0,-3.5)--(8.5,-3.5);
\draw(8.5,-5.0)--(8.0,-5.0);
\draw(11.0,-4.5)--(10.5,-4.5);
\draw(11.0,-4.5)--(11.0,-4.0);
\draw(8.5,-5.0)--(9.0,-5.0);
\draw(9.5,-7.0)--(9.0,-7.0);
\draw(8.5,-7.0)--(9.0,-7.0);
\draw(9.0,-5.5)--(9.0,-5.0);
\draw(10.0,-3.0)--(10.5,-3.0);
\draw(10.0,-2.5)--(10.5,-2.5);
\draw(8.0,-5.5)--(8.0,-5.0);
\draw(8.0,-4.5)--(8.0,-4.0);
\draw(8.0,-6.5)--(8.0,-6.0);
\draw(10.0,-3.0)--(10.0,-3.5);
\draw(10.5,-3.5)--(10.0,-3.5);
\draw(10.5,-3.5)--(10.5,-4.0);
\draw(9.5,-4.5)--(9.5,-4.0);
\draw(8.5,-5.0)--(8.5,-4.5);
\draw(11.0,-4.0)--(10.5,-4.0);
\draw(9.5,-3.5)--(10.0,-3.5);
\draw(10.0,-4.5)--(10.0,-4.0);
\draw(8.5,-7.0)--(8.5,-6.5);
\draw(9.0,-3.0)--(9.5,-3.0);
\draw(9.5,-3.5)--(9.5,-3.0);
\draw(10.5,-2.0)--(11.0,-2.0);
\draw(9.5,-3.0)--(9.5,-2.5);
\draw(9.0,-2.5)--(9.0,-3.0);
\draw(9.0,-3.5)--(9.0,-4.0);
\draw(9.5,-7.0)--(9.5,-6.5);
\draw(10.0,-4.5)--(10.5,-4.5);
\draw(8.0,-2.5)--(8.0,-3.0);
\draw(8.5,-2.5)--(8.5,-3.0);
\draw(8.5,-4.5)--(9.0,-4.5);
\draw(9.5,-6.5)--(9.5,-6.0);
\draw(8.5,-2.0)--(8.0,-2.0);
\draw(10.0,-2.0)--(10.5,-2.0);
\draw(8.5,-3.5)--(8.5,-3.0);
\draw(8.0,-4.0)--(8.0,-3.5);
\draw(9.5,-4.5)--(9.5,-5.0);
\draw(8.0,-7.0)--(8.0,-6.5);
\draw(9.0,-2.0)--(9.5,-2.0);
\draw(8.5,-6.5)--(9.0,-6.5);
\draw(9.0,-5.0)--(9.5,-5.0);
\draw(8.5,-2.0)--(9.0,-2.0);
\draw(9.5,-5.5)--(9.5,-5.0);
\draw(8.0,-4.0)--(8.5,-4.0);
\draw(9.0,-3.5)--(9.5,-3.5);
\draw(9.5,-6.0)--(9.5,-5.5);
\draw(9.5,-4.0)--(9.0,-4.0);
\draw(10.0,-4.0)--(10.0,-3.5);
\draw(9.5,-4.5)--(10.0,-4.5);
\draw(8.5,-4.0)--(9.0,-4.0);
\draw(9.0,-4.5)--(9.0,-4.0);
\draw(9.0,-2.5)--(9.0,-2.0);
\draw(9.0,-5.5)--(8.5,-5.5);
\draw(10.5,-3.0)--(11.0,-3.0);
\draw(8.0,-3.0)--(8.5,-3.0);
\draw(9.0,-6.0)--(9.0,-6.5);
\draw(10.0,-2.5)--(10.0,-2.0);
\draw(11.0,-4.0)--(11.0,-3.5);
\draw(8.0,-2.5)--(8.0,-2.0);
\draw(10.0,-2.5)--(9.5,-2.5);
\draw(9.0,-6.0)--(9.5,-6.0);
\draw(10.5,-2.5)--(10.5,-3.0);
\draw(9.0,-6.0)--(8.5,-6.0);
\draw(11.0,-2.5)--(11.0,-3.0);
\draw(11.0,-2.5)--(11.0,-2.0);
\end{tikzpicture}
\caption{Substitution rule for a variant of the chair tiling with affine expansion \(\diagonal(3,5)\).}
\label{figure:chair35-variant-2-affine}
\end{figure}%

\begin{example}[Chair \(3\times 5\) variant \(2\) (affine)]
The substitution rule is a variant of the chair tiling with \emph{affine} expansion \(\diagonal(3,5)\), and is shown in \cref{figure:chair35-variant-2-affine}. A supertile is shown in \cref{figure:chair35-2-patch3}. Each of the induced substitution matrices have integral eigenvalues. Due to the sizes of the induced matrices, we assume each of the direct limits splits completely.\\
\indent The six-term sequence is
\[
\begin{tikzcd}
\begin{tabular}{@{}c@{}}\(\mathbb{Z}[1/5]^2\oplus\mathbb{Z}[1/3]^7\)\\\(\oplus\mathbb{Z}[1/2]^2\oplus\mathbb{Z}^{94}\)\end{tabular}\arrow[r,"\evaluation"]&\begin{tabular}{@{}c@{}}\(\mathbb{Z}[1/15]\oplus\mathbb{Z}[1/5]^2\)\\\(\oplus\mathbb{Z}[1/3]^{13}\oplus\mathbb{Z}[1/2]^4\oplus\mathbb{Z}^{128}\)\end{tabular}\arrow[r,"\iota_\ast"]&\begin{tabular}{@{}c@{}}\(\mathbb{Z}[1/15]\oplus\mathbb{Z}[1/5]\)\\\(\oplus\mathbb{Z}[1/3]^7\oplus\mathbb{Z}[1/2]^2\oplus\mathbb{Z}^{35}\)\end{tabular}\arrow[d,two heads]\\
\mathbb{Z}[1/5]\oplus\mathbb{Z}[1/3]\arrow[u,hook]&0\arrow[l]&\mathbb{Z}\arrow[l]
\end{tikzcd}
\]
\noindent where the evaluation map sends \(\mathbb{Z}[1/5]\oplus\mathbb{Z}[1/3]^6\oplus\mathbb{Z}[1/2]^2\oplus\mathbb{Z}^{94}\leq K_0(C_r^\ast(\dot{G}_{AF});C_r^\ast(\dot{G}_u))\) isomorphically onto its image and the rest to \(0\). We have that
\begin{align*}
K_0(C_r^\ast(\dot{G}_u))&\cong\left.\left(\begin{tabular}{@{}c@{}}\(\mathbb{Z}[1/15]\oplus\mathbb{Z}[1/5]^2\)\\\(\oplus\mathbb{Z}[1/3]^{13}\oplus\mathbb{Z}[1/2]^4\oplus\mathbb{Z}^{128}\)\end{tabular}\right)\middle/(\mathbb{Z}[1/5]\oplus\mathbb{Z}[1/3]^6\oplus\mathbb{Z}[1/2]^2\oplus\mathbb{Z}^{94})\right.\oplus\mathbb{Z}\\
&\cong K_0(C_r^\ast(\dot{G}_{AF}))/(\mathbb{Z}[1/5]\oplus\mathbb{Z}[1/3]^6\oplus\mathbb{Z}[1/2]^2\oplus\mathbb{Z}^{94})\oplus\mathbb{Z}
\end{align*}
\noindent and
\begin{align*}
K_1(C_r^\ast(\dot{G}_u))&\cong(\mathbb{Z}[1/5]\oplus\mathbb{Z}[1/3]\oplus\mathbb{Z}^{41})/\mathbb{Z}^{41}\\
&\leq(\mathbb{Z}[1/5]^2\oplus\mathbb{Z}[1/3]^7\oplus\mathbb{Z}[1/2]^2\oplus\mathbb{Z}^{135})/\mathbb{Z}^{41}\\
&\cong C^1/\mathbb{Z}^{41}\\
&\cong K_0(C_r^\ast(\dot{G}_{AF}^{(1)}))/\mathbb{Z}^{41}\tag*{(\cref{proposition:af-cochain-nontrivial})}
\end{align*}
\noindent where we have to apply a slightly nontrivial isomorphism in the last step since this substitution does not satisfy the boundary hyperplane condition. The proposition applies since we assumed that \(C^1\) splits completely.
\end{example}
\begin{example}[Hexagon (ps)]
This is the \emph{pseudo}substitution rule in \cite[Figure 4]{clarksadun06}. Each of the induced substitution matrices have integral eigenvalues. One easily checks that the direct limits of the induced substitution matrices split completely.\\
\indent The six-term sequence is
\[
\begin{tikzcd}
\mathbb{Z}[1/2]^2\oplus\mathbb{Z}^6\arrow[r,"\evaluation"]&\mathbb{Z}[1/4]\oplus\mathbb{Z}^8\arrow[r,"\iota_\ast"]&\mathbb{Z}[1/4]\oplus\mathbb{Z}^3\arrow[d,two heads]\\
\mathbb{Z}[1/2]^2\arrow[u,hook]&0\arrow[l]&\mathbb{Z}\arrow[l]
\end{tikzcd}
\]
\noindent where the evaluation map sends \(\mathbb{Z}^6\leq K_0(C_r^\ast(\dot{G}_{AF});C_r^\ast(\dot{G}_u))\) isomorphically onto its image and the rest to \(0\). We have that
\begin{align*}
K_0(C_r^\ast(\dot{G}_u))&\cong(\mathbb{Z}[1/4]\oplus\mathbb{Z}^8)/\mathbb{Z}^6\oplus\mathbb{Z}\\
&\cong K_0(C_r^\ast(\dot{G}_{AF}))/\mathbb{Z}^6\oplus\mathbb{Z}
\end{align*}
\noindent and, since the substitution rule satisfies the boundary hyperplane condition,
\begin{align*}
K_1(C_r^\ast(\dot{G}_u))&\cong(\mathbb{Z}[1/2]^2\oplus\mathbb{Z}^5)/\mathbb{Z}^5\\
&\leq(\mathbb{Z}[1/2]^2\oplus\mathbb{Z}^{11})/\mathbb{Z}^5\\
&\cong K_0(C_r^\ast(\dot{G}_{AF}^{(1)}))/\mathbb{Z}^5.
\end{align*}
\end{example}
\begin{figure}[p]
\centering
\includegraphics[width=\textwidth,height=0.9\textheight,keepaspectratio]{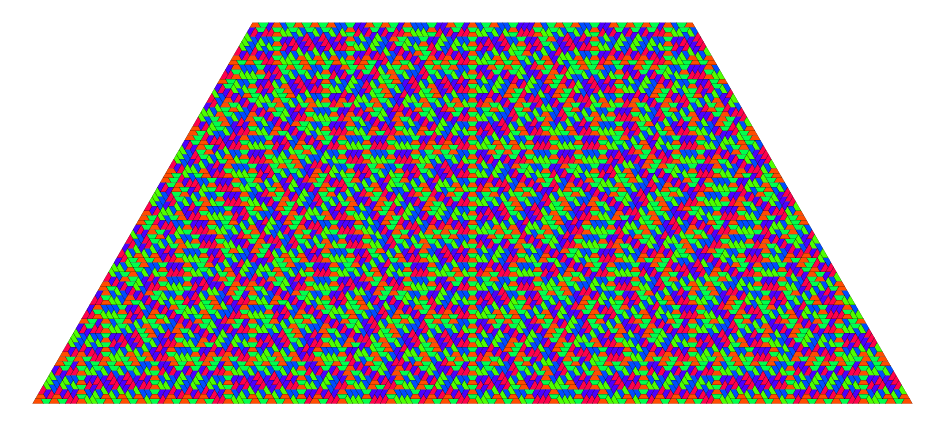}
\caption{A level-\(4\) supertile of the halfhex \(3\times 3\) tiling.}
\label{figure:halfhex3-patch4}
\end{figure}%

\begin{figure}[p]
\centering
\includegraphics[width=\textwidth,height=0.9\textheight,keepaspectratio]{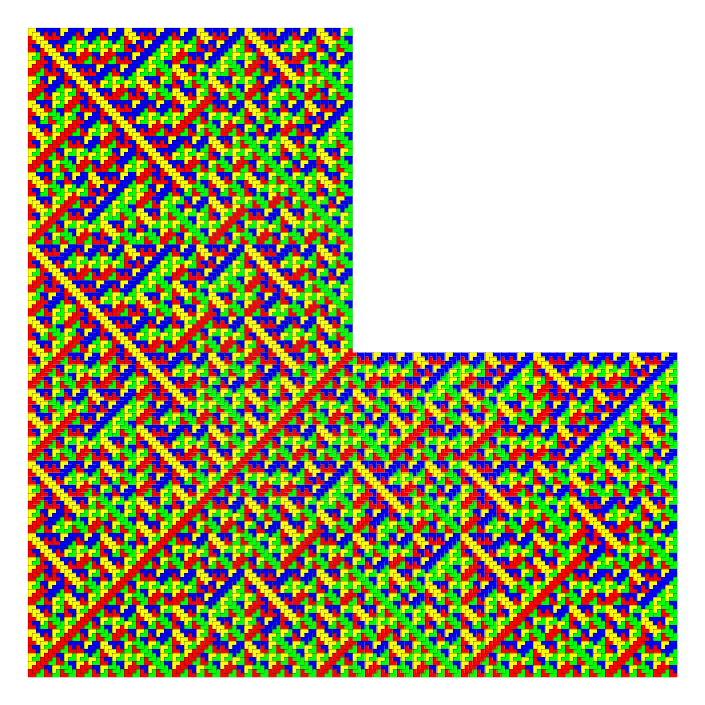}
\caption{A level-\(4\) supertile of the chair \(3\times 3\) tiling.}
\label{figure:chair3-patch4}
\end{figure}%

\begin{figure}[p]
\centering
\includegraphics[width=\textwidth,height=0.9\textheight,keepaspectratio]{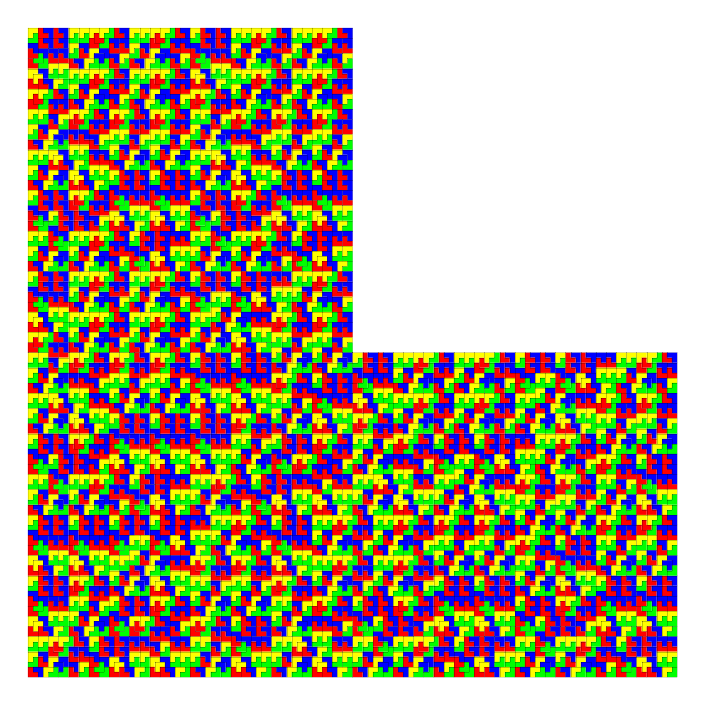}
\caption{A level-\(3\) supertile of the chair \(4\times 4\) tiling.}
\label{figure:chair4-patch3}
\end{figure}%

\begin{figure}[p]
\centering
\includegraphics[width=\textwidth,height=0.9\textheight,keepaspectratio]{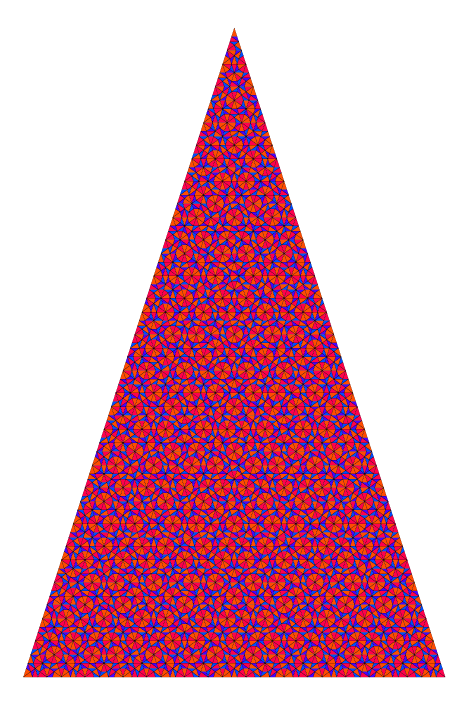}
\caption{A level-\(9\) supertile of the T\"ubingen triangle tiling.}
\label{figure:tubingen-patch9}
\end{figure}%

\begin{figure}[p]
\centering
\includegraphics[width=\textwidth,height=0.9\textheight,keepaspectratio]{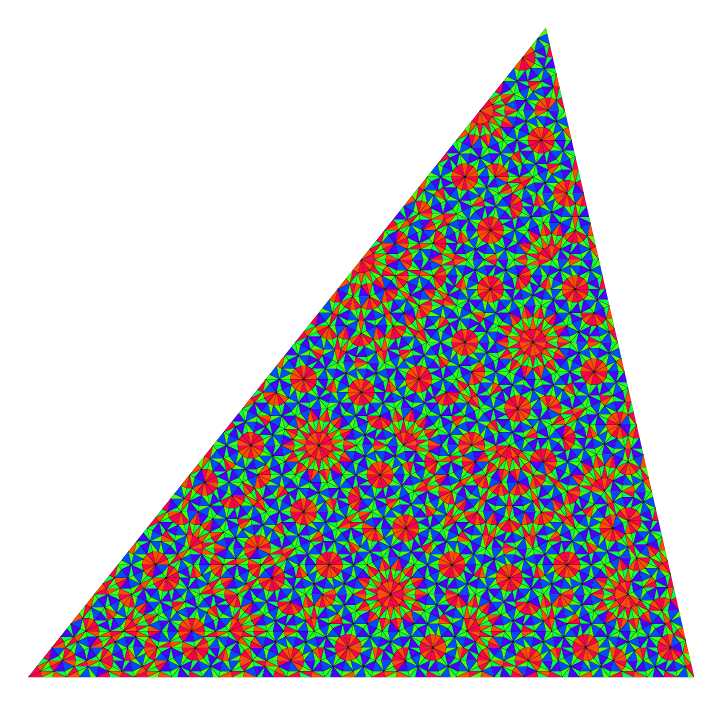}
\caption{A level-\(4\) supertile of the Danzer 7-fold tiling.}
\label{figure:danzer7-patch4}
\end{figure}%

\begin{figure}[p]
\centering
\includegraphics[width=\textwidth,height=0.9\textheight,keepaspectratio]{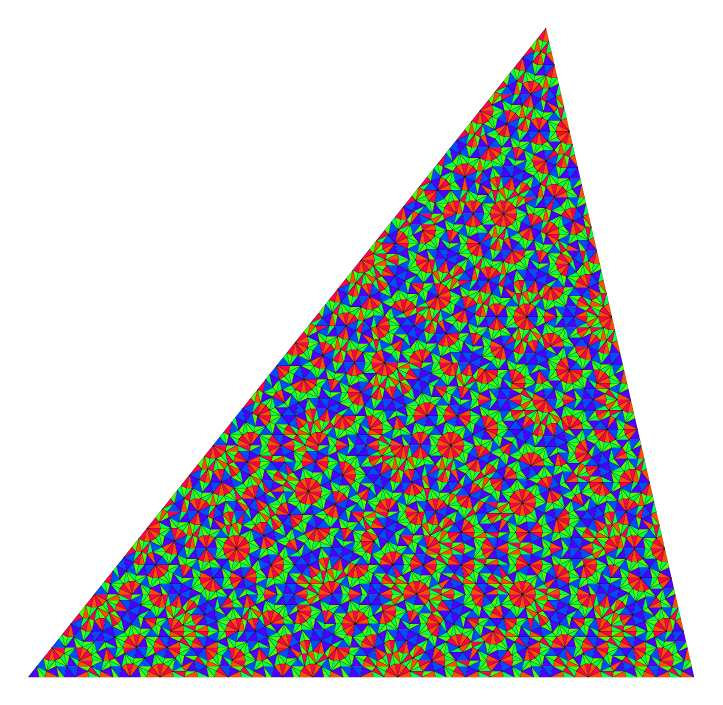}
\caption{A level-\(4\) supertile of the GKM 9.1.1.1 tiling.}
\label{figure:gkm9.1.1.1-patch4}
\end{figure}%

\begin{figure}[p]
\centering
\includegraphics[width=\textwidth,height=0.9\textheight,keepaspectratio]{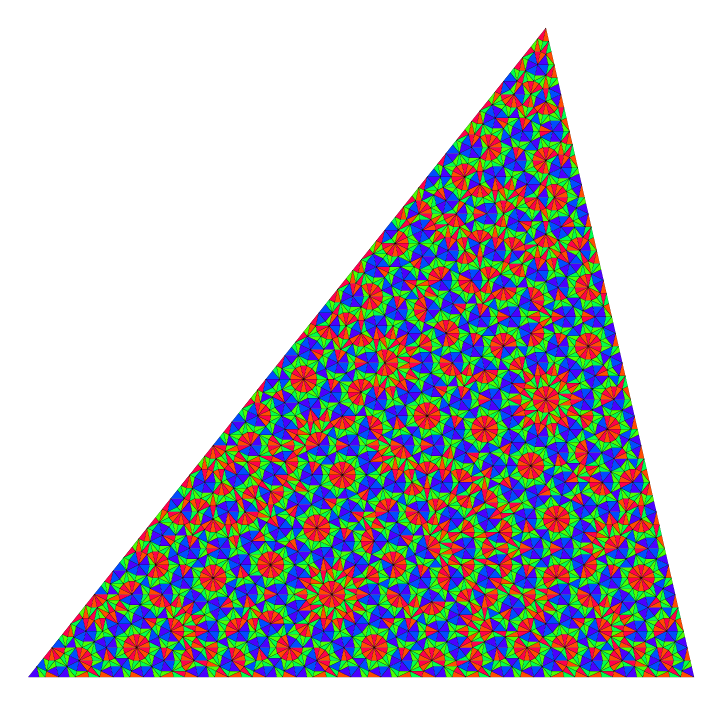}
\caption{A level-\(4\) supertile of the GKM 9.2.1.1 tiling.}
\label{figure:gkm9.2.1.1-patch4}
\end{figure}%

%
%
%

\begin{figure}[p]
\begin{minipage}{0.48\textwidth}
\centering
\includegraphics[width=\textwidth,height=0.9\textheight,keepaspectratio]{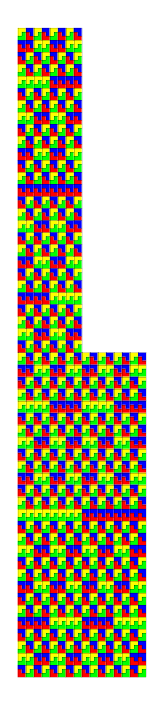}
\caption{A level-\(4\) supertile of the chair \(2\times 3\) tiling.}
\label{figure:chair23-patch4}
\end{minipage}\hfill
\begin{minipage}{0.48\textwidth}
\centering
\includegraphics[width=\textwidth,height=0.9\textheight,keepaspectratio]{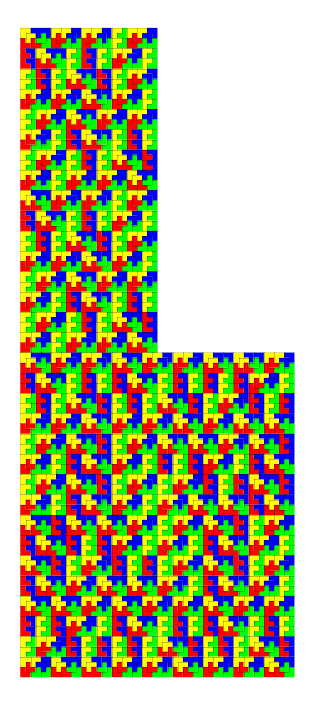}
\caption{A level-\(3\) supertile of the chair \(3\times 4\) tiling.}
\label{figure:chair34-patch3}
\end{minipage}
\end{figure}%

\begin{figure}[p]
\begin{minipage}{0.48\textwidth}
\centering
\includegraphics[width=\textwidth,height=0.9\textheight,keepaspectratio]{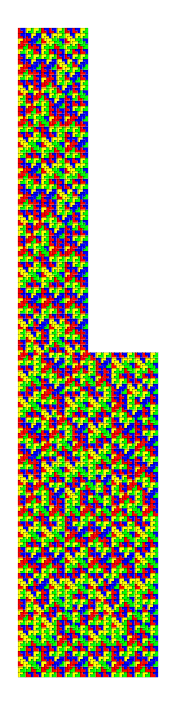}
\caption{A level-\(3\) supertile of the chair \(3\times 5\) variant \(1\) tiling.}
\label{figure:chair35-1-patch3}
\end{minipage}\hfill
\begin{minipage}{0.48\textwidth}
\centering
\includegraphics[width=\textwidth,height=0.9\textheight,keepaspectratio]{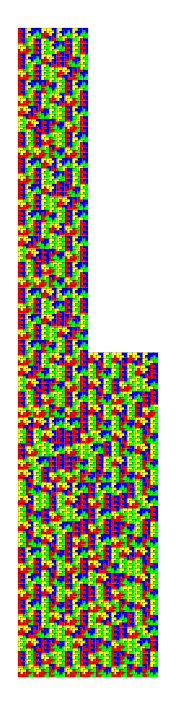}
\caption{A level-\(3\) supertile of the chair \(3\times 5\) variant \(2\) tiling.}
\label{figure:chair35-2-patch3}
\end{minipage}
\end{figure}%

\begin{landscape}
\footnotesize
\centering
\begin{tabularx}{\textwidth}{|l|c|c|c|c|c|c|c|}
\hline
&CP&PI&L(PI)&FT&L(FT)&Exact sequence&\begin{tabular}{@{}c@{}}\(K_0(C_r^\ast(\dot{G}_{AF});C_r^\ast(\dot{G}_u))\)\\\(\cong K_0(C_r^\ast(\dot{G}_{AF}^{(1)}))/\sim\)\end{tabular}\\
\hline
\begin{tabular}{@{}l@{}}Dyadic\\solenoid\end{tabular}&\(1\)&\(2\)&\(1\)&\(2\)&\(1\)&
\begin{tikzcd}[row sep=small,column sep=small,scale cd=0.75,ampersand replacement=\&]
\mathbb{Z}[1/2]^2\arrow[r,"\evaluation"]\&\mathbb{Z}[1/4]\arrow[r,"\iota_\ast"]\&\mathbb{Z}[1/4]\oplus\mathbb{Z}\arrow[d,two heads]\\
\mathbb{Z}[1/2]^2\arrow[u,hook]\&0\arrow[l]\&\mathbb{Z}\arrow[l]
\end{tikzcd}
&\(\mathbb{Z}[1/2]^2\)\\
\hline
Half-hex&\(24\)&\(60\)&\(38\)&\(42\)&\(20\)&
\begin{tikzcd}[row sep=small,column sep=small,scale cd=0.75,ampersand replacement=\&]
\mathbb{Z}[1/2]^3\oplus\mathbb{Z}^2\arrow[r,"\evaluation"]\&\mathbb{Z}[1/4]\oplus\mathbb{Z}[1/2]\oplus\mathbb{Z}^4\arrow[r,"\iota_\ast"]\&\mathbb{Z}[1/4]\oplus\mathbb{Z}^3\arrow[d,two heads]\\
\mathbb{Z}[1/2]^2\arrow[u,hook]\&0\arrow[l]\&\mathbb{Z}\arrow[l]
\end{tikzcd}
&\((\mathbb{Z}[1/2]^3\oplus\mathbb{Z}^3)/\mathbb{Z}\)\\
\hline
Half-hex \(3\times3\)&\(618\)&\(2832\)&\(2490\)&\(864\)&\(324\)&
\begin{tikzcd}[row sep=small,column sep=small,scale cd=0.75,ampersand replacement=\&]
\begin{tabular}{@{}c@{}}\(\mathbb{Z}[1/3]^6\oplus\mathbb{Z}[1/2]^{15}\)\\\(\oplus\mathbb{Z}^{259}\oplus\mathbb{Q}^6\)\end{tabular}\arrow[r,"\evaluation"]\&\begin{tabular}{@{}c@{}}\(\mathbb{Z}[1/9]\oplus\mathbb{Z}[1/3]^7\)\\\(\oplus\mathbb{Z}[1/2]^{24}\oplus\mathbb{Z}^{306}\)\\\(\oplus\mathbb{Q}^{12}\)\end{tabular}\arrow[r,"\iota_\ast"]\&\begin{tabular}{@{}c@{}}\(\mathbb{Z}[1/9]\oplus\mathbb{Z}[1/3]^3\)\\\(\oplus\mathbb{Z}[1/2]^9\oplus\mathbb{Z}^{48}\)\\\(\oplus\mathbb{Z}_2\oplus\mathbb{Q}^6\)\end{tabular}
\arrow[d,two heads]\\
\mathbb{Z}[1/3]^2\arrow[u,hook]\&0\arrow[l]\&\mathbb{Z}\arrow[l]
\end{tikzcd}
&\(\left.\left(\begin{tabular}{@{}c@{}}\((\mathbb{Z}[1/3]^6\oplus\mathbb{Z}[1/2]^{15}\)\\\(\oplus\mathbb{Z}^{402}\oplus\mathbb{Q}^6\)\end{tabular}\right)\middle/\mathbb{Z}^{143}\right.\)\\
\hline
Chair&\(56\)&\(224\)&\(175\)&\(84\)&\(34\)&
\begin{tikzcd}[row sep=small,column sep=small,scale cd=0.75,ampersand replacement=\&]
\mathbb{Z}[1/2]^2\oplus\mathbb{Z}^4\arrow[r,"\evaluation"]\&\mathbb{Z}[1/4]\oplus\mathbb{Z}[1/2]^2\oplus\mathbb{Z}^4\arrow[r,"\iota_\ast"]\&\mathbb{Z}[1/4]\oplus\mathbb{Z}[1/2]^2\oplus\mathbb{Z}\arrow[d,two heads]\\
\mathbb{Z}[1/2]^2\arrow[u,hook]\&0\arrow[l]\&\mathbb{Z}\arrow[l]
\end{tikzcd}
&\(\mathbb{Z}[1/2]^2\oplus\mathbb{Z}^4\)\\
\hline
Chair \(3\times 3\)&\(364\)&\(1616\)&\(1333\)&\(479\)&\(174\)&
\begin{tikzcd}[row sep=small,column sep=small,scale cd=0.75,ampersand replacement=\&]
\begin{tabular}{@{}c@{}}\(\mathbb{Z}[1/3]^4\oplus\mathbb{Z}[1/2]^2\)\\\(\oplus\mathbb{Z}^{76}\)\end{tabular}\arrow[r,"\evaluation"]\&\begin{tabular}{@{}c@{}}\(\mathbb{Z}[1/9]\oplus\mathbb{Z}[1/3]^6\)\\\(\oplus\mathbb{Z}[1/2]^4\oplus\mathbb{Z}^{101}\)\end{tabular}\arrow[r,"\iota_\ast"]\&\begin{tabular}{@{}c@{}}\(\mathbb{Z}[1/9]\oplus\mathbb{Z}[1/3]^4\)\\\(\oplus\mathbb{Z}[1/2]^2\oplus\mathbb{Z}^{26}\)\end{tabular}\arrow[d,two heads]\\
\mathbb{Z}[1/3]^2\arrow[u,hook]\&0\arrow[l]\&\mathbb{Z}\arrow[l]
\end{tikzcd}
&\(\left.\left(\begin{tabular}{@{}c@{}}\(\mathbb{Z}[1/3]^4\oplus\mathbb{Z}[1/2]^2\)\\\(\oplus\mathbb{Z}^{110}\)\end{tabular}\right)\middle/\mathbb{Z}^{34}\right.\)\\
\hline
Chair \(4\times 4\)&\(435\)&\(1782\)&\(1464\)&\(569\)&\(197\)&
\begin{tikzcd}[row sep=small,column sep=small,scale cd=0.75,ampersand replacement=\&]
\begin{tabular}{@{}c@{}}\(\mathbb{Z}[1/4]^2\oplus\mathbb{Z}[1/2]\)\\\(\oplus\mathbb{Z}^{19}\)\end{tabular}\arrow[r,"\evaluation"]\&\begin{tabular}{@{}c@{}}\(\mathbb{Z}[1/16]\oplus\mathbb{Z}[1/2]^4\)\\\(\oplus\mathbb{Z}^{30}\)\end{tabular}\arrow[r,"\iota_\ast"]\&\begin{tabular}{@{}c@{}}\(\mathbb{Z}[1/16]\oplus\mathbb{Z}[1/2]^3\)\\\(\oplus\mathbb{Z}^{12}\)\end{tabular}\arrow[d,two heads]\\
\mathbb{Z}[1/4]^2\arrow[u,hook]\&0\arrow[l]\&\mathbb{Z}\arrow[l]
\end{tikzcd}
&\(\mathbb{Z}[1/4]^2\oplus\mathbb{Z}[1/2]^2\oplus\mathbb{Z}^{19}\)\\
\hline
Table&\(80\)&\(398\)&\(334\)&\(102\)&\(42\)&
\begin{tikzcd}[row sep=small,column sep=small,scale cd=0.75,ampersand replacement=\&]
\mathbb{Z}[1/2]^6\oplus\mathbb{Z}^{29}\arrow[r,"\evaluation"]\&\begin{tabular}{@{}c@{}}\(\mathbb{Z}[1/4]\oplus\mathbb{Z}[1/2]^8\)\\\(\oplus\mathbb{Z}^{32}\)\end{tabular}\arrow[r,"\iota_\ast"]\&\begin{tabular}{@{}c@{}}\(\mathbb{Z}[1/4]\oplus\mathbb{Z}[1/2]^4\)\\\(\oplus\mathbb{Z}^4\oplus\mathbb{Z}_2\)\end{tabular}\arrow[d,two heads]\\
\mathbb{Z}[1/2]^2\arrow[u,hook]\&0\arrow[l]\&\mathbb{Z}\arrow[l]
\end{tikzcd}
&\((\mathbb{Z}[1/2]^6\oplus\mathbb{Z}^{50})/\mathbb{Z}^{21}\)\\
\hline
\begin{tabular}{@{}l@{}}Robinson\\triangle\end{tabular}&\(300\)&\(540\)&\(264\)&\(360\)&\(64\)&
\begin{tikzcd}[row sep=small,column sep=small,scale cd=0.75,ampersand replacement=\&]
\mathbb{Z}^{37}\arrow[r,"\evaluation"]\&\mathbb{Z}^{40}\arrow[r,"\iota_\ast"]\&\mathbb{Z}^9\arrow[d,two heads]\\
\mathbb{Z}^5\arrow[u,hook]\&0\arrow[l]\&\mathbb{Z}\arrow[l]
\end{tikzcd}
&\(\mathbb{Z}^{40}/\mathbb{Z}^3\)\\
\hline
\begin{tabular}{@{}l@{}}T\"ubingen\\triangle\end{tabular}&\(860\)&\(1710\)&\(880\)&\(950\)&\(130\)&
\begin{tikzcd}[row sep=small,column sep=small,scale cd=0.75,ampersand replacement=\&]
\mathbb{Z}^{101}\arrow[r,"\evaluation"]\&\mathbb{Z}^{120}\arrow[r,"\iota_\ast"]\&\mathbb{Z}^{25}\oplus\mathbb{Z}_5^2\arrow[d,two heads]\\
\mathbb{Z}^5\arrow[u,hook]\&0\arrow[l]\&\mathbb{Z}\arrow[l]
\end{tikzcd}
&\(\mathbb{Z}^{110}/\mathbb{Z}^9\)\\
\hline
Danzer \(7\)-fold&\(5628\)&\(11270\)&\(5922\)&\(6118\)&\(812\)&
\begin{tikzcd}[row sep=small,column sep=small,scale cd=0.75,ampersand replacement=\&]
\mathbb{Z}^{2017}\arrow[r,"\evaluation"]\&\mathbb{Z}^{\geq 2222}\oplus\mathbb{Q}^{\leq 18}\arrow[r,"\iota_\ast"]\&\mathbb{Z}^{\geq 212}\oplus\mathbb{Q}^{\leq 18}\arrow[d,two heads]\\
\mathbb{Z}^6\arrow[u,hook]\&0\arrow[l]\&\mathbb{Z}\arrow[l]
\end{tikzcd}
&\(\mathbb{Z}^{2226}/\mathbb{Z}^{209}\)\\
\hline
Ammann A\(2\)&\(68\)&\(228\)&\(163\)&\(110\)&\(45\)&
\begin{tikzcd}[row sep=small,column sep=small,scale cd=0.75,ampersand replacement=\&]
\mathbb{Z}^{6}\arrow[r,"\evaluation"]\&\mathbb{Z}^8\arrow[r,"\iota_\ast"]\&\mathbb{Z}^{7}\arrow[d,two heads]\\
\mathbb{Z}^4\arrow[u,hook]\&0\arrow[l]\&\mathbb{Z}\arrow[l]
\end{tikzcd}
&\(\mathbb{Z}^8/\mathbb{Z}^2\)\\
\hline
\begin{tabular}{@{}l@{}}Ammann A\(5\)\\decorated\end{tabular}&\(672\)&\(1608\)&\(960\)&\(824\)&\(184\)&
\begin{tikzcd}[row sep=small,column sep=small,scale cd=0.75,ampersand replacement=\&]
\mathbb{Z}^{65}\arrow[r,"\evaluation"]\&\mathbb{Z}^{80}\arrow[r,"\iota_\ast"]\&\mathbb{Z}^{24}\arrow[d,two heads]\\
\mathbb{Z}^8\arrow[u,hook]\&0\arrow[l]\&\mathbb{Z}\arrow[l]
\end{tikzcd}
&\(\mathbb{Z}^{72}/\mathbb{Z}^7\)\\
\hline
\begin{tabular}{@{}l@{}}Ammann A\(5\)\\undecorated\end{tabular}&\(220\)&\(496\)&\(281\)&\(272\)&\(65\)&
\begin{tikzcd}[row sep=small,column sep=small,scale cd=0.75,ampersand replacement=\&]
\mathbb{Z}^{16}\arrow[r,"\evaluation"]\&\mathbb{Z}^{20}\arrow[r,"\iota_\ast"]\&\mathbb{Z}^{10}\arrow[d,two heads]\\
\mathbb{Z}^5\arrow[u,hook]\&0\arrow[l]\&\mathbb{Z}\arrow[l]
\end{tikzcd}
&\(\mathbb{Z}^{16}\)\\
\hline
GKM \(9.1.1.1\)&\(4704\)&\(9604\)&\(5154\)&\(4984\)&\(660\)&
\begin{tikzcd}[row sep=small,column sep=small,scale cd=0.75,ampersand replacement=\&]
\mathbb{Z}^{602}\arrow[r,"\evaluation"]\&\mathbb{Z}^{672}\arrow[r,"\iota_\ast"]\&\mathbb{Z}^{84}\arrow[d,two heads]\\
\mathbb{Z}^{12}\arrow[u,hook]\&0\arrow[l]\&\mathbb{Z}^2\arrow[l]
\end{tikzcd}
&\(\mathbb{Z}^{686}/\mathbb{Z}^{84}\)\\
\hline
GKM \(9.2.1.1\)&\(3864\)&\(8008\)&\(4230\)&\(4060\)&\(522\)&
\begin{tikzcd}[row sep=small,column sep=small,scale cd=0.75,ampersand replacement=\&]
\mathbb{Z}^{1022}\arrow[r,"\evaluation"]\&\mathbb{Z}^{1148}\arrow[r,"\iota_\ast"]\&\mathbb{Z}^{140}\arrow[d,two heads]\\
\mathbb{Z}^{12}\arrow[u,hook]\&0\arrow[l]\&\mathbb{Z}^2\arrow[l]
\end{tikzcd}
&\(\mathbb{Z}^{1134}/\mathbb{Z}^{112}\)\\
\hline
GKM \(10.1.1.1\)&\(4872\)&\(9772\)&\(5096\)&\(5362\)&\(742\)&
\begin{tikzcd}[row sep=small,column sep=small,scale cd=0.75,ampersand replacement=\&]
\mathbb{Z}^{1331}\arrow[r,"\evaluation"]\&\mathbb{Z}^{\geq 1457}\oplus\mathbb{Q}^{\leq 27}\arrow[r,"\iota_\ast"]\&\mathbb{Z}^{\geq 133}\oplus\mathbb{Q}^{\leq 27}\arrow[d,two heads]\\
\mathbb{Z}^6\arrow[u,hook]\&0\arrow[l]\&\mathbb{Z}\arrow[l]
\end{tikzcd}
&\(\mathbb{Z}^{1470}/\mathbb{Z}^{139}\)\\
\hline
\begin{tabular}{@{}l@{}}Chair \(2\times 3\)\\(affine)\end{tabular}&\(36\)&\(104\)&\(72\)&\(62\)&\(28\)&
\begin{tikzcd}[row sep=small,column sep=small,scale cd=0.75,ampersand replacement=\&]
\begin{tabular}{@{}c@{}}\(\mathbb{Z}[1/3]\oplus\mathbb{Z}[1/2]\)\\\(\oplus\mathbb{Z}^3\)\end{tabular}\arrow[r,"\evaluation"]\&\begin{tabular}{@{}c@{}}\(\mathbb{Z}[1/6]\oplus\mathbb{Z}[1/2]\)\\\(\oplus\mathbb{Z}^4\)\end{tabular}\arrow[r,"\iota_\ast"]\&\begin{tabular}{@{}c@{}}\(\mathbb{Z}[1/6]\oplus\mathbb{Z}[1/2]\)\\\(\oplus\mathbb{Z}^2\)\end{tabular}\arrow[d,two heads]\\
\mathbb{Z}[1/3]\oplus\mathbb{Z}[1/2]\arrow[u,hook]\&0\arrow[l]\&\mathbb{Z}\arrow[l]
\end{tikzcd}
&\(\mathbb{Z}[1/3]\oplus\mathbb{Z}[1/2]\oplus\mathbb{Z}^3\)\\
\hline
\begin{tabular}{@{}l@{}}Chair \(3\times 4\)\\(affine)\end{tabular}&\(160\)&\(714\)&\(587\)&\(194\)&\(70\)&
\begin{tikzcd}[row sep=small,column sep=small,scale cd=0.75,ampersand replacement=\&]
\begin{tabular}{@{}c@{}}\(\mathbb{Z}[1/4]\oplus\mathbb{Z}[1/3]^2\)\\\(\oplus\mathbb{Z}[1/2]^3\oplus\mathbb{Z}^{16}\)\end{tabular}\arrow[r,"\evaluation"]\&\begin{tabular}{@{}c@{}}\(\mathbb{Z}[1/12]\oplus\mathbb{Z}[1/4]\)\\\(\oplus\mathbb{Z}[1/3]^2\oplus\mathbb{Z}[1/2]^4\)\\\(\oplus\mathbb{Z}^{22}\)\end{tabular}\arrow[r,"\iota_\ast"]\&\begin{tabular}{@{}c@{}}\(\mathbb{Z}[1/12]\oplus\mathbb{Z}[1/4]\)\\\(\oplus\mathbb{Z}[1/3]\oplus\mathbb{Z}[1/2]\)\\\(\oplus\mathbb{Z}^7\)\end{tabular}\arrow[d,two heads]\\
\mathbb{Z}[1/4]\oplus\mathbb{Z}[1/3]\arrow[u,hook]\&0\arrow[l]\&\mathbb{Z}\arrow[l]
\end{tikzcd}
&\(\left.\left(\begin{tabular}{@{}c@{}}\(\mathbb{Z}[1/4]\oplus\mathbb{Z}[1/3]^2\)\\\(\oplus\mathbb{Z}[1/2]^3\oplus\mathbb{Z}^{20}\)\end{tabular}\right)\middle/\mathbb{Z}^4\right.\)\\
\hline
\begin{tabular}{@{}l@{}}Chair \(3\times 5\)\\variant \(1\)\\(affine)\end{tabular}&\(404\)&\(1706\)&\(1368\)&\(548\)&\(192\)&
\begin{tikzcd}[row sep=small,column sep=small,scale cd=0.75,ampersand replacement=\&]
\begin{tabular}{@{}c@{}}\(\mathbb{Z}[1/5]^2\oplus\mathbb{Z}[1/3]^2\)\\\(\oplus\mathbb{Z}^{32}\)\end{tabular}\arrow[r,"\evaluation"]\&\begin{tabular}{@{}c@{}}\(\mathbb{Z}[1/15]\oplus\mathbb{Z}[1/5]^2\)\\\(\oplus\mathbb{Z}[1/3]^5\oplus\mathbb{Z}^{48}\)\end{tabular}\arrow[r,"\iota_\ast"]\&\begin{tabular}{@{}c@{}}\(\mathbb{Z}[1/15]\oplus\mathbb{Z}[1/5]\)\\\(\oplus\mathbb{Z}[1/3]^4\oplus\mathbb{Z}^{17}\)\end{tabular}\arrow[d,two heads]\\
\mathbb{Z}[1/5]\oplus\mathbb{Z}[1/3]\arrow[u,hook]\&0\arrow[l]\&\mathbb{Z}\arrow[l]
\end{tikzcd}
&\(\left.\left(\begin{tabular}{@{}c@{}}\(\mathbb{Z}[1/5]^2\oplus\mathbb{Z}[1/3]^2\)\\\(\oplus\mathbb{Z}^{36}\)\end{tabular}\right)\middle/\mathbb{Z}^4\right.\)\\
\hline
\begin{tabular}{@{}l@{}}Chair \(3\times 5\)\\variant \(2\)\\(affine)\end{tabular}&\(620\)&\(2958\)&\(2497\)&\(760\)&\(254\)&
\begin{tikzcd}[row sep=small,column sep=small,scale cd=0.75,ampersand replacement=\&]
\begin{tabular}{@{}c@{}}\(\mathbb{Z}[1/5]^2\oplus\mathbb{Z}[1/3]^7\)\\\(\oplus\mathbb{Z}[1/2]^2\oplus\mathbb{Z}^{94}\)\end{tabular}\arrow[r,"\evaluation"]\&\begin{tabular}{@{}c@{}}\(\mathbb{Z}[1/15]\oplus\mathbb{Z}[1/5]^2\)\\\(\oplus\mathbb{Z}[1/3]^{13}\oplus\mathbb{Z}[1/2]^4\)\\\(\oplus\mathbb{Z}^{128}\)\end{tabular}\arrow[r,"\iota_\ast"]\&\begin{tabular}{@{}c@{}}\(\mathbb{Z}[1/15]\oplus\mathbb{Z}[1/5]\)\\\(\oplus\mathbb{Z}[1/3]^7\oplus\mathbb{Z}[1/2]^2\)\\\(\oplus\mathbb{Z}^{35}\)\end{tabular}\arrow[d,two heads]\\
\mathbb{Z}[1/5]\oplus\mathbb{Z}[1/3]\arrow[u,hook]\&0\arrow[l]\&\mathbb{Z}\arrow[l]
\end{tikzcd}
&\(\left.\left(\begin{tabular}{@{}c@{}}\(\mathbb{Z}[1/5]^2\oplus\mathbb{Z}[1/3]^7\)\\\(\oplus\mathbb{Z}[1/2]^2\oplus\mathbb{Z}^{135}\)\end{tabular}\right)\middle/\mathbb{Z}^{41}\right.\)\\
\hline
Hexagon (ps)&\(24\)&\(102\)&\(80\)&\(42\)&\(22\)&
\begin{tikzcd}[row sep=small,column sep=small,scale cd=0.75,ampersand replacement=\&]
\mathbb{Z}[1/2]^2\oplus\mathbb{Z}^6\arrow[r,"\evaluation"]\&\mathbb{Z}[1/4]\oplus\mathbb{Z}^8\arrow[r,"\iota_\ast"]\&\mathbb{Z}[1/4]\oplus\mathbb{Z}^3\arrow[d,two heads]\\
\mathbb{Z}[1/2]^2\arrow[u,hook]\&0\arrow[l]\&\mathbb{Z}\arrow[l]
\end{tikzcd}
&\((\mathbb{Z}[1/2]^2\oplus\mathbb{Z}^{11})/\mathbb{Z}^5\)\\
\hline
\caption{CP is the number of collared prototiles, PI is the number of partial isometries, FT is the number of full translations, L(\(\cdot\)) is the number of loops, the second to last column is the corresponding six-term sequence in relative \(K\)-theory, and the last column is the contribution of \(K_0(C_r^\ast(\dot{G}_{AF}^{(1)}))\) towards \(K_1(C_r^\ast(\dot{G}_u))\). If an entry in the last column is not written as a quotient, then \(K_1(C_r^\ast(\dot{G}_u))\leq K_0(C_r^\ast(\dot{G}_{AF}^{(1)}))\).}
\label{table:computation-d-2}
\end{tabularx}
\end{landscape}

\clearpage

\nocite{putnam96,sims17,goncalvesramirezsolano17,goncalves11,renault06,rordamlarsenlaustsen00,davidson96,putnam18,weggeolsen93,sadun08,karoubi08,grunbaumshephard87,putnam15,lindmarcus21,clarksadun06,trevino20,solomyaktrevino22,langford40,savinienbellissard09,savinien08,forresthunton99}
\printbibliography

\end{document}